\documentclass[a4paper,10pt,reqno]{amsart}
\usepackage{amsmath}%
\usepackage{amsfonts}%
\usepackage{amscd}
\usepackage{amssymb}%
\usepackage{graphicx}
\usepackage{enumerate}
\usepackage[all]{xy}
\usepackage{hyperref}
\usepackage{mathtools}
\usepackage{ stmaryrd }
\usepackage{xcolor}
\usepackage[pageref]{backref}

\newtheorem{theorem}{Theorem}[section]
\theoremstyle{plain}
\newtheorem{proposition}[theorem]{Proposition}
\newtheorem{lemma}[theorem]{Lemma}
\newtheorem{corollary}[theorem]{Corollary}

 \theoremstyle{definition}

 \theoremstyle{remark}

 \newtheorem{remark}[theorem]{Remark}

\numberwithin{equation}{section}

\DeclareMathOperator{\im}{Im}
\DeclareMathOperator{\h}{h}\DeclareMathOperator{\rank}{r}
\DeclareMathOperator{\Spec}{Spec}
\DeclareMathOperator{\supp}{supp}
\DeclareMathOperator{\drank}{rank}
\DeclareMathOperator{\DC}{DC}
\DeclareMathOperator{\di}{d}
\DeclareMathOperator{\coker}{coker}
\DeclareMathOperator{\Aut}{Aut}

\begin{document}
\title[On Graded Division Rings]{On Graded Division Rings}
\author{Daniel E. N. Kawai}
\address[Daniel E. N. Kawai]{Department of Mathematics - IME, University of S\~ao Paulo,
Rua do Mat\~ao 1010, S\~ao Paulo, SP, 05508-090, Brazil}
\email[Daniel E. N. Kawai]{daniel.kawai@usp.br}%
\author{Javier S\'anchez}
\address[Javier S\'anchez]{Department of Mathematics - IME, University of S\~ao Paulo,
Rua do Mat\~ao 1010, S\~ao Paulo, SP, 05508-090, Brazil}
\email[Javier S\'anchez]{jsanchez@ime.usp.br}%
\thanks{The first author was supported by grant \#2017/06800-0, S\~ao Paulo Research Foundation (FAPESP), Brazil.}
\thanks{The second author was partially supported by grant \#2015/09162-9, S\~ao Paulo Research Foundation (FAPESP), Brazil,   by Grant CNPq 305773/2018-6
and by the DGI-MINECO and European Regional Development
Fund, jointly, through the grant MTM2017-83487-P.}

\date{\today}

\begin{abstract}
We develop the theory of group graded division rings parallel to the one by P. Cohn for (ungraded) division rings.
\end{abstract}

\maketitle

\setcounter{tocdepth}{1}
\tableofcontents

\section*{Introduction}

Let $R$ be a commutative ring. It is well known that the prime ideals of $R$ classify the homomorphisms from
$R$ to division rings. Indeed, for any prime ideal $P$ of $R$, we obtain a 
homomorphism from $R$ to a division ring  via the natural homomorphism
$R\rightarrow Q(R/P)$, where $Q(R/P)$ denotes the field of fractions of $R/P$. Conversely, if
$\varphi\colon R\rightarrow D$ is a homomorphism from $R$ to a division ring $D$, then $P=\ker \varphi$
is a prime ideal of $R$, $\varphi$ factors through $R\rightarrow Q(R/P)$ and therefore
the division subring of $D$ generated by the image of $R$ is $R$-isomorphic to $Q(R/P)$. Moreover,
let $P\subseteq P'$ be prime ideals of $R$. The localization of $R/P$ at the prime ideal
$P'/P$ yields a local subring of $Q(R/P)$ with residue field isomorphic to $Q(R/P')$.
This implies that any fraction 
 $ab^{-1}\in Q(R)$ which is defined in  $Q(R/P')$ it is also defined in $Q(R/P)$. Also, looking at the determinants
of matrices, one sees that any matrix with entries in $R$ that becomes  invertible in $Q(R/P')$  also  becomes invertible in $Q(R/P)$.

If the ring $R$ is not commutative, prime ideals no longer classify the homomorphisms to division rings. It may even be possible that $R$ has infinitely many different ``fields of fractions'', see 
for example \cite[Section~9]{Lam2}.

Let $R$ be any ring. An epic $R$-division ring is a ring homomorphism $R\rightarrow K$ where 
$K$ is a division ring generated by the image of $R$. 
In  \cite{CohnfreeringsI}, P. M. Cohn showed that the epic $R$-division rings are characterized
up to $R$-isomorphism by the collection of square matrices over $R$ which are carried to matrices singular
over $K$. This set of matrices is called the singular kernel of $R\rightarrow K$. He also gave the
precise conditions for a set of square matrices over $R$ to be a singular kernel, calling
such a collection a prime matrix ideal of $R$. The name comes from the fact that, if we endow the
set of square matrices over $R$ with certain two operations of sum and product (the sum being a partial operation), those sets have a similar behaviour to prime ideals. These operations are defined so that, when defined on square matrices
over a commutative ring, the determinant of the sum of matrices equals the sum of the determinants and
the determinant of a product of matrices equals the product of the determinants.
Also in \cite{CohnfreeringsI}, Cohn showed that if $\mathcal{P}$, $\mathcal{P}'$ are prime matrix ideals of $R$  and $R\rightarrow K_\mathcal{P}$, $R\rightarrow K_{\mathcal{P}'}$ are 
the corresponding epic $R$-division rings, then $\mathcal{P}\subseteq\mathcal{P}'$ if and only if there exists a local subring of $K_\mathcal{P}$
containing the image of $R$ with residue class division ring $R$-isomorphic to 
$K_{\mathcal{P}'}$.  We say that there exists a specialization from $K_\mathcal{P}$
to $K_{\mathcal{P}'}$. Furthermore, if a rational expression built up from elements of $R$
 makes sense in $K_{\mathcal{P}'}$, then it can also be evaluated in $K_\mathcal{P}$. P. M. Cohn also provided conditions on
square matrices over $R$ equivalent to the existence of (injective) homomorphisms from $R$ to division rings
and to the  existence of a \emph{best} epic $R$-division ring in the sense that a rational expression
that makes sense in some epic $R$-division ring, makes sense in it.

In \cite{Malcolmsondetermining}, P. Malcolmson described several alternative ways of determining
epic $R$-division rings. One of them is induced from the notion of rank of a matrix over a division ring.
If $\varphi\colon R\rightarrow K$ is an epic $R$-division ring, we can associate to each
matrix over $R$ the rank of this matrix when considered over $K$ via $\varphi$. He determined
which functions from the set of matrices over $R$ with values in $\mathbb{N}$ are
rank functions induced from epic $R$-division rings. Another alternative way of determining epic $R$-division rings described by P. Malcolmson is induced from the notion of dimension over a division ring.
More precisely, if $R\rightarrow K$ is an epic $R$-division ring, we can associate with each finitely
presented right $R$-module $M$ the number $\dim_K(M\otimes_R K)\in\mathbb{N}$. He described
which functions from the class of finitely presented right $R$-modules with values in
$\mathbb{N}$ are induced from epic $R$-division rings as dimensions. Another important
feature of rank functions is that,
theoretically speaking, it is easy to know when there exists a specialization from
an epic $R$-division ring to another in terms of rank functions as defined by P. Malcolmson.
In \cite{Schofieldbook}, A. Schofield gave another equivalent notion to that of epic $R$-division rings in terms
of a rank function that satisfies certain natural conditions. This time it is a function from the
class of homomorphisms between finitely generated projective right $R$-modules  with values in
$\mathbb{N}$.
We would like to remark that Sylvester rank functions with values in $\mathbb{R}_+$ have proved useful
in many different situations \cite{AraClaramuntUniqueness,AraClaramuntSylvester,Eleklamplighter,Elekinfinitedimensional,GoodearlvonNeumann,JaikinZapirainbasechangeatiyah,JaikinLopezstrongatiyah,Schofieldbook}.

The theory of group graded rings has played an important role in Ring Theory (see for example \cite{Hazrat_2016}, \cite{NastasescuvanOystaeyenMethodsgraded}) and many results in classical
ring theory have a mirrored version for group graded rings. Furthermore, if $R$ is a filtered ring, it has proved fruitful to study the associated graded ring, which usually is a simpler object, in order to obtain information  about the original ring.

The main aim of this article is to develop Cohn's theory on division rings in the context of group graded rings.
More precisely, let $\Gamma$ be a group and $R=\bigoplus\limits_{\gamma\in \Gamma}R_\gamma$
be a $\Gamma$-graded ring. A $\Gamma$-graded epic $R$-division ring is a 
homomorphism of $\Gamma$-graded rings $R\rightarrow K$ where $K$ is a
$\Gamma$-graded division ring generated by the image of $R$.
Matrices over $R$ represent homomorphisms between finitely
generated free $R$-modules. Homomorphisms of $\Gamma$-graded modules between $\Gamma$-graded free $R$-modules
are given by (what we call) \emph{homogeneous matrices}. These are $m\times n$
matrices $A$ for which there exist $\alpha_1,\dotsc,\alpha_m,\beta_1,\dotsc,\beta_n\in \Gamma$ such that
each $(i,j)$ entry of $A$ belongs to $R_{\alpha_i\beta_j^{-1}}$. We show that 
$\Gamma$-graded epic $R$-division rings $R\rightarrow K$ are characterized, up to $R$-isomorphism
of $\Gamma$-graded rings, by the collection of homogeneous matrices which are carried to singular
matrices over $K$. These sets are called the gr-singular kernel of $R\rightarrow K$. 
We give the precise conditions under which a collection of homogeneous matrices over $R$
is a gr-singular kernel and thus defining the concept of gr-prime matrix ideal.
If $\mathcal{P},\mathcal{P}'$ are gr-prime matrix ideals of $R$
and $R\rightarrow K_\mathcal{P}$, $R\rightarrow K_{{\mathcal{P}'}}$ are the corresponding $\Gamma$-graded
epic $R$-division rings, then $\mathcal{P}\subseteq\mathcal{P}'$ if and only if there exists a $\Gamma$-graded local subring of
$K_\mathcal{P}$ that contains the image of $R$ with residue class $\Gamma$-graded division ring
$R$-isomorphic to $K_{\mathcal{P}'}$ as $\Gamma$-graded rings. 
Furthermore,
if a homogeneous rational expression obtained from elements of $R$ make sense in
$K_{\mathcal{P}'}$ then it can also be evaluated in $K_{\mathcal{P}}$. We then provide conditions on the set of square
homogeneous matrices over $R$ that characterize when there exists an
(injective) homomorphism of $\Gamma$-graded rings from $R$ to a $\Gamma$-graded division ring and
when there exists a \emph{best} $\Gamma$-graded epic $R$-division ring.
We also provide the graded concepts corresponding to the different rank functions defined by 
Malcolmson and Schofield. We show they give alternative ways of determining
$\Gamma$-graded epic $R$-division rings in terms of rank functions from the set of homogeneous matrices,
from the class of $\Gamma$-graded finitely presented modules and from the class of $\Gamma$-graded
homomorphisms between $\Gamma$-graded projective $R$-modules, respectively, all of them with 
values in $\mathbb{N}$.

In the study of division rings, one of the pioneering works carrying the information  
from  the associated graded ring to the original  filtered ring
was \cite{Cohnembeddingringsskewfields}. P. M. Cohn showed that if a ring $R$ endowed with a valuation with values in
$\mathbb{Z}$ is such that its associated graded ring is a (graded) Ore domain, then $R$
can be embedded in a division ring. Other proofs of this result can be found in  
 \cite{Lichtmanvaluationmethods} and in  \cite{AsensioVandenBerghVanOystaeyen} together with
\cite{Limicrolocalizationskewfields}. More recently, a generalization of the result by Cohn has been given
by A. I. Valitskas \cite{Valitskas2014filtered}. We believe that our work could be helpful in order to generalize the result by Cohn to a greater extent than has been done by Valitskas. 

An elementary application of our theory is as follows. Suppose that  $R$ is a ring graded by a group
$\Gamma$.  As an immediate  consequence of \cite[Proposition~1.2.2]{NastasescuvanOystaeyenMethodsgraded},
one obtains that if there exists an (injective) homomorphism 
from $R$ to a division ring, then there exists an (injective) homomorphism of $\Gamma$-graded rings
from $R$ to a $\Gamma$-graded division ring. Thus if one shows that there do not exist (injective)
homomorphisms of $\Gamma$-graded rings from $R$ to  $\Gamma$-graded division rings, then there
do not exist (injective) homomorphisms 
from $R$ to  division rings. See section~\ref{sec:grprimespectrum} for other similar results.

It is also interesting to remark that the existence of an (injective) homomorphism from a $\Gamma$-graded ring $R$
to a division ring is not equivalent to the existence of a homomorphism of $\Gamma$-graded rings from $R$ to a $\Gamma$-graded division ring.
For that we  produce an easy example of a graded ring
for which there does not exist a homomorphism to a division ring but it
is embeddable in a graded division ring. 
Let $T$ be the ring obtained as localization of $\mathbb{Z}$ at the prime ideal $3\mathbb{Z}$.
Let $R$ be the ring $T[i]\subseteq \mathbb{C}$. Let
$C_2=\langle x\rangle$ be the cyclic group of order two, and let 
$\sigma\colon C_2\rightarrow \Aut(R)$ be the homomorphism of groups which sends $x$
to the automorphism induced by the complex conjugation. Set now $S=R[C_2;\sigma]$. That is,
$S$ is the skew group ring of $G$ over $R$ induced by $\sigma$. Hence $S$ is a
$C_2$-graded ring, $S=S_e+S_x$
where $S_e=R$ and $S_x=Rx$ and the product is determined by $xr=\overline{r}x$ for all $r\in R$.
Clearly $S$ is embeddable in the $C_2$ graded division ring $\mathbb{Q}[C_2;\sigma]$.
Suppose that there exists a homomorphism of rings from $S$ to a division ring $K$.
Let $\varphi\colon S\rightarrow K$ be such homomorphism. 
Since $(1-x)(1+x)=0$, then either $\varphi(1+x)=0$ or $\varphi(1-x)=0$.
If $\varphi(1+x)=0$, then $0=\varphi(1+x)=1+\varphi(x)$. Thus $\varphi(x)=-1$. But then
$(-1)\varphi(i)=\varphi(xi)=\varphi(-ix)=-\varphi(i)(-1)=\varphi(i)$. Since $\varphi(i)\neq 0$,
then $K$ has characteristic two. This is a contradiction because $\varphi$ induces a homomorphism
from $R=S_e$ to $K$ and $2$ is invertible in $R$. In the same way, it can be shown that
if $\varphi(1-x)=0$, then $\varphi(x)=1$ and, again, it  implies  that the characteristic of $K$
is $2$, a contradiction.

\medskip

In Section~\ref{sec:basicdefinitions}, we introduce some of the notation that will be used throughout the paper and provide a short survey about the results on graded rings that will be used.

Let $\Gamma$ be a group. 
A $\Gamma$-almost graded division ring is a (not necessarily graded) homomorphic image of a 
$\Gamma$-graded division ring. For example, let $K$ be a field and consider the group ring $K[\Gamma]$.
It is a $\Gamma$-graded division ring, and the augmentation map 
$K[\Gamma]\rightarrow K$, which is not a homomorphism of $\Gamma$-graded rings, endows $K$ with
as structure of $\Gamma$-almost graded division ring. 
In the nongraded context, this concept is not necessary because
a nontrivial image of a division ring is again a division ring.    
In Section~\ref{sec:almostgradeddivisionrings}, we show that if $R$ is a $\Gamma$-graded ring, $\varphi\colon R\rightarrow D$ is a homomorphism of $\Gamma$-graded rings with $D$ a $\Gamma$-graded division ring and $\psi\colon D\rightarrow E$ is a ring homomorphism where $E$ is a nonzero ring,
then the  homogeneous matrices over $R$
that become invertible via $\varphi$ and via $\psi\varphi$ are the same. Thus (a posteriori)
$\psi(D)$ determines a $\Gamma$-graded epic $R$-division ring.

The main results in Section~\ref{sec:gradedrationalclosure} are as follows.
Let  $\varphi\colon R\rightarrow D$ be a homomorphism of $\Gamma$-graded rings and
let $\Sigma$ be a set of square homogeneous matrices with entries in $R$. Suppose that the matrices
of $\Sigma$ become invertible in $D$ via $\varphi$. Then, under certain natural conditions on $\Sigma$,
the entries of the inverses of the matrices in $\Sigma$ are the homogeneous elements of a $\Gamma$-graded
subring of $R$. Moreover, if $D$ is a $\Gamma$-graded division ring generated by the
image of $\varphi$  and $\Sigma$ the set of homogeneous
matrices that become invertible under $\varphi$, then any homogeneous element of $D$ is an entry of
the inverse of some matrix in $\Sigma$.

Section~\ref{sec:categoryspecializations} begins showing that the universal localization $R_\Sigma$
of the $\Gamma$-graded ring $R$ at a set of homogeneous matrices is again a $\Gamma$-graded ring.
Then it is shown that a homomorphism of $\Gamma$-graded rings $\varphi\colon R\rightarrow D$, where $D$ is $\Gamma$-graded division ring,
is an epimorphism in the category of $\Gamma$-graded rings if and only if
$D$ is generated by the image of $\varphi$.  If this is the case, we say that $(D,\varphi)$
is a $\Gamma$-graded epic $R$-division ring and we prove that if $\Sigma$ is the
set of square homogeneous matrices that become invertible in $D$ via $\varphi$, then
$R_\Sigma$ is a $\Gamma$-graded local ring with $\Gamma$-graded residue division ring $R$-isomorphic
to $D$. 
Then the concept of gr-specialization between $\Gamma$-graded epic $R$-division rings is defined. 
The section ends showing 
that the existence of a gr-specialization from $(D,\varphi)$ to another $\Gamma$-graded epic $R$-division ring $(D',\varphi')$ is equivalent to say that all the homogeneous rational expressions (from elements of $R$) that make sense in 
$(D',\varphi')$ make sense in $(D,\varphi)$ too, and that it is also equivalent to the fact that
any homogeneous matrix over $R$ that becomes invertible in $(D',\varphi')$ becomes invertible in $(D,\varphi)$ too.

Section~\ref{sec:Malcolmsonscriterion} is devoted to the proof of the graded version of the so called Malcolmson's criterion \cite{Malcolmsonscriterion} and an important consequence. This criterion determines the kernel of the natural homomorphism from $R$
to the universal localization $R_\Sigma$ of $R$ at certain sets $\Sigma$ of homogeneous matrices. As a corollary one obtains a sufficient condition for the ring $R_\Sigma$ not to be the zero ring. 
The long and technical proof of Malcolmson's criterion consists on an element-wise construction of the ring $R_\Sigma$. This construction will also be used in the next section.

The concept of gr-prime matrix ideal is given in Section~\ref{sec:grprimematrixidealyields} and it is shown
that the different  $\Gamma$-graded epic $R$-division rings are determined by the gr-prime matrix ideals
up to $R$-isomorphism of $\Gamma$-graded rings.

In Section~\ref{sec:grmatrixideals}, 
the concepts of a gr-matrix ideal and of the radical of a gr-matrix ideal are defined and it is characterized
how is the gr-matrix ideal generated by a set of homogeneous square matrices. Then it is proved
that gr-prime matrix ideals behave like prime ideals in a commutative ring. All these concepts are used to
provide necessary and sufficient condition for the existence of  homomorphisms (embeddings) of $\Gamma$-graded rings to $\Gamma$-graded division rings.

The basic theory of Sylvester rank functions in the graded context with values in $\mathbb{N}$  is developed in Section~\ref{sec:grSylvesterrankfunctions}.
The main difference with the ungraded case stems from the  fact that, in the graded case, the same
homogeneous  matrix  can define more than one homomorphism between $\Gamma$-graded free modules.
As far as we know, this is the first paper where Sylvester rank functions are considered for graded objects.

In Section~\ref{sec:grprimespectrum}, we deal with a new situation that appears in the graded context.
If $\Gamma$ is a group and $R$ is a $\Gamma$-graded ring, then the ring $R$ can be considered as a $\Gamma/\Omega$-graded ring for any normal subgroup $\Omega$ of $\Gamma$.
Thus there are $\Gamma$-graded and $\Gamma/\Omega$-graded versions  of the concepts studied before. In this section, we try to relate them. Note that when $\Omega=\Gamma$, a $\Gamma/\Omega$-graded
epic $R$-division ring is simply  an $R$-division ring and thus one can  relate
 the theory of $\Gamma$-graded division rings and the theory of division rings as developed by Cohn.

The last section is devoted to identify inverse limits in the category of $\Gamma$-graded
epic $R$-division rings with specializations as morphisms with certain ultraproducts of $\Gamma$-graded epic
$R$-division rings. In the context of division rings, a similar result was given in 
\cite[Section~7]{Lichtmanvaluationmethodsingroupringsand}, but  our proof is more direct
and general even when specialized to the ungraded case.

A second paper is in the works where we deal, among other topics, with the graded versions
of weak algorithm, (semi)firs and (pseudo-)Sylvester domains.

\medskip

We would like to finish this introduction by pointing out that most of the techniques used in this paper are
adaptations of the ones from
the works by P. M. Cohn and P. Malcolmson. We just take credit for realizing that they can be applied in the more general setting of group graded rings.

\section{Basic definitions and notation}\label{sec:basicdefinitions}

Rings are supposed to be associative with $1$. We recall that a \emph{domain} is a nonzero ring such that
for elements $x,y$ of the ring, the equality $xy=0$ implies that either $x=0$ or $y=0$. A \emph{division ring}
is a nonzero ring such that every nonzero element is invertible.
For a ring $R$, we define $\mathbb{M}(R)$ to be the set of all square matrices of any size.
Also, for each $i$ with $1\leq i\leq n$,  let $e_i$ denote the column
$$\left(\begin{smallmatrix}
0\\\vdots\\1\\\vdots\\0
\end{smallmatrix}\right)$$
in which the $i$-th entry is $1$ and the other entries are zero.

Let $A\in M_n(R)$. We say that  $A$ is \emph{full} if whenever $A=PQ$, with $P\in M_{n\times r}(R)$ and 
$Q\in M_{r\times n}(R)$, then $r\geq n$. If we think of $A$ as an endomorphism of the free (right) $R$-module $R^n$,
it means that $A$ does not factor through $R^r$ with $r<n$.
We say that $A$ is \emph{hollow} if it has an $r\times s$ block of zeros where $r+s>n$. It is well known that a hollow matrix is not full.

Let $S$ be a ring and $f\colon R\rightarrow S$ be a
ring homomorphism. For each matrix $M$ with entries in $R$, we denote by $M^f$ the
matrix whose entries are the images  of the entries of $M$ by $f$, that is, if $a_{ij}\in R$
is the $(i,j)$-entry of $M$, then the $(i,j)$-entry of $M^f$ is $f(a_{ij})$. 
Given a set of matrices $\Sigma$, we denote $\Sigma^f=\{M^f\colon M\in\Sigma\}$. 
We say that the ring homomorphism $f\colon R\rightarrow S$ is \emph{$\Sigma$-inverting}
if the matrix $M^f$ is invertible in $S$ for each $M\in\Sigma$.

\bigskip

We proceed to give some basics on group graded rings that can be found in
\cite{NastasescuvanOystaeyenMethodsgraded} and  \cite{Hazrat_2016}, for example.

If $\Gamma$ is a group, the identity element of $\Gamma$ will be denoted by $e$.

Let $\Gamma$ be a group. A ring $R$ is called a \emph{$\Gamma$-graded ring} if 
$R=\bigoplus\limits_{\gamma\in\Gamma}R_\gamma$ where each $R_\gamma$ is an additive subgroup of $R$
and $R_\gamma R_\delta\subseteq R_{\gamma\delta}$ for all $\gamma,\delta\in\Gamma$. 
The \emph{support of $R$} is defined as the set $\supp R=\{\gamma\in\Gamma\colon R_\gamma\neq \{0\}\}$.  The set
$\h(R)=\bigcup\limits_{\gamma\in\Gamma}R_\gamma$ is called the set of \emph{homogeneous elements} of $R$.
It is well known that the
identity element $1\in R$ belongs to $R_e$,  that $R_e$ is a subring of $R$ and
that if $x\in R_\gamma$ is invertible in $R$, then $x^{-1}\in
R_{\gamma^{-1}}$. A (two-sided) ideal $I$ of $R$ is called a \emph{graded ideal} if
$I=\bigoplus\limits_{\gamma\in\Gamma}(I\cap R_\gamma)$. Thus $I$ is a graded ideal
if and only if for any $x\in I$, $x=\sum x_i$, where $x_i\in \h(R)$, implies that $x_i\in I$.
Observe that if $X\subseteq \h(R)$, then the ideal of $R$ generated by $X$ is a
graded ideal.
If $I$ is a graded ideal, then the quotient ring $R/I$ is a $\Gamma$-graded ring with $R/I=\bigoplus_{\gamma\in\Gamma}(R/I)_\gamma$,
where $(R/I)_\gamma=(R_\gamma+I)/I$.

A \emph{$\Gamma$-graded domain} is a nonzero $\Gamma$-graded ring such that if $x,y\in\h(R)$, the
equality $xy=0$ implies that either $x=0$ or $y=0$.
A \emph{$\Gamma$-graded division ring} is a nonzero $\Gamma$-graded ring such that every nonzero homogeneous element is invertible. A commutative $\Gamma$-graded division ring is a \emph{$\Gamma$-graded field}.
Clearly, any $\Gamma$-graded division ring is a $\Gamma$-graded domain.

A $\Gamma$-graded ring $R$ is called a $\Gamma$-\emph{graded local ring} if the two-sided ideal 
$\mathfrak{m}$ generated
by the noninvertible homogeneous elements is a proper ideal. In this case, the $\Gamma$-graded ring $R/\mathfrak{m}$ is a $\Gamma$-graded division ring  and it will be called the \emph{residue class $\Gamma$-graded division ring} of $R$.

For $\Gamma$-graded rings $R$ and $S$, a \emph{homomorphism of $\Gamma$-graded rings} 
$f\colon R\rightarrow S$ is a ring homomorphism such that $f(R_\gamma)\subseteq S_\gamma$ for all $\gamma\in\Gamma$.
An \emph{isomorphism of $\Gamma$-graded rings} is a homomorphism of $\Gamma$-graded rings which
is bijective. Notice that the inverse is also an isomorphism of $\Gamma$-graded rings. 

Let $\Omega$ be a normal subgroup of $\Gamma$. Consider the $\Gamma$-graded ring 
$R=\bigoplus_{\gamma\in\Gamma}R_\gamma$. It can be regarded as a $\Gamma/\Omega$-graded ring as follows
$$R=\bigoplus_{\alpha\in\Gamma/\Omega} R_\alpha,\quad \textrm{where } 
R_\alpha=\bigoplus_{\gamma\in\alpha}R_\gamma.$$

\medskip

Let $R$ be a $\Gamma$-graded ring. A $\Gamma$-\emph{graded (right) $R$-module} $M$ is defined to be
a right $R$-module with a direct sum decomposition $M=\bigoplus_{\gamma\in\Gamma} M_\gamma$, where
each $M_\gamma$ is an additive subgroup of $M$ such that $M_\lambda R_\gamma\subseteq M_{\lambda\gamma}$
for all $\lambda,\gamma\in\Gamma$. A submodule $N$ of $M$ is called a graded submodule if 
$N=\bigoplus_{\gamma\in\Gamma}(N\cap M_\gamma)$. In this case, the factor module $M/N$ forms 
a $\Gamma$-graded $R$-module with $M/N=\bigoplus_{\gamma\in\Gamma} (M/N)_\gamma$,
where $(M/N)_{\gamma}=(M_\gamma+N)/N$.

For $\Gamma$-graded $R$-modules $M$ and $N$, a \emph{homomorphism of $\Gamma$-graded $R$-modules} 
$f\colon M\rightarrow N$ is a homomorphism  of $R$-modules such that $f(M_\gamma)\subseteq N_\gamma$
for all $\gamma\in\Gamma$. In this case, $\ker f$ is a graded submodule of $M$ and $\im f$ is a
graded submodule of $N$. 

If $\Omega$ is a normal subgroup of $\Gamma$, then a $\Gamma$-graded $R$-module $M=\bigoplus_{\gamma\in\Gamma}M_\gamma$ can be regarded as a $\Gamma/\Omega$-graded over the $\Gamma/\Omega$-graded ring
$R$ as follows $$M=\bigoplus_{\alpha\in\Gamma/\Omega}M_\alpha,\quad \textrm{where } M_\alpha=\bigoplus_{\gamma\in \alpha} M_\gamma.$$
Moreover, a homomorphism of $\Gamma$-graded $R$-modules is also a homomorphism of $\Gamma/\Omega$-graded
$R$-modules.

Let $\{M_i\colon i\in I\}$ be a set of $\Gamma$-graded $R$-modules. Then $\bigoplus_{i\in I}M_i$
has a natural structure of $\Gamma$-graded $R$-module given by
$(\bigoplus_{\gamma\in \Gamma}M_i)_\gamma=\bigoplus_{i\in I}({M_i})_\gamma$.

Let $M$ be a $\Gamma$-graded right $R$-module and $N$ be a $\Gamma$-graded left $R$-module. 
Then the tensor product $M\otimes_R N$ has a natural structure of $\Gamma$-graded $\mathbb{Z}$-module 
where $(M\otimes_R N)_\gamma=\{\sum_i m_i\otimes n_i\colon m_i\in M_{\gamma'}, n_i\in N_{\gamma''},
\gamma'\gamma''=\gamma\}$.

Let $M$ be a $\Gamma$-graded $R$-module. For $\delta\in\Gamma$, we define the \emph{$\delta$-shifted}
$\Gamma$-graded $R$-module $M(\delta)$ as
$$M(\delta)=\bigoplus_{\gamma\in\Gamma} M(\delta)_\gamma,\quad \textrm{where } M(\delta)_\gamma=M_{\delta\gamma}.$$

A $\Gamma$-graded $R$-module $F$ is called a \emph{$\Gamma$-graded free $R$-module} if $F$ is a free
$R$-module with a homogeneous basis. It is well known that the $\Gamma$-graded free $R$-modules are of the form $$\bigoplus_{i\in I}R(\delta_i),\quad \textrm{where } I \textrm{ is an indexing set and } \delta_i\in 
\Gamma.$$
If  $I=\{1,\dotsc,n\}$, then $\bigoplus_{i\in I}R(\delta_i)=R(\delta_1)\oplus\dotsb\oplus R(\delta_n)$,
will also be denoted by $R^n(\overline{\delta})$ where $\overline{\delta}=(\delta_1,\dotsc,\delta_n)\in\Gamma^n$.

A $\Gamma$-graded $R$-module $P$ is called a \emph{$\Gamma$-graded projective module}
if for any diagram of $\Gamma$-graded $R$-modules and homomorphisms of $\Gamma$-graded modules
$$\xymatrix{ & P\ar[d]^u \ar@{-->}[dl]_h & \\
 M\ar[r]^g & N \ar[r] & 0 }  $$
there is a graded $R$-module homomorphism $h\colon P\rightarrow M$ with $gh=u$. As in the ungraded case, the following statements are equivalent ways of saying that  $P$ is a $\Gamma$-graded projective module 
\begin{enumerate}
	\item $P$ is $\Gamma$-graded and projective as an $R$-module.
	\item Every short exact sequence of homomorphisms of $\Gamma$-graded $R$-modules
	$0\rightarrow L\rightarrow M\rightarrow P\rightarrow 0$.
	splits via a homomorphism of $\Gamma$-graded $R$-modules.
	\item $P$ is isomorphic, as $\Gamma$-graded $R$-module, to a direct summand of a 
	$\Gamma$-graded free $R$-module.
\end{enumerate}

Let $P$ be a $\Gamma$-graded projective $R$-module and let $\Omega$
be a normal subgroup of $\Gamma$. If we regard $P$ as a $\Gamma/\Omega$-graded $R$-module, then
$P$ is also projective as a $\Gamma/\Omega$-graded $R$-module.

\bigskip

Let $\Gamma$ be a group and $R=\bigoplus\limits_{\gamma\in\Gamma}R_\gamma$ be a $\Gamma$-graded ring.
Following \cite{Hazrat_2016}, for $\overline{\alpha}=(\alpha_1,\dotsc,\alpha_m)\in\Gamma^m$ and
$\overline{\beta}=(\beta_1,\dotsc,\beta_n)$, set
$$M_{m\times n}(R)[\overline{\alpha}][\overline{\beta}]=\begin{pmatrix}
	R_{\alpha_1\beta_1^{-1}} & R_{\alpha_1\beta_2^{-1}} & \dotsb & R_{\alpha_1\beta_n^{-1}} \\
	R_{\alpha_2\beta_1^{-1}} & R_{\alpha_2\beta_2^{-1}} & \dotsb & R_{\alpha_2\beta_n^{-1}} \\
	\vdots & \vdots & \ddots & \vdots \\
	R_{\alpha_m\beta_1^{-1}} & R_{\alpha_m\beta_2^{-1}} & \dotsb & R_{\alpha_m\beta_n^{-1}}
\end{pmatrix}.$$
That is $M_{m\times n}(R)[\overline{\alpha}][\overline{\beta}]$ consists of the matrices whose
$(i,j)$-entry belongs to $R_{\alpha_i\beta_j^{-1}}$. 
Such a matrix $A\in M_{m\times n}(R)[\overline{\alpha}][\overline{\beta}]$ gives a homomorphism of $\Gamma$-graded $R$-modules $$R^n(\overline{\beta})\rightarrow
R^m(\overline{\alpha}),\quad \left(\begin{smallmatrix} x_1 \\ \vdots \\ x_n \end{smallmatrix}\right)\mapsto
A \left(\begin{smallmatrix} x_1 \\ \vdots \\ x_n \end{smallmatrix}\right),$$
and in this way $M_{m\times n}(R)[\overline{\alpha}][\overline{\beta}]$ can be identified
with the set of all homomorphisms of $\Gamma$-graded $R$-modules $R^n(\overline{\beta})\rightarrow R^m(\overline{\alpha})$.

By $A\in\mathfrak{M}_{m\times n}(R)$, we mean
that $A\in M_{m\times n}[\overline{\alpha}][\overline{\beta}]$ of some
$\overline{\alpha}\in\Gamma^m$ and $\overline{\beta}\in\Gamma^{n}$. It is
important to note that, for a matrix $A\in \mathfrak{M}_{m\times n}(R)$, it is possible that
$A\in M_{m\times n}[\overline{\alpha}][\overline{\beta}]\cap M_{m\times n}[\overline{\alpha'}][\overline{\beta'}]$ even if $\overline{\alpha}\neq \overline{\alpha'}$ or $\overline{\beta}\neq \overline{\beta'}$.
The matrix $A$ belongs to that intersection if whenever the $(i,j)$-entry of $A$ is not zero,
then  $\alpha_i\beta_j^{-1}= \alpha'_i{\beta_j'}^{-1}$.

We set
$$\mathfrak{M}_\bullet(R)=\bigcup_{m,n}\mathfrak{M}_{m\times n}(R).$$

We remark that if $A\in M_{m\times n}(R)[\overline{\alpha}][\overline{\beta}]$ and 
$B\in M_{n\times p}(R)[\overline{\beta}][\overline{\varepsilon}]$ then
$AB\in M_{m\times p}(R)[\overline{\alpha}][\overline{\varepsilon}]$.
We will say that $A,B$ are \emph{compatible}.

When $m=n$, we will write
$M_{n}(R)[\overline{\alpha}][\overline{\beta}]$ and $\mathfrak{M}_n(R)$. The set of all such matrices will be
denoted by $\mathfrak{M}(R)$, that is,
$$\mathfrak{M}(R)=\bigcup\limits_n \mathfrak{M}_n(R).$$

If $A\in M_{n}(R)[\overline{\alpha}][\overline{\beta}]$ is an invertible matrix, then 
$A^{-1}\in M_{n}(R)[\overline{\beta}][\overline{\alpha}]$.

If $\Sigma\subseteq \mathfrak{M}(R)$, we will write $\Sigma_n[\overline{\alpha}][\overline{\beta}]$
to denote the set $\Sigma\cap M_n(R)[\overline{\alpha}][\overline{\beta}]$.

\medskip

A matrix $A\in \mathfrak{M}_n(R)$ is \emph{gr-full} if every time that $A=PQ$ 
for some  matrices
$P\in M_{n\times r}[\overline{\alpha}][\overline{\lambda}]$, 
$Q\in M_{r\times n}(R)[\overline{\lambda}][\overline{\beta}]$, 
then $r\geq n$. If we think of $A$ as a homomorphism of $\Gamma$-graded modules between
two  $\Gamma$-graded free $R$-modules, it means that for all 
$\overline{\alpha},\overline{\beta}\in\Gamma^n$,  such that
$A$ defines a graded homomorphism $R^n(\overline{\beta})\rightarrow R^n(\overline{\alpha})$, then
it never factors by any graded homommorphism $R^n(\overline{\beta})\rightarrow R^r(\overline{\lambda})$ with $r< n$.

Suppose that $A\in M_n(R)[\overline{\alpha}][\overline{\beta}]$, $E\in M_n(R)$ is a permutation matrix obtained permuting the rows of $I_n$ according to the permutation $\sigma\in S_n$. Then
$E\in M_n(R)[\overline{\alpha'}][\overline{\alpha}]$,
where $\overline{\alpha'}=(\alpha_{\sigma(1)},\dotsc,\alpha_{\sigma(n)})$, and $EA\in M_n(R)[\overline{\alpha'}][\overline{\beta}]$. Similarly the matrix  $E\in M_{n}(R)[\overline{\beta}][{\overline{\beta'}}]$,
where $\overline{\beta'}=(\beta_{\sigma(1)},\dotsc,\beta_{\sigma(n)})$, and
$AE\in M_n(R)[\overline{\alpha}][\overline{\beta'}]$. 
Hence, for permutation matrices $E,F$ of appropriate size, a matrix $A\in \mathfrak{M}(R)$ is gr-full if, and only if, $EAF$ is gr-full.

A hollow matrix $A\in\mathfrak{M}(R)$ is not gr-full. Indeed, suppose that $A$ 
has an $r\times s$ block of zeros. There exist permutation matrices $E,F$ such that
$EAF=\left(\begin{smallmatrix}T & 0\\U & V \end{smallmatrix}\right)$, that is, the block of $r\times s$ zeros is in the north-east corner. Then \[
\begin{pmatrix}
T & 0\\U & V
\end{pmatrix}=
\begin{pmatrix}
T & 0\\0 & I
\end{pmatrix}
\begin{pmatrix}
I & 0\\U & V
\end{pmatrix},
\]
where $T\in M_{r\times(n-s)}(R)[\overline{\alpha}][\overline{\beta}]$, 
$U\in M_{(n-r)\times(n-s)}(R)[\overline{\delta}][\overline{\beta}]$,  $V\in M_{(n-r)\times s}(R)[\overline{\delta}][\overline{\varepsilon}]$ for some sequences $\overline{\alpha},\overline{\beta},\overline{\delta},\overline{\varepsilon}$
of elements of $\Gamma$.
The result now follows because $\left(\begin{smallmatrix}
T & 0\\0 & I
\end{smallmatrix}\right)\in M_{n\times (2n-r-s)}(R)[\overline{\alpha}*\overline{\delta}][\overline{\beta}*\overline{\delta}]$ and
$\left(\begin{smallmatrix}
I & 0\\U & V
\end{smallmatrix}\right)\in M_{(2n-r-s)\times n}(R)[\overline{\beta}*\overline{\delta}][\overline{\beta}*\overline{\varepsilon}]$.
\medskip

Let $D$ be a $\Gamma$-graded division ring and $M$ be a $\Gamma$-graded $D$-module.
As in the ungraded case, the following assertions hold true
\begin{enumerate}[(1)]
	\item Any $\Gamma$-graded $D$-module is graded free.
	\item Any $D$-linearly independent subset of $M$
	consisting of homogeneous elements can be extended to a homogeneous basis of $M$.
	\item Any two homogeneous bases of $M$ over $D$ have the same cardinality.
	\item If $N$ is a $\Gamma$-graded submodule of $M$, then $\dim_D(N)+\dim_D(M/N)=\dim_D(M)$.
\end{enumerate}
We remark that, over a $\Gamma$-graded division ring, the concepts of gr-full matrix  and
of invertible matrix coincide.

\medskip

Let $D$ be a $\Gamma$-graded division ring. Let $A\in \mathfrak{M}_{m\times n}(D)$. 

The \emph{elementary homogeneous row (column) operations} on $A$ are 
\begin{enumerate}[(1)]
	\item Interchange two rows (columns) of $A$.
	\item Multiply a row on the left (a column on the right) by a nonzero homogeneous element.
	\item Suppose that $A\in M_{m\times n}(D)[\overline{\alpha}][\overline{\beta}]$. 
	Multiply row $i$ on the left
	by an element of $R_{\alpha_j\alpha_i^{-1}}$ and add the result to row $j$ 
	(multiply column $i$ on the right by an element of $R_{\beta_i\beta_j^{-1}}$ and add the result
	to column $j$).
\end{enumerate}  
Notice that those three operations on the rows (columns) can be obtained multiplying $A$
on the left (right) by an invertible matrix in $M_m(D)[\overline{\alpha'}][\overline{\alpha}]$ 
(in $M_n(D)[\overline{\beta}][\overline{\beta'}]$).

The \emph{rank} of $A$ is the dimension of the right $D$-module spanned by its columns. The matrix $A$
can be regarded as a $D$-linear map of right $D$-modules 
$R(\beta_1)\oplus\dotsb\oplus R(\beta_n)\rightarrow R(\alpha_1)\oplus\dotsb\oplus R(\alpha_n)$. The rank of $A$ coincides with the dimension of the image of $A$. The rank of $A$ can also be computed
reducing the matrix $A$ to column echelon form by homogeneous column operations. It is the number of nonzero
columns of the column echelon form. 

The rank of $A$ equals also the dimension of the left free $D$-module spanned by its rows. The matrix
$A$ can be regarded as a $D$-linear map of left $D$-modules 
$R(\alpha^{-1})\oplus\dotsb\oplus R(\alpha_m^{-1})\rightarrow R(\beta_1^{-1})\oplus\dotsb\oplus R(\beta_n^{-1})$. The rank of $A$ coincides with the dimension of the image of $A$. The rank of $A$ can also be computed
reducing the matrix $A$ to row echelon form by homogeneous row operations. It is the number of nonzero
rows of the row echelon form.

Furthermore, the rank of $A$ coincides with the size of a largest invertible square submatrix  
(obtained by eliminating rows and/or columns). We will denote the rank of $A$ by $\drank(A)$.

\bigskip



\section{Almost graded division rings}\label{sec:almostgradeddivisionrings}

\emph{Throughout this section, let $\Gamma$ be a group}. 

We say that a ring $R$ is a \emph{$\Gamma$-almost graded ring}
if there is  a family $\{R_\gamma\colon \gamma\in\Gamma\}$ of additive subgroups $R_\gamma$ of $R$
such that $1\in R_e$, $R=\sum\limits_{\gamma\in\Gamma}R_\gamma$ and $R_\gamma R_{\gamma'}\subseteq
R_{\gamma\gamma'}$ for all $\gamma,\gamma'\in\Gamma$. The name of almost graded rings was chosen to
be compatible with the definition of almost strongly graded rings given in \cite[p.14]{NastasescuvanOystaeyenMethodsgraded}. We define $\supp{R}=\{\gamma\in\Gamma\colon R_\gamma\neq \{0\}\}$.
Given two $\Gamma$-almost graded rings $R$ and $S$, a ring homomorphism $f\colon R\rightarrow S$
is a \emph{homomorphism of $\Gamma$-almost graded rings} if $f(R_\gamma)\subseteq S_\gamma$
for all $\gamma\in\Gamma$.
Clearly, any $\Gamma$-graded ring $R=\bigoplus\limits_{\gamma\in\Gamma}R_\gamma$ is a
	$\Gamma$-almost graded ring in the natural way. Given two $\Gamma$-graded rings
	$R,S$, a homomorphism of $\Gamma$-almost graded rings is in fact a homomorphism of 
	$\Gamma$-graded rings.

 Let $R$ be a $\Gamma$-graded ring, $S$ be ring and $f\colon R\rightarrow S$
	be a ring homomorphism. Then $\im f$ is a $\Gamma$-almost graded ring with 
	$(\im f)_\gamma=f(R_\gamma)$ and the restriction $f\colon R\rightarrow \im f$ is a
	homomorphism of $\Gamma$-almost graded rings. Furthermore, any $\Gamma$-almost graded
	ring can be regarded in this way. More precisely, suppose that $R=\sum\limits_{\gamma\in\Gamma}R_\gamma$
	is a $\Gamma$-almost graded ring. Set $\widetilde{R}_\gamma$ to be a
	disjoint copy of $R_\gamma$. If $a\in R_\gamma$, denote by $\tilde{a}\in \widetilde{R}_\gamma$
	the disjoint copy of $a\in R_\gamma$. Consider the $\Gamma$-graded additive group
	$\widetilde{R}=\bigoplus\limits_{\gamma\in\Gamma} \widetilde{R}_\gamma$.
	Define $\widetilde{R}_\gamma\times \widetilde{R}_{\gamma'}\rightarrow \widetilde{R}_{\gamma\gamma'}$
	by $(\tilde{a},\tilde{b})\mapsto \widetilde{ab}$, and extend it by distributivity to
	$\widetilde{R}\times\widetilde{R}\rightarrow \widetilde{R}$. This
	endows $\widetilde{R}$ with a structure of $\Gamma$-graded ring such that  $\supp\widetilde{R}=\supp R$
	and 	$\varphi\colon\widetilde{R}=\bigoplus\limits_{\gamma\in\Gamma}\widetilde{R}_\gamma\rightarrow 
	R=\sum\limits_{\gamma\in\Gamma}R_\gamma$ determined by $\tilde{a}\mapsto a$ for all
	$a\in R_\gamma,\ \gamma\in \Gamma$, is a homomorphism of $\Gamma$-almost graded rings such that
	$\varphi(\widetilde{R}_\gamma)=R_\gamma$, $\im\varphi=R$. 

Another important example  is as follows. If $S=\bigoplus_{\alpha\in\Gamma/\Omega} S_\alpha$ is a $\Gamma/\Omega$-graded ring, then $S$ can be endowed with a structure of $\Gamma$-almost graded ring defining
$S_\gamma=S_\alpha$ for all $\gamma\in\alpha$.
Suppose that $R=\bigoplus_{\gamma\in\Gamma}R_\gamma$ is a $\Gamma$-graded ring and that	that $\Omega$ is a normal subgroup of $\Gamma$. The ring $R$
is a $\Gamma/\Omega$-graded  ring defining $R_\alpha=\bigoplus_{\gamma\in\alpha} R_\gamma$
for each $\alpha\in \Gamma/\Omega$.  If $f\colon R\rightarrow S$ is a homomorphism
of $\Gamma/\Omega$-graded rings, then it is a homomorphism of $\Gamma$-almost graded rings. 
Let $\Gamma'$ be the subgroup of $\Gamma$ generated by $\supp R$.
Observe that if $\Omega$ is a normal subgroup of $\Gamma'$ (instead of $\Gamma$), $S$ is a
$\Gamma'/\Omega$-graded ring and $f\colon R\rightarrow S$ a homomorphism of $\Gamma'/\Omega$-graded rings,
then $f$ is a homomorphism of $\Gamma'$-almost graded rings.

	\bigskip
	
	We say that a nonzero ring $E$ is a \emph{$\Gamma$-almost graded division ring} if $E$ is a $\Gamma$-almost graded
ring such that every nonzero  element $x\in E_\gamma$, $\gamma\in\Gamma$, is invertible
with inverse $x^{-1}\in E_{\gamma^{-1}}$.	
	Note that if $E$ is a $\Gamma$-almost graded division ring, then $\widetilde{E}$ is a 
	$\Gamma$-graded division ring. 
	
	The following easy result tells us that $\Gamma$-almost graded
	division rings are graded division rings although not necessarily of type $\Gamma$.
	
	\begin{lemma}\label{lem:almostgradeddivisionrings}
	Let $E$ be a $\Gamma$-almost graded division ring. The following assertions hold true.
	\begin{enumerate}[\rm(1)]
	\item If $0\neq b\in E_\gamma$, then $bE_{\gamma'}=E_{\gamma\gamma'}$ and
	$E_{\gamma'}b=E_{\gamma'\gamma}$ for $\gamma\in\Gamma$.
		\item $E_\gamma\cdot E_{\gamma'}=E_{\gamma\gamma'}$ for all $\gamma,\gamma'\in\Gamma$.
		\item $\supp E$ is a subgroup of $\Gamma$.
		\item There exists a normal subgroup $N$ of $\supp E$ such that $E$ is a 
		$\frac{\supp E}{N}$-graded division ring.
	\end{enumerate}
		\end{lemma}
	
	\begin{proof}
	If $u\in E_{\gamma\gamma'}$, then $b\cdot b^{-1}u=u$ where $b^{-1}u\in E_{\gamma'}$. 
	The other part is analogous. Thus (1) is proved.
	
	(2) is a consequence of (1).
	
	Since $1\in E_e$, then (3) follows from (2).
	
	(4) First note that, for each $\gamma\in\Gamma$,
	the condition $E_\gamma\cap E_e\neq\{0\}$ implies that $E_\gamma=E_e$. 
	Indeed, if $0\neq b\in E_\gamma\cap E_e$, (1) implies that $E_\gamma=bE_e=E_e$. 
	
	Define $N=\{\gamma\in\Gamma\colon E_\gamma=E_e\}$. We show that $N$ is a normal subgroup of $\supp\Gamma$.
	
	Clearly $e\in N$. If $\gamma,\gamma'\in N$, then $E_{\gamma\gamma'}=E_\gamma E_{\gamma'}=E_eE_e=E_e$. 
	Now, if $\gamma\in N$, then 
	$E_e=E_\gamma E_{\gamma^{-1}}=E_eE_{\gamma^{-1}}=E_{\gamma^{-1}}$. Thus $\gamma^{-1}\in N$. 
	Suppose $\gamma\in N$ and $\sigma\in\supp\Gamma$. Then 
	$E_{\sigma\gamma\sigma^{-1}}=E_\sigma E_\gamma E_{\sigma^{-1}}=E_\sigma E_e E_{\sigma^{-1}}
	= E_e$. Thus $\sigma\gamma\sigma^{-1}\in N$.

Now let $\gamma,\gamma'\in \supp \Gamma$. Then
$$\gamma^{-1}\gamma'\in N \Leftrightarrow E_{\gamma^{-1}\gamma'}=E_e
\Leftrightarrow  E_\gamma=E_{\gamma'}.$$	
And the result is proved.	
	\end{proof}
	

	Let $R$ be a $\Gamma$-graded ring, $S$ be a ring and $f\colon R\rightarrow S$ be 
	a ring homomorphism. For each $\gamma\in\Gamma$, define 
	$$(S_0)_\gamma=f(R_\gamma).$$
	If $n\geq 0$, and $(S_n)_\gamma$ has been defined for each $\gamma\in\Gamma$, define
	\[(T_{n+1})_\gamma=\{y^{-1}\colon y\in (S_n)_{\gamma^{-1}} \textrm{ and $y$ is invertible in }S \},\]
	$$(S_{n+1})_\gamma=\left. \begin{array}{c} \textrm{Additive subgroup of $S$ generated by } \\
	\{x_1x_2\dotsm x_r\colon r\in\mathbb{N},\, x_i\in (S_n)_{\gamma_i}\cup(T_{n+1})_{\gamma_i},\,
	\gamma_1\gamma_2\dotsm\gamma_n=\gamma\}\end{array}\right.$$
	Now set $(\DC(f))_\gamma=\textrm{Subgroup generated by } \bigcup_{n\geq 0}(S_n)_\gamma.$ Then
	the subring of $S$ defined by 
	$$\DC(f)=\textrm{Additive subgroup generated by }\bigcup_{\gamma\in\Gamma}(\DC(f))_\gamma$$
	is the \emph{almost graded division closure of $f\colon R\rightarrow S$}. 
	Note that $\DC(f)$ is a $\Gamma$-almost graded ring such that  
	if $x\in(\DC(f))_\gamma$ and $x$ is invertible in $S$, then $x^{-1}\in (\DC(f))_{\gamma^{-1}}$.
	It is the least subring of $S$ that contains $\im f$ and is closed under inversion
	of almost homogeneous elements. 
	
	If $\DC(f)=S$ and $\DC(f)$ is a $\Gamma$-almost graded divison ring, we say that
	$S$ is the  \emph{$\Gamma$-almost graded division ring generated by $\im f$}.

	Notice also that if $S$ is a division ring,
	then $\DC(f)$ is a $\Gamma$-almost graded division ring.

	Note that if $S$ is a $\Gamma$-graded ring, and $f\colon R\rightarrow S$ is a homomorphism
	of $\Gamma$-graded rings, then $(S_n)_\gamma\subseteq S_\gamma$ for each $n\geq 0$. Therefore 
	$(\DC(f))_\gamma\subseteq S_\gamma$ and $\DC(f)$ is a $\Gamma$-graded subring of $S$.
	It is the least subring of $S$ that contains $\im f$ and is closed under inversion
	of homogeneous elements. Moreover  if $S$ is a $\Gamma$-graded division ring,
	then $\DC(f)$ is a $\Gamma$-graded division subring of $S$. In this case, if $S=\DC(f)$
	we say that $S$ is the \emph{$\Gamma$-graded division ring generated by $\im f$}.

\begin{proposition}\label{prop:almostdivisionring}
Let $\Gamma$ be a group, $D=\bigoplus_{\gamma\in\Gamma}D_\gamma$ be a $\Gamma$-graded division ring,
and let $f\colon D\rightarrow S$ be a ring homomorphism with $S$ a nonzero ring. The following
assertions hold true.
\begin{enumerate}[\rm (1)]
\item $\DC(f)$ is a $\Gamma$-almost graded division ring with 
$$\DC(f)_\gamma=(\im f)_\gamma=\{f(x)\colon x\in D_\gamma\}$$
 and $D=\widetilde{\DC(f)}$.
\item The sets 
$$\Upsilon=\{A\in\mathfrak{M}(D)\colon A \textrm{ is invertible over }D\},$$
$$\Sigma=\{A\in\mathfrak{M}(D)\colon A^f \textrm{ is invertible over } S\}$$
coincide.
\item If $R$ is a $\Gamma$-graded ring and $\varphi\colon R\rightarrow D$
is a homomorphism of $\Gamma$-graded rings then the sets
$$\Upsilon_\varphi=\{A\in\mathfrak{M}(R)\colon A^{\varphi} \textrm{ is invertible over }D\},$$
$$\Sigma_\varphi=\{A\in\mathfrak{M}(R)\colon A^{(f\varphi)} \textrm{ is invertible over } S\}$$
coincide.
\end{enumerate}
\end{proposition}

\begin{proof}
(1) has already been proved.

(2)
Clearly, if $A\in \Upsilon$, then $A\in\Sigma$. Suppose now that 
$A\in M_n(R)[\overline{\alpha}][\overline{\beta}]$ such that
$A\notin\Upsilon.$ Then there exists a nonzero homogeneous column
$\left(\begin{smallmatrix}x_1\\ \vdots \\ x_n \end{smallmatrix}\right)\in M_{n\times 1}(R)[\overline{\beta}][\delta]$ such that $A\left(\begin{smallmatrix}x_1\\ \vdots \\ x_n \end{smallmatrix}\right)=0$.
Note that $\left(\begin{smallmatrix}x_1\\ \vdots \\ x_n \end{smallmatrix}\right)^f\neq 0$ because
$D$ is a graded division ring and $S$ is not the zero ring.
Thus $A^f\left(\begin{smallmatrix}x_1\\ \vdots \\ x_n \end{smallmatrix}\right)^f=0$, which implies that
$A\notin\Sigma$.

(3) follows from (2) because $\Sigma_\varphi=\{A\in\mathfrak{M}(R)\colon A^\varphi\in\Sigma\}$ and
$\Upsilon_\varphi=\{A\in\mathfrak{M}(R)\colon A^\varphi\in\Upsilon\}.$
\end{proof}


\section{Graded rational  closure}\label{sec:gradedrationalclosure}

\emph{Throughout this section, let $\Gamma$ be a group.}

\medskip

We begin this section introducing some important notation that will be used throughout. 

Let $\overline{\alpha}=(\alpha_1,\dotsc,\alpha_n)\in\Gamma^n$,
$\overline{\alpha'}=(\alpha_1',\dotsc,\alpha_m')\in\Gamma^m$ and $\delta\in\Gamma$, then we define
$$\overline{\alpha}*\overline{\alpha'}\coloneqq(\alpha_1,\dotsc,\alpha_n,\alpha_1',\dotsc,\alpha_m')\in\Gamma^{n+m}.$$
$$\overline{\alpha}\cdot \delta\coloneqq(\alpha_1\delta,\dotsc,\alpha_n\delta)\in\Gamma^n$$

Let  $R$ be a $\Gamma$-graded ring and $S$ be a ring.

For each  $A\in\mathfrak{M}_n(R)$,
the last column will be called $A_\infty$ and the matrix consisting
of the remaining $n-1$ columns will be called $A_\bullet$. We will
write $A=(A_\bullet\ A_\infty)$.

For each sequence $\overline{\alpha}=(\alpha_1,\dots,\alpha_n)\in\Gamma^n$, the last 
 element $\alpha_n$ will be denoted $\alpha_\infty$,
and $(\alpha_1,\dotsc,\alpha_{n-1})$ will be denoted by $\alpha_\bullet$. Thus $\overline{\alpha}=\alpha_\bullet*\alpha_\infty$

For  $u\in M_{n\times 1}(S)$, the last entry of $u$ will be denoted by $u_\infty$
and the $(n-1)\times 1$ column consisting of the remaining entries will be denoted by $u_\bullet$. Hence
$u=\left(\begin{smallmatrix}
u_\bullet\\u_\infty
\end{smallmatrix}\right).$
We remark that if $n=1$, then $A_\bullet,\alpha_\bullet,u_\bullet$ are empty and thus
$A=A_\infty$, $\overline{\alpha}=\alpha_\infty$ and $u=u_\infty$.

If $A\in\mathfrak{M}_{n\times (n+1)}(R)$, we will
denote by $A_0$ its first column, by $A_\infty$ its last column and by $A_\bullet$
the matrix consisting of the other $n-1$ columns, that is, we will write
$A=(A_0\ A_\bullet\ A_\infty)$. We will call the matrix $(A_0\ A_\bullet)$
the \emph{numerator} of $A$ and the matrix $(A_\bullet\ A_\infty)$ the \emph{denominator} of $A$.
If $A\in M_{n\times (n+1)}(R)[\overline{\alpha}][\overline{\beta}]$, we suppose $\overline{\beta}$ is divided as $\beta_0*\beta_\bullet*\beta_\infty$. 
If $u\in M_{(n+1)\times 1}(S)$, we will write 
$u=\left(\begin{smallmatrix}
	u_0 \\ u_\bullet \\ u_\infty
\end{smallmatrix}\right).$
Again, we remark that if $n=1$, then $A_\bullet$, $\beta_\bullet$, $u_\bullet$ are empty and thus
$A=(A_0\ A_\infty)$, $\overline{\beta}=(\beta_0, \beta_\infty)$ and $u=\left(\begin{smallmatrix}
	u_0 \\ u_\infty
\end{smallmatrix}\right)$.

\bigskip

Let $R=\bigoplus\limits_{\gamma\in\Gamma} R_\gamma$ be a $\Gamma$-graded ring and
$\Sigma\subseteq \mathfrak{M}(R)$.

We say that the subset $\Sigma$ of $\mathfrak{M}(R)$ is \emph{gr-lower semimultiplicative} if
it satisfies the following two conditions: 
\begin{enumerate}[(i)]
	\item $(1)\in \Sigma$, i.e. the identity matrix of size $1\times 1$ belongs to $\Sigma$.
	\item If $A\in\Sigma_n[\overline{\alpha}][\overline{\beta}]$ and  $B\in\Sigma_m[\overline{\alpha'}]\overline{\beta'}]$, then
	the matrix $\begin{pmatrix}
		A & 0 \\ C & B
	\end{pmatrix}\in\Sigma$ for any  $C\in M_{m\times n}(R)[\overline{\alpha'}][\overline{\beta}]$. 
	Notice that the matrix $\begin{pmatrix}
		A & 0 \\ C & B
	\end{pmatrix}\in M_{(n+m)}(R)[\overline{\alpha}*\overline{\alpha'}][\overline{\beta}*\overline{\beta'}]$.
\end{enumerate}
An gr-upper semimultiplicative subset of $\mathfrak{M}(R)$ is defined analogously.

A subset $\Sigma$ of $\mathfrak{M}(R)$ is \emph{gr-multiplicative} if it satisfies the following
two conditions
\begin{enumerate}[(i)]
	\item $\Sigma$ is lower gr-semimultiplicative.
	\item If $A\in\Sigma$, then $EAF\in\Sigma$ for any permutation matrices $E,F$ of appropriate size.
\end{enumerate}
\begin{remark}\label{rem:grmultiplicative}
We remark that if  $\Sigma$ is gr-multiplicative then it is also an upper gr-semimultiplicative subset of $\mathfrak{M}(R)$. Indeed,
suppose that $A\in\Sigma_n[\overline{\alpha}][\overline{\beta}]$, $B\in\Sigma_m[\overline{\alpha'}][\overline{\beta'}]$ and $C\in M_{n\times m}(R)[\overline{\alpha}][\overline{\beta'}]$. Then, since $\Sigma$
is lower gr-semimultiplicative, $\left(\begin{smallmatrix}
B & 0 \\ C & A  \end{smallmatrix}\right)\in\Sigma$. But now
$\left(\begin{smallmatrix}
A & C \\ 0 & B  \end{smallmatrix}\right)=E^{-1}\left(\begin{smallmatrix}
B & 0 \\ C & A  \end{smallmatrix}\right)E\in\Sigma$ for some permutation matrix $E$, as desired.
\end{remark}

\begin{proposition}
Let $R$ be a $\Gamma$-graded ring, $S$ be a ring and $f\colon R\rightarrow S$ be a ring homomorphism.
Then the set $$\Sigma=\{M\in\mathfrak{M}(R)\colon M^f \textrm{ is invertible over } S \}$$ is gr-multiplicative.
\end{proposition}
\begin{proof}
Clearly the $1\times 1$ matrix $(1)\in \Sigma.$. 

Let  $A\in\Sigma_n[\overline{\alpha}][\overline{\beta}]$,  $B\in\Sigma_m[\overline{\alpha'}][\overline{\beta'}]$ and $C\in M_{m\times n}[\overline{\alpha'}][\overline{\beta}]$. Then
	the matrix $\begin{pmatrix}
		A & 0 \\ C & B
	\end{pmatrix}^f$ belongs to $\Sigma$ because 
 it is invertible with inverse
\[
\begin{pmatrix}{(A^f)}^{-1}&0\\-(B^f)^{-1}(C^f)(A^f)^{-1}&(B^f)^{-1}\end{pmatrix}.
\]
Notice that if $E,F$ are permutation matrices, then $E^f, F^f$ are also permutation matrices.
Hence if $A\in \Sigma$, then $(EAF)^f$ is invertible with inverse $(F^f)^{-1}(A^f)^{-1}(E^f)^{-1}$.
\end{proof}

Note that if $S$ is a $\Gamma$-graded ring, $f\colon R\rightarrow S$ 
is a graded homomorphism and $A\in M_n(R)[\overline{\alpha}][\overline{\beta}]$,  then $A^f
\in M_n(S)[\overline{\alpha}][\overline{\beta}]$. Moreover, if $A^f$ is invertible, then $(A^f)^{-1}
\in M_n(S)[\overline{\beta}][\overline{\alpha}]$, 
and the $(j,i)$-entry of $(A^f)^{-1}$ belongs to $R_{\beta_j\alpha_i^{-1}}$.
With this in mind, we make the following definition.

Let $R=\bigoplus\limits_{\gamma\in\Gamma} R_\gamma$ be a $\Gamma$-graded ring and
$\Sigma\subseteq \mathfrak{M}(R)$.
Let $S$ be a ring (not necessarily graded) and 
 $f\colon R\rightarrow S$ be a $\Sigma$-inverting ring homomorphism. For $\gamma\in\Gamma$,
 we define the \emph{homogeneous rational closure of degree $\gamma$} as the set
$(Q_f(\Sigma))_\gamma$ consisting of all $x\in S$ such that there exist
$\overline{\alpha},\overline{\beta}\in\Gamma^n$ and
$A\in\Sigma_n[\overline{\alpha}][\overline{\beta}]$    such that $\gamma=(\alpha_i\beta_j^{-1})^{-1}=\beta_j\alpha_i^{-1}$ and $x$ is the $(j,i)$-entry of $(A^f)^{-1}$ (for some
positive integer $n$ and $i,j\in\{1,\dotsc,n\}$).
The \emph{homogeneous rational closure} is the set $$Q_f(\Sigma)=\bigcup\limits_{\gamma\in\Gamma}(Q_f(\Sigma))_\gamma.$$
The \emph{graded rational closure}, denoted by $R_f(\Sigma)$, is the additive
subgroup of $S$ generated by $Q_f(\Sigma)$.

When the set $\Sigma$ is gr-lower semimultiplicative, the graded rational closure $R_f(\Sigma)$
is a subring of $S$ as the following results show.

\begin{lemma}\label{lem:homogeneousclosure}
Let $R$ be a $\Gamma$-graded ring and 
 $\Sigma$ be a gr-lower semimultiplicative subset of $\mathfrak{M}(R)$.  
Let $S$ be a ring and $f\colon R\rightarrow S$ be a $\Sigma$-inverting ring homomorphism.
Fix $\gamma\in \Gamma$. For $x\in S$, the following conditions are equivalent.
\begin{enumerate}[\rm(1)]
	\item $x\in (Q_f(\Sigma))_\gamma$.
	\item There exist $\overline{\alpha},\overline{\beta}\in\Gamma^n$ and $A\in\Sigma_n[\overline{\alpha}][\overline{\beta}]$ 
	such that $\alpha_i=e$, $\beta_j=\gamma$ and $x$ is the $(j,i)$-entry of $(A^f)^{-1}$.

\item There exist $\overline{\alpha},\overline{\beta}\in\Gamma^n$, $A\in\Sigma_n[\overline{\alpha}][\overline{\beta}]$  and $u\in M_{n\times 1}(S)$ 
such that $\alpha_i=e$, $\beta_j=\gamma$, $u_j=x$ and $A^fu=e_i$.

	\item There exist  $\overline{\alpha},\overline{\beta}\in\Gamma^n$, $A\in\Sigma_n[\overline{\alpha}][\overline{\beta}]$, 
	$a\in M_{n\times 1}(R)[\overline{\alpha}][{e}]$ and $u\in M_{n\times 1}(S)$
	such that	$\beta_j=\gamma$, $u_j=x$ and $A^fu=a^f$.

	\item There exist $\overline{\alpha},\overline{\beta}\in\Gamma^n$, $A\in\Sigma_n[\overline{\alpha}][\overline{\beta}]$, 
	$a\in M_{n\times 1}(R)[\overline{\alpha}][{e}]$ and $u\in M_{n\times 1}(S)$
 such that
	$\beta_\infty=\gamma$, $u_\infty=x$ and $A^fu=a^f$.
	
	\item There exist $\overline{\alpha},\overline{\beta}\in\Gamma^n$, $A\in\Sigma_n[\overline{\alpha}][\overline{\beta}]$, 
	$b\in M_{1\times n}(R)[\gamma][\overline{\beta}]$ and
	$c\in M_{n\times 1}(R)[\overline{\alpha}][{e}]$ 
 such that $x=b^f(A^f)^{-1}c^f$.

\item There exist $\overline{\alpha}\in\Gamma^n$, $\overline{\beta}\in\Gamma^{n+1}$ , $A\in M_{n\times (n+1)}(R)[\overline{\alpha}][\overline{\beta}]$ and 
$u\in M_{(n+1)\times 1} (S)$ 
such that $\beta_0=e$, $\beta_\infty=\gamma$, $u_0=1$, $u_\infty=x$,
$(A_\bullet\ A_\infty)\in \Sigma$ and $A^fu=0$.
\end{enumerate}
\end{lemma}

\begin{proof}  
	(1)$\Rightarrow$(2) Let $A\in\Sigma_n[\overline{\alpha}][\overline{\beta}]$ 
such that  $x$ is the $(j,i)$-entry of $(A^f)^{-1}$ and 
  $\gamma=(\alpha_i\beta_j^{-1})^{-1}=\beta_j\alpha_i^{-1}$ for some $i,j$. 
Then $A$ can be regarded as a matrix 
in $A\in\Sigma_n[\overline{\alpha}\cdot \alpha_i^{-1}][\overline{\beta}\cdot\alpha_i^{-1}]$ 
and thus (2) follows.

(2)$\Rightarrow$(3) Suppose that (2) holds. Let $u$ be the $i$th column of 
$(A^f)^{-1}$. Then $A^fu=e_i$, as desired.


(3)$\Rightarrow$(4) It is clear because $e_i\in M_{n\times 1}(R)[\overline{\alpha}][e]$ and $e_i^f=e_i$.

(4)$\Rightarrow$(5)
Let $A\in\Sigma$, $i$, $j$, $a$ and $u$ be as in (4). Suppose that $A^fu=a^f$ 
with $u_j=x$. The matrix 
$\begin{pmatrix}A&0\\-e_j^t&1\end{pmatrix}\in \Sigma_n[\overline{\alpha}*\beta_j][\overline{\beta}*\beta_j]$. Notice  that it belongs to $\Sigma$
 because $\Sigma$ is gr-lower semimultiplicative. 
The matrix $\begin{pmatrix}a\\0\end{pmatrix}\in M_{(n+1)\times 1}(R)[\overline{\alpha}*\beta_j][e]$.
Now (5) follows from the following equality
\[
\begin{pmatrix}A^f&0\\-e_j^t&1\end{pmatrix}\begin{pmatrix}u\\x\end{pmatrix}=\begin{pmatrix}a^f\\0\end{pmatrix}
=\begin{pmatrix}a\\0\end{pmatrix}^f.
\]

(5)$\Rightarrow$(6) From (5) we obtain that $u=(A^f)^{-1}a^f$. Hence
 $$x=(e_n^t)^fu=(e_n^t)^f(A^f)^{-1}a^f.$$ Now (6) follows because $e_n^t\in M_{1\times n}(R)[\gamma][\overline{\beta}]$.

(6)$\Rightarrow$(1)
Let $A$, $b$ and $c$ as in (6). Then
\[
\begin{pmatrix}1&0&0\\c&A&0\\0&b&1\end{pmatrix}\in\Sigma_{(n+2)\times(n+2)}[e*\overline{\alpha}*\gamma][e*\overline{\beta}*\gamma].
\] Moreover
\[
\begin{pmatrix}1&0&0\\c^f&A^f&0\\0&b^f&1\end{pmatrix}^{-1}=\begin{pmatrix}1&0&0\\-(A^f)^{-1}c^f&(A^f)^{-1}&0\\b^f(A^f)^{-1}c^f&-b^f(A^f)^{-1}&1\end{pmatrix}
\]
Thus $x=b^f(A^f)^{-1}c^f$ belongs to  $(Q_f(\Sigma))_\gamma$.

(5)$\Leftrightarrow$(7)
Suppose $A\in M_n(R)[\overline{\alpha}][\overline{\beta}]$ with $\beta_\infty=\gamma$, 
$a\in M_{n\times 1}(R)[\overline{\alpha}][e]$ and $u\in M_{n\times 1}(S)$ with $u_\infty=x$. Then
the equality $A^fu=a^f$ is equivalent to the equality 
 \[
\begin{pmatrix}-a^f&A^f\end{pmatrix}\begin{pmatrix}1\\u\end{pmatrix}=0.
\]
Notice that  $(-a\ A)\in M_{n\times (n+1)}(R)[\overline{\alpha}][e*\overline{\beta}]$.
\end{proof}

\begin{theorem}\label{theo:epimorphism}
Let  $R$ be a $\Gamma$-graded ring and $\Sigma$ be
a gr-lower semimultiplicative subset of $\mathfrak{M}(R)$. Let $S$ be
a ring and $f\colon R\rightarrow S$ a $\Sigma$-inverting ring homomorphism. Then
\begin{enumerate}[\rm(1)]
\item For each $\gamma\in \Gamma$,  $f(R_\gamma)\subseteq(Q_f(\Sigma))_\gamma$.
\item If $\gamma\in \Gamma$ and $x,y\in(Q_f(\Sigma))_\gamma$, then $x+y\in(Q_f(\Sigma))_\gamma$.
\item If $\gamma,\delta\in\Gamma$ and $x\in(Q_f(\Sigma))_\gamma$, $y\in(Q_f(\Sigma))_\delta$, then $xy\in(Q_f(\Sigma))_{\gamma\delta}$.
\end{enumerate}
Hence $R_f(\Sigma)$ is a $\Gamma$-almost graded ring (which is a subring of $S$) that contains $\im(f)$. Furthermore
\begin{enumerate}[\rm(1)]
 \setcounter{enumi}{3}
\item The restriction $f\colon R\rightarrow R_f(\Sigma)$ is a ring epimorphism.
	\item If $S$ is a $\Gamma$-graded ring and $f\colon R\rightarrow S$ is 
	a homomorphism of $\Gamma$-graded rings, then
$(Q_f(\Sigma))_\gamma\subseteq S_\gamma$ for each $\gamma\in\Gamma$ and 
$R_f(\Sigma)=\bigoplus\limits_{\gamma\in\Gamma} (Q_f(\Sigma))_\gamma$ is a $\Gamma$-graded subring of $S$ such that $\h(R_f(\Sigma))=Q_f(\Sigma)$.
\end{enumerate}
\end{theorem}

\begin{proof}
(1) Let $r\in R_\gamma$. Then $f(1)f(r)=f(r)$ where $1\in M_1(R)[\gamma][\gamma]$ and $r
\in M_1(R)[\gamma][e]$. Then Lemma~\ref{lem:homogeneousclosure}(5) implies that $f(r)\in (Q_f(\Sigma))_\gamma$.

(2) Let $x,y\in(Q_f(\Sigma))_\gamma$. By Lemma~\ref{lem:homogeneousclosure}(5), 
 there exist $\overline{\alpha},\overline{\beta}\in\Gamma^n$, $A\in\Sigma_n[\overline{\alpha}][\overline{\beta}]$, 
	$a\in M_{n\times 1}(R)[\overline{\alpha}][{e}]$ and $u\in M_{n\times 1}(S)$
 such that
	$\beta_\infty=\gamma$, $u_\infty=x$ and 
	$$A^fu=(A_\bullet^f\ A_\infty^f)\begin{pmatrix} u_\bullet \\ x
\end{pmatrix}=a^f.$$ 
There also exist   $B\in\Sigma_{n'}[\overline{\alpha'}][\overline{\beta'}]$, 
	$b\in M_{n'\times 1}(R)[\overline{\alpha'}][{e}]$ and $v\in M_{n'\times 1}(S)$
 such that
	$\beta'_\infty=\gamma$, $v_\infty=y$ and 
 $$B^fv=(B_\bullet^f\ B_\infty^f)\begin{pmatrix} v_\bullet \\ y
\end{pmatrix}=b^f.$$
Then the matrix $\left(\begin{array}{cc|c}A_\bullet&A_\infty&0\\\hline0&-B_\infty&B\end{array}\right)\in\Sigma_{n+n'}[\overline{\alpha}*\overline{\alpha'}][\overline{\beta}*\overline{\beta'}]$, the column
$\begin{pmatrix}a\\b\end{pmatrix}\in M_{(n+n')\times 1}[\overline{\alpha}*\overline{\alpha'}][e]$ and we have the following equality
\[
\left(\begin{array}{cc|c}A_\bullet^f&A_\infty^f&0\\\hline0&-B_\infty^f&B^f\end{array}\right)\begin{pmatrix}u_\bullet\\x\\ \hline v_\bullet\\x+y\end{pmatrix}=\begin{pmatrix}a^f\\b^f\end{pmatrix}=\begin{pmatrix}a\\b\end{pmatrix}^f.
\]
Hence $x+y\in (Q_f(\Sigma))_\gamma$.

(3) Let $x\in(Q_f(\Sigma))_\gamma$ and $y\in(Q_f(\Sigma))_\delta$. 
 There exist $A\in\Sigma_n[\overline{\alpha}][\overline{\beta}]$, 
	$a\in M_{n\times 1}(R)[\overline{\alpha}][{e}]$ and $u\in M_{n\times 1}(S)$
 such that
	$\beta_\infty=\gamma$, $u_\infty=x$ and $A^fu=a^f$. 
 There also exist $B\in\Sigma_{n'}[\overline{\alpha'}][\overline{\beta'}]$, 
	$b\in M_{n'\times 1}(R)[\overline{\alpha'}][{e}]$ and $v\in M_{n'\times 1}(S)$
 such that
	$\beta_\infty=\delta$, $v_\infty=y$ and $B^fv=b^f$.
Now $\left(\begin{array}{cc|c}B_\bullet^f&B_\infty^f&0\\\hline0&-a^f&A^f\end{array}\right)\in
\Sigma_{n'+n}[\overline{\alpha'}*\overline{\alpha}\beta'_\infty][\overline{\beta'}*\overline{\beta}\beta'_\infty]$ 
with $(\overline{\beta'}*\overline{\beta}\beta'_\infty)_\infty=\gamma\delta$, $\begin{pmatrix}b\\0\end{pmatrix}
\in M_{(n'+n)\times 1}(R)[\overline{\alpha'}*\overline{\alpha}\beta'_\infty][e]$ and
we have the equality 
\[
\left(\begin{array}{cc|c}B_\bullet^f&B_\infty^f&0\\\hline0&-a^f&A^f\end{array}\right)\begin{pmatrix}v_\bullet\\  y\\\hline u_\bullet y\\xy\end{pmatrix}=\begin{pmatrix}b^f\\0\end{pmatrix}=\begin{pmatrix}b\\0\end{pmatrix}^f.
\]
Hence $xy\in(Q_f(\Sigma))_{\gamma\delta}$.

From (1)--(3), it is easy to show that $R_f(\Sigma)$  is a $\Gamma$-almost graded ring and a
subring of $S$. 

(4) Let $g,h\colon R_f(\Sigma)\rightarrow T$ be ring homomorphisms such that
$gf=hf$. If $x\in (Q_f(\Sigma))_\gamma$,
then $x$ is an entry of a square matrix $B$ which is the inverse of $A^f$ for some $A\in\Sigma$.
From $A^fB=BA^f=I$, it follows that $A^{gf}B^g=B^gA^{gf}=I$ and $A^{hf}B^h=B^hA^{hf}=I$.
Thus $B^g=B^h$, and $g(x)=h(x)$. Since $R_f(\Sigma)$ is generated by $(Q_f(\Sigma))_\gamma$, $\gamma\in\Gamma$,
then $f\colon R\rightarrow R_f(\Sigma)$ is a ring epimorphism.

(5) Now suppose that $S$ is a $\Gamma$-graded ring and $f\colon R\rightarrow S$ is homomorphism
of  $\Gamma$-graded rings. 
 Let $x\in(Q_f(\Sigma))_\gamma$. 
 There exist $A\in\Sigma_n[\overline{\alpha}][\overline{\beta}]$, 
	$a\in M_{n\times 1}(R)[\overline{\alpha}][{e}]$ and $u\in M_{n\times 1}(S)$
 such that
	$\beta_\infty=\gamma$, $u_\infty=x$ and $A^fu=a^f$.
Notice that $A^f\in M_n(S)[\overline{\alpha}][\overline{\beta}]$ is
an invertible  matrix 
and that $a^f\in M_{n\times 1}(S)[\overline{\alpha}][e]$. The
matrix $(A^f)^{-1}\in M_n(S)[\overline{\beta}][\overline{\alpha}]$. Now $(A^f)^{-1}$ and $a^f$ are compatible
and $u=(A^f)^{-1}a^f$. Then $x=u_\infty\in S_{\beta_\infty}$, that is $x\in S_\gamma$.

By (1)--(3), it is easy to prove that $R_f(\Sigma)$ 
is a graded subring of $S$ whose set of homogeneous elements equals $Q_f(\Sigma)$.
\end{proof}

\begin{lemma}[Cramer's rule]\label{lem:Cramersrule} 
Let $R$ be a $\Gamma$-graded ring and $\Sigma$ be
a  subset of $\mathfrak{M}(R)$. Let $S$ be
a ring and $f\colon R\rightarrow S$ be a $\Sigma$-inverting ring homomorphism.

Let $\gamma\in\Gamma$ and $x\in (Q_f(\Sigma))_\gamma$. 
Suppose that $\overline{\alpha}\in\Gamma^n$,
$\overline{\beta}\in\Gamma^{n+1}$,  $A\in M_{n\times (n+1)}(R)[\overline{\alpha}][\overline{\beta}]$ 
and 
$u\in M_{(n+1)\times 1} (S)$ 
such that $\beta_0=e$, $\beta_\infty=\gamma$, $u_0=1$, $u_\infty=x$,
$(A_\bullet\ A_\infty)\in \Sigma$ and $A^fu=0$.
 Then the following assertions hold true.
\begin{enumerate}[\rm(1)]
	\item $x$ is invertible in $S$ if, and only if, the matrix $(A_0\ A_\bullet)^f$ is invertible
	in $M_n(S)$.
	\item $x$ is a regular element of $S$ if, and only if, the matrix $(A_0\ A_\bullet)^f$ is a
	regular element of $M_n(S)$.
	\item If $x=0$, then the matrix $(A_0\ A_\bullet)^f$ is not full over $S$. Furthermore,
	if $S$ is a $\Gamma$-graded ring and $f\colon R\rightarrow S$ is a homomorphism
	of graded rings, then matrix $(A_0\ A_\bullet)^f\in M_n(S)[\overline{\alpha}][\overline{\beta}]$ is not gr-full over $S$.
\end{enumerate}
\end{lemma}

\begin{proof}
First note the equality
\begin{equation}\label{eq:Cramer}
\begin{pmatrix}A_\bullet^f&-A_0^f\end{pmatrix}=\begin{pmatrix}A_\bullet^f&A_\infty^f\end{pmatrix}\begin{pmatrix}I&u_\bullet\\0&x\end{pmatrix}.
\end{equation}
Also notice that the homogeneous matrix $(A_\bullet^f\ A_\infty^f)$ is invertible in $M_n(S)$ because $f$ is $\Sigma$-inverting.

(1) Suppose that $x$ is invertible in $S$. Then
\[
\begin{pmatrix}I&u_\bullet\\0&x\end{pmatrix}
\]
is invertible in $M_n(S)$. Hence 
$(A_\bullet^f\ -A_0^f)$ is invertible, and therefore $(A_0^f\ A_\bullet^f)$ is invertible in $M_n(S)$. 

Conversely, suppose that $(A_0^f\ A_\bullet^f)$ is invertible in
$M_n(S)$. Hence the fact that
$(A_\bullet^f\ -A_0^f)$ is invertible and  \eqref{eq:Cramer} imply that 
\[
\begin{pmatrix}I&u_\bullet\\0&x\end{pmatrix}
\]
is invertible in $M_n(S)$. Thus there exists $\left(\begin{smallmatrix}v&w\\y&z\end{smallmatrix}\right)\in M_n(S)$
such that
\[
\begin{array}{cc}
\begin{pmatrix}I&u_\bullet\\0&x\end{pmatrix}\begin{pmatrix}v&w\\y&z\end{pmatrix}=I,\quad\begin{pmatrix}v&w\\y&z\end{pmatrix}\begin{pmatrix}I&u_\bullet\\0&x\end{pmatrix}=I.
\end{array}
\]
Thus $xz=1$, $yI=0$ and $yu_\bullet+zx=1$. Therefore $x$ is invertible in $S$.

(2) Follows easily from \eqref{eq:Cramer}.

(3) Suppose that $x=0$. Then \eqref{eq:Cramer} can be expressed as
\[
\begin{array}{rcl}
\begin{pmatrix}A_\bullet^f&-A_0^f\end{pmatrix}&=&\begin{pmatrix}A_\bullet^f&A_\infty^f\end{pmatrix}\begin{pmatrix}I&u_\bullet\\0&x\end{pmatrix}\\&=&\begin{pmatrix}A_\bullet^f&A_\infty^f\end{pmatrix}\begin{pmatrix}I&u_\bullet\\0&0\end{pmatrix}\\&=&\begin{pmatrix}A_\bullet^f&A_\infty^f\end{pmatrix}\begin{pmatrix}I\\0\end{pmatrix}\begin{pmatrix}I&u_\bullet\end{pmatrix}
\end{array}
\]
which implies that $(A_\bullet^f\ -A_0^f)$ is not full and therefore
 $(A_0^f\ A_\bullet^f)$ is not full. 

If moreover $f\colon R\rightarrow S$ is a homomorphism of $\Gamma$-graded rings, then 
$(A_\bullet\ -A_0)^f\in M_n(S)[\overline{\alpha}][\beta_\bullet*e]$,
$(A_\bullet\ A_\infty)^f\in M_n(S)[\overline{\alpha}][\beta_\bullet*\beta_\infty]$,
$\begin{pmatrix}I\\0\end{pmatrix}\in M_{n\times(n-1)}(S)[\beta_\bullet*\beta_\infty][\beta_\bullet]$ and
$(I\ u_\bullet)\in M_{(n-1)\times n}[\beta_\bullet][\beta_\bullet*e]$.
It implies that the matrix $(A_\bullet\ -A_0)^f$ is not gr-full which in turn implies that 
$(A_0\ A_\bullet)^f$ is not gr-full, as desired. 
\end{proof}

Given $A$ and $x$ as in Lemma~\ref{lem:Cramersrule}, we say that $(A_0\ A_\bullet)$ is the \emph{numerator
of $x$} and $(A_\bullet\ A_\infty)$ is the \emph{denominator of $x$}. Thus, $x$ is invertible in $S$
if and only if its numerator  is invertible in $M_n(S)$.

\begin{theorem}
Let $R$ be a $\Gamma$-graded ring. Let $S$ be
a ring and $f\colon R\rightarrow S$ be a  ring homomorphism.
Set $$\Sigma=\{A\in\mathfrak{M}(R)\colon A^f \textrm{ is invertible over }S \}.$$
If $x\in (Q_f(\Sigma))_\gamma$ is invertible in $S$, then $x^{-1}\in (Q_f(\Sigma))_{\gamma^{-1}}$.

Moreover, if $S$ is a $\Gamma$-almost graded division ring and $f\colon R\rightarrow S$ is a homomorphism
of $\Gamma$-almost graded rings, then $R_f(\Sigma)$ is a $\Gamma$-almost graded division subring of $S$.
\end{theorem}

\begin{proof}
Let $x\in(Q_f(\Sigma))_\gamma$. By Lemma~\ref{lem:homogeneousclosure}(7),
there exist $A\in M_{n\times (n+1)}(R)[\overline{\alpha}][\overline{\beta}]$ and 
$u\in M_{(n+1)\times 1} (S)$ 
such that $\beta_0=e$, $\beta_\infty=\gamma$, $u_0=1$, $u_\infty=x$,
$(A_\bullet\ A_\infty)\in \Sigma$ and  $A^fu=(A_0^f\ A_\bullet^f\ A_\infty^f)
\left(\begin{smallmatrix} 1 \\ u_\bullet \\ x \end{smallmatrix}\right)=0$.
Equivalently $A_0^f+A_\bullet^f u_\bullet +A_\infty^f u_\infty=0$. Hence
$A_0^fx^{-1}+A_\bullet^f u_\bullet x^{-1}+A_\infty^f=0$, or equivalently
\[
\begin{pmatrix}A_\infty^f&A_\bullet^f&A_0^f\end{pmatrix}\begin{pmatrix}1\\u_\bullet x^{-1}\\x^{-1}\end{pmatrix}=0.
\]
Since $x$ is invertible, Cramer's rule implies that  the matrix $(A_0^f\ A_\bullet^f)$
 is invertible over $S$. Thus $(A_\bullet^f\ A_0^f)$ is invertible over $S$, and therefore
 $(A_\bullet\ A_0)\in \Sigma$.  
Moreover, notice that  $(A_\infty\ A_\bullet\ A_0)\in
M_{n\times(n+1)}(R)[\overline{\alpha}][\beta_\infty*\beta_\bullet*\beta_0]$. 
This can also be expressed as $A\in M_{n\times (n+1)}(R)[\overline{\alpha}\beta_\infty^{-1}][\beta_\infty\beta_\infty^{-1}*\beta_\bullet\beta_\infty^{-1}*\beta_0\beta_\infty^{-1}]$. 
By Lemma~\ref{lem:homogeneousclosure}(7), and observing the equality  
$\beta_\infty\beta_\infty^{-1}*\beta_\bullet\beta_\infty^{-1}*\beta_0\beta_\infty^{-1}=
 e*\beta_\bullet\beta_\infty^{-1}*\gamma^{-1}$ in $\Gamma^{n+1}$, we get that $x^{-1}\in (Q_f(\Sigma))_{\gamma^{-1}}$.

The second part follows because, by Theorem~\ref{theo:epimorphism}, $R_f(\Sigma)$ is 
a $\Gamma$-almost graded subring of $S$ which, by the foregoing, is closed under inverses of homogeneous elements.
\end{proof}

\begin{corollary}
Let  $R$ be a $\Gamma$-graded ring, $K$ be a $\Gamma$-graded division ring and
$f\colon R\rightarrow K$ be a homomorphism of $\Gamma$-graded rings. If 
$$\Sigma=\{A\in\mathfrak{M}(R)\colon A^f \textrm{ is invertible over }K \}$$
and $K$ is generated as a $\Gamma$-graded division ring by the image of $f$, then $K=R_f(\Sigma)$. \qed
\end{corollary}

We end this section with an interesting result, but that will not be used in later sections.
We show that two  elements (and by induction any finite number of elements) can be brought to a common denominator.

\begin{lemma}
Let  $R$ be a $\Gamma$-graded ring and $\Sigma$ be
a gr-lower semimultiplicative subset of $\mathfrak{M}(R)$. Let $S$ be
a ring and $f\colon R\rightarrow S$ a $\Sigma$-inverting ring homomorphism.

If $x\in(Q_f(\Sigma))_\gamma$ and $y\in(Q_f(\Sigma))_\delta$ for some $\gamma,\delta\in \Gamma$, then
they can be brought to a common denominator. 
\end{lemma}

\begin{proof}
Let $x\in(Q_f(\Sigma))_\gamma$ and $y\in(Q_f(\Sigma))_\delta$. 
 There exist $A\in M_{n\times (n+1)}(R)[\overline{\alpha}][\overline{\beta}]$ and 
$u\in M_{(n+1)\times 1} (S)$ 
such that $\beta_0=e$, $\beta_\infty=\gamma$, $u_0=1$, $u_\infty=x$,
$(A_\bullet\ A_\infty)\in \Sigma$ and  $A^fu=(A_0^f\ A_\bullet^f\ A_\infty^f)
\left(\begin{smallmatrix} 1 \\ u_\bullet \\ x \end{smallmatrix}\right)=0$.
There also exist $B\in M_{n'\times (n'+1)}(R)[\overline{\alpha'}][\overline{\beta'}]$ and 
$v\in M_{(n'+1)\times 1} (S)$ 
such that $\beta_0'=e$, $\beta_\infty'=\delta$, $v_0=1$, $v_\infty=y$,
$(B_\bullet\ B_\infty)\in \Sigma$ and  $B^fv=(B_0^f\ B_\bullet^f\ B_\infty^f)
\left(\begin{smallmatrix} 1 \\ v_\bullet \\ y \end{smallmatrix}\right)=0$. Then
 $$\left(\begin{array}{c|cc|cc}A_0&A_\bullet&A_\infty&0&0\\\hline0&0&-B_\infty&B_\bullet&B_\infty\end{array}\right)^f\left(\begin{array}{c} 1 \\\hline u_\bullet \\ x \\\hline 0\\ x \end{array}\right)=0 $$

$$\left(\begin{array}{c|cc|cc}0&A_\bullet&A_\infty&0&0\\\hline B_0&0&-B_\infty&B_\bullet&B_\infty\end{array}\right)^f\left(\begin{array}{c} 1\\ \hline 0 \\0\\ \hline v_\bullet\\ y\end{array}\right)=0$$
Now 
$\left(\begin{array}{c|cc|cc}A_0&A_\bullet&A_\infty&0&0\\\hline0&0&-B_\infty&B_\bullet&B_\infty\end{array}\right)\in M_{{(n+n')}\times (n+n'+1)}(R)[\overline{\alpha}*\overline{\alpha'}{\beta_{\infty}'}^{-1}\beta_\infty]
[\overline{\beta}*\overline{\varepsilon}]$ where $\overline{\varepsilon}=\beta_\bullet'{\beta_{\infty}'}^{-1}\beta_\infty* \beta_\infty$. 
The matrix $\left(\begin{array}{c|cc|cc}0&A_\bullet&A_\infty&0&0\\\hline B_0&0&-B_\infty&B_\bullet&B_\infty\end{array}\right)$ belongs to $M_{{(n+n')}\times (n+n'+1)}(R)[\overline{\alpha}\beta_\infty^{-1}\beta_\infty'*\overline{\alpha'}][\overline{\nu}]$ where $\overline{\nu}=
\beta_0'* \beta_\bullet\beta_\infty^{-1}\beta'_\infty* \beta_\infty'* \beta_\bullet'* \beta_\infty'$.
\end{proof}

\subsection{Generalization of Cramer's Rule}

In this little appendix, we shall show that one can generalize the Cramer's rule for certain matrices over $R_f(\Sigma)$ that are like homogeneous matrices. Let us fix the positive integers $n$ and $m$.

If $r\geq 0$ and $A\in\mathfrak{M}_{(r+m)\times (n+r+m)}(R)$, we will
denote by $A_{\overline{0}}$ the submatrix consisting of its first $n$ columns, by $A_{\overline{\infty}}$ the submatrix consisting of the last $m$ columns and by $A_{\overline{\bullet}}$
the matrix consisting of the other $r$ columns, that is, we will write
$A=(A_{\overline{0}}\ A_{\overline{\bullet}}\ A_{\overline{\infty}})$.
If $A\in M_{(r+m)\times (n+r+m)}(R)[\overline{\alpha}][\overline{\beta}]$, we suppose $\overline{\beta}$ is divided as $\beta_{\overline{0}}*\beta_{\overline{\bullet}}*\beta_{\overline{\infty}}$. 
If $U\in M_{(n+r+m)\times n}(S)$, we will write 
$U=\left(\begin{smallmatrix}
	U_{\overline{0}} \\ U_{\overline{\bullet}} \\ U_{\overline{\infty}}
\end{smallmatrix}\right).$
We remark that if $r=0$, then $A_{\overline{\bullet}}$, $\beta_{\overline{\bullet}}$, $U_{\overline{\bullet}}$ are empty and thus
$A=(A_{\overline{0}}\ A_{\overline{\infty}})$, $\overline{\beta}=(\beta_{\overline{0}}, \beta_{\overline{\infty}})$ and $U=\left(\begin{smallmatrix}
	U_{\overline{0}} \\ U_{\overline{\infty}}
\end{smallmatrix}\right)$.

\begin{proposition}
Let $R$ be a $\Gamma$-graded ring and $\Sigma$ be a gr-multiplicative subset of $\mathfrak{M}(R)$. Let $S$ be a ring and $f:R\rightarrow S$ be a $\Sigma$-inverting ring homomorphism. Fix a matrix $P=(p_{ij})\in M_{m\times n}(S)$ such that for some $\overline{\gamma}\in\Gamma^m$ and $\overline{\delta}\in\Gamma^n$ we have $p_{ij}\in(Q_f(\Sigma))_{\gamma_i\delta_j^{-1}}$. Then there exist $r\geq 0$, $\overline{\alpha}\in\Gamma^{r+m}$, $\overline{\beta}\in\Gamma^{n+r+m}$ , $A\in M_{(r+m)\times (n+r+m)}(R)[\overline{\alpha}][\overline{\beta}]$ and 
$U\in M_{(n+r+m)\times n} (S)$ 
such that $\beta_{\overline{0}}=\overline{\delta}$, $\beta_{\overline{\infty}}=\overline{\gamma}$, $U_{\overline{0}}=I_n$, $U_{\overline{\infty}}=P$,
$(A_{\overline{\bullet}}\ A_{\overline{\infty}})\in \Sigma$ and $A^fU=0$.
\end{proposition}
\begin{proof}
First, suppose that $P'=(p'_{ij})\in M_{m\times n}(S)$, where $p'_{ij}\in(Q_f(\Sigma))_{\gamma_i\delta_j^{-1}}$ and there exist $r'\geq 0$, $\overline{\alpha'}\in\Gamma^{r'+m}$, $\overline{\beta'}\in\Gamma^{n+r'+m}$ , $A'\in M_{(r'+m)\times (n+r'+m)}(R)[\overline{\alpha'}][\overline{\beta'}]$ and 
$U'\in M_{(n+r'+m)\times n} (S)$ 
such that $\beta'_{\overline{0}}=\overline{\delta}$, $\beta'_{\overline{\infty}}=\overline{\gamma}$, $U'_{\overline{0}}=I_n$, $U'_{\overline{\infty}}=P'$,
$(A'_{\overline{\bullet}}\ A'_{\overline{\infty}})\in \Sigma$ and $(A')^fU'=0$, and suppose that $P''=(p''_{ij})\in M_{m\times n}(S)$,  where $p''_{ij}\in(Q_f(\Sigma))_{\gamma_i\delta_j^{-1}}$ and there exist $r''\geq 0$, $\overline{\alpha''}\in\Gamma^{r''+m}$, $\overline{\beta''}\in\Gamma^{n+r''+m}$ , $A''\in M_{(r''+m)\times (n+r''+m)}(R)[\overline{\alpha''}][\overline{\beta''}]$ and 
$u''\in M_{(n+r''+m)\times n} (S)$ 
such that $\beta''_{\overline{0}}=\overline{\delta}$, $\beta''_{\overline{\infty}}=\overline{\gamma}$, $u''_{\overline{0}}=I_n$, $u''_{\overline{\infty}}=P'$,
$(A''_{\overline{\bullet}}\ A''_{\overline{\infty}})\in \Sigma$ and $(A'')^fU''=0$. Then:
\[
\begin{pmatrix}
A'_{\overline{0}}&A'_{\overline{\bullet}}&A'_{\overline{\infty}}&0&0\\
A''_{\overline{0}}&0&-A''_{\overline{\infty}}&A''_{\overline{\bullet}}&A''_{\overline{\infty}}
\end{pmatrix}^f
\begin{pmatrix}
I\\
U'\\
P'\\
U''\\
P'+P''
\end{pmatrix}
=0
\]
and:
\[
\begin{pmatrix}
A'_{\overline{0}}&A'_{\overline{\bullet}}&A'_{\overline{\infty}}&0&0\\
A''_{\overline{0}}&0&-A''_{\overline{\infty}}&A''_{\overline{\bullet}}&A''_{\overline{\infty}}
\end{pmatrix}\in M_{(r'+m+r''+m)\times(n+r'+m+r''+m)}(R)[\overline{\alpha}][\overline{\beta}],
\]
where $\overline{\alpha}=\overline{\alpha'}*\overline{\alpha''}$ and $\overline{\beta}=\overline{\delta}*\beta'_{\overline{\bullet}}*\overline{\gamma}*\beta''_{\overline{\bullet}}*\overline{\gamma}$.

Hence it is enough, by induction on the number of entries of $A$, to consider a matrix:
\[
X=\begin{pmatrix}
0&0&0\\
0&x&0\\
0&0&0
\end{pmatrix}\in M_{m\times n}(S)
\]
where the $(i,j)$-entry is some $x\in(Q_f(\Sigma))_{\gamma_i\delta_j^{-1}}$ and the other entries are zero. Then there exist $r\geq 0$, $\overline{\zeta}\in\Gamma^{r+1}$, $\overline{\eta}\in\Gamma^{r+2}$ , $C\in M_{(r+1)\times (r+2)}(R)[\overline{\zeta}][\overline{\eta}]$ and 
$w\in M_{(r+2)\times 1} (S)$ 
such that $\eta_0=e$, $\eta_\infty=\gamma_i\delta_j^{-1}$, $w_0=1$, $w_\infty=x$,
$(C_\bullet\ C_\infty)\in \Sigma$ and $C^fw=0$. Now we have:
\[
\begin{pmatrix}
0&0&0&0&I_{i-1}&0&0\\
0&C_0&0&C_\bullet&0&C_\infty&0\\
0&0&0&0&0&0&I_{m-i}
\end{pmatrix}^f
\begin{pmatrix}
I_{j-1}&0&0\\
0&1&0\\
0&0&I_{n-j}\\
0&w_\bullet&0\\
0&0&0\\
0&x&0\\
0&0&0
\end{pmatrix}
=0
\]
and:
\[
\begin{pmatrix}
0&0&0&0&I_{i-1}&0&0\\
0&C_0&0&C_\bullet&0&C_\infty&0\\
0&0&0&0&0&0&I_{m-i}
\end{pmatrix}\in M_{(r+m)\times(n+r+m)}(R)[\overline{\alpha}][\overline{\beta}],
\]
where:
\[
\overline{\alpha}=(\gamma_1,\dots,\gamma_{i-1})*\zeta\delta_j*(\gamma_{i+1},\dots,\gamma_m)
\]
and:
\[
\overline{\beta}=(\delta_1,\dots,\delta_{j-1})*\delta_j*(\delta_{j+1},\dots,\delta_n)*\eta_\bullet\delta_j*(\gamma_1,\dots,\gamma_{i-1})*\gamma_i*(\gamma_{i+1},\dots,\gamma_m).
\]
\end{proof}

Now we state the generalization of Cramer's rule.

\begin{lemma}[Generalized Cramer's rule]\label{lem:GenCramersrule} 
Let $R$ be a $\Gamma$-graded ring and $\Sigma$ be
a  subset of $\mathfrak{M}(R)$. Let $S$ be
a ring and $f\colon R\rightarrow S$ be a $\Sigma$-inverting ring homomorphism.

Fix a matrix $P=(p_{ij})\in M_{m\times n}(S)$ such that for some $\overline{\gamma}\in\Gamma^m$ and $\overline{\delta}\in\Gamma^n$ we have $p_{ij}\in(Q_f(\Sigma))_{\gamma_i\delta_j^{-1}}$.
Suppose that $r\geq 0$, $\overline{\alpha}\in\Gamma^{r+m}$, $\overline{\beta}\in\Gamma^{n+r+m}$ , $A\in M_{(r+m)\times (n+r+m)}(R)[\overline{\alpha}][\overline{\beta}]$ and 
$U\in M_{(n+r+m)\times n} (S)$ 
such that $\beta_{\overline{0}}=\overline{\delta}$, $\beta_{\overline{\infty}}=\overline{\gamma}$, $U_{\overline{0}}=I_n$, $U_{\overline{\infty}}=P$,
$(A_{\overline{\bullet}}\ A_{\overline{\infty}})\in \Sigma$ and $A^fU=0$.
 Then the following assertions hold true.
\begin{enumerate}[\rm(1)]
	\item $P$ is invertible over $S$ if, and only if, the $(r+m)\times (n+r)$ matrix $(A_{\overline{0}}\ A_{\overline{\bullet}})^f$ is invertible
	 over $S$.
	\item $P$ is a regular matrix over $S$ if, and only if, the $(r+m)\times (n+r)$ matrix $(A_{\overline{0}}\ A_{\overline{\bullet}})^f$ is 
	regular.
	\item If $P$ is not full over $S$, then the matrix $(A_{\overline{0}}\ A_{\overline{\bullet}})^f$ is not full over $S$. Furthermore,
	if $S$ is a $\Gamma$-graded ring, $f\colon R\rightarrow S$ is a homomorphism
	of graded rings and $P$ has some factorization $P=XY$ where $X\in M_{m\times s}(S)[\overline{\gamma}][\overline{\varepsilon}]$ and $Y\in M_{s\times n}[\overline{\varepsilon}][\overline{\delta}]$ for some $\overline{\varepsilon}\in\Gamma^s$, $s<\mathrm{min}(m,n)$, then the matrix $(A_{\overline{0}}\ A_{\overline{\bullet}})^f\in M_{(r+m)\times (n+r)}(S)[\overline{\alpha}][\beta_{\overline{0}}*\beta_{\overline{\bullet}}]$ is not gr-full over $S$.
\end{enumerate}
\end{lemma}

\begin{proof}
First note the equality
\begin{equation}\label{eq:GenCramer}
\begin{pmatrix}A_{\overline{\bullet}}^f&-A_{\overline{0}}^f\end{pmatrix}=\begin{pmatrix}A_{\overline{\bullet}}^f&A_{\overline{\infty}}^f\end{pmatrix}\begin{pmatrix}I&U_{\overline{\bullet}}\\0&P\end{pmatrix}.
\end{equation}
Also notice that the homogeneous matrix $(A_{\overline{\bullet}}^f\ A_{\overline{\infty}}^f)$ is invertible over $S$ because $f$ is $\Sigma$-inverting.

(1) Suppose that $P$ is invertible over $S$. Then
\[
\begin{pmatrix}I&U_{\overline{\bullet}}\\0&P\end{pmatrix}
\]
is invertible over $S$. Hence 
$(A_{\overline{\bullet}}^f\ -A_{\overline{0}}^f)$ is invertible, and therefore $(A_{\overline{0}}^f\ A_{\overline{\bullet}}^f)$ is invertible over $S$. 

Conversely, suppose that $(A_{\overline{0}}^f\ A_{\overline{\bullet}}^f)$ is invertible over
$S$. Hence the fact that
$(A_{\overline{\bullet}}^f\ -A_{\overline{0}}^f)$ is invertible and  \eqref{eq:GenCramer} imply that 
\[
\begin{pmatrix}I&U_{\overline{\bullet}}\\0&P\end{pmatrix}
\]
is invertible over $S$. Thus there exists $\left(\begin{smallmatrix}V&W\\Y&Z\end{smallmatrix}\right)\in M_{(r+n)\times(r+m)}(S)$
such that
\[
\begin{array}{cc}
\begin{pmatrix}I&U_{\overline{\bullet}}\\0&P\end{pmatrix}\begin{pmatrix}V&W\\Y&Z\end{pmatrix}=I_{r+m},\quad\begin{pmatrix}V&W\\Y&Z\end{pmatrix}\begin{pmatrix}I&U_{\overline{\bullet}}\\0&P\end{pmatrix}=I_{r+n}.
\end{array}
\]
Thus $PZ=I$, $YI=0$ and $YU_{\overline{\bullet}}+ZP=I$. Therefore $P$ is invertible over $S$.

(2) Follows easily from \eqref{eq:GenCramer}.

(3) Suppose that $P$ has a non-full factorization $P=XY$, where $X=M_{m\times s}(S)$, $Y\in M_{s\times n}(S)$ and $s<\mathrm{min}(m,n)$. Then \eqref{eq:GenCramer} can be expressed as
\[
\begin{array}{rcl}
\begin{pmatrix}A_{\overline{\bullet}}^f&-A_{\overline{0}}^f\end{pmatrix}&=&\begin{pmatrix}A_{\overline{\bullet}}^f&A_{\overline{\infty}}^f\end{pmatrix}\begin{pmatrix}I&U_{\overline{\bullet}}\\0&P\end{pmatrix}\\\\
&=&\begin{pmatrix}A_{\overline{\bullet}}^f&A_{\overline{\infty}}^f\end{pmatrix}\begin{pmatrix}I&U_{\overline{\bullet}}\\0&XY\end{pmatrix}\\\\
&=&\begin{pmatrix}A_{\overline{\bullet}}^f&A_{\overline{\infty}}^f\end{pmatrix}\begin{pmatrix}I&0\\0&X\end{pmatrix}\begin{pmatrix}I&U_{\overline{\bullet}}\\0&Y\end{pmatrix}
\end{array}
\]
which implies that $(A_{\overline{\bullet}}^f\ -A_{\overline{0}}^f)$ is not full and therefore
 $(A_{\overline{0}}^f\ A_{\overline{\bullet}}^f)$ is not full. 

If moreover $S$ is a $\Gamma$-graded ring, $f\colon R\rightarrow S$ is a homomorphism of $\Gamma$-graded rings and $P$ has some factorization $P=XY$ where $X\in M_{m\times s}(S)[\overline{\gamma}][\overline{\varepsilon}]$ and $Y\in M_{s\times n}[\overline{\varepsilon}][\overline{\delta}]$ for some $\overline{\varepsilon}\in\Gamma^s$, $s<\mathrm{min}(m,n)$, then
$(A_{\overline{\bullet}}\ -A_{\overline{0}})^f\in M_n(S)[\overline{\alpha}][\beta_{\overline{\bullet}}*\overline{\delta}]$,
$(A_{\overline{\bullet}}\ A_{\overline{\infty}})^f\in M_n(S)[\overline{\alpha}][\beta_{\overline{\bullet}}*\overline{\gamma}]$,
$\begin{pmatrix}I&0\\0&X\end{pmatrix}\in M_{(r+m)\times(r+s)}(S)[\beta_{\overline{\bullet}}*\overline{\gamma}][\beta_{\overline{\bullet}}*\overline{\varepsilon}]$ and
$\begin{pmatrix}I& U_{\overline{\bullet}}\\0&Y\end{pmatrix}\in M_{(r+s)\times (r+n)}[\beta_{\overline{\bullet}}*\overline{\varepsilon}][\beta_{\overline{\bullet}}*\overline{\delta}]$.
It implies that the matrix $(A_{\overline{\bullet}}\ -A_{\overline{0}})^f$ is not gr-full which in turn implies that 
$(A_{\overline{0}}\ A_{\overline{\bullet}})^f$ is not gr-full, as desired. 
\end{proof}


\section{The category of graded $R$-division rings and gr-specializations}\label{sec:categoryspecializations}

This section is an adaptation of \cite[Section~7.2]{Cohnfreeeidealringslocalization} to the graded situation.

\medskip

\emph{Throughout this section, let $\Gamma$ be a group}

Let $R=\bigoplus\limits_{\gamma\in\Gamma}R_\gamma$ be a $\Gamma$-graded ring.

A \emph{$\Gamma$-graded $R$-ring} is a pair $(K,\varphi)$ where $K$ is a $\Gamma$-graded ring and $\varphi\colon R\rightarrow K$ is a homomorphism of graded rings.
A \emph{graded $R$-subring} of $(K,\varphi)$ is a graded subring $L$ of $K$ such that $\varphi(R)\subseteq L$.

A \emph{$\Gamma$-graded $R$-division ring} is a $\Gamma$-graded $R$-ring $(K,\varphi)$ such that
$K$ is a $\Gamma$-graded division ring.  If $K=\DC(\varphi)$, that is, $K$ is the $\Gamma$-graded
division ring generated by the image of $\varphi$, we say that $(K,\varphi)$ is a \emph{$\Gamma$-graded epic
$R$-field} 

A \emph{homomorphism of $\Gamma$-graded $R$-rings} between $\Gamma$-graded $R$-rings $(K,\varphi)$
and $(K',\varphi')$ is a homomorphism of graded rings $f\colon K\rightarrow K'$
such that $\varphi'=f\circ \varphi$. If, moreover, $f\colon K\rightarrow K'$ is an isomorphism of
$\Gamma$-graded rings, we say that $f$ is an \emph{isomorphism of $\Gamma$-graded $R$-rings}.

\medskip

Let now $\Sigma\subseteq \mathfrak{M}(R)$. The \emph{universal localization of
$R$ at $\Sigma$} is a pair $(R_\Sigma,\lambda)$ where
$R_\Sigma$ is a ring and $\lambda\colon R\rightarrow R_\Sigma$
is a $\Sigma$-inverting homomorphism such that
for any other $\Sigma$-inverting ring homomorphism $f\colon R\rightarrow S$
there exists a unique ring homomorphism $F\colon R_\Sigma \rightarrow S$ with $f=F\lambda$. 

Now we give some important properties of $R_\Sigma$.

\begin{proposition}\label{prop:basicsonuniversallocalization}
Let $R$ be a $\Gamma$-graded ring and let $\Sigma\subseteq \mathfrak{M}(R)$. Then the following 
statements hold true
\begin{enumerate}[\rm(1)]
	\item There exists the universal localization $(R_\Sigma,\lambda)$ of $R$ at $\Sigma$.
	\item $\lambda\colon R\rightarrow R_\Sigma$ is a ring epimorphism.
	\item The ring $R_\Sigma$ is a $\Gamma$-graded ring, $\lambda \colon R\rightarrow R_\Sigma$
	is a homomorphism of $\Gamma$-graded rings, and $(R_\Sigma,\lambda)$ is unique up to
	isomorphism of $\Gamma$-graded $R$-rings.
	\item Suppose that $S$ is a $\Gamma$-graded ring, $f\colon R\rightarrow S$ is a $\Sigma$-inverting homomorphism 
	of $\Gamma$-graded rings and   $F\colon R_\Sigma\rightarrow S$ is the unique homomorphism of rings such that $f=F\lambda$.
	Then $F\colon R_\Sigma\rightarrow S$ is a homomorphism of $\Gamma$-graded $R$-rings. Moreover, if
	$\Sigma$ is gr-lower semimultiplicative, then  
$\im F=R_f(\Sigma)$. 
\end{enumerate} 
\end{proposition}

\begin{proof}First we construct a free ring $\mathbb{Z}\langle X\rangle$ where $X$ is constructed as follows.
For each $\gamma\in\Gamma$ and $r\in R_\gamma$, consider a symbol $x_r^\gamma$. For
each matrix $A=(a_{ij})\in\Sigma$, fix $(\overline{\alpha},\overline{\beta})$ such that $A\in M_{n}(R)[\overline{\alpha}][\overline{\beta}]$ and consider a matrix $A^*$
whose entries are symbols $A^*=(a_{ij}^*)$. Then let $X$ be the disjoint union
$$X=\{x_{r}^{\gamma}\colon r\in R_\gamma,\gamma\in\Gamma\}\cup \{a_{ij}^*\colon a_{ij} \textrm{ is the $(i,j)$-entry of }  A\in\Sigma\}. $$

Now we turn $\mathbb{Z}\langle X\rangle$ into a $\Gamma$-graded ring by giving  degrees to the elements
of $X$. If $r\in R_\gamma$, we set $x_r^\gamma$ to be of degree $\gamma$. If $A=(a_{ij})\in\Sigma$
with fixed $(\overline{\alpha},\overline{\beta})$, then $a_{ij}\in R_{\alpha_i\beta_j^{-1}}$, thus we
let  $a_{ij}^*$ of degree $\beta_i\alpha_j^{-1}$. Notice that $A^*\in M_n(\mathbb{Z}\langle X\rangle)
[\overline{\beta}][\overline{\alpha}]$. 

Let $I$ be the ideal of $\mathbb{Z}\langle X\rangle$ generated by the 
homogeneous elements of any of the following forms
\begin{itemize}
\item $x_{r+s}^{\gamma}-x_{r}^{\gamma}-x_{s}^{\gamma}$ for $r,s\in R_\gamma$.
\item $x_{rs}^{\gamma\delta}-x_{r}^{\gamma}x_{s}^{\delta}$ for $r\in R_\gamma$ and $s\in R_\delta$. 
\item $x_1^e - 1$.
\item $\sum_k x_{a_{ik}}^{\alpha_i\beta_k^{-1}}a_{kj}^*-\delta_{i,j}$ for $A\in\Sigma$.
\item $\sum_k a_{ik}^*x_{a_{kj}}^{\alpha_k\beta_j^{-1}}-\delta_{i,j}$ for $A\in\Sigma$.
\end{itemize}

Set $R_\Sigma=\mathbb{Z}\langle X\rangle /I$ and $\lambda \colon R\rightarrow R_\Sigma$
be the homomorphism of $\Gamma$-graded rings determined by $\lambda (r)=\overline{x_r^\gamma}$ for each $r\in R_\gamma$, $\gamma\in\Gamma$. 
Since $I$ is a graded ideal of 
$\mathbb{Z}\langle X\rangle$, then $R_\Sigma$ is a $\Gamma$-graded ring and 
$\lambda$ is a homomorphism of graded rings.

Suppose that $f\colon R\rightarrow S$ is a $\Sigma$-inverting ring homomorphism. 
For each $A=(a_{ij})\in\Sigma[\overline{\alpha}][\overline{\beta}]$, suppose that $(A^f)^{-1}=(b_{ij})$. 
Then there
exists a unique homomorphism of  rings $F'\colon \mathbb{Z}\langle X\rangle\rightarrow S$ such that
$F'(x_r^\gamma)=f(r)$ for each $r\in R_\gamma$, $\gamma\in\Gamma$, and $F'(a_{ij}^*)=b_{ij}$. 
Note that $I\subseteq \ker F$, and let $F\colon R_\Sigma\rightarrow S$
be the induced homomorphism. Hence $F\lambda=f$, as desired. 
To prove the uniqueness and the fact that $\lambda\colon R\rightarrow R_\Sigma$ is a ring epimorphism,
notice that  from $ F\lambda=f$, we obtain that $F(\overline{x_r^\gamma})=f(r)$, and now the same argument
of Theorem~\ref{theo:epimorphism}(4) shows that $F\big(\overline{a_{ij}^*}\big)=b_{ij}$.

Now we proceed to show (4). Suppose that $S$ is a $\Gamma$-graded ring and $f\colon R\rightarrow S$ is a $\Sigma$-inverting homomorphism
of graded rings. Notice that $A^f\in M_n(S)[\overline{\alpha}][\overline{\beta}]$ and
$(A^f)^{-1}\in M_n(S)[\overline{\beta}][\overline{\alpha}]$.
 Then
$f(r)\in S_\gamma$ for each $r\in R_\gamma$, $\gamma\in\Gamma$, and
$b_{ji}\in R_{\beta_{j}\alpha_{i}^{-1}}$ for each $A=(a_{ij})\in \Sigma_n[\overline{\alpha}][\overline{\beta}]$. Hence $F'$ and $F$ are homomorphisms of $\Gamma$-graded rings. 
Now, if $\Sigma$ is gr-lower semimultiplicative, then $R_f(\Sigma)$ is a subring of $S$ generated
by $\im f$ and the entries of the inverses of the matrices in $\Sigma^f$, and that is exactly the image of $F$.
\end{proof}

Now our aim is to show that if $(K,\varphi)$ is a \emph{$\Gamma$-graded epic
$R$-field}, then $\varphi\colon R\rightarrow K$ is in fact an epimorphism of ($\Gamma$-graded) rings. 
For the sake of completion, we preferred to give the proof of the following lemma, but this could be shown as a direct consequence of \cite[Proposition~7.2.1]{Cohnfreeeidealringslocalization} and the fact that if $f\colon R\rightarrow S$ is a homomorphism of $\Gamma$-graded rings that is
an epimorphism in the category of $\Gamma$-graded rings, then it is an epimorphism in the category of rings. 
The proof of this fact is as follows, if $g_1,g_2\colon S\rightarrow T$ are homomorphisms of rings
such that $g_1f=g_2f$, there exist homomorphisms of $\Gamma$-graded rings $\widetilde{g_1}\colon S\rightarrow \widetilde{\im g_1f}$, $\widetilde{g_2}\colon S\rightarrow \widetilde{\im g_2f}$
and homomorphism of rings $\pi\colon\widetilde{\im g_1f}\rightarrow T$ such that $\widetilde{g_1}f=\widetilde{g_2}f$ and
$g_1=\pi\widetilde{g_1}$, $g_2=\pi\widetilde{g_2}$. Since
$f$ is an epimorphism of $\Gamma$-graded rings, then $\widetilde{g_1}=\widetilde{g_2}$. Thus $g_1=g_2$.

\begin{lemma}
Let $R$, $S$ be  $\Gamma$-graded rings and $f\colon R\rightarrow S$ be a homomorphism of $\Gamma$-graded rings. 
The following statements are equivalent.
\begin{enumerate}[\rm(1)]
	\item $f$ is an epimorphism of $\Gamma$-graded rings.
	\item In the $\Gamma$-graded $S$-bimodule $S\otimes_RS$, $x\otimes 1=1\otimes x$
	for all $x\in S$.
	\item The natural map $f\colon S\otimes_R S\rightarrow S$ determined by $f(x\otimes y)=xy$
	is an isomorphism of graded $S$-bimodules.
\end{enumerate}
\end{lemma}

\begin{proof}
$(1)\Rightarrow (2)$ Consider the $\Gamma$-graded 
additive group $M=S\oplus (S\otimes_R S)$. It can be
 endowed with a structure of $\Gamma$-graded ring via the multiplication
$(x,u)(y,v)=(xy,xv+uy)$. Notice that if $(x,u)\in M_\gamma$ and $(y,v)\in M_\delta$,
then $x,u$ have degree $\gamma$ and $y,v$ have degree $\delta.$ Hence 
$xy$ and $xv+uy$ have degree $\gamma\delta$. 

Consider the homomorphisms of $\Gamma$-graded rings $g,h\colon S\rightarrow M$ defined by
$g(x)=(x,0)$ and $h(x)=(x,x\otimes 1-1\otimes x)$. Since $gf=hf$ and $f$ is an epimorphism
of graded rings, then $x\otimes 1=1\otimes x$.

$(2)\Rightarrow (1)$ Let $g,h\colon S\rightarrow T$ be homomorphisms of $\Gamma$-graded rings such
that $gf=hf$. Then there exists a well defined 
map $F\colon S\otimes_R S\rightarrow T$, $x\otimes y\mapsto g(x)h(y)$. 
For each $x\in S$, since $x\otimes 1=1\otimes x$, we obtain that $g(x)=F(x\otimes 1)=F(1\otimes x)=h(x)$.
Thus $g=h$, as desired.

$(2)\Rightarrow (3)$ First note that $f$ is a homomorphism of $\Gamma$-graded $S$-bimodules.
Clearly $f$ is surjective. Now,
since $f\left(\sum_i x_i\otimes y_i\right)=\sum_i x_iy_i$, injectivity follows from the fact
that $\sum_i x_i\otimes y_i=\sum_ix_i(1\otimes y_i)=\sum_ix_i(y_i\otimes 1)=
\sum_ix_iy_i\otimes 1=\left(\sum_i x_iy_i\right)\otimes 1$.

$(3)\Rightarrow (2)$ Since for each $x\in S$, $f(x\otimes 1)=x=f(1\otimes x)$ and $f$ is
an isomorphism, the result follows.
\end{proof}

\begin{proposition}\label{prop:epimorphism=divisionclosure}
Let $R$ be a $\Gamma$-graded ring, $K$ be a $\Gamma$-graded division ring and $f\colon R\rightarrow K$
be a homomorphism of $\Gamma$-graded rings. Then $f$ is an epimorphism of graded rings if, and only if,
$K=\DC(f)$.
\end{proposition}

\begin{proof}
Suppose that $f\colon R\rightarrow K$ is an epimorphism of $\Gamma$-graded rings. Consider the
graded division subring $\DC(f)$ of $K$. Let $\mathcal{B}$ be a set of homogeneous elements
of $K$ that is a basis of $K$ as a right $\DC(f)$-module. Then we have the following isomorphisms
of graded right $\DC(f)$-modules
\begin{eqnarray*} K\cong K\otimes_{DC(f)}K & \cong & \left(\bigoplus_{b\in\mathcal{B}}b\DC(f) \right)\otimes_{\DC(f)}K
\cong\bigoplus_{b\in\mathcal{B}}\left(b\DC(f)\otimes_{\DC(f)}K\right) \\ 
&\cong&\bigoplus_{b\in\mathcal{B}}b\otimes_{\DC(f)} K
\cong \bigoplus_{b\in\mathcal{B}}K(\gamma_b),\end{eqnarray*}
for some $\gamma_b\in\Gamma$. 
Hence $\mathcal{B}$ must consist of just one element.

Conversely, suppose that $\DC(f)=K$. Let $$\Sigma=\{A\in\mathfrak{M}(R)\colon A^f \textrm{ is invertible over }K\}.$$
Then $K=\DC(f)=R_f(\Sigma)$. By Theorem~\ref{theo:epimorphism}(4), $f\colon R\rightarrow K$ is a ring epimorphism,
and therefore an epimorphism of $\Gamma$-graded rings. 
\end{proof}

\begin{theorem}\label{theo:gradedlocal}
Let $R$ be a $\Gamma$-graded ring.
\begin{enumerate}[\rm(1)]
	\item If $\Sigma\subseteq\mathfrak{M}(R)$ is such that
	the universal localization $R_\Sigma$ is a $\Gamma$-graded local ring with maximal graded ideal $\mathfrak{m}$,
		then $R_\Sigma/\mathfrak{m}$ is a $\Gamma$-graded epic $R$-division ring. 
		
		\item Let $K$ be a $\Gamma$-almost graded division ring and
$f\colon R\rightarrow K$ be a homomorphism of $\Gamma$-almost graded rings such that 
$\DC(f)=K$. Let $$\Sigma=\{A\in\mathfrak{M}(R)\colon A^f \textrm{ is invertible over }K \}.$$
The following assertions hold true.
\begin{enumerate}[\rm(a)]
	\item $R_\Sigma$ is a $\Gamma$-graded local ring.
	\item If $\mathfrak{m}$ is the maximal graded ideal of $R_\Sigma$, then
	$R_\Sigma/\mathfrak{m}$ is a $\Gamma$-graded epic $R$-division ring satisfying the following statements.
	\begin{enumerate}[\rm(i)]
		\item There exists a surjective homomorphism of $\Gamma$-almost graded
		rings $\widetilde{F}\colon R_\Sigma/\mathfrak{m}\rightarrow K$ such that the following diagram is commutative
		$$\xymatrix{R\ar[r]^\lambda\ar[rd]_f & R_\Sigma \ar[d]^F\ar[r]^\pi & R_\Sigma/\mathfrak{m}\ar[ld]^{\widetilde{F}} \\ 
		& K & }$$
		\item If $K$ is a $\Gamma$-graded division ring, then $\widetilde{F}\colon R_\Sigma/\mathfrak{m}
		\rightarrow K$ is an isomorphism of $\Gamma$-graded epic $R$-division rings.
	\end{enumerate}
\end{enumerate}
\end{enumerate}

\end{theorem}

\begin{proof}
(1) 
The homomorphism $\lambda \colon R\rightarrow R_\Sigma$ is a ring epimorphism by Proposition~\ref{prop:basicsonuniversallocalization}(2).
The natural homomorphism $\pi\colon R_\Sigma\rightarrow R_\Sigma/\mathfrak{m}$ is surjective. Therefore
$\pi\lambda\colon R\rightarrow R_\Sigma/\mathfrak{m}$ is a ring epimorphism, thus a $\Gamma$-graded epic $R$-division ring.

(2) Let $\lambda\colon R\rightarrow R_\Sigma$ be the canonical homomorphism. Hence there exists a unique homomorphism
of $\Gamma$-almost graded $R$-rings $F\colon R_\Sigma\rightarrow K$ 
such that $F\lambda=f$. Set $\mathfrak{m}=(\ker F)_g$, in other words,
$\mathfrak{m}=\bigoplus\limits_{\gamma\in\Gamma}(\ker F\cap (R_\Sigma)_\gamma)\subseteq \ker F$. 
Let $x\in (R_\Sigma)_\gamma\setminus\mathfrak{m}$. Then $F(x)\neq 0$ and then $F(x)\in K_\gamma$ is invertible
in $K$. By Proposition~\ref{prop:basicsonuniversallocalization}(4), $R_\Sigma=R_\lambda(\Sigma)$. Thus, 
there exist $\overline{\alpha}\in\Gamma^n$, $\overline{\beta}\in\Gamma^{n+1}$ , $A\in M_{n\times (n+1)}(R)[\overline{\alpha}][\overline{\beta}]$ and 
$u\in M_{(n+1)\times 1} (S)$ 
such that $\beta_0=e$, $\beta_\infty=\gamma$, $u_0=1$, $u_\infty=x$,
$(A_\bullet\ A_\infty)\in \Sigma$ and  
$$(A_0^\lambda\ A_\bullet^\lambda\ A_\infty^\lambda)\left(\begin{smallmatrix}1 \\ u_\bullet \\ x \end{smallmatrix}\right)=0.$$
Applying $F$ to the entries of the  matrices involved we obtain 
$$(A_0^f\ A_\bullet^f\ A_\infty^f)\left(\begin{smallmatrix}1 \\ u_\bullet^F \\ F(x) \end{smallmatrix}\right)=0.$$
Since $F(x)$ is invertible, by Cramer's rule~\ref{lem:Cramersrule}, $(A_0^f\ A_\bullet^f)$ is invertible in $K$. Therefore
$(A_0\ A_\bullet)\in\Sigma$ and $(A_0^\lambda\ A_\bullet^\lambda)$ is invertible in $R_\Sigma.$ Again by Cramer's rule,
$x$ is invertible in $R_\Sigma$. Hence $R_\Sigma$ is a $\Gamma$-graded local ring where $\mathfrak{m}$ is the ideal
generated by the nonivertible homogeneous elements of $R_\Sigma$, and (a) is proved.
The ring $R_\Sigma/\mathfrak{m}$ is a $\Gamma$-graded
division ring and, by (1), (b) follows.

(i) and (ii) follow, respectively, because $\mathfrak{m}\subseteq \ker F$ and $\mathfrak{m}=\ker F$ if $K$ is a $\Gamma$-graded
division ring.
\end{proof}

We proceed to give a result that characterizes when a universal localization $R_\Sigma$ is a graded local ring. Its proof  follows 
the one given in the ungraded result in \cite[Proposition~7.2.6]{Cohnfreeeidealringslocalization}. But before that, we need following result is well known and can be found, for example, in 
\cite[Proposition~1.1.31]{Hazrat_2016}.

\begin{lemma}\label{lem:gradedlocal}
Let $R$ be a $\Gamma$-graded ring. Then $R$ is a $\Gamma$-graded local ring
if and only if $R_e$ is a local ring.\qed 
\end{lemma}

\begin{proposition}\label{lem:locuniversalgradedlocal}
Let  $R$ be a $\Gamma$-graded ring, $\Sigma$ be 
a gr-multiplicative subset of $\mathfrak{M}(R)$ and 
$\lambda\colon R\rightarrow R_\Sigma$ be the natural homomorphism of $\Gamma$-graded rings.
Then $R_\Sigma$ is a $\Gamma$-graded local ring if and only if 
it satisfies the following two conditions:
\begin{enumerate}[\rm(1)]
\item $R_\Sigma\neq\{0\}$;
	\item For a matrix $A\in\Sigma_n[\overline{\alpha'}*e][\overline{\beta'}*e]$, if
	$B$, the
	$(n,n)$-minor of $A$, is such that $B^\lambda$ is not invertible over $R_\Sigma,$ then $(A-e_{nn})^\lambda$ is invertible over $R_\Sigma$, where $e_{nn}$ denotes the matrix with $1$ in the $(n,n)$ entry and
	zeros everywhere else. 
\end{enumerate}
\end{proposition}

\begin{proof}
Consider the canonical homomorphism of $\Gamma$-graded local rings
$\lambda\colon R\rightarrow R_\Sigma$. 

Suppose that $R_\Sigma$ is a $\Gamma$-graded local ring with maximal graded ideal $\mathfrak{m}$ and
canonical homomorphism $\pi\colon R_\Sigma\rightarrow R_\Sigma/\mathfrak{m}$. 
Since $R_\Sigma$ is graded local, by definition, $R_\Sigma\neq \{0\}$. 
Recall that any matrix $C\in \mathfrak{M}(R_\Sigma)$ is invertible if and only if $C^\pi$
is invertible over $R_\Sigma/\mathfrak{m}$. Let $A\in\Sigma_n[\overline{\alpha'}*e][\overline{\beta'}*e]$
such that its $(n,n)$-minor $B$ is not invertible over $R_\Sigma$. It is enough to show that $(A-e_{nn})^\pi$
is invertible. Some non-trivial left linear combination (over the graded division ring $R/\mathfrak{m}$)
with homogeneous coefficients of
the rows of $B^\pi$ is zero. If we take the corresponding left linear combination of the first $n-1$
rows of $A^\pi$, we obtain $(0,0,\dotsc,0,c)$ where $c$ is homogeneous and $c\neq 0$, because $A^\pi$ is
invertible. We now subtract from the last row of $A$, $c^{-1}$ times this combination of the other rows
and obtain the matrix $A-e_{nn}$, which is therefore invertible in $R_\Sigma/\mathfrak{m}$ because it is
the product of the matrix corresponding to those elementary operations on $A^\pi$ times $A^\pi$.

Conversely, suppose now that conditions (1) and (2) are satisfied. By Lemma~\ref{lem:gradedlocal},
it is enough to prove that $(R_\Sigma)_e$ is a local ring. Let $x\in (R_\Sigma)_e$. 
By Lemma~\ref{lem:homogeneousclosure}(3),
there exist $\overline{\alpha},\overline{\beta}\in\Gamma^n$,
$A\in\Sigma_n[\overline{\alpha}][\overline{\beta}]$ and
$u\in M_{n\times 1}(R_\Sigma)$ such that $\alpha_i=e$, $\beta_j=e$, $u_j=x$ and $A^\lambda u=e_i$.
Since $\Sigma$ is gr-multiplicative, we may suppose that $A\in \Sigma_n[\overline{\alpha'}*e][\overline{\beta'}*e]$, $u_n=x$ and $A^\lambda u=e_n$. Suppose $x$ is not invertible in $R_\Sigma$. Equivalently, by Lemma~\ref{lem:Cramersrule},
the matrix $(A_\bullet^\lambda\, e_n^\lambda)$ is not invertible in $R_\Sigma$. This implies that
the $(n,n)$-minor of $(A_\bullet^\lambda\,\, e_n^\lambda)$, which is the $(n,n)$-minor of
$A$ is not invertible in $R_\Sigma.$ Hence, $(A-e_{nn})^\lambda$ is invertible over $R_\Sigma$.
Then the matrix $(A^\lambda)^{-1}(A-e_{nn})^\lambda=I-(A^\lambda)^{-1}e_{nn}$ is invertible in
$R_\Sigma$. Since this matrix is of the form
$$\begin{pmatrix} 1&0&\cdots&0&* \\
0 & 1 & \cdots &0 &   *\\
\vdots &  \dots & \ddots &  1&*\\
 0 & \cdots & \cdots &0            &1-x\end{pmatrix},$$
we obtain that $1-x$ is invertible in $R_\Sigma$, as desired.
\end{proof}

Now we proceed to define the category of graded epic $R$-division rings and gr-specializations.

\medskip

Let $R$ be a $\Gamma$-graded ring. 

Suppose that $(K,\varphi)$, $(L,\psi)$ are $\Gamma$-graded epic $R$-division rings and 
set $$\Sigma=\{A\in\mathfrak{M}(R)\colon A^\psi \textrm{ is invertible over } L\}.$$
If there exists a  homomorphism of $\Gamma$-graded $R$-rings
$\Phi\colon R_\Sigma\rightarrow K$, we define the
\emph{core of $L$ in $K$} as $\mathfrak{C}_{L}(K)=\Phi(R_\Sigma)$. We remark that, if it exists,
it is unique and observe that, by Proposition~\ref{prop:basicsonuniversallocalization}(4), 
$\mathfrak{C}_L(K)=R_\varphi(\Sigma)$.
By Theorem~\ref{theo:gradedlocal}(2)(a),
$R_\Sigma$ is a $\Gamma$-graded local ring. Therefore $\mathfrak{C}_L(K)$ is a $\Gamma$-graded 
local subring of $K$ that contains $R$. Moreover, the natural homomorphism
of $\Gamma$-graded $R$-rings $\Psi\colon R_\Sigma\rightarrow L$ factors through $\mathfrak{C}_L(K)$ in
a unique way, because $L\cong R_\Sigma/\mathfrak{m}$ where $\mathfrak{m}$ is the maximal graded ideal
of $R_\Sigma$.

A  \emph{gr-subhomomorphism} is a homomorphism of $\Gamma$-graded $R$-rings $f\colon K_f\rightarrow L$
where $K_f$ is a graded $R$-subring of $K$ such that $x^{-1}\in K_f$ for each $x\in \h(K_f)\setminus \ker f$.
Note that $K_f$ is a graded local subring of $K$ because any homogeneous element not in the graded ideal
$\ker f$ is invertible. Hence $K_f/\ker f$ is a $\Gamma$-graded  $R$-division ring 
contained in $L$. This implies that
$f$ is a surjective homomorphism of $\Gamma$-graded $R$-rings and that $K_f/\ker f\cong L$ is 
a   $\Gamma$-graded epic  $R$-division ring. 
For each
$A\in\Sigma$, consider $A^\varphi$ which belong to $\mathfrak{M}(K)$. Since
$K_f$ is a $\Gamma$-graded local $R$-ring whose residue graded
division ring is $L$, we get that $A^\varphi$ is invertible in $K_f$.
Thus there exists a unique homomorphism of graded $R$-rings 
$\Phi\colon R_\Sigma\rightarrow K_f\subseteq K$ and 
a commutative diagram of homomorphisms of $\Gamma$-graded $R$-rings
\begin{equation}\label{eq:commutativityspecializations}
\xymatrixrowsep{3pt}
\xymatrixcolsep{20pt}
\xymatrix{   & K_f\ar[dd]^f \\
R_\Sigma\ar[rd]_\Psi\ar[ru]^\Phi & \\ 
		&  L }
\end{equation}		
Thus $\mathfrak{C}_L(K)$
is contained in the domain of any subhomomorphism from $K$ to $L$, 
it is a $\Gamma$-graded local $R$-subring of $K$,  the restriction of any 
subhomomorphism to $\mathfrak{C}_L(K)$ is a subhomomorphism and and all
such restrictions coincide in $\mathfrak{C}_L(K)$, because of the commutativity of \eqref{eq:commutativityspecializations}.

Now we give another description of $\mathfrak{C}_L(K)$.
Let  $f\colon K_f\rightarrow L$ be a gr-subhomomorphism between the $\Gamma$-graded epic $R$-fields $(K,\varphi)$,
$(L,\psi)$. For each $\gamma\in \Gamma$
define $(c(f)_0)_\gamma=\varphi(R_\gamma)$, and if $n\geq 0$, set
$$(c(f)_{n+1})_\gamma=\left. \begin{array}{cc}\textrm{ Additive subgroup of $K$ generated by } \\ 
\{x_1\dotsm x_r\colon r\geq 1, x_i\in (c(f)_n)_{\gamma_i} \textrm{ or } 
x_i=y_i^{-1} \textrm{ where } \\ y_i\in (c(f)_n)_{\gamma_i^{-1}}\setminus\ker f,\ \gamma_1\dotsc\gamma_r=\gamma   \}\end{array}  \right.$$
Then define $c(f)_\gamma=\bigcup\limits_{n\geq 0} (c(f)_n)_\gamma$, and $C_L(K)=\bigoplus_{\gamma\in\Gamma} c(f)_\gamma.$
Note that $C_L(K)$ is a $\Gamma$-graded local $R$-subring of $K_f$ with maximal graded ideal 
$C_L(K)\cap\ker f$ and
such that the restriction $f\colon C_L(K)\rightarrow L$ is a gr-subhomomorphism. 
If we take  $K_f=\mathfrak{C}_L(K)$, then we obtain that
$C_L(K)\subseteq \mathfrak{C}_L(K)$, but since $\mathfrak{C}_L(K)$ is contained in the domain
of any gr-subhomomorphism, we get that $C_L(K)=\mathfrak{C}_L(K)$.
Roughly speaking, this equality means that any rational homogeneous 
expression obtained from the elements of (the image of) $R$
in $L$ makes sense in $K$ and the elements obtained with those rational expressions from the
elements of (the image of) $R$ in $K$  form $\mathfrak{C}_L(K)$.


Because if there exist gr-subhomomorphisms
between the $\Gamma$-graded epic $R$-division rings 
$(K,\varphi)$ and $(L,\psi)$, then they all coincide in the core, we make
the following definition. A \emph{gr-specialization} 
is the unique homomorphism of $\Gamma$-graded $R$-rings $f\colon \mathfrak{C}_L(K)\rightarrow L$.

Suppose that $(K,\varphi)$, $(L,\psi)$ and $(M,\phi)$ are $\Gamma$-graded epic $R$-division rings.
If $f\colon K_f\rightarrow L$ and $g\colon L_g\rightarrow M$ are gr-subhomomorphisms, then the restriction
$gf\colon P=f^{-1}(L_g)\rightarrow M$ is a gr-subhomomorphism which
will be called the \emph{composition gr-subhomomorphism} of $f$ and $g$. Indeed, suppose that $z\in\h(P)\setminus \ker(gf)$.
Since $g(f(z))\neq 0$, then $f(z)^{-1}\in L_g$. As $f(z)\neq 0$, and thus $z^{-1}\in K_f$, then $z^{-1}\in P$. 
We define the composition
of the corresponding gr-specializations, as the gr-specialization corresponding to the  composition gr-subhomomorphism of $f$ and $g$. In other words, the unique homomorphism of $\Gamma$-graded $R$-rings
$\mathfrak{C}_M(K)\rightarrow M$.
It follows that the composition of gr-specializations is associative.

Note that the only subhomomorphism from the $\Gamma$-graded epic $R$-division ring $(K,\varphi)$ to $(K,\varphi)$ is the identity map on $K$. Therefore
$\mathfrak{C}_K(K)=K$ and the corresponding specialization is the identity map.

We define the category $\mathcal{E}_R$ as the category whose objects are the $\Gamma$-graded
epic $R$-division rings and whose morphisms are the gr-specializations. We remark that  there is at most one
morphism between two objects in this category and that isomorphisms correspond to
isomorphisms of $\Gamma$-graded $R$-rings. Indeed, if the composition
of two gr-specializations $f$ and $g$ is the identity gr-specialization, then they have to be isomorphisms of $\Gamma$-graded $R$-rings.

An initial object $(K,\varphi)$ in the category $\mathcal{E}_R$ is a \emph{universal $\Gamma$-graded
 epic $R$-division ring}. In other words, there exists a gr-specialization
from $(K,\varphi)$ to any other $\Gamma$-graded epic $R$-division ring $(L,\psi)$. If
moreover, $\varphi\colon R\rightarrow K$ is injective, we say that this initial object is 
a \emph{universal $\Gamma$-graded epic $R$-division ring of fractions of $R$}.

Now we give the following important result.

\begin{theorem}\label{theo:specialization}
Let $R$ be a $\Gamma$-graded ring and let $(K_1,\varphi_1)$, $(K_2,\varphi_2)$ be $\Gamma$-graded epic
$R$-division rings. Set 
$$\Sigma_i=\{A\in\mathfrak{M}(R)\colon A^{\varphi_i} \textrm{ is invertible over } K_i \},\ i=1,2.$$
The following statements are equivalent.
\begin{enumerate}[\rm(1)]
	\item There exists a gr-specialization from $(K_1,\varphi_1)$ to $(K_2,\varphi_2)$.
	\item $\Sigma_2\subseteq \Sigma_1$.
	\item There exists a homomorphism $R_{\Sigma_2}\rightarrow R_{\Sigma_1}$
	of $\Gamma$-graded $R$-rings.
\end{enumerate}
Furthermore, if there exists a gr-specialization from $(K_1,\varphi_1)$ to $(K_2,\varphi_2)$ and another
gr-specialization from $(K_2,\varphi_2)$ to $(K_1,\varphi_1)$, then $K_1$ and $K_2$ are isomorphic $\Gamma$-graded $R$-rings.
\end{theorem}

\begin{proof}
$(1)\Rightarrow(2)$ By definition, there exists a homomorphism of $\Gamma$-graded $R$ rings
$\mathfrak{C}_{K_2}(K_1)\rightarrow K_2$. By definition of $\mathfrak{C}_{K_2}(K_1)$,
any matrix in $\Sigma_2$ is invertible over $\mathfrak{C}_{K_2}(K_1)\subseteq K_1$.
Thus $\Sigma_2\subseteq \Sigma_1$


$(2)\Rightarrow(3)$ If $\Sigma_2\subseteq \Sigma_1$, the universal property of $R_{\Sigma_2}$
implies the existence of a homomorphism of $\Gamma$-graded $R$-rings $R_{\Sigma_2}\rightarrow R_{\Sigma_1}$.

$(3)\Rightarrow(1)$ Consider the unique homomorphisms of $\Gamma$-graded $R$-rings
$\Phi_i\colon R_{\Sigma_i}\rightarrow K_i$, $i=1,2$. 
Let $h\colon R_{\Sigma_2}\rightarrow R_{\Sigma_1}$
be a homomorphism of $\Gamma$-graded $R$-rings. Then there exists 
the homomorphism of graded $R$-rings $\Phi_1h\colon R_{\Sigma_2}\rightarrow K_1$. Then
by what has been explained above, $\Phi_2$ factors through $\mathfrak{C}_{K_2}(K_1)$,
and gives the desired specialization.

Now suppose that there 
exist gr-specializations $f\colon \mathfrak{C}_{K_2}(K_1)\rightarrow K_2$ and 
$g\colon \mathfrak{C}_{K_1}(K_2)\rightarrow K_1$. Then the composition $gf$ gives a gr-specialization
from $K_1$ in itself. Thus it has to be the identity. Similarly the composition $fg$
gives a gr-specialization from $K_2$ in itself. Hence, $f$ is an isomorphism in the category $\mathcal{E}_R$
of $\Gamma$-graded epic $R$-division rings. Therefore, $f$ is an isomorphism of graded $R$-rings.
\end{proof}

\begin{corollary}\label{coro:divisionringuniversallocalization}
Let $R$ be a $\Gamma$-graded ring.  Suppose
that there exists $\Omega\subseteq\mathfrak{M}(R)$ such that $(R_\Omega,\lambda)$,
where $\lambda\colon R\rightarrow R_\Omega$ is the canonical homomorphism, is a 
$\Gamma$-graded (epic) $R$-division ring. Then the only gr-specializations to $(R_\Omega,\lambda)$
are isomorphisms of $\Gamma$-graded $R$-rings.
\end{corollary}

\begin{proof}
Suppose there exists a gr-specialization from the $\Gamma$-graded epic $R$-division ring
$(K,\varphi)$ to $(R_\Omega,\lambda)$.  By Theorem~\ref{theo:specialization}(3), then there exists a 
(unique) homomorphism of $\Gamma$-graded $R$-rings $R_\Omega\rightarrow R_\Sigma\rightarrow K$,
where $$\Sigma=\{A\in\mathfrak{M}(R)\colon A^\varphi \textrm{ is invertible over } K\}.$$
Now, since $R_\Omega$ and $K$ are  $\Gamma$-graded epic $R$-division rings, the image of $R_\Omega$
must be $K$ and therefore they are isomorphic as $\Gamma$-graded $R$-rings.  
\end{proof}

\begin{corollary}
Let $R$ be a $\Gamma$-graded ring with a universal $\Gamma$-graded epic $R$-division ring $(U,\rho)$. Suppose that
$\Sigma\subseteq\mathfrak{M}(R)$ is such that there exists a homomorphism of $\Gamma$-graded
rings $R_\Sigma\rightarrow L$ for some $\Gamma$-graded division ring $L$. Then $(U,\rho)$
is a universal $\Gamma$-graded epic $R_\Sigma$-division ring. 
\end{corollary}

\begin{proof}
Consider the canonical homomorphism $\lambda\colon R\rightarrow R_\Sigma$.
Let $f\colon R_\Sigma\rightarrow L$ be a homomorphism of $\Gamma$-graded rings with $L$ a $\Gamma$-graded
division ring. Then $(\DC(f),f\lambda)$ is a $\Gamma$-graded epic $R$-division ring such that the matrices
in $\Sigma$ become invertible. Hence, by Theorem~\ref{theo:specialization},  $\Sigma^\rho$ consists of invertible matrices in $U$. Thus there exists a unique homomorphism of $\Gamma$-graded rings $\psi\colon R_\Sigma\rightarrow U$ and $(U,\psi)$
is a $\Gamma$-graded epic $R_\Sigma$-division ring.

Consider a $\Gamma$-graded epic $R_\Sigma$-division ring $(K,\varphi)$. The composition
$\varphi\lambda\colon R\rightarrow K$ is an epimorphism of $\Gamma$-graded rings, because
$\lambda$ and $\varphi$ are. Hence $(K,\varphi\lambda)$ is a $\Gamma$-graded epic $R$-division ring and
therefore there exists a specialization from $(U,\rho)$ to $(K,\varphi\lambda)$.
\end{proof}

Adapting \cite[p.426]{Cohnfreeeidealringslocalization} to the graded context, we give some examples to illustrate the concepts of universal graded division ring and graded division rings that are universal localizations.

Let $R=\bigoplus\limits_{\gamma\in\Gamma}R_\gamma$ be a commutative $\Gamma$-graded domain. Then the localization of $R$ at the set $\h(R)\setminus\{0\}$ of nonzero homogeneous elements yields
a $\Gamma$-graded epic $R$-field $(F,\varphi)$. We point out that $F=\bigoplus\limits_{\gamma\in\Gamma}F_\gamma$ is a $\Gamma$-graded field with $$F_\gamma=\{ab^{-1}\mid a\in R_\delta, b\in R_\varepsilon, \delta\varepsilon^{-1}=\gamma\}$$ for each $\gamma\in\Gamma$.
Furthermore, if $(K,\psi)$ is a
$\Gamma$-graded epic $R$-division ring, then $\ker \psi$ is a graded prime ideal of $R$. That is,
$\ker\varphi\neq R$ and if $x,y\in \h(R)\cap\ker \psi$ with $xy\in\ker\psi$, then $x\in\ker\psi$ of
$y\in \ker\psi$. Hence $\h(R)\setminus\ker \psi$ is a multiplicative subset of $R$. Then
the localization of $R$ at $\h(R)\setminus\ker \psi$ is a $\Gamma$-graded local subring of $F$
with $\Gamma$-graded residue division ring $R$-isomorphic to $K$. Therefore $(F,\varphi)$
is a $\Gamma$-graded universal $R$-division ring of fractions that is a universal localization.

Let $S=E\times F$ be the direct product of two $\Gamma$-graded fields
$E=\bigoplus\limits_{\gamma\in\Gamma} E_\gamma$ and $F=\bigoplus\limits_{\gamma\in\Gamma} F_\gamma$.
Then $S=\bigoplus\limits_{\gamma\in\Gamma} S_\gamma$ is a $\Gamma$-graded ring
with $S_\gamma=E_\gamma\times F_\gamma$.
Suppose  $(D,\rho)$ is a $\Gamma$-graded epic $S$-division ring. Since
$(1,1)=(1,0)+(0,1)$ and $(1,0)(0,1)=(0,0)$, then either $\rho(1,0)=0$ or $\rho(0,1)=0$.
If $\rho(1,0)=0$, then $\rho(E\times\{0\})=0$ and if $\rho(0,1)=0$, then
$\rho(\{0\}\times F)=0$. Hence $S$ has only two epic $S$-division rings which are $E$ and $F$. Note
that none of them is a universal $\Gamma$-graded epic $R$-division ring. On the other hand,
both are universal localizations. For example $E$ is the universal localization of
$S$ at $\{(1,1)\}\cup\{(a,0)\mid a\in \h(E)\setminus\{0\}\}$.

Let now $E$ be a $\Gamma$-graded field. Then the polynomial ring $E[x]=\bigoplus\limits_{\gamma\in\Gamma}
E[x]_\gamma$ is a $\Gamma$-graded ring with 
$$E[x]_\gamma=E_\gamma[x]=\{a_0+a_1x+\dotsb+a_nx^n\mid a_i\in E_\gamma, n\in\mathbb{N}\}.$$
The ideal $(x^2)$ is a graded ideal of $E[x]$. Hence $T=E[x]/(x^2)$ is a $\Gamma$-graded local ring
with maximal graded ideal $(x)/(x^2)$. Then $E$ is the unique
$\Gamma$-graded epic $T$-division ring, and thus $E$ is a universal $\Gamma$-graded epic $T$-division ring. Notice that $E$ is not a 
universal localization at matrices in $\mathfrak{M}(R)$ because the matrices which
become invertible in $E$ are already invertible in $T$ since $E$ is the $\Gamma$-graded residue division ring
of $T$.

The ring $U=T\times F$ with $T$ as before and $F$ a $\Gamma$-graded field has $E$ and $F$ as 
$\Gamma$-graded epic $U$-division rings, but only $F$ is a universal localization. 


\section{Malcolmson's construction of the universal localization}\label{sec:Malcolmsonscriterion}

\emph{Throughout this section, let $\Gamma$ be a group.}

\medskip

This section is devoted to show that the natural extension of the construction of the ring $R_\Sigma$ given  by
P. Malcolmson in \cite{Malcolmsonscriterion} works in the context of graded rings.
Although technical, this construction will be important for us in the next section. 
In Section~\ref{subsec:Malcolmsonscriterion}, we give the graded version of the main results in \cite{Malcolmsonscriterion}, the so called Malcolmson's criterion and a sufficient condition for the universal localization $R_\Sigma$ not to be the zero ring.

Let $R$ be a $\Gamma$-graded ring and $\lambda\colon R\rightarrow R_\Sigma$ the universal localilzation
at a gr-lower semimultiplicative subset
of $\mathfrak{M}(R)$.
By Lemma~\ref{lem:homogeneousclosure}(6) and Proposition~\ref{prop:basicsonuniversallocalization}(4),
every homogeneous element of $(R_\Sigma)_\gamma$ is of the form
$F^\lambda(A^{\lambda})^{-1}X^\lambda$ where $A\in \Sigma_n[\overline{\alpha}][\overline{\beta}]$,
$F\in M_{1\times n}(R)[\gamma][\overline{\beta}]$, $X\in M_{n\times 1}[\overline{\alpha}][e]$.
For each $\gamma\in\Gamma$,
$(R_\Sigma)_\gamma$ is constructed as a set of equivalent classes of 5-tuples
$(F,A,X,\alpha,\beta)$. The equivalence class $[(F,A,X,\alpha,\beta)]$ of $(F,A,X,\alpha,\beta)$ is interpreted as the element
$F^\lambda(A^{\lambda})^{-1}X^\lambda$ of $R_\Sigma$ and addition and product are defined according to this interpretation. Thus
for $[F',A',X',\alpha',\beta']+[F,A,X,\alpha,\beta]\in (R_\Sigma)_\gamma$
$$[F',A',X',\alpha',\beta']+[F,A,X,\alpha,\beta]=
\left[\begin{pmatrix}F'&F\end{pmatrix},\begin{pmatrix}A'&0\\0&A\end{pmatrix},\begin{pmatrix}X'\\X\end{pmatrix},
\alpha'*\alpha,\beta'*\beta\right],$$
and for  $(F',A',X',\alpha',\beta')\in (R_\Sigma)_{\gamma'}$ and $(F,A,X,\alpha,\beta)\in (R_\Sigma)_\gamma$
$$
[F',A',X',\alpha',\beta']\cdot[F,A,X,\alpha,\beta]= \setlength{\arraycolsep}{1.2pt}
\left[\begin{pmatrix}0&F'\end{pmatrix},\begin{pmatrix}A&0\\-X'F& A'\end{pmatrix},\begin{pmatrix}X\\0\end{pmatrix},
\alpha*\alpha'\gamma, \beta*\beta'\gamma\!\right].
$$
Also $-[(F,A,X,\alpha,\beta)]=
[(-F,A,X,\alpha,\beta)]$ and
$\lambda(r)=[(r,1,1,e,e)]$ for $r\in R_\gamma$.

\medskip 

In this section, for the ease of exposition, we use the following notation. By the expression
\emph{$A$ is a homogeneous matrix}, we mean $A\in \mathfrak{M}_\bullet (R)$. We will
also use the terms \emph{homogeneous row}, \emph{homogeneous column} to emphasize that the matrix in
question is a row or a column, respectively. If 
$A\in M_{m\times n}(R)[\overline{\alpha}][\overline{\beta}]$, but we do not want to make reference to the
size of $A$, we will say \emph{$A$ is a homogenous matrix of distribution $(\alpha,\beta)$}. 
Also, the sequence $\overline{\alpha}\gamma$ will be denoted by $\alpha\gamma$ for each $\overline{\alpha}\in\Gamma^n$ and $\gamma\in \Gamma$.

\subsection{Equivalence relation}\label{subsec:equivalencerelation}

Let  $R$ be a $\Gamma$-graded ring and $\Sigma$ be a gr-lower semimultiplicative subset
of $\mathfrak{M}(R)$.

For $\gamma\in \Gamma$, let $(T_\Sigma)_\gamma$ be the set of $5$-tuples $(F,A,X,\alpha,\beta)$
where $A\in \Sigma$ of distribution $(\alpha,\beta)$, $F$ is a homogeneous row of distribution 
$(\gamma,\beta)$, and $X$ is a homogeneous column of distribution $(\alpha,e)$.

Let $(F,A,X,\alpha,\beta), (G,B,Y,\delta,\varepsilon)\in (T_\Sigma)_\gamma$. We say that
$$(F,A,X,\alpha,\beta)\sim (G,B,Y,\delta,\varepsilon),$$
if and only if there exist $L,M,P,Q\in \Sigma$, homogeneous rows $J,U$ and
homogeneous columns $W,V$ such that
\begin{equation}\label{eq:equivalencerelation}
\left(\begin{array}{cccc|c}
A&0&0&0&X\\0&B&0&0&Y\\0&0&L&0&W\\0&0&0&M&0\\\hline F&-G&0&J&0
\end{array}\right)=
\left(\begin{array}{c}
P\\\hline U
\end{array}\right)
\left(\begin{array}{c|c}
Q&V
\end{array}\right)
\end{equation}
where $P,U,Q,V$ have distributions $(\pi,\omega),(\gamma,\omega),(\omega, \theta),(\omega,e)$, 
respectively, and, 
if we think of $$\pi=\pi_1*\pi_2*\pi_3*\pi_4\ \textrm{ and } \
\theta=\theta_1*\theta_2*\theta_3*\theta_4$$ 
then $\pi_1=\alpha$,
 $\pi_2=\delta$, $\theta_1=\beta$, $\theta_2=\varepsilon$.

The right hand side of \eqref{eq:equivalencerelation} will also be denoted by
\[
\left(\begin{array}{cccc}
P_{11}&P_{12}&P_{13}&P_{14}\\P_{21}&P_{22}&P_{23}&P_{24}\\P_{31}&P_{32}&P_{33}&P_{34}\\P_{41}&P_{42}&P_{43}&P_{44}\\\hline U_{1}&U_{2}&U_{3}&U_{4}
\end{array}\right)
\left(\begin{array}{cccc|c}
Q_{11}&Q_{12}&Q_{13}&Q_{14}&V_{1}\\Q_{21}&Q_{22}&Q_{23}&Q_{24}&V_{2}\\Q_{31}&Q_{32}&Q_{33}&Q_{34}&V_{3}\\Q_{41}&Q_{42}&Q_{43}&Q_{44}&V_{4}
\end{array}\right).
\]

\begin{lemma}
Let $(F,A,X,\alpha,\beta),(G,B,Y,\delta,\varepsilon)\in (T_\Sigma)_\gamma$ such that there is
a factorization as
a product of homogeneous matrices of any of these forms with $L,M,P,Q\in \Sigma$ and with the corresponding
distributions
\begin{enumerate}[\rm(1)]
\item
$
\left(\begin{array}{cc|c}
A&0&X\\0&B&Y\\\hline F&-G&0
\end{array}\right)=
\left(\begin{array}{c}
P\\\hline U
\end{array}\right)
\left(\begin{array}{c|c}
Q&V
\end{array}\right)
$

\item
$
\left(\begin{array}{ccc|c}
A&0&0&X\\0&B&0&Y\\0&0&M&W\\\hline F&-G&0&0
\end{array}\right)=
\left(\begin{array}{c}
P\\\hline U
\end{array}\right)
\left(\begin{array}{c|c}
Q&V
\end{array}\right)
$

\item
$
\left(\begin{array}{ccc|c}
A&0&0&X\\0&B&0&Y\\0&0&L&0\\\hline F&-G&J&0
\end{array}\right)=
\left(\begin{array}{c}
P\\\hline U
\end{array}\right)
\left(\begin{array}{c|c}
Q&V
\end{array}\right)
$

\end{enumerate}
\noindent
then $(F,A,X,\alpha,\beta)\sim (G,B,Y,\delta,\varepsilon)$. 
\end{lemma}

\begin{proof}
(1) Suppose 
$P,U,Q,V$ have distributions $(\pi,\omega)$, $(\gamma,\omega)$, $(\omega,\theta)$, $(\omega,e)$ where 
$\pi_1=\alpha$, $\pi_2=\delta$,  $\theta_1=\beta$ and $\theta_2=\varepsilon$, and that we have the factorization
\[
\left(\begin{array}{cc|c}
A&0&X\\0&B&Y\\\hline F&-G&0
\end{array}\right)=
\left(\begin{array}{cc}
P_{11}&P_{12}\\P_{21}&P_{22}\\\hline U_1&U_2
\end{array}\right)
\left(\begin{array}{cc|c}
Q_{11}&Q_{12}&V_1\\Q_{21}&Q_{22}&V_2
\end{array}\right)
\]
Thus, we have the factorization
\[
\left(\begin{array}{cccc|c}
A&0&0&0&X\\0&B&0&0&Y\\0&0&1&0&1\\0&0&0&1&0\\\hline F&-G&0&1&0
\end{array}\right) = 
\left(\begin{array}{cccc}
P_{11}&P_{12}&0&0\\P_{21}&P_{22}&0&0\\0&0&1&0\\0&0&0&1\\\hline U_{1}&U_{2}&0&1
\end{array}\right)
\left(\begin{array}{cccc|c}
Q_{11}&Q_{12}&0&0&V_{1}\\Q_{21}&Q_{22}&0&0&V_{2}\\0&0&1&0&1\\0&0&0&1&0
\end{array}\right)
\]
where the factors of the right hand side have distributions
\[
\begin{array}{cccc}
(\alpha*\delta*e*\gamma*\gamma,&\omega_1*\omega_2*e*\gamma),&(\omega_1*\omega_2*e*\gamma,&\beta*\varepsilon*e*\gamma*e)
\end{array}
\]

\noindent (2) 
Suppose $P,U,Q,V$ have distributions $(\pi,\omega)$, $(\gamma,\omega)$, $(\omega,\theta)$, $(\omega,e)$ where
$\pi_1=\alpha$, $\pi_2=\delta$, $\theta_1=\beta$ and $\theta_2=\varepsilon$ and we have the factorization:
\[
\left(\begin{array}{ccc|c}
A&0&0&X\\0&B&0&Y\\0&0&L&W\\\hline F&-G&0&0
\end{array}\right)=
\left(\begin{array}{ccc}
P_{11}&P_{12}&P_{13}\\P_{21}&P_{22}&P_{23}\\P_{31}&P_{32}&P_{33}\\\hline U_1&U_2&U_3
\end{array}\right)
\left(\begin{array}{ccc|c}
Q_{11}&Q_{12}&Q_{13}&V_1\\Q_{21}&Q_{22}&Q_{23}&V_2\\Q_{31}&Q_{32}&Q_{33}&V_3
\end{array}\right)
\]
Thus, we have the following equality
\[\setlength{\arraycolsep}{4pt}
\left(\begin{array}{cccc|c}
A&0&0&0&X\\0&B&0&0&Y\\0&0&L&0&W\\0&0&0&1&0\\\hline F&-G&0&1&0
\end{array}\right)
=
\left(\begin{array}{cccc}
P_{11}&P_{12}&P_{13}&0\\P_{21}&P_{22}&P_{23}&0\\P_{31}&P_{32}&P_{33}&0\\0&0&0&1\\\hline U_{1}&U_{2}&U_{3}&1
\end{array}\right)
\left(\begin{array}{cccc|c}
Q_{11}&Q_{12}&Q_{13}&0&V_{1}\\Q_{21}&Q_{22}&Q_{23}&0&V_{2}\\Q_{31}&Q_{32}&Q_{33}&0&V_{3}\\0&0&0&1&0
\end{array}\right)
\]
where the factors of the right hand side have distributions
\[
\begin{array}{cccc}
(\alpha*\delta*\pi_3*\gamma*\gamma,&\omega_1*\omega_2*\omega_3*\gamma),&(\omega_1*\omega_2*\omega_3*\gamma,&\beta*\varepsilon*\theta_3*\gamma*e)
\end{array}
\]

\noindent(3) 
Suppose $P,U,Q,V$ have distributions $(\pi,\omega)$, $(\gamma,\omega)$, $(\omega,\theta)$, $(\omega,e)$ where
$\pi_1=\alpha$, $\pi_2=\delta$, $\theta_1=\beta$ and $\theta_2=\varepsilon$ and that we have the factorization
\[
\left(\begin{array}{ccc|c}
A&0&0&X\\0&B&0&Y\\0&0&M&0\\\hline F&-G&J&0
\end{array}\right)=
\left(\begin{array}{ccc}
P_{11}&P_{12}&P_{13}\\P_{21}&P_{22}&P_{23}\\P_{31}&P_{32}&P_{33}\\\hline U_1&U_2&U_3
\end{array}\right)
\left(\begin{array}{ccc|c}
Q_{11}&Q_{12}&Q_{13}&V_1\\Q_{21}&Q_{22}&Q_{23}&V_2\\Q_{31}&Q_{32}&Q_{33}&V_3
\end{array}\right)
\]
Thus, we have the factorization
\[
\setlength{\arraycolsep}{2.5pt}\left(\begin{array}{cccc|c}
A&0&0&0&X\\0&B&0&0&Y\\0&0&M&0&0\\0&0&0&M&0\\\hline F&-G&0&J&0
\end{array}\right)
=
\left(\begin{array}{cccc}
P_{11}&P_{12}&P_{13}&0\\P_{21}&P_{22}&P_{23}&0\\P_{31}&P_{32}&P_{33}&0\\P_{31}&P_{32}&P_{33}&M\\\hline U_{1}&U_{2}&U_{3}&J
\end{array}\right)
\left(\begin{array}{cccc|c}
Q_{11}&Q_{12}&Q_{13}&0&V_{1}\\Q_{21}&Q_{22}&Q_{23}&0&V_{2}\\Q_{31}&Q_{32}&Q_{33}&0&V_{3}
\\0&0&-1&1&0
\end{array}\right)
\]
where the factors of the right hand side have distribution
\[
\begin{array}{cccc}
(\alpha*\delta*\pi_3*\pi_3*\gamma,&\omega_1*\omega_2*\omega_3*\theta_3),&(\omega_1*\omega_2*\omega_3*\theta_3,&\beta*\varepsilon*\theta_3*\theta_3*e)
\end{array},
\]
respectively.
\end{proof}

\begin{lemma}
For each $\gamma\in\Gamma$, the relation $\sim$ defined in $(T_\Sigma)_\gamma$ is an equivalence relation.
\end{lemma}

\begin{proof}
Let $(F,A,X,\alpha,\beta), (G,B,Y,\delta,\varepsilon), (H,C,Z,\zeta,\eta)\in (T_\Sigma)_\gamma$.

The relation $\sim$ is reflexive. Indeed, we have the 
factorization
\[
\left(\begin{array}{cc|c}
A&0&X\\0&A&X\\\hline F&-F&0
\end{array}\right)=
\left(\begin{array}{cc}
I&0\\I&-A\\\hline 0&F
\end{array}\right)
\left(\begin{array}{cc|c}
A&0&X\\I&-I&0
\end{array}\right)
\]
where the factors are homogeneous matrices that have distributions $(\alpha*\alpha*\gamma\,,\,\alpha*\beta)$ and
$(\alpha*\beta\,,\,\beta*\beta*e)$, respectively. This shows that $(F,A,X,\alpha,\beta)\sim(F,A,X,\alpha,\beta)$.

The relation $\sim$ is symmetric. Indeed, suppose that $(F,A,X,\alpha,\beta)\sim (G,B,Y,\delta,\varepsilon)$. 
There exist $L,M,P,Q\in \Sigma$, homogeneous rows $J,U$ and
homogeneous columns $W,V$ such that
\begin{equation}\label{eq:equivalencerelation2}
\setlength{\arraycolsep}{1.5pt}\left(\begin{array}{cccc|c}
A&0&0&0&X\\0&B&0&0&Y\\0&0&L&0&W\\0&0&0&M&0\\\hline F&-G&0&J&0
\end{array}\right)=
\left(\begin{array}{cccc}
P_{11}&P_{12}&P_{13}&P_{14}\\P_{21}&P_{22}&P_{23}&P_{24}\\P_{31}&P_{32}&P_{33}&P_{34}\\P_{41}&P_{42}&P_{43}&P_{44}\\\hline U_{1}&U_{2}&U_{3}&U_{4}
\end{array}\right)
\left(\begin{array}{cccc|c}
Q_{11}&Q_{12}&Q_{13}&Q_{14}&V_{1}\\Q_{21}&Q_{22}&Q_{23}&Q_{24}&V_{2}\\Q_{31}&Q_{32}&Q_{33}&Q_{34}&V_{3}\\Q_{41}&Q_{42}&Q_{43}&Q_{44}&V_{4}
\end{array}\right),
\end{equation}
where $P,U,Q,V$ have distributions $(\pi,\omega),(\gamma,\omega),(\omega, \theta),(\omega,e)$, 
respectively, and \linebreak $\pi_1=\alpha$,
 $\pi_2=\delta$, $\theta_1=\beta$, $\theta_2=\varepsilon$.
Then we have the factorization
$$\scalebox{0.85}{$
\setlength{\arraycolsep}{2pt}
\left(\begin{array}{cccccc|c}
B&0&0&0&0&0&Y\\0&A&0&0&0&0&X\\0&0&B&0&0&0&0\\0&0&0&L&0&0&W\\0&0&0&0&M&0&0\\0&0&0&0&0&B&0\\\hline G&-F&0&0&-J&G&0
\end{array}\right)=  
\left(\begin{array}{cccccc}
I&0&0&0&0&0\\0&P_{11}&P_{12}&P_{13}&P_{14}&0\\-I&P_{21}&P_{22}&P_{23}&P_{24}&0\\0&P_{31}&P_{32}&P_{33}&P_{34}&0\\0&P_{41}&P_{42}&P_{43}&P_{44}&0\\-I&P_{21}&P_{22}&P_{23}&P_{24}&B\\\hline 0&-U_1&-U_2&-U_3&-U_4&G
\end{array}\right)
\left(\begin{array}{cccccc|c}
B&0&0&0&0&0&Y\\Q_{12}&Q_{11}&Q_{12}&Q_{13}&Q_{14}&0&V_1\\Q_{22}&Q_{21}&Q_{22}&Q_{23}&Q_{24}&0&V_2\\
Q_{32}&Q_{31}&Q_{32}&Q_{33}&Q_{34}&0&V_3\\Q_{42}&Q_{41}&Q_{42}&Q_{43}&Q_{44}&0&V_4\\0&0&-I&0&0&I&0
\end{array}\right),
$}$$
where the factors  have distributions
$(\delta*\alpha*\delta*\pi_3*\pi_4*\delta*\gamma\,,\, \delta*\omega_1*\omega_2*\omega_3*\omega_4*\varepsilon)$ and
$(\delta*\omega_1*\omega_2*\omega_3*\omega_4*\varepsilon\ , \ \varepsilon*\beta*\varepsilon*\theta_3*
\theta_4*\varepsilon*e)$ respectively. Hence, $(G,B,Y,\delta,\varepsilon)
\sim(F,A,X,\alpha,\beta)$, and the symmetric property  of the relation $\sim$ is proved.

Now we proceed to prove that $\sim$ satisfies the transitive property. 
Suppose that $(F,A,X,\alpha,\beta)\sim (G,B,Y,\delta,\varepsilon)$ and
$(G,B,Y,\delta,\varepsilon)\sim(H,C,Z,\zeta,\eta)$. Hence, 
there exist $L,M,P,Q\in \Sigma$, homogeneous rows $J,U$ and
homogeneous columns $W,V$ as in \eqref{eq:equivalencerelation2},
and there exist $L',M',P',Q'\in \Sigma$, homogeneous rows $J',U'$ and
homogeneous columns $W',V'$ such that
\[
\setlength{\arraycolsep}{1.5pt}\left(\begin{array}{cccc|c}
B&0&0&0&Y\\0&C&0&0&Z\\0&0&L'&0&W'\\0&0&0&M'&0\\\hline G&-H&0&J'&0
\end{array}\right)=
\left(\begin{array}{cccc}
P_{11}'&P_{12}'&P_{13}'&P_{14}'\\P_{21}'&P_{22}'&P_{23}'&P_{24}'\\P_{31}'&P_{32}'&P_{33}'&P_{34}'\\P_{41}'&P_{42}'&P_{43}'&P_{44}'\\\hline U_{1}'&U_{2}'&U_{3}'&U_{4}'
\end{array}\right)
\left(\begin{array}{cccc|c}
Q_{11}'&Q_{12}'&Q_{13}'&Q_{14}'&V_{1}'\\Q_{21}'&Q_{22}'&Q_{23}'&Q_{24}'&V_{2}'\\Q_{31}'&Q_{32}'&Q_{33}'&Q_{34}'&V_{3}'\\Q_{41}'&Q_{42}'&Q_{43}'&Q_{44}'&V_{4}'
\end{array}\right), \]
where $P',U',Q',V'$ have distributions $(\pi',\omega'),(\gamma,\omega'),(\omega', \theta'),(\omega',e)$, 
respectively, and $\pi_1'=\delta$,
 $\pi_2'=\zeta$, $\theta_1'=\varepsilon$, $\theta_2'=\eta$. Then we have the factorization of the matrix
\[
\left(\begin{array}{ccccccccccc|c}
C & 0 & 0 & 0 & 0 & 0 & 0 & 0 & 0 & 0 & 0 & Z \\
0 & A & 0 & 0 & 0 & 0 & 0 & 0 & 0 & 0 & 0 & X \\ 
0 & 0 & B & 0 & 0 & 0 & 0 & 0 & 0 & 0 & 0 & Y \\
0 & 0 & 0 & L & 0 & 0 & 0 & 0 & 0 & 0 & 0 & W \\
0 & 0 & 0 & 0 & M & 0 & 0 & 0 & 0 & 0 & 0 & 0 \\
0 & 0 & 0 & 0 & 0 & L' & 0 & 0 & 0 & 0 & 0 & W' \\
0 & 0 & 0 & 0 & 0 & 0 & B & 0 & 0 & 0 & 0 & 0 \\
0 & 0 & 0 & 0 & 0 & 0 & 0 & C & 0 & 0 & 0 & 0 \\
0 & 0 & 0 & 0 & 0 & 0 & 0 & 0 & L' & 0 & 0 & 0 \\
0 & 0 & 0 & 0 & 0 & 0 & 0 & 0 & 0 & M' & 0 & 0 \\
0 & 0 & 0 & 0 & 0 & 0 & 0 & 0 & 0 & 0 & M & 0 \\ \hline 
H &-F & 0 & 0 & 0 & 0 &-G & H & 0 &-J' &J & 0  
\end{array}\right)
\]
as a product of the matrices
\[
\scalebox{0.78}{$\setlength{\arraycolsep}{0.4pt}
\left(\begin{array}{ccccccccccc}
I & 0 & 0 & 0 & 0 & 0 & 0 & 0 & 0 & 0 & 0 \\
0 & P_{11} & P_{12} & P_{13} & P_{14} & 0 & 0 & 0 & 0 & 0 & 0 \\
0 & P_{21} & P_{22} & P_{23} & P_{24} & 0 & 0 & 0 & 0 & 0 & 0 \\
0 & P_{31} & P_{32} & P_{33} & P_{34} & 0 & 0 & 0 & 0 & 0 & 0 \\
0 & P_{41} & P_{42} & P_{43} & P_{44} & 0 & 0 & 0 & 0 & 0 & 0  \\
0 & 0 & 0 & 0 & 0 &                     I & 0 & 0 & 0 & 0 & 0 \\
0 & -P_{21} & -P_{22} & -P_{23} & -P_{24} & 0 & P_{11}' & P_{12}' & P_{13}' & P_{14}' & 0 \\
-I & 0 & 0 & 0 & 0 & 0                        & P_{21}' & P_{22}' & P_{23}' & P_{24}' & 0 \\
0 & 0 & 0 & 0 & 0 &                        -I & P_{31}' & P_{32}' & P_{33}' & P_{34}' & 0 \\
0 & 0 & 0 & 0 & 0 & 0                         & P_{41}' & P_{42}' & P_{43}' & P_{44}' & 0 \\
0 & -P_{41} & -P_{42} & -P_{43} & -P_{44} & 0 & 0 & 0 & 0 & 0 & M \\ \hline
0 & -U_1 & -U_2 & -U_3 & -U_4 & 0 & -U_1' & -U_2' & -U_3' & -U_4' & J 
\end{array}\right)\!\!
\left(\begin{array}{ccccccccccc|c}
C & 0 & 0 & 0 & 0 & 0 & 0 & 0 & 0 & 0 & 0 & Z \\
0 & Q_{11} & Q_{12} & Q_{13} & Q_{14} & 0 & 0 & 0 & 0 & 0 & 0 & V_1 \\
0 & Q_{21} & Q_{22} & Q_{23} & Q_{24} & 0 & 0 & 0 & 0 & 0 & 0 & V_2 \\
0 & Q_{31} & Q_{32} & Q_{33} & Q_{34} & 0 & 0 & 0 & 0 & 0 & 0 & V_3 \\
0 & Q_{41} & Q_{42} & Q_{43} & Q_{44} & 0 & 0 & 0 & 0 & 0 &  0  &V_4  \\
0 & 0 & 0 & 0 & 0 &                     L' & 0 & 0 & 0 & 0 &   0 & W' \\
Q_{12}' & 0 & Q_{11}' & 0 & 0 & Q_{13}' & Q_{11}' & Q_{12}' & Q_{13}' & Q_{14}' & 0 & V_1'\\
Q_{22}' & 0 & Q_{21}' & 0 & 0 & Q_{23}' & Q_{21}' & Q_{22}' & Q_{23}' & Q_{24}' & 0 & V_2'\\
Q_{32}' & 0 & Q_{31}' & 0 & 0 & Q_{33}' & Q_{31}' & Q_{32}' & Q_{33}' & Q_{34}' & 0 & V_3'\\
Q_{42}' & 0 & Q_{41}' & 0 & 0 & Q_{43}' & Q_{41}' & Q_{42}' & Q_{43}' & Q_{44}' & 0  & V_4'\\
0 & 0 & 0 & 0 & I & 0 & 0 & 0 & 0 & 0 & I & 0 \\ 
 \end{array}\right)$},
\]
where the factors have distributions
$$(\zeta*\alpha*\delta*\pi_3*\pi_4*\pi_3'*\delta*\zeta*\pi_3'*\pi_4'*\pi_4*\gamma \ , \ 
\zeta*\omega_1*\omega_2*\omega_3*\omega_4*\pi_3'*\omega_1'*\omega_2'*\omega_3'*\omega_4'*\theta_4),$$
$$(\zeta*\omega_1*\omega_2*\omega_3*\omega_4*\pi_3'*\omega_1'*\omega_2'*\omega_3'*\omega_4'*\theta_4 \ , \  
\eta*\beta*\varepsilon*\theta_3*\theta_4*\theta_3'*\varepsilon*\eta*\theta_3'*\theta_4'*\theta_4*e),$$
respectively. This factorization implies that $(F,A,X,\alpha,\beta)\sim (H,C,Z,\zeta,\eta)$, as desired. 
\end{proof}


\subsection{Operations}\label{subsec:operations}

Let $\gamma\in \Gamma$. 

If $(F',A',X',\alpha',\beta'), (F,A,X,\alpha,\beta)\in (T_\Sigma)_\gamma$, then we define
$$(F',A',X',\alpha',\beta') + (F,A,X,\alpha,\beta)=
\left(\begin{pmatrix}F'&F\end{pmatrix},\begin{pmatrix}A'&0\\0&A\end{pmatrix},\begin{pmatrix}X'\\X\end{pmatrix},
\alpha'*\alpha,\beta'*\beta\right).$$ Note that it belongs to $(T_\Sigma)_\gamma$.

If  $(F',A',X',\alpha',\beta')\in (T_\Sigma)_{\gamma'}$ and $(F,A,X,\alpha,\beta)\in (T_\Sigma)_\gamma$, then we define
\[
(F',A',X',\alpha',\beta')\cdot(F,A,X,\alpha,\beta)= \setlength{\arraycolsep}{1.2pt}
\left(\begin{pmatrix}0&F'\end{pmatrix},\begin{pmatrix}A&0\\-X'F& A'\end{pmatrix},\begin{pmatrix}X\\0\end{pmatrix},
\alpha*\alpha'\gamma, \beta*\beta'\gamma\!\right).
\]
Note that this element belongs to $(T_\Sigma)_{\gamma'\gamma}$ because
the homogeneous matrix $(0\ F')$ has distribution $(\gamma'\gamma,\beta*\beta'\gamma)$ and 
the homogeneous matrix $\left(\begin{smallmatrix}
X \\ 0	
\end{smallmatrix}\right)$ has distribution $(\alpha*\alpha'\gamma,e)$.

If $(F,A,X,\alpha,\beta)\in (T_\Sigma)_\gamma$, we define 
$$-(F,A,X,\alpha,\beta)=(-F,A,X,\alpha,\beta)\in (T_\Sigma)_\gamma.$$

Finally, if $r\in R_\gamma$, we define $$\mu(r)=(r,1,1,e,e)\in (T_\Sigma)_\gamma.$$

Now we prove a series of lemmas that show the compatibility of the operations just defined and the equivalence relation $\sim$.

\begin{lemma}\label{lem:easywelldefined}
The following assertions hold true.
\begin{enumerate}[\rm(1)]
\item If $(F',A,X,\alpha,\beta),(F,A,X,\alpha,\beta)\in(T_\Sigma)_\gamma$, then $$(F',A,X,\alpha,\beta)+(F,A,X,\alpha,\beta)\sim(F'+F,A,X,\alpha,\beta).$$
\item If $(F,A,X',\alpha,\beta),(F,A,X,\alpha,\beta)\in(T_\Sigma)_\gamma$, then $$(F,A,X',\alpha,\beta)+(F,A,X,\alpha,\beta)\sim
(F,A,X'+X,\alpha,\beta).$$
\item If $r\in R_{\gamma'}$ and $(F,A,X,\alpha,\beta)\in(T_\Sigma)_\gamma$, then 
$$\mu(r)\cdot(F,A,X,\alpha,\beta)\sim(rF,A,X,\alpha,\beta)\in (T_\Sigma)_{\gamma'\gamma}.$$
\item If $(F',A',X',\alpha',\beta')\in(T_\Sigma)_{\gamma'}$ and $r\in R_\gamma$, then $$(F',A',X',\alpha',\beta')\cdot\mu(r)\sim(F',A',X'r,\alpha'\gamma,\beta'\gamma)\in (T_\Sigma)_{\gamma'\gamma}.$$
\end{enumerate}

\end{lemma}

\begin{proof}
(1) It follows from the following factorization
\[
\left(\begin{array}{ccc|c}
A&0&0&X\\0&A&0&X\\0&0&A&X\\\hline F'&F&-F'-F&0
\end{array}\right)=
\left(\begin{array}{ccc}
I&0&0\\I&A&0\\I&0&A\\\hline 0&F&-F'-F
\end{array}\right)
\left(\begin{array}{ccc|c}
A&0&0&X\\-I&I&0&0\\-I&0&I&0
\end{array}\right)
\]
where the factors of the right hand side have distributions
$(\alpha*\alpha*\alpha*\gamma\, , \, \alpha*\beta*\beta)$ and
$(\alpha*\beta*\beta\, , \, \beta*\beta*\beta*e)$, respectively.

(2) It follows from the equality
\[
\left(\begin{array}{ccc|c}
A&0&0&X'\\0&A&0&X\\0&0&A&X'+X\\\hline F&F&-F&0
\end{array}\right)=
\left(\begin{array}{ccc}
I&0&0\\0&I&0\\I&I&A\\\hline 0&0&-F
\end{array}\right)
\left(\begin{array}{ccc|c}
A&0&0&X'\\0&A&0&X\\-I&-I&I&0
\end{array}\right)
\]
where the factors of the right hand side have distributions
$(\alpha*\alpha*\alpha*\gamma\, , \, \alpha*\alpha*\beta)$ and
$(\alpha*\alpha*\beta\, , \, \beta*\beta*\beta*e)$, respectively.

(3) It follows from the factorization
\[
\left(\begin{array}{ccc|c}
A&0&0&X\\-F&1&0&0\\0&0&A&X\\\hline0&r&-rF&0
\end{array}\right)=
\left(\begin{array}{ccc}
I&0&0\\0&1&0\\I&0&A\\\hline 0&r&-rF
\end{array}\right)
\left(\begin{array}{ccc|c}
A&0&0&X\\-F&1&0&0\\-I&0&I&0
\end{array}\right)
\]
where the factors of the right hand side have distributions $(\alpha*\gamma*\alpha*\gamma'\gamma \, , \,
\alpha*\gamma*\beta)$ and $(\alpha*\gamma*\beta \, , \, \beta*\gamma*\beta*e)$, respectively.

(4) It follows from the factorization
\[
\left(\begin{array}{ccc|c}
1&0&0&1\\-X'r&A'&0&0\\0&0&A'&X'r\\\hline0&F'&-F'&0
\end{array}\right)=
\left(\begin{array}{ccc}
1&0&0\\0&I&0\\X'r&I&A'\\\hline 0&0&-F'
\end{array}\right)
\left(\begin{array}{ccc|c}
1&0&0&1\\-X'r&A'&0&0\\0&-I&I&0
\end{array}\right)
\]
where the factors of the right hand side have distributions $(e*\alpha'\gamma*\alpha'\gamma*\gamma'\gamma \, , \,
e*\alpha'\gamma*\beta'\gamma)$ and $(e*\alpha'\gamma*\beta'\gamma \, , \, e*\beta'\gamma*\beta'\gamma*e)$, respectively.
\end{proof}

\begin{lemma}\label{lem:welldefined}
The relation $\sim$ is compatible with the operations defined on the $(T_\Sigma)_\gamma$'s. More precisely, the following
assertions hold true.
\begin{enumerate}[\rm(1)]
\item For $x',x\in(T_\Sigma)_\gamma$, then $x+x'\sim x'+x$.
\item For $x',x,y\in(T_\Sigma)_\gamma$ such that $x\sim y$, then $x'+x\sim x'+y$ and 
$x+x'\sim y+x'$.
\item For $x,y\in(T_\Sigma)_\gamma$ and $x'\in(T_\Sigma)_{\gamma'}$ such that $x\sim y$, then 
$x'x\sim x'y$ and $xx'\sim yx'$.
\item For $x,y\in(T_\Sigma)_\gamma$ such that $x\sim y$, then $-x\sim -y$.
\end{enumerate}
\end{lemma}

\begin{proof}
(1) Let $(F',A',X',\alpha',\beta'),(F,A,X,\alpha,\beta)\in (T_\Sigma)_\gamma$.
The equality 
\[
\setlength{\arraycolsep}{4pt}\left(\begin{array}{cccc|c}
A'&0&0&0&X'\\0&A&0&0&X\\0&0&A&0&X\\0&0&0&A'&X'\\\hline F'&F&-F&-F'&0
\end{array}\right)=
\left(\begin{array}{cccc}
I&0&0&0\\0&I&0&0\\0&I&A&0\\I&0&0&A'\\\hline 0&0&-F&-F'
\end{array}\right)
\left(\begin{array}{cccc|c}
A'&0&0&0&X'\\0&A&0&0&X\\0&-I&I&0&0\\-I&0&0&I&0
\end{array}\right)
\]
where the factors of the right hand side have distributions
$(\alpha'*\alpha*\alpha*\alpha'*\gamma\,,\,\alpha'*\alpha*\beta*\beta')$ and
$(\alpha'*\alpha*\beta*\beta'\,,\,\beta'*\beta*\beta*\beta'*e)$, respectively,
shows (1).

(2) First note that, by (1), it is enough to prove that $x'+ x\sim x'+y.$ Now
let $\gamma\in\Gamma$, let $(F',A',X',\alpha',\beta')\in (T_\Sigma)_{\gamma}$ and
let 
$(F,A,X,\alpha,\beta),(G,B,Y,\delta,\varepsilon)\in(T_\Sigma)_\gamma$ be such that
$(F,A,X,\alpha,\beta)\sim(G,B,Y,\delta,\varepsilon)$. Thus, there exist $L,M,P,Q\in \Sigma$,
homogenous rows $J,U$, and homogeneous columns $W,V$ as in \eqref{eq:equivalencerelation}. The result
follows because the matrix
\[
\left(\begin{array}{ccccccccc|c}
A'&0&0&0&0&0&0&0&0&X'\\0&A&0&0&0&0&0&0&0&X\\0&0&A'&0&0&0&0&0&0&X'\\0&0&0&B&0&0&0&0&0&Y\\0&0&0&0&L&0&0&0&0&W\\
0&0&0&0&0&A&0&0&0&0\\0&0&0&0&0&0&B&0&0&0\\0&0&0&0&0&0&0&L&0&0\\0&0&0&0&0&0&0&0&M&0\\\hline F'&F&-F'&-G&0&F&-G&0&J&0
\end{array}\right)
\]
can be expressed as the product of the homogeneous matrices
\[
\setlength{\arraycolsep}{1.5pt}\left(\begin{array}{ccccccccccc}
I&0&0&0&0&0&0&0&0\\
0&I&0&0&0&0&0&0&0\\
I&0&A'&0&0&0&0&0&0\\
0&0&0&I&0&0&0&0&0\\
0&0&0&0&I&0&0&0&0\\
0&-I&0&0&0&P_{11} &P_{12}&P_{13}&P_{14}\\
0&0&0&-I&0&P_{21}&P_{22}&P_{23}&P_{24}\\
0&0&0&0&-I&P_{31}&P_{32}&P_{33}&P_{34}\\
0&0&0&0&0&P_{41}&P_{42}&P_{43}&P_{44}\\\hline 
0&0&-F'&0&0&U_1&U_2&U_3&U_4
\end{array}\right)
\left(\begin{array}{ccccccccc|c}
A'&0&0&0&0&0&0&0&0&X'\\
0&A&0&0&0&0&0&0&0&X\\
-I&0&I&0&0&0&0&0&0&0\\
0&0&0&B&0&0&0&0&0&Y\\
0&0&0&0&L&0&0&0&0&W\\
0&Q_{11}&0&Q_{12}&Q_{13}&Q_{11}&Q_{12}&Q_{13}&Q_{14}&V_1\\
0&Q_{21}&0&Q_{22}&Q_{23}&Q_{21}&Q_{22}&Q_{23}&Q_{24}&V_2\\
0&Q_{31}&0&Q_{32}&Q_{33}&Q_{31}&Q_{32}&Q_{33}&Q_{34}&V_3\\
0&Q_{41}&0&Q_{42}&Q_{43}&Q_{41}&Q_{42}&Q_{43}&Q_{44}&V_4\\
\end{array}\right)
\]
that have distributions
$(\alpha'*\alpha*\alpha'*\delta*\pi_3*\alpha*\delta*\pi_3*\pi_4*\gamma\,,\,
\alpha'*\alpha*\beta'*\delta*\pi_3*\omega_1*\omega_2*\omega_3*\omega_4)$ and
$(\alpha'*\alpha*\beta'*\delta*\pi_3*\omega_1*\omega_2*\omega_3*\omega_4\,,\,
\beta'*\beta*\beta'*\varepsilon*\theta_3*\beta*\varepsilon*\theta_3*\theta_4*e)$, respectively.

(3) Let $\gamma,\gamma'\in\Gamma$, let $(F',A',X',\alpha',\beta')\in (T_\Sigma)_{\gamma'}$ and
let 
$(F,A,X,\alpha,\beta),\linebreak
(G,B,Y,\delta,\varepsilon)\in(T_\Sigma)_\gamma$ be such that
$(F,A,X,\alpha,\beta)\sim(G,B,Y,\delta,\varepsilon)$. Thus, there exist $L,M,P,Q\in \Sigma$,
homogenous rows $J,U$, and homogeneous columns $W,V$ as in \eqref{eq:equivalencerelation}. 

We prove first that 
$$ (F',A',X',\alpha',\beta')\cdot(F,A,X,\alpha,\beta)\sim 
(F',A',X',\alpha',\beta')\cdot(G,B,Y,\delta,\varepsilon).$$
It follows because the following homogeneous matrix
\[
\left(\begin{array}{cccccccccc|c}
A&0&0&0&0&0&0&0&0&0&X\\
-X'F&A'&0&0&0&0&0&0&0&0&0\\
0&0&B&0&0&0&0&0&0&0&Y\\
0&0&-X'G&A'&0&0&0&0&0&0&0\\
0&0&0&0&L&0&0&0&0&0&W\\
0&0&0&0&0&A&0&0&0&0&0\\
0&0&0&0&0&0&B&0&0&0&0\\
0&0&0&0&0&0&0&L&0&0&0\\
0&0&0&0&0&0&0&0&M&0&0\\
0&0&0&0&0&-X'F&X'G&0&-XJ&A'&0\\\hline 
0&F'&0&-F'&0&0&0&0&0&F'&0
\end{array}\right)
\]
has the factorization, as a product of homogeneous matrices,
\[
\scalebox{0.8}{$\setlength{\arraycolsep}{1pt}
\left(\begin{array}{cccccccccc}
I&0&0&0&0&0&0&0&0&0\\
0&I&0&0&0&0&0&0&0&0\\
0&0&I&0&0&0&0&0&0&0\\
0&0&0&I&0&0&0&0&0&0\\
0&0&0&0&I&0&0&0&0&0\\
-I&0&0&0&0&P_{11}&P_{12}&P_{13}&P_{14}&0\\
0&0&-I&0&0&P_{21}&P_{22}&P_{23}&P_{24}&0\\
0&0&0&0&-I&P_{31}&P_{32}&P_{33}&P_{34}&0\\
0&0&0&0&0&P_{41}&P_{42}&P_{43}&P_{44}&0\\
0&-I&0&I&0&-X'U_1&-X'U_2&-X'U_3&-X'U_4&A'\\\hline 
0&0&0&0&0&0&0&0&0&F'
\end{array}\right)
\left(\begin{array}{cccccccccc|c}
A&0&0&0&0&0&0&0&0&0&X\\
-X'F&A'&0&0&0&0&0&0&0&0&0\\
0&0&B&0&0&0&0&0&0&0&Y\\
0&0&-X'G&A'&0&0&0&0&0&0&0\\
0&0&0&0&L&0&0&0&0&0&W\\
Q_{11}&0&Q_{12}&0&Q_{13}&Q_{11}&Q_{12}&Q_{13}&Q_{14}&0&V_1\\
Q_{21}&0&Q_{22}&0&Q_{23}&Q_{21}&Q_{22}&Q_{23}&Q_{24}&0&V_2\\
Q_{31}&0&Q_{32}&0&Q_{33}&Q_{31}&Q_{32}&Q_{33}&Q_{34}&0&V_3\\
Q_{41}&0&Q_{42}&0&Q_{43}&Q_{41}&Q_{42}&Q_{43}&Q_{44}&0&V_4\\
0&I&0&-I&0&0&0&0&0&I&0
\end{array}\right)$}
\]
where the factors have distributions
$$(\alpha*\alpha'\gamma*\delta*\alpha'\gamma*\pi_3*\alpha*\delta*\pi_3*\pi_4*\alpha'\gamma*\gamma'\gamma \, , \,
\alpha*\alpha'\gamma*\delta*\alpha'\gamma*\pi_3*\omega_1*\omega_2*\omega_3*\omega_4*\beta'\gamma),$$
$$(\alpha*\alpha'\gamma*\delta*\alpha'\gamma*\pi_3*\omega_1*\omega_2*\omega_3*\omega_4*\beta'\gamma  \, , \, 
\beta*\beta'\gamma*\varepsilon*\beta'\gamma*\theta_3*\beta*\varepsilon*\theta_3*\theta_4*\beta'\gamma*e),$$
respectively.

Now let $\gamma,\gamma'\in\Gamma$, let $(F,A,X,\alpha,\beta),\in(T_\Sigma)_\gamma$ and 
let $(F',A',X',\alpha',\beta'),\linebreak (G',B',Y',\delta',\varepsilon')\in (T_\Sigma)_{\gamma'}$ 
 be such that $(F',A',X',\alpha',\beta')\sim(G',B',Y',\delta',\varepsilon')$. 
Thus, there exist $L',M',P,Q\in \Sigma$,
homogenous rows $J',U$, and homogeneous columns $W',V$ such that
\begin{equation}
\setlength{\arraycolsep}{1.5pt}\left(\begin{array}{cccc|c}
A'&0&0&0&X'\\0&B'&0&0&Y'\\0&0&L'&0&W'\\0&0&0&M'&0\\\hline F'&-G'&0&J'&0
\end{array}\right)=
\left(\begin{array}{cccc}
P_{11}&P_{12}&P_{13}&P_{14}\\P_{21}&P_{22}&P_{23}&P_{24}\\P_{31}&P_{32}&P_{33}&P_{34}\\P_{41}&P_{42}&P_{43}&P_{44}\\\hline U_{1}&U_{2}&U_{3}&U_{4}
\end{array}\right)
\left(\begin{array}{cccc|c}
Q_{11}&Q_{12}&Q_{13}&Q_{14}&V_{1}\\Q_{21}&Q_{22}&Q_{23}&Q_{24}&V_{2}\\Q_{31}&Q_{32}&Q_{33}&Q_{34}&V_{3}\\Q_{41}&Q_{42}&Q_{43}&Q_{44}&V_{4}
\end{array}\right),
\end{equation}
where $P,U,Q,V$ have distributions $(\pi',\omega'),(\gamma',\omega'),(\omega', \theta'),(\omega',e)$, 
respectively, and \linebreak $\pi_1'=\alpha'$,
 $\pi_2'=\delta'$, $\theta_1'=\beta'$, $\theta_2'=\varepsilon'$.
We show that 
$$(F',A',X',\alpha',\beta')\cdot(F,A,X,\alpha,\beta)\sim 
(G',B',Y',\delta',\varepsilon')\cdot (F,A,X.\alpha,\delta).$$
It follows because the following homogeneous matrix
\[
\left(\begin{array}{ccccccccccc|c}
A&0&0&0&0&0&0&0&0&0&0&X\\
-X'F&A'&0&0&0&0&0&0&0&0&0&0\\
0&0&A&0&0&0&0&0&0&0&0&X\\
0&0&-Y'F&B'&0&0&0&0&0&0&0&0\\
0&0&0&0&A&0&0&0&0&0&0&X\\
0&0&0&0&-W'F&L'&0&0&0&0&0&0\\
0&0&0&0&0&0&A&0&0&0&0&0\\
0&0&0&0&0&0&0&A&0&0&0&X\\
0&0&0&0&0&0&0&-X'F&A'&0&0&0\\
0&0&0&0&0&0&0&0&0&L'&0&0\\
0&0&0&0&0&0&0&0&0&0&M'&0\\\hline 
0&F'&0&-G'&0&0&0&F'&-G'&0&J'&0
\end{array}\right)
\]
can be expressed as the following product of homogeneous matrices

\[
\scalebox{0.75}{$\setlength{\arraycolsep}{1pt}
\left(\begin{array}{ccccccccccc}
I&0&0&0&0&0&0&0&0&0&0\\
0&I&0&0&0&0&0&0&0&0&0\\
0&0&I&0&0&0&0&0&0&0&0\\
0&0&0&I&0&0&0&0&0&0&0\\
0&0&I&0&A&0&0&0&0&0&0\\
0&0&0&0&0&I&0&0&0&0&0\\
I&0&-I&0&0&0&A&0&0&0&0\\
0&-I&0&0&0&0&-X'F&P_{11}&P_{12}&P_{13}&P_{14}\\
0&0&0&-I&0&0&0&P_{21}&P_{22}&P_{23}&P_{24}\\
0&0&0&0&-W'F&-I&0&P_{31}&P_{32}&P_{33}&P_{34}\\
0&0&0&0&0&0&0&P_{41}&P_{42}&P_{43}&P_{44}\\\hline 
0&0&0&0&0&0&0&U_1&U_2&0U_3&U_4
\end{array}\right)
\left(\begin{array}{ccccccccccc|c}
A&0&0&0&0&0&0&0&0&0&0&X\\
-X'F&A'&0&0&0&0&0&0&0&0&0&0\\
0&0&A&0&0&0&0&0&0&0&0&X\\
0&0&-Y'F&B'&0&0&0&0&0&0&0&0\\
0&0&-I&0&I&0&0&0&0&0&0&0\\
0&0&0&0&-W'F&L'&0&0&0&0&0&0\\
-I&0&-I&0&0&0&I&0&0&0&0&0\\
0&Q_{11}&-V_1F&Q_{12}&0&Q_{13}&0&Q_{11}&Q_{12}&Q_{13}&Q_{14}&0\\
0&Q_{21}&-V_2F&Q_{22}&0&Q_{23}&0&Q_{21}&Q_{22}&Q_{23}&Q_{24}&0\\
0&Q_{31}&-V_3F&Q_{32}&0&Q_{33}&0&Q_{31}&Q_{32}&Q_{33}&Q_{34}&0\\
0&Q_{41}&-V_4F&Q_{34}&0&Q_{43}&0&Q_{41}&Q_{42}&Q_{43}&Q_{44}&0
\end{array}\right)$}
\]
where the factors have distributions
$(\alpha*\alpha'\gamma*\alpha*\delta'\gamma*\alpha*\pi_3'\gamma*\alpha*\alpha'\gamma*\delta'\gamma*
\pi_3'\gamma*\pi_4'\gamma*\gamma'\gamma \, , \,
\alpha*\alpha'\gamma*\alpha*\delta'\gamma*\beta*\pi_3'\gamma*\beta*\omega_1'\gamma*\omega_2'\gamma*
\omega_3'\gamma*\omega_4'\gamma)$ and
$(\alpha*\alpha'\gamma*\alpha*\delta'\gamma*\beta*\pi_3'\gamma*\beta*\omega_1'\gamma*\omega_2'\gamma*
\omega_3'\gamma*\omega_4'\gamma\,,\,
\beta*\beta'\gamma*\beta*\varepsilon'\gamma*\beta*\theta_3'\gamma*\beta*\beta'\gamma*\varepsilon'\gamma
*\theta_3'\gamma*\theta_4'\gamma*e)$, respectively.

(4) Let $(F,A,X,\alpha,\beta),(G,B,Y,\delta,\varepsilon)\in(T_\Sigma)_\gamma$ be such that
$(F,A,X,\alpha,\beta)\sim(G,B,Y,\delta,\varepsilon)$. Thus, there exist $L,M,P,Q\in \Sigma$,
homogenous rows $J,U$, and homogeneous columns $W,V$ as in \eqref{eq:equivalencerelation}. The result
follows because we have the factorization
\[
\setlength{\arraycolsep}{3pt}\left(\begin{array}{cccc|c}
A&0&0&0&X\\0&B&0&0&Y\\0&0&L&0&W\\0&0&0&M&0\\\hline -F&G&0&-J&0
\end{array}\right)
=
\left(\begin{array}{cccc}
P_{11}&P_{12}&P_{13}&P_{14}\\
P_{21}&P_{22}&P_{23}&P_{24}\\
P_{31}&P_{32}&P_{33}&P_{34}\\
P_{41}&P_{42}&P_{43}&P_{44}\\\hline -U_{1}&-U_{2}&-U_{3}&-U_{4}
\end{array}\right)
\left(\begin{array}{cccc|c}
Q_{11}&Q_{12}&Q_{13}&Q_{14}&V_{1}\\
Q_{21}&Q_{22}&Q_{23}&Q_{24}&V_{2}\\
Q_{31}&Q_{32}&Q_{33}&Q_{34}&V_{3}\\
Q_{41}&Q_{42}&Q_{43}&Q_{44}&V_{4}
\end{array}\right)
\]
where the factors have distributions
$
(\alpha*\delta*\pi_3*\pi_4*\gamma\,,\,\omega_1*\omega_2*\omega_3*\omega_4)$ and $(\omega_1*\omega_2*\omega_3*\omega_4\,,\,
\beta*\varepsilon*\theta_3*\theta_4*e)$, respectively.

\end{proof}


\subsection{Graded ring structure}\label{subsec:gradedringstructure}
We define $(\mathcal{R}_\Sigma)_\gamma$ as the set of equivalence classes
in $(T_\Sigma)_\gamma$ under the equivalence relation $\sim$. The equivalent class of 
$(F,A,X,\alpha,\beta)\in (T_\Sigma)_\gamma$ will be denoted by $[F,A,X,\alpha,\beta]$.

In Section~\ref{subsec:operations}, we proved that the operation + is well defined in
$(\mathcal{R}_\Sigma)_\gamma$ for each $\gamma\in\Gamma$. 

\begin{lemma}
Let  $\gamma\in\Gamma$. Then $(\mathcal{R}_\Sigma)_\gamma$ is an abelian group with sum defined by
$$[F',A',X',\alpha',\beta']+[F,A,X,\alpha,\beta]=
\left[\begin{pmatrix}F'&F\end{pmatrix},\begin{pmatrix}A'&0\\0&A\end{pmatrix},\begin{pmatrix}X'\\X\end{pmatrix},
\alpha'*\alpha,\beta'*\beta\right]$$
\end{lemma}
\begin{proof}
The operation is well defined and commutative by Lemma~\ref{lem:welldefined}(2) and (1).
 
Now we show that the operation is associative.
Let $[F'',A'',X'',\alpha'',\beta'']$,\linebreak $[F',A',X',\alpha',\beta']$,$[F,A,X,\alpha,\beta]\in (T_\Sigma)_\gamma$. Then
\[
\begin{array}{cl}
&[F'',A'',X'',\alpha'',\beta'']+\left([F',A',X',\alpha',\beta']+[F,A,X,\alpha,\beta]\right)\\\\
=&[F'',A'',X'',\alpha'',\beta'']+\left[\begin{pmatrix}F'&F\end{pmatrix},\begin{pmatrix}A'&0\\0&A\end{pmatrix},\begin{pmatrix}X'\\X\end{pmatrix},\alpha'*\alpha,\beta''*\beta\right]\\\\
=&\left[\begin{pmatrix}F''&F'&F\end{pmatrix},\begin{pmatrix}A''&0&0\\0&A'&0\\0&0&A\end{pmatrix},\begin{pmatrix}X''\\X'\\X\end{pmatrix},\alpha''*\alpha'*\alpha,\beta''*\beta'*\beta\right]\\\\
=&\left[\begin{pmatrix}F''&F'\end{pmatrix},\begin{pmatrix}A''&0\\0&A'\end{pmatrix},\begin{pmatrix}X''\\X'\end{pmatrix},\alpha''*\alpha',\beta''*\beta'\right]+[F,A,X,\alpha,\beta]\\\\
=&\left([F'',A'',X'',\alpha'',\beta'']+[F',A',X',\alpha',\beta']\right)+[F,A,X,\alpha,\beta],
\end{array}
\]
as desired.

The element $\mu(0)=[0,1,1,e,e]$ is the zero element. Indeed,
let $(F,A,X,\alpha,\beta)\in (T_\Sigma)_\gamma$. Then we have
the following factorization
\[
\left(\begin{array}{ccc|c}
A&0&0&X\\0&1&0&1\\0&0&A&X\\\hline F&0&-F&0
\end{array}\right)=
\left(\begin{array}{ccc}
I&0&0\\0&1&0\\I&0&A\\\hline 0&0&-F
\end{array}\right)
\left(\begin{array}{ccc|c}
A&0&0&X\\0&1&0&1\\-I&0&I&0 
\end{array}\right)
\] by Lemma~\ref{lem:easywelldefined}
where the factors have distributions $(\alpha*e*\alpha*\gamma\,,\,\alpha*e*\beta)$ and
$(\alpha*e*\beta\, ,\, \beta*e*\beta*e)$, respectively. Thus,
$[F,A,X,\alpha,\beta]+[0,1,1,e,e]=[F,A,X,\alpha,\beta]$.

Given  $(F,A,X,\alpha,\beta)\in (T_\Sigma)_\gamma$, 
the element $[-F,A,X,\alpha,\beta]$ is well defined by Lemma~\ref{lem:welldefined}(4). We claim
that  it is
the additive inverse of $[F,A,X,\alpha,\beta]$ in $\mathcal{R}_\gamma$.  
Thus, consider the following factorization\[
\left(\begin{array}{ccc|c}
A&0&0&X\\0&A&0&X\\0&0&1&1\\\hline F&-F&0&0
\end{array}\right)=
\left(\begin{array}{ccc}
I&0&0\\I&A&0\\0&0&1\\\hline 0&-F&0
\end{array}\right)
\left(\begin{array}{ccc|c}
A&0&0&X\\-I&I&0&0\\0&0&1&0
\end{array}\right)
\]
where the factors have distributions $(\alpha*\alpha*e*\gamma\,,\,\alpha*\beta*e)$ and
$(\alpha*\beta*e\,,\, \beta*\beta*e*e)$, respectively. It shows that
$[F,A,X,\alpha,\beta]+[-F,A,X,\alpha,\beta]=[0,1,1,e,e]$, as claimed.
\end{proof}

In Section~\ref{subsec:operations}, we showed that the product
functions $(\mathcal{R}_\Sigma)_{\gamma'}\times (\mathcal{R}_\Sigma)_\gamma\rightarrow (\mathcal{R}_\Sigma)_{\gamma'\gamma}$ are well defined. 
Now we define $\mathcal{R}_\Sigma=\bigoplus_{\gamma\in\Gamma} (\mathcal{R}_\Sigma)_\gamma$.
By the foregoing lemma, it
is an additive group. We now prove that it is a $\Gamma$-graded ring with the induced product.

\begin{lemma}
$\mathcal{R}_\Sigma$ is a $\Gamma$-graded ring with the product determined by the rule
\[
[F',A',X',\alpha',\beta']\cdot[F,A,X,\alpha,\beta]= \setlength{\arraycolsep}{1.2pt}
\left[\begin{pmatrix}0&F'\end{pmatrix},\begin{pmatrix}A&0\\-X'F& A'\end{pmatrix},\begin{pmatrix}X\\0\end{pmatrix},
\alpha*\alpha'\gamma, \beta*\beta'\gamma\!\right],
\]
for any $(F',A',X',\alpha',\beta')\in (T_\Sigma)_{\gamma'}$ and $(F,A,X,\alpha,\beta)\in (T_\Sigma)_\gamma$.
\end{lemma}

\begin{proof}
By Lemma~\ref{lem:welldefined}(3), the product is well defined.

By Lemma~\ref{lem:easywelldefined}(3) and (4), the identity element is 
$[1,1,1,e,e]$.

Now we proceed to show that the product is associative. Let\linebreak
 $(F'',A'',X'',\alpha'',\beta'')\in(T_\Sigma)_{\gamma''},$
$(F',A',X',\alpha',\beta')\in(T_\Sigma)_{\gamma'}$ and 
$(F,A,X,\alpha,\beta)\in(T_\Sigma)_{\gamma}$. Then
\[
\begin{array}{cl}
&\left([F'',A'',X'',\alpha'',\beta'']\cdot[F',A',X',\alpha',\beta']\right)\cdot
[F,A,X,\alpha,\beta]\\\\
=&\left[\begin{pmatrix}0&F'\end{pmatrix},\begin{pmatrix}A'& 0\\-X''F'&A'\end{pmatrix},\begin{pmatrix}X'\\0\end{pmatrix},
\alpha'*\alpha''\gamma^{'},\beta'*\beta''\gamma'\right]\cdot[F,A,X,\alpha,\beta]\\\\
=&\setlength{\arraycolsep}{2pt}\left[\begin{pmatrix}0&0&F''\end{pmatrix},\begin{pmatrix}A&0&0\\-X'F&A'&0\\0&-X''F'&A''\end{pmatrix},\begin{pmatrix}X\\0\\0\end{pmatrix},\alpha*\alpha'\gamma*\alpha''\gamma'\gamma,\beta*\beta'\gamma*\beta''\gamma'\gamma\right]\\\\
=&[F'',A'',X'',\alpha'',\beta'']\cdot\left[\begin{pmatrix}0&F'\end{pmatrix},\begin{pmatrix}A&0\\-X'F&A'\end{pmatrix},\begin{pmatrix}X\\0\end{pmatrix},\alpha*\alpha'\gamma,\beta*\beta'\gamma\right]\\\\
=&[F'',A'',X'',\alpha'',\beta'']\cdot\left([F',A',X',\alpha',\beta']\cdot[F,A,X,\alpha,\beta]\right),
\end{array}
\]
which shows that the product is associative.

It  remains to show that the distributive laws are satisfied.
Let
 $(F',A',X',\alpha',\beta'),$
$(G',B',Y',\delta',\varepsilon')\in(T_\Sigma)_{\gamma'}$ and 
$(F,A,X,\alpha,\beta)\in(T_\Sigma)_{\gamma}$. 
First note that 
\begin{align}
\left( [F',A',X',\alpha',\beta']+[G',B',Y',\delta',\varepsilon'] \right)\cdot [F,A,X,\alpha,\beta] =\phantom{aaaaaaaaaaaaaa}\nonumber\\
\left[\begin{pmatrix}0&F'&G'\end{pmatrix},\begin{pmatrix}A&0&0\\-X'F&A'&0 \\ -Y'F&0&B' \end{pmatrix},\begin{pmatrix}X\\0\\0\end{pmatrix},\alpha*\alpha'\gamma*\delta'\gamma,\beta*\beta'\gamma*\varepsilon'\gamma\right]. \label{eq:distributive1}
\end{align}
Second observe that
\begin{align}
[F',A',X',\alpha',\beta']\cdot[F,A,X,\alpha,\beta]+[G',B',Y',\delta',\varepsilon']\cdot[F',A',X',\alpha',\beta']=\phantom{aaa}\nonumber\\  \setlength{\arraycolsep}{1.8pt}
\left[\begin{pmatrix}0&F'&0&G'\end{pmatrix},\begin{pmatrix}A&0&0&0\\-X'F&A'&0&0 \\ 
0&0&A&0\\
0&0&-Y'F&B' \end{pmatrix},\begin{pmatrix}X\\0\\X\\0\end{pmatrix},\alpha*\alpha'\gamma*\alpha*\delta'\gamma,\beta*\beta'\gamma*\beta*\varepsilon'\gamma\right].\label{eq:distributive2}
\end{align}
The fact that \eqref{eq:distributive1} equals \eqref{eq:distributive2} follows because the homogeneous matrix
\[
\left(
\begin{array}{ccccccc|c}
A&0&0&0&0&0&0&X\\
-X'F&A'&0&0&0&0&0&0\\
-Y'F&0&B'&0&0&0&0&0\\
0&0&0&A&0&0&0&X\\
0&0&0&-X'F&A'&0&0&0\\
0&0&0&0&0&A&0&X\\
0&0&0&0&0&-Y'F&B'&0\\\hline 
0&F'&G'&0&-F'&0&-G'&0
\end{array}
\right)
\]
factorizes as the product of homogeneous matrices
\[ \setlength{\arraycolsep}{3pt}
\left(\begin{array}{ccccccc}
I&0&0&0&0&0&0\\
0&I&0&0&0&0&0\\
0&0&I&0&0&0&0\\
I&0&0&A&0&0&0\\
0&I&0&-X'F&A'&0&0\\
I&0&0&0&0&A&0\\
0&0&I&0&0&-Y'F&B'\\\hline 
0&0&0&0&-F'&0&-G'
\end{array}\right)
\left(\begin{array}{ccccccc|c}
A&0&0&0&0&0&0&X\\
-X'F&A'&0&0&0&0&0&0\\
-Y'F&0&B'&0&0&0&0&0\\
-I&0&0&I&0&0&0&0\\
0&-I&0&0&I&0&0&0\\
-I&0&0&0&0&I&0&0\\
0&0&-I&0&0&0&I&0
\end{array}\right)
\]
where the factors have distributions
$$(\alpha*\alpha'\gamma*\delta'\gamma*\alpha*\alpha'\gamma*\alpha*\delta'\gamma*\gamma'\gamma \, ,\,
\alpha*\alpha'\gamma*\delta'\gamma*\beta*\beta'\gamma*\beta*\varepsilon'\gamma),$$
$$(\alpha*\alpha'\gamma*\delta'\gamma*\beta*\beta'\gamma*\beta*\varepsilon'\gamma \, , 
\, \beta*\beta'\gamma*\varepsilon'\gamma*\beta*\beta'\gamma*\beta*\varepsilon'\gamma*e),$$
respectively.

Let now $(F',A',X',\alpha',\beta')\in (T_\Sigma)_{\gamma'}$ and
$(F,A,X,\alpha,\beta)$,$(G,B,Y,\delta,\varepsilon)\in (T_\Sigma)_\gamma$. 
First note that
\begin{align}
 [F',A',X',\alpha',\beta']\cdot\left( [F,A,X,\alpha,\beta]+[G,B,Y,\delta,\varepsilon] \right) =\phantom{aaaaaaaaaaaaaa}\nonumber\\
\left[\begin{pmatrix}0&0&F'\end{pmatrix},\begin{pmatrix}A&0&0\\0&B&0 \\ -X'F&-X'G&A' \end{pmatrix},\begin{pmatrix}X\\Y\\0\end{pmatrix},\alpha*\delta*\alpha'\gamma,\beta*\varepsilon*\beta'\gamma\right]. \label{eq:distributive3}
\end{align}
Second observe that
\begin{align}
[F',A',X',\alpha',\beta']\cdot[F,A,X,\alpha,\beta]+[F',A',X',\alpha',\beta']\cdot[G,B,Y,\delta,\varepsilon]=\phantom{aaaaa}\nonumber\\  \setlength{\arraycolsep}{1.8pt}
\left[\begin{pmatrix}0&F'&0&F'\end{pmatrix},\begin{pmatrix}A&0&0&0\\-X'F&A'&0&0 \\ 
0&0&B&0\\
0&0&-X'G&A' \end{pmatrix},\begin{pmatrix}X\\0\\Y\\0\end{pmatrix},\alpha*\alpha'\gamma*\delta*\alpha'\gamma,\beta*\beta'\gamma*\varepsilon*\beta'\gamma\right].\label{eq:distributive4}
\end{align}
The fact that \eqref{eq:distributive3} equals \eqref{eq:distributive4} follows because the homogeneous matrix
\[
\left(
\begin{array}{ccccccc|c}
A&0&0&0&0&0&0&X\\
-X'F&A'&0&0&0&0&0&0\\
0&0&B&0&0&0&0&Y\\
0&0&-X'G&A'&0&0&0&0\\
0&0&0&0&A&0&0&X\\
0&0&0&0&0&B&0&Y\\
0&0&0&0&-X'F&-X'G&A'&0\\\hline 
0&F'&0&F'&0&0&-F'&0
\end{array}
\right)
\]
factorizes as the product of homogeneous matrices
\[ \setlength{\arraycolsep}{3pt}
\left(\begin{array}{ccccccc}
I&0&0&0&0&0&0\\
0&I&0&0&0&0&0\\
0&0&I&0&0&0&0\\
0&0&0&I&0&0&0\\
I&0&0&0&A&0&0\\
0&0&I&0&0&B&0\\
0&I&0&I&-X'F&-X'G&A'\\\hline 
0&0&0&0&0&0&-F'
\end{array}\right)
\left(\begin{array}{ccccccc|c}
A&0&0&0&0&0&0&X\\
-X'F&A'&0&0&0&0&0&0\\
0&0&B&0&0&0&0&Y\\
0&0&-X'G&A'&0&0&0&0\\
-I&0&0&0&I&0&0&0\\
0&0&-I&0&0&I&0&0\\
0&-I&0&-I&0&0&I&0
\end{array}\right)
\]
where the factors have distributions
$$(\alpha*\alpha'\gamma*\delta*\alpha'\gamma*\alpha*\delta*\alpha'\gamma*\gamma'\gamma \, ,\,
\alpha*\alpha'\gamma*\delta*\alpha'\gamma*\beta*\varepsilon*\beta'\gamma),$$
$$(\alpha*\alpha'\gamma*\delta*\alpha'\gamma*\beta*\varepsilon*\beta'\gamma \, , 
\, \beta*\beta'\gamma*\varepsilon*\beta'\gamma*\beta*\varepsilon*\beta'\gamma*e),$$
respectively.
\end{proof}


\subsection{Universal localization property}\label{subsec:universallocalizationproperty}

\begin{proposition}\label{prop:universallocalization}
Consider the map $\mu\colon R\rightarrow \mathcal{R}_\Sigma$ determined by
$\mu(r)=[r,1,1,e,e]$ for all $r\in R_\gamma$, $\gamma\in\Gamma$. Then
the pair $(\mathcal{R}_\Sigma,\mu)$ is the universal localization of $R$ at $\Sigma$.
\end{proposition}

\begin{proof}
By Lemma~\ref{lem:easywelldefined}(1) and (3), $\mu$ is a homomorphism of $\Gamma$-graded rings.

By $E_i$ we will denote the column matrix consisting of 1 as its $i$-entry and all the other entries are zero,
and by $E_i^T$ its transpose, the row matrix consisting of 1 as its $i$-entry and all other entries are zero. 

Let $A=(a_{ij})\in\Sigma$ an $n\times n$ homogeneous matrix of distribution $(\alpha,\beta)$.

We claim that the $n\times n$  matrix $B=([E_i^T, A, E_j,\alpha\alpha_j^{-1},\beta\alpha_j^{-1}])_{ij}$ is the inverse of $A^\mu$.

First observe
that $[E_i^T, A, E_j,\alpha\alpha_j^{-1},\beta\alpha_j^{-1}]\in \mathcal{R}_{\beta_i\alpha_j^{-1}}$ because
$E_i^T$ has distribution \linebreak $(\beta_i\alpha_j^{-1},\beta\alpha_j^{-1})$,
$A$ has distribution $(\alpha\alpha_j^{-1}, \beta\alpha_j^{-1})$ and $E_j$ has distribution
$(\alpha\alpha_j^{-1},e)$. Thus, $B$
is homogeneous of distribution $(\beta,\alpha)$.

Second, using Lemma~\ref{lem:easywelldefined}(3) and (1), we obtain that the
product of the $i$-th line of $A^\mu$ with the $j$-th column of $B$ equals
\[
\begin{array}{cl}
&\sum_k\mu(a_{i,k})[E_k^T,A,E_j,\alpha\alpha_j^{-1},\beta\alpha_j^{-1}]\\\\
=&\sum_k[a_{i,k}E_k^T,A,E_j,\alpha\alpha_j^{-1},\beta\alpha_j^{-1}]\\\\
=&[\sum_ka_{i,k}E_k^T,A,E_j,\alpha\alpha_j^{-1},\beta\alpha_j^{-1}]\\\\
=&[E_i^TA,A,E_j,\alpha\alpha_j^{-1},\beta\alpha_j^{-1}]\in \mathcal{R}_{\alpha_i\alpha_j^{-1}}.
\end{array}
\]

Third, we show that $$[E_i^TA,A,E_j,\alpha\alpha_j^{-1},\beta\alpha_j^{-1}]=\mu(\delta_{ij})=
[\delta_{ij},1,1,e,e]=\left\{\begin{smallmatrix} 
[1,1,1,e,e] \textrm{ if } i=j, \\ [0,1,1,e,e] \textrm{ if } i\neq j.
\end{smallmatrix}\right.$$
It follows from the following factorization
\[
\left(\begin{array}{cc|c}
A&0&E_j\\0&1&1\\\hline E_i^TA&-\delta_{ij}&0
\end{array}\right)=
\left(\begin{array}{cc}
I&0\\0&1\\\hline E_i^T&-\delta_{ij}
\end{array}\right)
\left(\begin{array}{cc|c}
A&0&E_j\\0&1&1
\end{array}\right)
\]
where the factors have distributions 
$(\alpha\alpha_j^{-1}*e*\alpha_i\alpha_j^{-1}\, , \, \alpha\alpha_j^{-1}*e)$ and 
\linebreak$(\alpha\alpha_j^{-1}*e,\beta\alpha_j^{-1}*e*e)$, respectively. 
Therefore $B$ is the right inverse of $A^\mu$.

Now we proceed to prove that $B$ is the left inverse of $A^\mu$. 
Using Lemma~\ref{lem:easywelldefined}(4) and (2) and considering that $a_{kj}\in R_{\alpha_k\beta_j^{-1}}$, we obtain that the
product of the $i$-th line of $B$ with the $j$-th column of $A^\mu$  equals
\[
\begin{array}{cl}
&\sum_k[E_i^T,A,E_k,\alpha\alpha_k^{-1},\beta\alpha_k^{-1}]\mu(a_{k,j})\\\\
=&\sum_k[E_i^T,A,E_ka_{k,j},\alpha\alpha_k^{-1}\cdot\alpha_k\beta_j^{-1},\beta\alpha_k^{-1}\cdot\alpha_k\beta_j^{-1}]\\\\
=&[E_i^T,A,\sum_kE_ka_{k,j},\alpha\beta_j^{-1},\beta\beta_j^{-1}]\\\\
=&[E_i^T,A,AE_j,\alpha\beta_j^{-1},\beta\beta_j^{-1}]\in\mathcal{R}_{\beta_i\beta_j^{-1}}.
\end{array}
\]
As before, we show that $[E_i^T,A,AE_j,\alpha\beta_j^{-1},\beta\beta_j^{-1}]=
\mu(\delta_{ij})$. It follows from
\[
\left(\begin{array}{cc|c}
A&0&AE_j\\0&1&1\\\hline E_i^T&\delta_{ij}&0
\end{array}\right)=
\left(\begin{array}{cc}
A&0\\0&1\\\hline E_i^T&-\delta_{ij}
\end{array}\right)
\left(\begin{array}{cc|c}
I&0&E_j\\0&1&1
\end{array}\right)
\]
where the factors have distributions $(\alpha\beta_j^{-1}*e*\beta_i\beta_j^{-1} \,,\,
\beta\beta_j^{-1}*e)$ and $(\beta\beta_j^{-1}*e\,,\,\beta\beta_j^{-1}*e*e)$, respectively.

Therefore, the claim is proved.

It remains to prove that $\mu\colon R\rightarrow \mathcal{R}_\Sigma$ is universal.

Note that if $(F,A,X,\alpha,\beta)\in (T_\Sigma)_\gamma$, and we suppose that 
$F=(f_1,\dotsc,f_n)$ and $X=\left(\begin{smallmatrix} x_1\\ \vdots \\ x_n \end{smallmatrix}\right)$,
then
\begin{eqnarray}
F^\mu(A^\mu)^{-1}X^\mu
&=&\sum_{i,j}\mu(f_i)[E_i^T,A,E_j,\alpha\alpha_j^{-1},\beta\alpha_j^{-1}]\mu(x_j) \nonumber \\
&=&\sum_{i,j}[f_iE_i^T,A,E_jx_j,\alpha\alpha_j^{-1}\cdot\alpha_j,\beta\alpha_j^{-1}\cdot\alpha_j] \nonumber \\
&=&[\sum_if_iE_i^T,A,\sum_jE_jx_j,\alpha,\beta] \nonumber \\
&=&[F,A,X,\alpha,\beta]. \label{eq:uniqueness}
\end{eqnarray}

Let now $S$ be a $\Gamma$-graded ring and $\varphi\colon R\rightarrow S$ be  a
$\Sigma$-inverting homomorphism of graded rings. We define 
$\Phi\colon \mathcal{R}_\Sigma \rightarrow S$ as follows. Let 
$(F,A,X,\alpha,\beta)\in (T_\Sigma)_\gamma$, then
\[
\Phi([F,A,X,\alpha,\beta])
=F^\varphi(A^\varphi)^{-1}X^\varphi\in S_\gamma.
\]
Now we show that $\Phi$ is well defined. Let
$(F,A,X,\alpha,\beta),(G,B,Y,\delta,\varepsilon)\in (T_\Sigma)_\gamma$ 
be such that $(F,A,X,\alpha,\beta)\sim(G,B,Y,\delta,\varepsilon)$. Then
there exist $L,M,P,Q\in\Sigma$, homogeneous rows $J,U$
and homogeneous columns $W,V$ such that
\[
\left(\begin{array}{cccc|c}
A&0&0&0&X\\0&B&0&0&Y\\0&0&L&0&W\\0&0&0&M&0\\\hline F&-G&0&J&0
\end{array}\right)=
\left(\begin{array}{c}
P\\\hline U
\end{array}\right)
\left(\begin{array}{c|c}
Q&V
\end{array}\right)
\]
where $P,U,Q,V$ have distributions $(\pi,\omega),(\gamma,\omega),(\omega, \theta),(\omega,e)$, 
respectively, and \linebreak $\pi_1=\alpha$,
 $\pi_2=\delta$, $\theta_1=\beta$, $\theta_2=\varepsilon$.
Then
\begin{eqnarray*}
0 & = & U^\varphi V^\varphi \\\\
  & = &U^\varphi Q^\varphi(Q^\varphi)^{-1}(P^\varphi)^{-1}P^\varphi V^\varphi \\\\
	& = &(U Q)^\varphi((PQ)^\varphi)^{-1}(PV)^\varphi \\\\
	& = & \left(\begin{array}{cccc}F^\varphi&-G^\varphi&0&J^\varphi\end{array}\right)
	\left(\begin{array}{cccc}A^\varphi&0&0&0\\0&B^\varphi&0&0\\0&0&L^\varphi&0\\0&0&0&M^\varphi\end{array}\right)^{-1}
	\left(\begin{array}{c}X^\varphi\\Y^\varphi\\W^\varphi\\0\end{array}\right)\\\\
	& = &\setlength{\arraycolsep}{3pt}\left(\begin{array}{cccc}F^\varphi&-G^\varphi&0&J^\varphi\end{array}\right)
	\left(\begin{array}{cccc}(A^\varphi)^{-1}&0&0&0\\0&(B^\varphi)^{-1}&0&0\\0&0&(L^\varphi)^{-1}&0\\0&0&0&(M^\varphi)^{-1}\end{array}\right)\left(\begin{array}{c}X^\varphi\\Y^\varphi\\W^\varphi\\0\end{array}\right)\\\\
	& = & F^\varphi (A^\varphi)^{-1}X^\varphi-G^\varphi (B^\varphi)^{-1}Y^\varphi+
	 0(L^\varphi)^{-1}W^\varphi+J^\varphi(M^\varphi)^{-1}0\\\\
	& = &F^\varphi (A^\varphi)^{-1}X^\varphi-G^\varphi (B^\varphi)^{-1}Y^\varphi,
\end{eqnarray*}
which shows that $\Phi$ is well defined.

Let now $(F',A',X',\alpha',\beta')$, $(F,A,X,\alpha,\beta)\in (T_\Sigma)_\gamma$. Then
\begin{flalign*}
\Phi([F',A',X',\alpha',\beta']+[F,A,X,\alpha,\beta]) &\setlength{\arraycolsep}{2pt} =  
\left(\begin{array}{cc}F'&F\end{array}\right)^\varphi\left(\left(\begin{array}{cc}A'&0\\0&A\end{array}\right)^\varphi\right)^{-1}\left(\begin{array}{c}X'\\X\end{array}\right)^\varphi\\\\&=
\setlength{\arraycolsep}{2pt}\left(\begin{array}{cc}{F'}^\varphi&F\end{array}\right)\left(\begin{array}{cc}({A'}^{\varphi})^{-1}&0\\0&(A^\varphi)^{-1}\end{array}\right)\left(\begin{array}{c}{X'}^\varphi\\X^\varphi\end{array}\right)\\\\&={F'}^\varphi 
({A'}^\varphi)^{-1}{X'}^\varphi+F^\varphi (A^{\varphi})^{-1}X^\varphi\\\\&=\Phi([F',A',X',\alpha',\beta'])+
\Phi([F,A,X,\alpha,\beta]).
\end{flalign*}
Thus, $\Phi$ is an additive map.

Let $(F',A',X',\alpha',\beta')\in(T_\Sigma)_{\gamma'}$ and $(F,A,X,\alpha,\beta)\in(T_\Sigma)_\gamma$. Then
\begin{flalign*} 
\Phi([F',A',X',\alpha',\beta']&\cdot[F,A,X,\alpha,\beta]) =  \\\\
  & =    \setlength{\arraycolsep}{2pt}\left(\begin{array}{cc}0&F'\end{array}\right)^\varphi \left(\left(\begin{array}{cc}A&0\\-X'F&A'\end{array}\right)^\varphi\right)^{-1} \left(\begin{array}{c}X\\0\end{array}\right)^\varphi \\\\
  & =   \setlength{\arraycolsep}{2pt} \left(\begin{array}{cc} 0  &{F'}^\varphi \end{array}\right)\left(\begin{array}{cc} A^{\varphi} & 0\\
-{X'}^\varphi F^\varphi &{A'}^\varphi \end{array}\right)^{-1} \left(\begin{array}{c}X^\varphi\\0\end{array}\right)
 \\\\
  & =  \setlength{\arraycolsep}{3pt} \left(\begin{array}{cc}0&{F'}^{\varphi}\end{array}\right)\left(\begin{array}{cc}(A^\varphi)^{-1}& 0\\
-({A'}^\varphi)^{-1}{X'}^\varphi F^\varphi (A^\varphi)^{-1} &({A'}^\varphi)^{-1}\end{array}\right)\left(\begin{array}{c}X^\varphi\\0\end{array}\right)\\\\
& =  {F'}^\varphi ({A'}^\varphi)^{-1}{X'}^\varphi F^\varphi (A^{\varphi})^{-1}X^\varphi\\\\
& =  \Phi([F',A',X',\alpha',\beta'])\cdot \Phi([F,A,X,\alpha,\beta]).
\end{flalign*}
Hence, $\Phi$ is a homomorphism of graded rings. Clearly $\Phi\mu=\varphi$. The uniqueness of $\Phi$
now follows from \eqref{eq:uniqueness}.
\end{proof}


\subsection{Malcolmson's criterion}\label{subsec:Malcolmsonscriterion}

Now we proceed to prove two results that determine the kernel of the natural homomorphism $R\rightarrow R_\Sigma$ and an important case in which the ring $R_\Sigma$ is not zero. The following theorem is known as Malcolmson's Criterion.

\begin{theorem}\label{theo:Malcolmsoncriterion}
Let  $R$ be a $\Gamma$-graded ring and $\Sigma$ be a gr-lower semimultiplicative
subset of $\mathfrak{M}(R)$.
Consider the canonical homomorphism of $\Gamma$-graded rings $\lambda\colon R\rightarrow R_\Sigma$.
For $\gamma\in\Gamma$, a homogeneous element $r\in R_\gamma$ belongs to 
$\ker \lambda$ if and only if there exist $L,M,P,Q\in\Sigma$, homogeneous rows $J,U$ and
homogeneous columns $W,V$ such that  
\[
\left(\begin{array}{cc|c}
L&0&W\\0&M&0\\\hline 0&J&r
\end{array}\right)=
\left(\begin{array}{c}
P\\\hline U
\end{array}\right)
\left(\begin{array}{c|c}
Q&V
\end{array}\right),
\]
where $P,U,Q,V$ have distributions $(\pi,\omega),(\gamma,\omega),(\omega,\theta),(\omega,e)$, respectively.
\end{theorem}

\begin{proof}

By Proposition~\ref{prop:universallocalization}, $\mu\colon R\rightarrow \mathcal{R}_\Sigma$ is
the universal localization of $R$ at $\Sigma$. Thus $\lambda(r)=0$ if and only if $\mu(r)=0$.

Hence, suppose that $r\in R_\gamma$ is such that $\mu(r)=0$. It means that $[r,1,1,e,e]\sim [0,1,1,e,e]$.
Thus there exist $L,M,P,Q\in \Sigma$, homogeneous lines $J,U$ and
homogeneous columns $W,V$ such that 
\[\setlength{\arraycolsep}{3pt}
\left(\begin{array}{cccc|c}
1&0&0&0&1\\0&1&0&0&1\\0&0&L&0&W\\0&0&0&M&0\\\hline r&0&0&J&0
\end{array}\right) 
=
\left(\begin{array}{cccc}
P_{11}&P_{12}&P_{13}&P_{14}\\P_{21}&P_{22}&P_{23}&P_{24}\\P_{31}&P_{32}&P_{33}&P_{34}\\P_{41}&P_{42}&P_{43}&P_{44}\\\hline U_{1}&U_{2}&U_{3}&U_{4}
\end{array}\right)
\left(\begin{array}{cccc|c}
Q_{11}&Q_{12}&Q_{13}&Q_{14}&V_{1}\\Q_{21}&Q_{22}&Q_{23}&Q_{24}&V_{2}\\Q_{31}&Q_{32}&Q_{33}&Q_{34}&V_{3}\\Q_{41}&Q_{42}&Q_{43}&Q_{44}&V_{4}
\end{array}\right)
\]
where $P$ has distribution $(\pi,\omega)$, $U$ has distribution $(\gamma,\omega)$,
$Q$ has distribution $(\omega,\theta)$ and $V$ has distribution $(\omega,e)$. 
Now the following equality 
\[\setlength{\arraycolsep}{2pt}
\left(\begin{array}{ccccc|c}
1&0&0&0&0&1\\0&1&0&0&0&1\\0&0&L&0&0&W\\0&0&0&M&0&0\\0&0&0&0&1&0\\\hline 0&0&0&-J&r&r
\end{array}\right)
=
\left(\begin{array}{ccccc}
P_{11}&P_{12}&P_{13}&P_{14}&0\\P_{21}&P_{22}&P_{23}&P_{24}&0\\P_{31}&P_{32}&P_{33}&P_{34}&0\\P_{41}&P_{42}&P_{43}&P_{44}&0\\-P_{11}&-P_{12}&-P_{13}&-P_{14}&1\\\hline -U_{1}&-U_{2}&-U_{3}&-U_{4}&r
\end{array}\right)
\left(\begin{array}{ccccc|c}
Q_{11}&Q_{12}&Q_{13}&Q_{14}&0&V_{1}\\Q_{21}&Q_{22}&Q_{23}&Q_{24}&0&V_{2}\\Q_{31}&Q_{32}&Q_{33}&Q_{34}&0&V_{3}\\
Q_{41}&Q_{42}&Q_{43}&Q_{44}&0&V_{4}\\1&0&0&0&1&1
\end{array}\right)
\]
where the homogeneous matrices of the right hand side have distributions
$(e*e*\pi_3*\pi_4*e\gamma\,,\,\omega_1*\omega_2*\omega_3*\omega_4*e)$ and
$(\omega_1*\omega_2*\omega_3*\omega_4*e\,,\,
e*e*\theta_3*\theta_4*e*e )$, respectively, shows the result.

Conversely, suppose 
there exist $L,M,P,Q\in\Sigma$, homogeneous rows $J,U$ and
homogeneous columns $W,V$ such that  
\[
\left(\begin{array}{cc|c}
L&0&W\\0&M&0\\\hline 0&J&r
\end{array}\right)=
\left(\begin{array}{c}
P\\\hline U
\end{array}\right)
\left(\begin{array}{c|c}
Q&V
\end{array}\right),
\]
where $P,U,Q,V$ have distributions $(\pi,\omega),(\gamma,\omega),(\omega,\theta),(\omega,e)$, respectively.
It follows that $[0,1,1,e,e]\sim[r,1,1,e,e]$ because
\[
\left(\begin{array}{cccc|c}
1&0&0&0&1\\0&1&0&0&1\\0&0&L&0&W\\0&0&0&M&0\\\hline 0&-r&0&J&0
\end{array}\right)
=
\left(\begin{array}{ccc}
1&0&0\\0&1&0\\0&0&P\\\hline 0&-r&U
\end{array}\right)
\left(\begin{array}{ccc|c}
1&0&0&1\\0&1&0&1\\0&0&Q&V
\end{array}\right)
\]
where the factors of the right hand side have distributions
$(e*e*\pi*\gamma\,,\, e*e*\omega)$ and $(e*e*\omega\,,\, e*e*\theta*e)$, respectively.

\end{proof}

\begin{corollary}\label{coro:Malcolmsoncriterion}
Let  $R$ be a  $\Gamma$-graded ring and $\Sigma$ be a gr-multiplicative
subset of $\mathfrak{M}(R)$ consisting of gr-full matrices. Then $R_\Sigma$ is a nonzero
$\Gamma$-graded ring.
\end{corollary}

\begin{proof}
It is enough to prove that $1\in R_e$ is not in the kernel of the canonical homomorphism of
graded rings $\lambda\colon R\rightarrow R_\Sigma$. Suppose that $1\in \ker \lambda$. Then,
by Theorem~\ref{theo:Malcolmsoncriterion},
there exist $L,M,P,Q\in\Sigma$, homogeneous rows $J,U$ and
homogeneous columns $W,V$ such that  
\[
\left(\begin{array}{cc|c}
L&0&W\\0&M&0\\\hline 0&J&1
\end{array}\right)=
\left(\begin{array}{c}
P\\\hline U
\end{array}\right)
\left(\begin{array}{c|c}
Q&V
\end{array}\right),
\]
where $P,U,Q,V$ have distributions $(\pi,\omega),(e,\omega),(\omega,\theta),(\omega,e)$, respectively.
Making elementary column operations, we obtain 
\[
\left(\begin{array}{cc|c}
L&-WJ&W\\0&M&0\\\hline 0&0&1
\end{array}\right)=
\left(\begin{array}{c}
P\\\hline U
\end{array}\right)
\left(\begin{array}{c|c}
Q'&V
\end{array}\right),
\]
where $P,U,Q',V$ have distributions $(\pi,\omega),(e,\omega),(\omega,\theta),(\omega,e)$, 
respectively.
Since $\Sigma$ is gr-multiplicative, it is also upper gr-semimultiplicative
 by Remark~\ref{rem:grmultiplicative}. Thus, the matrix $\left(\begin{smallmatrix}
L&-WJ&W\\0&M&0\\ 0&0&1
\end{smallmatrix}\right)\in \Sigma$ but it is not gr-full, a contradiction. Therefore, $1\notin\ker\lambda$.
\end{proof}

\section{A gr-prime matrix ideal yields a graded division ring, and vice versa}\label{sec:grprimematrixidealyields}

The first part of this section is the adaptation to the graded context of the first part of \cite[Section~7.3]{Cohnfreeeidealringslocalization}. The proof
of the main result Theorem~\ref{theo:primematrixequalsdivisionring} is the graded version of
\cite{Malcolmsonprimematrixideal} using the construction of Section~\ref{sec:Malcolmsonscriterion}.
It could also have been
proved without using the results  in Section~\ref{sec:Malcolmsonscriterion} via a graded version of \cite{Malcolmsonprimematrixideal} that can be found in \cite{DanielMAT0148}.

\medskip

\emph{Throughout this section, let $\Gamma$ be a group.}

\medskip

Let $R$ be a $\Gamma$-graded ring. 
If $(K,\varphi)$ is
a $\Gamma$-graded epic $R$-division ring, the set 
$$\{A\in\mathfrak{M}(R)\colon A^\varphi \textrm{ is not invertible over } K\}$$
will be called the \emph{singular kernel of $(K,\varphi)$}. Now we show that
gr-singular kernels are gr-prime matrix ideals. 
The aim of this section is to show that gr-singular kernels determine graded epic $R$-division rings
in a similar way as commutative $R$-fields are determined by prime ideals of $R$.

Given an $n\times n$ matrix $A$
with entries in $R$, if we write $A=(A_1\ A_2\, \dotsc\, A_n)$
we understand that $A_1,\dotsc,A_n$ are the columns of $A$. And if we write
$A=\left(\begin{smallmatrix} A_1\\ \vdots \\ A_n \end{smallmatrix}\right)$ 
we understand that $A_1,\dotsc,A_n$ are the rows of $A$.

Given two  matrices $A,B\in\mathfrak{M}(R)$, 
we define the \emph{diagonal sum} of $A$ and $B$ as
$$A\oplus B = \left(\begin{array}{cc}
	A&0 \\ 0 & B
\end{array}\right).$$ 
Notice that if $A\in M_m(R)[\overline{\alpha}][\overline{\beta}]$
and $B\in M_n(R)[\overline{\alpha'}][\overline{\beta'}]$, then $A\oplus B\in M_{m+n}(R)[\overline{\alpha}*\overline{\alpha'}][\overline{\beta}*\overline{\beta'}]$.

Let $A,B\in M_n(R)[\overline{\alpha}][\overline{\beta}]$. If
they differ at most in the $i$-th column, then we define the \emph{determinantal sum}
of $A$ and $B$ with respect to the $i$-th column as 
$$A\nabla B=(A_1\ \dotsc\  A_i+B_i\ \dotsc\ A_n).$$ Similarly, if they differ
at most in the $i$-th row we define the determinantal sum of $A$ and $B$
with respect to the $i$-th row as
$$A\nabla B=\left(\begin{array}{c}  A_1 \\ \vdots \\ A_i+B_i \\ \vdots \\A_n\end{array}\right).$$
The matrix $A\nabla B$, when defined, has the same distribution as $A$ and $B$.

Note that the operation $\oplus$ is associative. On the other
hand, the operation $\nabla$ is not always defined, and as
a consequence it is not associative.

Notice  that  distributive laws are satisfied. More precisely,
 if $C$ is another homogeneous matrix, then
$C\oplus(A\nabla B)=(C\oplus A) \nabla (C\oplus B)$ and
$(A\nabla B)\oplus C=(A\oplus C)\nabla (B\oplus C)$ whenever $A\nabla B$ is defined.

On the other hand, if $B,C\in M_n(R)[\overline{\alpha}][\overline{\beta}]$ which differ on at most one column, $A\in M_n(R)[\overline{\alpha'}][\overline{\alpha}]$, $D\in M_n(R)[\overline{\beta}][\overline{\beta'}]$, then may happen that
$$A(B\nabla C)\neq AB\nabla AC,\qquad (B\nabla C)D\neq BD\nabla CD,$$
because, for example, $AB$ and $AC$ ($BD,BC$) may differ on more than $1$ row/column. But in some cases we one
can apply the distributive law. Let $X\in \mathfrak{M}(R)$ and suppose that either
$X$ is a diagonal matrix, or $X$ is a permutation matrix, 
then $$X(B\nabla C)= XB\nabla XC,\qquad (B\nabla C)X = BX\nabla CX.$$ Moreover,
we can regard $X\in M_n(R)[\overline{\alpha'}][\overline{\alpha}]
\cap M_n(R)[\overline{\beta}][{\overline{\beta'}}]$ for some $\alpha',\beta'\in \Gamma^n$. Thus
$X(B\nabla C)\in M_n(R)[\overline{\alpha'}][\overline{\beta}]$ and
$(B\nabla C)X\in M_n(R)[\overline{\alpha}][\overline{\beta'}]$.

\medskip

Let $\Gamma$ be a group and let $R$ be a $\Gamma$-graded ring. 
A subset $\mathcal{P}$ of $\mathfrak{M}(R)$ is a gr-prime matrix ideal if 
the following conditions are satisfied.
\begin{enumerate}[(PM1)]\label{def:grprimematrixideal}
	\item $\mathcal{P}$ contains all the homogeneous matrices that are not gr-full;
	\item If $A,B\in\mathcal{P}$ and their determinantal sum (with respect to
	a row or column) exists, then
	$A\nabla B\in\mathcal{P}$;
	\item If $A\in\mathcal{P}$, then $A\oplus B\in\mathcal{P}$ for all
	$B\in\mathfrak{M}(R)$;
	\item For $A,B\in\mathfrak{M}(R)$,  $A\oplus B\in\mathcal{P}$ implies that
	$A\in\mathcal{P}$ or $B\in\mathcal{P}$;
	\item $1\notin \mathcal{P}$;
	\item If $A\in\mathcal{P}$ and $E,F$ are permutation matrices of appropriate size, then $EAF\in\mathcal{P}$.
\end{enumerate}

We remark that when $\Gamma=\{1\}$, that is, the ungraded case, (PM6) is a consequence of
(PM1)--(PM5) as shown in \cite[(g) in p.431]{Cohnfreeeidealringslocalization}. 
We have not been able to obtain (PM6) from the others in the general graded case.

\begin{proposition}\label{prop:viceversa}
Let  $R$ be a $\Gamma$-graded ring. Let
$K$ be a $\Gamma$-almost graded division ring and 
$\varphi\colon R\rightarrow K$ be a homomorphism of $\Gamma$-almost graded rings.
Then 
$$\mathcal{P}=\{A\in\mathfrak{M}(R)\colon A^\varphi \textrm{ is not invertible }\}$$
is a gr-prime matrix ideal.
Therefore, the following assertions hold true.
\begin{enumerate}[\rm(1)]
	\item If $(K,\varphi)$ is a $\Gamma$-graded epic $R$-division ring, then the gr-singular
	kernel of $(K,\varphi)$ is a $\Gamma$-gr-prime matrix ideal.
	\item Let $N$ be a normal subgroup of $\Gamma$ and consider $R$ as a $\Gamma/N$-graded
	ring. Let $(K,\varphi)$ be a $\Gamma/N$-graded epic $R$-division ring. Then
	$$\mathcal{P}=\{A\in\mathfrak{M}_\Gamma(R)\colon A^\varphi \textrm{ is not invertible }\}$$
	is a $\Gamma$-gr-prime matrix ideal.
\end{enumerate}
\end{proposition}

\begin{proof}
 Let
$K$ be a $\Gamma$-almost graded division ring and 
$\varphi\colon R\rightarrow K$ be a homomorphism of $\Gamma$-almost graded rings. 

First suppose that $K=\DC(\varphi)$ and let $$\Sigma=\mathfrak{M}(R)\setminus\mathcal{P}=\{A\in\mathfrak{M}(R)\colon A^\varphi \textrm{ is invertible over } K\}.$$
By Theorem~\ref{theo:gradedlocal}(2), $R_\Sigma$ is a local ring. If $\mathfrak{m}$ is
the maximal graded ideal of $R_\Sigma$, there exists a surjective homomorphism
of $\Gamma$-almost graded rings $\widetilde{\Phi}\colon R_\Sigma/\mathfrak{m}\rightarrow K$
such that the following diagram is commutative
$$\xymatrix{R\ar[r]^\lambda\ar[rd]_\varphi & R_\Sigma \ar[d]^\Phi\ar[r]^\pi & R_\Sigma/\mathfrak{m}\ar[ld]^{\widetilde{\Phi}} \\ 
		& K & }$$


By Proposition~\ref{prop:almostdivisionring}(3), the sets $\{A\in\mathfrak{M}(R)\colon
A^{(\pi\lambda)} \textrm{ is invertible over } R_\Sigma/\mathfrak{m}\}$ and
$\{A\in\mathfrak{M}(R)\colon A^{(\widetilde{\Phi}\pi\lambda)} \textrm{ is invertible over }K\}$
are equal. Because this last set equals $\Sigma$, we get that
$\Sigma=\{A\in\mathfrak{M}(R)\colon
A^{(\pi\lambda)} \textrm{ is invertible over } R_\Sigma/\mathfrak{m}\}.$
Now, since $(R_\Sigma/\mathfrak{m}, \pi\lambda)$ is a $\Gamma$-graded epic $R$-division ring,
it is enough to prove (1). Thus, suppose that $K$ is a $\Gamma$-graded division ring and 
$(K,\varphi)$ is a $\Gamma$-graded epic $R$-division ring. 
Let $\mathcal{P}=\{A\in\mathfrak{M}(R)\colon A^\varphi \textrm{ is not invertible over } K\}$. 

If $A\in\mathfrak{M}(R)$ is not gr-full, then $A^\varphi$ is not gr-full. Since $K$
is a $\Gamma$-graded division ring, $A^\varphi$ is not invertible over $K$. Thus, (PM1) is
satisfied.

Let now $A,B\in\mathcal{P}_n[\overline{\alpha}][\overline{\beta}]$ such that $A\nabla B$ is defined.  We may suppose that
$A,B$ differ on the first column. Hence $A=(A_1\ C_2\,\dotsc\,C_n)$ and
$B=(B_1\ C_2\,\dotsc\,C_n)$. 
Since $A^\varphi$ and $B^\varphi$ are not invertible over $K$, the columns of $A^\varphi$ and $B^\varphi$ are right linearly dependent over $K$.
If  the columns $C_2^\varphi,\dotsc,C_n^\varphi$ are right linearly dependent over $K$,
then the columns of $(A\nabla B)^\varphi$ are right linearly dependent over $K$ and thus $A\nabla B\in\mathcal{P}$.
Hence we can suppose that there exist homogeneous elements
$a_1,\dotsc,a_n, b_1,\dotsc,b_n\in K$, with $a_1,b_1\neq 0$, such that $$A_1a_1+C_2a_2+\dotsb+C_na_n=0,\quad
B_1b_1+C_2b_2+\dotsb+C_nb_n=0.$$
But then 
$$A_1+B_1+ C_2(a_2a_1^{-1}+b_2b_1^{-1})+\dotsb +C_n(a_na_1^{-1}+b_nb_1^{-1})=0,$$
which shows that $A\nabla B\in\mathcal{P}$. Thus (PM2) is proved.

Let $A\in\mathcal{P}$ and $B\in\mathfrak{M}(R)$, then $A^\varphi$ is not invertible over $K$,
but then $A^\varphi\oplus B^\varphi=(A\oplus B)^\varphi$ is not invertible over $K$. It implies (PM3).

Now suppose that $A,B\in\mathfrak{M}(R)$ are such that $A\oplus B\in\mathcal{P}$. It means that the homogeneous matrix
$A^\varphi\oplus B^\varphi$ is not invertible over $K$. It implies that either $A^\varphi$ or $B^\varphi$
is not invertible. That is, $A\in\mathcal{P}$ or $B\in\mathcal{P}$ and (PM4) follows.

Clearly, (PM5) is satisfied.

Let $A\in\mathcal{P}$ and $E,F$ be  permutation matrices with entries in $R$. Notice that $E^\varphi, F^\varphi$ are permutation matrices
with entries in $K$. Thus,  if $(EAF)^\varphi=E^\varphi A^\varphi F^\varphi$ were invertible over $K$,
then $A^\varphi=(E^\varphi)^{-1}(EAF)^\varphi (F^\varphi)^{-1}$ would be invertible over $K$, a contradiction. Thus $EAF\in\mathcal{P}$ and
(PM6) is shown.
\end{proof}

\begin{lemma}\label{lem:grprimematrixideal}
Let $R$ be a $\Gamma$-graded ring and $\mathcal{P}$ be a gr-prime matrix ideal. 
 Let $A,B\in\mathfrak{M}(R)$. The following assertions hold true.
\begin{enumerate}[\rm(1)]
\item  If $A$ and $B$ are such that
$C=A\nabla B$ exists and $B$ is not gr-full. Then $A\in\mathcal{P}$ if and only if $C\in\mathcal{P}$.
\item Let $A\in\mathcal{P}$. The result of adding a suitable  right multiple of one column of $A$ to another column again
lies in $\mathcal{P}$. More precisely, if $A\in M_n(R)[\overline{\alpha}][\overline{\beta}]$ and $a\in 
R_{\beta_i\beta_j^{-1}}$,
then $(A_1\,\dotsc\, A_{j-1}\ A_j+A_ia\ A_{j+1}\,\dotsc\, A_n)\,$ belongs to $\mathcal{P}$.
\item If $A\oplus B\in\mathcal{P}$, then $B\oplus A\in\mathcal{P}$.
\item Suppose that $A\in M_m(R)[\overline{\alpha}][\overline{\beta}]$ and $B\in M_n(R)[\overline{\delta}][\overline{\varepsilon}]$.
For $C\in M_{n\times m}(R)[\overline{\delta}][\overline{\beta}]$, 
\[
\begin{pmatrix}
A&0\\C&B
\end{pmatrix}\in\mathcal{P}\ \textrm{ if and only if } \
\begin{pmatrix}
A&0\\0&B
\end{pmatrix}\in\mathcal{P}
\]
Similarly, for $C\in M_{m\times n}(R)[\overline{\beta}][\overline{\varepsilon}]$,
\[
\begin{pmatrix}
A&C\\0&B
\end{pmatrix}\in\mathcal{P}\ \textrm{ if and only if } \
\begin{pmatrix}
A&0\\0&B
\end{pmatrix}\in\mathcal{P}
\]
\item The set $\mathfrak{M}(R)\setminus \mathcal{P}$ is gr-multiplicative.
\item No identity matrix belongs to $\mathcal{P}$
\item Suppose that $A\in M_n(R)[\overline{\alpha}][\overline{\beta}]$ and $B\in M_n(R)[\overline{\beta}][\overline{\delta}]$. Then $AB\in\mathcal{P}$ if, and only if, $A\oplus B\in\mathcal{P}$.
\item No invertible matrix in $\mathfrak{M}(R)$ belongs to $\mathcal{P}$.
\item Suppose that $A$ and $B$ are such that
$C=A\nabla B$ exists and $B\in\mathcal{P}$. Then $A\in\mathcal{P}$ if, and only if, $C\in\mathcal{P}$.
\end{enumerate}
\end{lemma}

\begin{proof}
(1) By (PM1) and (PM2), if $A\in\mathcal{P}$, then $C\in\mathcal{P}$. Conversely, suppose that $C\in\mathcal{P}$. 
Clearly $A=C\nabla B'$ where  $B'$ is obtained from $B$ changing the sign of a row or column. 
Now $A\in\mathcal{P}$ because $B'$ is not gr-full.

(2) Suppose that $\overline{\beta}=\beta_1*\overline{\beta'}$ with and $c\in R_{\beta_2\beta_1^{-1}}$. If $A=(A_1\ A_2\,\dotsc\,A_n)$,
then
 \[
\begin{array}{lcl}
(A_1+A_2c\ A_2\ \dots\ A_n)& = & (A_1\ A_2\ \dots\ A_n)\nabla(A_2c\ A_2\ \dots\ A_n)\\
& = & A\nabla (A_2\ A_3\ \dots\ A_n)
\begin{pmatrix}
c & 1 & & 0\\&&\ddots&\\0&0&&1
\end{pmatrix}
\end{array}.
\]
Thus, the right hand side is a determinantal sum of $A$ and $(A_2c\ A_2\ \dots\ A_n)$, which is
not gr-full. Indeed, it is 
the product of  $(A_2\ A_3\ \dots\ A_n)\in M_{n\times(n-1)}(R)[\overline{\alpha}][\overline{\beta'}]$ and 
$\left(\begin{smallmatrix}
c & 1 & & 0\\&&\ddots&\\0&0&&1
\end{smallmatrix}\right)\in M_{(n-1)\times n}(R)[\overline{\beta'}][\overline{\beta}]$.
respectively.

(3) It follows from (PM6).

(4) We show the first statement, the other can be proved analogously. 
If we write $A=(A_1\ A')$ and $C=(C_1\ C')$, then
\[
\begin{pmatrix}
A&0\\C&B
\end{pmatrix}=
\begin{pmatrix}
A_1&A'&0\\0&C'&B
\end{pmatrix}\nabla
\begin{pmatrix}
0&A'&0\\C_1&C'&B
\end{pmatrix}
\]
The second matrix of the right hand side is a matrix with a 
submatrix that is a block of zeros of size $m\times (n+1)$.
Since $m+n+1>m+n$, that matrix is hollow and therefore not gr-full. 
By (1), 
\[
\begin{pmatrix}
A&0\\C&B
\end{pmatrix}\in\mathcal{P}\ \textrm{ if and only if }\
\begin{pmatrix}
A_1&A'&0\\0&C'&B
\end{pmatrix}\in\mathcal{P}
\]
Similarly, one can repeat the argument applied to columns of $A'$ and $C'$ and so on, 
to obtain the desired result.

(5) Let $\Sigma=\mathfrak{M}(R)\setminus \mathcal{P}$. By (PM5), $1\in \Sigma$.
By (PM4), $A\oplus B\in\Sigma$ if $A,B\in\Sigma$. Now (4) implies that $\Sigma$
is lower gr-semimultiplicative. Finally, (PM6) shows that $\Sigma$ is gr-multiplicative.  

(6) It follows from (PM4) and (PM5).

(7) First notice that, by (6) and (PM4),  a matrix $C\in \mathfrak{M}(R)$ belongs to $\mathcal{P}$ if and only if $C\oplus I\in \mathcal{P}$ for the identity matrix $I$ of the same size as $C$. 

We claim that $C\in\mathcal{P}$ if and only if $-C\in\mathcal{P}$. Indeed,
\begin{flalign*}
\begin{pmatrix}
C&0\\0&I
\end{pmatrix}\in\mathcal{P}\stackrel{(2)}{\Leftrightarrow}
\begin{pmatrix}
C&-C\\0&I
\end{pmatrix}\in\mathcal{P} \stackrel{(2)}{\Leftrightarrow} 
\begin{pmatrix}
0&-C\\I&I
\end{pmatrix}\in\mathcal{P} & \\ \stackrel{\textrm{(PM6)}}{\Leftrightarrow} 
\begin{pmatrix}
-C&0\\I&I
\end{pmatrix}\in\mathcal{P}\stackrel{(4)}{\Leftrightarrow}
\begin{pmatrix}
-C&0\\0&I
\end{pmatrix}\in\mathcal{P},
\end{flalign*}
and the claim is proved. Then
\begin{flalign*}
\begin{pmatrix}
A&0\\0&B
\end{pmatrix}\in\mathcal{P} \stackrel{(4)}{\Leftrightarrow}
\begin{pmatrix}
A&0\\I&B
\end{pmatrix}\in\mathcal{P}\stackrel{(2)}{\Leftrightarrow}
\begin{pmatrix}
A&-AB\\I&0
\end{pmatrix}\in\mathcal{P}& \\ \stackrel{\textrm{(PM6)}}{\Leftrightarrow}
\begin{pmatrix}
-AB&A\\0&I
\end{pmatrix}\in\mathcal{P}\stackrel{(4)}{\Leftrightarrow}
\begin{pmatrix}
-AB&0\\0&I
\end{pmatrix},
\end{flalign*}
and, by the claim, the result follows.

(8) If $A\in M_n(R)[\overline{\alpha}][\overline{\beta}]$ is  invertible,
then $A^{-1}\in M_n(R)[\overline{\beta}][\overline{\alpha}]$. Since
$AA^{-1}=I\notin\mathcal{P}$, (7) implies that $A\oplus A^{-1}\notin\mathcal{P}$.
Now (PM3) shows that $A\notin\mathcal{P}$.

(9) By (PM2), if $A\in\mathcal{P}$, then $C\in\mathcal{P}$. 
Conversely, suppose that $C\in\mathcal{P}$. 
Clearly $A=C\nabla B'$ where  $B'$ is obtained from $B$ changing the sign of a row or column.
More precisely, $B'$ is the product of $B$ by a diagonal matrix $D$ whose diagonal elements
are $1$ or $-1$.  Now $B\oplus D \in\mathcal{P}$
because $B\in\mathcal{P}$. Thus  $B'\in\mathcal{P}$ by (7).
Therefore $A\in\mathcal{P}$ by (PM2).
\end{proof}

The proof of Lemma~\ref{lem:grprimematrixideal} is very similar to the one for the ungraded case, see for example \cite[p. 430--431]{Cohnfreeeidealringslocalization}. The main difference is that we were not
able to show \cite[(d) p. 430]{Cohnfreeeidealringslocalization} because not every multiple of a column
can be added to another column so that the matrix remains homogeneous. As a consequence the proof
of Lemma~\ref{lem:grprimematrixideal}(7) is also different.

\medskip

Let $R$ be a graded ring and let $\mathcal{P}$ be a gr-prime matrix ideal. 
The universal localization of $R$ at the set $\Sigma=\mathfrak{M}(R)\setminus \mathcal{P}$
will be denoted by $R_\mathcal{P}$ (instead of $R_\Sigma$).

\begin{theorem}\label{theo:primematrixequalsdivisionring}
Let  $R$ be a $\Gamma$-graded ring. The following assertions hold true.
\begin{enumerate}[\rm(1)]
	\item If $\mathcal{P}$ is any gr-prime matrix ideal of $R$,  then the localization
	$R_\mathcal{P}$ is a $\Gamma$-graded local ring. Moreover, its
	residue class $\Gamma$-graded division ring is a $\Gamma$-graded epic $R$-division ring
	such that its gr-singular kernel equals $\mathcal{P}$.
	\item If $(K,\varphi)$ is a $\Gamma$-graded epic $R$-division ring, with gr-singular kernel $\mathcal{P}$, then
	$\mathcal{P}$ is a gr-prime matrix ideal and the $\Gamma$-graded local ring $R_\mathcal{P}$
	has residue class graded division ring $R$-isomorphic to $K$.
\end{enumerate}
\end{theorem}

\begin{proof}
(2) By Proposition~\ref{prop:viceversa}, $\mathcal{P}$ is a gr-prime matrix ideal.

By (1), $R_\mathcal{P}$ is a $\Gamma$-graded local ring and its residue class graded division ring is a
graded epic $R$-division ring with singular kernel $\mathcal{P}$. Then, by 
Theorem~\ref{theo:gradedlocal}(b)(ii), $K$ and the residue class graded division ring of $R_\mathcal{P}$
are isomorphic $\Gamma$-graded $R$-rings.

(1) Let $\Sigma=\mathfrak{M}(R)\setminus \mathcal{P}$. We will use the identification of
$R_\mathcal{P}$ with the ring $\mathcal{R}_\Sigma=\bigoplus\limits_{\gamma\in\Gamma}(\mathcal{R}_\Sigma)_\gamma$ given in Sections~\ref{subsec:gradedringstructure},\ref{subsec:universallocalizationproperty}.
The elements of $(\mathcal{R}_\Sigma)_\gamma$ are equivalence classes $[F,A,X,\alpha,\beta]$ under the
equivalence relation $\sim$ defined in Section~\ref{subsec:equivalencerelation} of
5-tuples $(F,A,X,\alpha,\beta)$, where $A\in\Sigma$ of distribution $(\alpha,\beta)$, $F$ is
a homogeneous row of distribution $(\gamma,\beta)$ and $X$ is a homogeneous column of distribution
$(\alpha,e)$.

For each $\gamma\in\Gamma$, let $\mathfrak{P}_\gamma$ be the subset of $(\mathcal{R}_\Sigma)_\gamma$ consisting of the elements $[F,A,X,\alpha,\beta]\in(\mathcal{R}_\Sigma)_\gamma$ such that
$\begin{pmatrix}
A & X \\ F & 0
\end{pmatrix}\in\mathcal{P}$. Notice that this matrix is homogeneous of distribution $(\alpha*\gamma,\, \beta*e)$.

\underline{Step~1:} For each $\gamma\in\Gamma$, $\mathfrak{P}_\gamma$ is well defined. That is,
if $(F,A,X,\alpha,\beta)\sim(G,B,Y,\delta,\varepsilon)$ and
$\begin{pmatrix}
A & X \\ F & 0
\end{pmatrix}\in\mathcal{P}$, then $\begin{pmatrix}
B & Y \\ G & 0
\end{pmatrix}\in\mathcal{P}$.

Suppose  $(F,A,X,\alpha,\beta)\sim(G,B,Y,\delta,\varepsilon)$. There exist $L,M,P,Q\in \Sigma$, homogeneous rows $J,U$ and
homogeneous columns $W,V$ such that
\[
\left(\begin{array}{cccc|c}
A&0&0&0&X\\0&B&0&0&Y\\0&0&L&0&W\\0&0&0&M&0\\\hline F&-G&0&J&0
\end{array}\right)
\]
is not gr-full, by \eqref{eq:equivalencerelation}, and thus belongs to $\mathcal{P}$.
Applying permutations of rows and columns we obtain that 
\[ \left(\begin{array}{ccccc}
A&0&0&0&X\\0&B&0&0&Y\\F&-G&0&J&0\\0&0&L&0&W\\ 0&0&0&M&0
\end{array}\right)
\in\mathcal{P}, \quad \left(\begin{array}{ccccc}
A&0&0&X&0\\0&B&0&Y&0\\F&-G&0&0&J\\0&0&L&W&0\\ 0&0&0&0&M
\end{array}\right)\in\mathcal{P}.
\]
Since $M\notin\mathcal{P}$, then 
\[
\left(\begin{array}{cccc}
A&0&0&X\\0&B&0&Y\\F&-G&0&0\\0&0&L&W
\end{array}\right)\in\mathcal{P}
\]
by Lemma~\ref{lem:grprimematrixideal}(4) and
(PM4). Again, aplying column permutations, Lemma~\ref{lem:grprimematrixideal}(4),
(PM4) and the fact that $L\notin\mathcal{P}$ we obtain
\[
\left(\begin{array}{ccc}
A&0&X\\0&B&Y\\F&-G&0
\end{array}\right)\in\mathcal{P}.
\]
After permuting some rows and columns, we obtain
\[
\left(\begin{array}{ccc}
A&X&0\\0&Y&B\\F&0&-G
\end{array}\right)\in\mathcal{P},\quad 
\left(\begin{array}{ccc}
A&X&0\\F&0&-G\\0&Y&B
\end{array}\right)\in\mathcal{P}.
\]
Now observe that $\left(\begin{array}{ccc}
A&X&0\\F&0&-G\\0&0&B
\end{array}\right)\in\mathcal{P}$ because
$[F,A,X,\alpha,\beta]\in\mathfrak{P}_\gamma$, Lemma~\ref{lem:grprimematrixideal}(4) and (PM3).
Hence the equality 
\[\left(\begin{array}{ccc}
A&0&0\\F&0&-G\\0&Y&B
\end{array}\right)\nabla \left(\begin{array}{ccc}
A&X&0\\F&0&-G\\0&0&B
\end{array}\right)
=\left(\begin{array}{ccc}
A&X&0\\F&0&-G\\0&Y&B
\end{array}\right)\]
implies that $\left(\begin{array}{ccc}
A&0&0\\F&0&-G\\0&Y&B
\end{array}\right)\in\mathcal{P}$ by Lemma~\ref{lem:grprimematrixideal}(9). Then, since 
$A\notin\mathcal{P}$, $\left(\begin{array}{cc}
0&-G\\Y&B
\end{array}\right)\in\mathcal{P}$ by Lemma~\ref{lem:grprimematrixideal}(4) and (PM4). After
permuting some rows and columns, $\left(\begin{array}{cc}
B&Y\\-G&0
\end{array}\right)\in\mathcal{P}$.
Now
\begin{equation}\label{eq:step1}
\left(\begin{array}{cc}I&0\\0&-1
\end{array}\right)\left(\begin{array}{cc}
B&Y\\-G&0
\end{array}\right)=\left(\begin{array}{cc}
B&Y\\G&0
\end{array}\right)\in\mathcal{P}
\end{equation}
by (PM3) and Lemma~\ref{lem:grprimematrixideal}(7),(8).

\underline{Step~2:} For each $\gamma\in\Gamma$, $\mathfrak{P}_\gamma$ is an
additive subgroup of $(\mathcal{R}_\Sigma)_\gamma$.

The $2\times 2$ matrix $\left(\begin{array}{cc}
1&1\\0&0
\end{array}\right)\in\mathcal{P}$ is hollow and therefore not gr-full. Thus it belongs to $\mathcal{P}$,
and $[0,1,1,e,e]\in\mathfrak{P}_\gamma$. Thus, the zero element of $(\mathcal{R}_\Sigma)_\gamma$ belongs to $\mathfrak{P}_\gamma$.

If $[F,A,X,\alpha,\beta]\in\mathfrak{P}_\gamma$, then
$-[F,A,X,\alpha,\beta]=[-F,A,X,\alpha,\beta]\in\mathfrak{P}_\gamma$ because
if $\left(\begin{array}{cccc}
A&X\\F&0
\end{array}\right)\in\mathcal{P}$, then $\left(\begin{array}{cccc}
A&X\\-F&0
\end{array}\right)\in\mathcal{P}$ by the same argument as \eqref{eq:step1}.

So let now $[F',A',X',\alpha',\beta'],[F,A,X,\alpha,\beta]\in\mathfrak{P}_\gamma$. Then
\[
[F',A',X',\alpha',\beta'] + [F,A,X,\alpha,\beta]=
\left[\begin{pmatrix}F'&F\end{pmatrix},\begin{pmatrix}A'&0\\0&A\end{pmatrix},\begin{pmatrix}X'\\X\end{pmatrix},
\alpha'*\alpha,\beta'*\beta\right].
\]
To show that the previous sum belongs to $\mathfrak{P}_\gamma$, since 
\[\begin{pmatrix}
A' & 0 &0 \\ 0&A&X\\ F'&F&0
\end{pmatrix} \nabla \begin{pmatrix}
A' & 0 & X'\\ 0&A&0\\ F'&F&0
\end{pmatrix}=\begin{pmatrix}
A'&0&X'\\ 0&A&X\\ F'&F&0
\end{pmatrix},\]
it is enough to show that both summands belong to $\mathcal{P}$ by (PM2). 
The homogeneous matrix $\begin{pmatrix}
A' & 0 &0 \\ 0&A&X\\ F'&F&0
\end{pmatrix}\in\mathcal{P}$ by (PM3) and Lemma~\ref{lem:grprimematrixideal}(4), because $\begin{pmatrix}
A&X\\ F&0
\end{pmatrix}\in\mathcal{P}$. By a similar argument $\begin{pmatrix}
A' &  X' & 0\\ F'&0&F\\ 0&0&A
\end{pmatrix}\in \mathcal{P}$. Permuting rows and columns, we obtain
that $\begin{pmatrix}
A' & 0 & X'\\ 0&A&0\\ F'&F&0
\end{pmatrix}\in\mathcal{P}$.

\underline{Step~3:} If $[F',A',X',\alpha',\beta']\in (\mathcal{R}_\Sigma)_{\gamma'}$, and
$[F,A,X,\alpha,\beta]\in\mathfrak{P}_\gamma$, then
$$[F',A',X',\alpha',\beta']\cdot [F,A,X,\alpha,\beta]\in\mathfrak{P}_{\gamma'\gamma}.$$
Similarly, if $[F',A',X',\alpha',\beta']\in \mathfrak{P}_{\gamma'}$, and
$[F,A,X,\alpha,\beta]\in(\mathcal{R}_\Sigma)_\gamma$, then
$$[F',A',X',\alpha',\beta']\cdot [F,A,X,\alpha,\beta]\in\mathfrak{P}_{\gamma'\gamma}.$$

We prove both cases at the same time. Observe that
\[
[F',A',X',\alpha',\beta']\cdot[F,A,X,\alpha,\beta]= \setlength{\arraycolsep}{1.2pt}
\left[\begin{pmatrix}0&F'\end{pmatrix},\begin{pmatrix}A&0\\-X'F& A'\end{pmatrix},\begin{pmatrix}X\\0\end{pmatrix},
\alpha*\alpha'\gamma, \beta*\beta'\gamma\right].
\]
First note that the matrix
$\begin{pmatrix}
A & X & 0 & 0 \\
F & 0 & 0 & 1 \\
0 & 0 & A' & X' \\
0 & 0 & F' & 0
\end{pmatrix}\in\mathcal{P}$ by Lemma~\ref{lem:grprimematrixideal}(4) and (PM3), and that it is
of distribution $(\alpha*\gamma*\alpha'\gamma*\gamma'\gamma,\, \beta*e*\beta'\gamma*\gamma)$.
Now the matrix $\begin{pmatrix}
I&0&0&0\\ 0&1&0&0 \\ 0&-X'&I&0 \\ 0 &0&0&1
\end{pmatrix}$ is invertible, homogeneous and of distribution
$(\alpha*\gamma*\alpha'\gamma*\gamma^{-1}\gamma'^{-1},\, \alpha*\gamma*\alpha'\gamma*\gamma'\gamma).$
By Lemma~\ref{lem:grprimematrixideal}(7) and (PM3), we obtain
\[
\begin{pmatrix}
I&0&0&0\\ 0&1&0&0 \\ 0&-X'&I&0 \\ 0 &0&0&1
\end{pmatrix} \begin{pmatrix}
A & X & 0 & 0 \\
F & 0 & 0 & 1 \\
0 & 0 & A' & X' \\
0 & 0 & F' & 0
\end{pmatrix}= \begin{pmatrix}
A & X & 0 & 0 \\
F & 0 & 0 & 0 \\
-X'F & 0 & A' & 0 \\
0 & 0 & F' & 0
\end{pmatrix}\in\mathcal{P}.
\]
Permuting rows and columns, we get
$\begin{pmatrix}
A & 0 & X & 0 \\
-X'F & A' & 0 & 0 \\
0 & F' & 0 & 0 \\
F & 0 & 0 & 1
\end{pmatrix}\in\mathcal{P}.$ Now Lemma~\ref{lem:grprimematrixideal}(4), (PM5)
and (PM4) imply that $\begin{pmatrix}
A & 0 & X  \\
-X'F &  A' & 0 \\
0 &  F' & 0
\end{pmatrix}\in\mathcal{P}$, as desired.

\underline{Step~4:} Define $\mathfrak{P}=\bigoplus\limits_{\gamma\in\Gamma}\mathfrak{P}_\gamma$.
Then $\mathfrak{P}$ is a graded ideal of $\mathcal{R}_\Sigma$ by Steps~1--3. Moreover
$\mathfrak{P}\neq\mathcal{R}_\Sigma$ because $[1,1,1,e,e]$, the identity element of
$\mathcal{R}_\Sigma$, does not belong to $\mathfrak{P}$ by Lemma~\ref{lem:grprimematrixideal}(8), since the $2\times 2$ matrix
$\begin{pmatrix}
1&1\\1&0
\end{pmatrix}$ is invertible.

\underline{Step~5:} $\mathcal{R}_\Sigma$ is a $\Gamma$-graded local ring with
graded maximal ideal $\mathfrak{P}$.

Let $\varphi\colon R\rightarrow R_\mathcal{P}$ be the universal localization at $\Sigma$. By
(the proof of) Proposition~\ref{prop:universallocalization}, the isomorphism
$\Phi\colon \mathcal{R}_\Sigma\rightarrow R_\mathcal{P}$ sends
$[F,A,X,\alpha,\beta]\in(\mathcal{R}_\Sigma)_\gamma$ to
$F^\varphi (A^\varphi)^{-1}X^\varphi\in (R_\mathfrak{P})_\gamma$. 

Let $[F,A,X,\alpha,\beta]\in(\mathcal{R}_\Sigma)_\gamma\setminus \mathfrak{P}_\gamma$. Thus
$\begin{pmatrix}
F & X \\ A & 0
\end{pmatrix}\notin\mathcal{P}$ and $\begin{pmatrix}
F & X \\ A & 0
\end{pmatrix}^\varphi$ is invertible in $R_\mathcal{P}$. Also the matrices
$\begin{pmatrix}
A^\varphi & 0 \\ 0 & 1
\end{pmatrix}$ and $\begin{pmatrix}
I & 0 \\ -F^\varphi & 1
\end{pmatrix}$ are invertible in $R_\mathcal{P}$. Hence
\[ \begin{pmatrix}
A^\varphi & 0 \\ 0 & 1
\end{pmatrix} \begin{pmatrix}
I & 0 \\ -F^\varphi & 1
\end{pmatrix} \begin{pmatrix}
F & X \\ A & 0
\end{pmatrix}^\varphi = \begin{pmatrix}
I & (A^\varphi)^{-1}X \\ 0 & -F^\varphi (A^\varphi)^{-1} X^\varphi
\end{pmatrix} \]
is invertible in $R_\mathcal{P}$. Thus, the element
$F^\varphi(A^\varphi)^{-1}X^\varphi$ is invertible in $R_\mathcal{P}$, and therefore
$[F,A,X,\alpha,\beta]$ is invertible in $\mathcal{R}_\Sigma$.

\underline{Step~6:}
Set $\mathcal{K}=\mathcal{R}_\Sigma/\mathfrak{P}$ and let
 $\Psi\colon R\rightarrow \mathcal{K}$ be the composition of $\mu\colon R\rightarrow \mathcal{R}_\Sigma$
 with the natural projection $\mathcal{R}_\Sigma\rightarrow \mathcal{K}$, $[F,A,X,\alpha,\beta]\mapsto \overline{[F,A,X,\alpha,\beta]}$. Then
 $(\mathcal{K},\Psi)$ is a $\Gamma$-graded epic $R$-division ring by Proposition~\ref{prop:basicsonuniversallocalization}(2) and Proposition~\ref{prop:epimorphism=divisionclosure}.

\underline{Step~7:}
For $x\in\h(R)$, $\Psi(x)=0$ if, and only if, the $1\times 1$ matrix $x\in\mathcal{P}$.

Indeed, \[
\Psi(x)=\overline{[x,1,1,e,e]}=0 \Leftrightarrow [x,1,1,e,e]\in\mathfrak{P}
\Leftrightarrow \begin{pmatrix}
1&1\\ x&0
\end{pmatrix}\in \mathcal{P} \Leftrightarrow \begin{pmatrix}
1&1\\ 0&x
\end{pmatrix}\in\mathcal{P}.\]
By Lemma~\ref{lem:grprimematrixideal}(4) and (PM4) the last condition is equivalent to $x\in\mathcal{P}$.

\underline{Step~8:}
For $r\in R_\gamma$ and $[F,A,X,\alpha,\beta]\in (\mathcal{R}_\Sigma)_\gamma$,
$\Psi(r)=[F,A,X,\alpha,\beta]$ if and only if
$\begin{pmatrix}
A&X\\ F&r
\end{pmatrix}\in\mathcal{P}$.

First notice that $\Psi(r)=\overline{[F,A,X,\alpha,\beta]}$ if and only if
$\overline{[r,1,1,e,e]}=\overline{[F,A,X,\alpha,\beta]}$. Equivalently
\[
[F,A,X,\alpha,\beta]-[r,1,1,e,e]=\left[\begin{pmatrix}F&-r\end{pmatrix},\begin{pmatrix}A&0\\0&1\end{pmatrix},\begin{pmatrix}X\\1\end{pmatrix},
\alpha*e,\beta*e\right]\in\mathfrak{P}_\gamma,
\]
which means that $\begin{pmatrix}
A&0&X\\0&1&1\\F&-r&0
\end{pmatrix}\in\mathcal{P}$. This matrix belongs to $\mathcal{P}$ if and only if
\begin{equation}\label{eq:step8}
    \begin{pmatrix}
    A&0&X\\ 0&1&0\\ F&-r&r
    \end{pmatrix}\in\mathcal{P}
\end{equation}
by Lemma~\ref{lem:grprimematrixideal}(2), since the last matrix is obtained after substracting
the second-last column to the last one.
After permuting rows and columns, we get that the matrix \eqref{eq:step8} belongs to
$\mathcal{P}$ if and only if $\begin{pmatrix}
A & X & 0 \\ F & r & -r \\ 0 & 0 & 1
\end{pmatrix}\in\mathcal{P}$. By Lemma~\ref{lem:grprimematrixideal}(4), (PM3), (PM4) and
(PM5) this is equivalent to $\begin{pmatrix}
A & X \\ F&r
\end{pmatrix}\in\mathcal{P}$.

\underline{Step~9:} The singular kernel of $(\mathcal{K},\Psi)$ is $\mathcal{P}$.

Since $\mathcal{R}_\Sigma$ is $R$-isomorphic to $R_\mathcal{P}$, we get that any matrix in
$\mathfrak{M}(R)\setminus\mathcal{P}$ is inverted in $\mathcal{R}_\Sigma$ via $\mu$ and
therefore in $\mathcal{K}$ via $\Psi$.

It remains to show that the matrices in $\mathcal{P}$ do not become invertible via $\Psi$. Let 
$A\in\mathcal{P}$ of size $n\times n$. If all square submatrices (obtained by eliminating equal number of rows and columns)
of $A$ belong to $\mathcal{P}$, then, in particular, all entries of $A$ belong to $\mathcal{P}$.
By Step~7, $A^\Psi$ is the zero matrix and therefore not invertible. Hence suppose that 
$A_1$, of size $m\times m$ with $m>1$, is a submatrix of $A$ of largest size that does not belong to $\mathcal{P}$.
We will show that $A^\Psi$ is singular in $\mathcal{K}$ by expressing one column of $A^\Psi$
as a homogeneous linear combination of the others. 

Since rearrangement of rows and columns does not affect the singularity of $A$, we may suppose that
$$A=\begin{pmatrix}
A_1 & A_2 & A_3 \\ A_4 & A_5 & A_6
\end{pmatrix}$$ where $\begin{pmatrix}
A_2\\ A_5
\end{pmatrix}$ is a column of $A$ and $A\in M_n(R)[\overline{\alpha_1}*\overline{\alpha_2}][\overline{\beta_1}*e*\overline{\beta_2}]$. First observe
that for every $j\in\{1,\dotsc,n\}$
\[
\begin{pmatrix}
A_1 & A_2\\ E_j^T\begin{pmatrix}A_1\\ A_4\end{pmatrix} & E_j^T\begin{pmatrix}A_2\\ A_5\end{pmatrix}
\end{pmatrix}\in\mathcal{P}
\]
because if $j\leq m$, then it has a repeated row and if $j>m$ then it is a submatrix of $A$ of greater
size than $A$.
By Step~8, it means that
\[
\begin{pmatrix}
E_j^T\begin{pmatrix}A_2 \\ A_5 \end{pmatrix}
\end{pmatrix}^\Psi= \overline{\left[ E_j^T\begin{pmatrix}A_1\\A_4\end{pmatrix}, A_1,A_2,\alpha_1,\beta_1 \right]}
\textrm{ for } j=1,\dotsc,n\]
Hence, if we suppose that $\begin{pmatrix}A_1\\A_4\end{pmatrix}=(b_{kl})$,

\begin{eqnarray*}
\begin{pmatrix}
A_2 \\ A_5
\end{pmatrix}^\Psi &=& \begin{pmatrix} 
E_1^T\begin{pmatrix} A_2\\A_5\end{pmatrix} \vspace{0.1cm} \\
 E_2^T\begin{pmatrix} A_2\\ A_5 \end{pmatrix} \\
\vdots \\
 E_n^T\begin{pmatrix} A_2\\ A_5 \end{pmatrix}
\end{pmatrix}^\Psi = 
\begin{pmatrix}
\overline{\left[ E_1^T\begin{pmatrix}A_1\\A_4\end{pmatrix}, A_1,A_2,\alpha_1,\beta_1 \right]}\vspace{0.1cm} \\
\overline{\left[ E_2^T\begin{pmatrix}A_1\\A_4\end{pmatrix}, A_1,A_2,\alpha_1,\beta_1 \right]} \\
\vdots \\
\overline{\left[ E_n^T\begin{pmatrix}A_1\\A_4\end{pmatrix}, A_1,A_2,\alpha_1,\beta_1 \right]}
\end{pmatrix} \\
& = & 
\begin{pmatrix}
\overline{[b_{11},1,1,e,e]} & \dots & \overline{[b_{1m},1,1,e,e]} \\
\vdots & \ddots & \vdots \\
\overline{[b_{n1},1,1,e,e]} & \dots & \overline{[b_{nm},1,1,e,e]}
\end{pmatrix} \cdot
\begin{pmatrix}
\overline{\left[ E_1^T, A_1,A_2,\alpha_1,\beta_1 \right]}\vspace{0.1cm} \\
\overline{\left[ E_2^T, A_1,A_2,\alpha_1,\beta_1 \right]} \\
\vdots \\
\overline{\left[ E_m^T, A_1,A_2,\alpha_1,\beta_1 \right]}
\end{pmatrix} \\
&=&
\begin{pmatrix}
A_1 \\ A_4
\end{pmatrix}^\Psi \cdot  \begin{pmatrix}
\overline{\left[ E_1^T, A_1,A_2,\alpha_1,\beta_1 \right]}\vspace{0.1cm} \\
\overline{\left[ E_2^T, A_1,A_2,\alpha_1,\beta_1 \right]} \\
\vdots \\
\overline{\left[ E_m^T, A_1,A_2,\alpha_1,\beta_1 \right]}
\end{pmatrix} 
\end{eqnarray*}
where we have used Lemma~\ref{lem:easywelldefined}(3) in the second equality.
The result  follows noting that 
$\begin{pmatrix}
A_1 \\ A_4
\end{pmatrix}\in M_{n\times m}(R)[\overline{\alpha_1}*\overline{\alpha_2}][\overline{\beta_1}]$
and $\begin{pmatrix}
\overline{\left[ E_1^T, A_1,A_2,\alpha_1,\beta_1 \right]}\vspace{0.1cm} \\
\overline{\left[ E_2^T, A_1,A_2,\alpha_1,\beta_1 \right]} \\
\vdots \\
\overline{\left[ E_m^T, A_1,A_2,\alpha_1,\beta_1 \right]}
\end{pmatrix}\in M_{m\times 1}(\mathcal{K})[\overline{\beta_1}][e]$.
\end{proof}

The following is Theorem~\ref{theo:specialization}, but expressed in terms of gr-prime matrix ideals.

\begin{corollary}\label{coro:specialization}
Let $R$ be a $\Gamma$-graded ring, $(K_i,\varphi_i)$, $i=1,2$, be $\Gamma$-graded epic
$R$-division rings with singular kernels $\mathcal{P}_i$, respectively. 
The following statements are equivalent.
\begin{enumerate}[\rm(1)]
	\item There exists a gr-specialization from $K_1$ to $K_2$.
	\item $\mathcal{P}_1\subseteq \mathcal{P}_2$.
	\item There exists a homomorphism $R_{\mathcal{P}_2}\rightarrow R_{\mathcal{P}_1}$
	of $\Gamma$-graded $R$-rings.
\end{enumerate}
Furthermore, if there exists a gr-specialization from $K_1$ to $K_2$ and another
gr-specialization from $K_2$ to $K_1$, then $K_1$ and $K_2$ are isomorphic graded $R$-rings. \qed
\end{corollary}

The following corollaries are the graded version of the results in \cite[page~442]{Cohnfreeeidealringslocalization}.

\begin{corollary}
Let  $R$ be a $\Gamma$-graded ring and $(K,\varphi)$ be
a graded epic  $R$-division ring with singular kernel $\mathcal{P}$. Suppose that $\gamma\in\Gamma$. 
Consider the
universal localization $\lambda\colon R\rightarrow R_\mathcal{P}$ and let 
$\Phi\colon R_\mathcal{P} \rightarrow K$ be the homomorphism of $\Gamma$-graded rings such that
$\varphi=\Phi\lambda$. 
\begin{enumerate}[\rm(1)]
	\item Let $x\in K_\gamma$. Then $x=0$ if and only if its numerator belongs to $\mathcal{P}$.
	\item Let $x\in (R_\mathcal{P})_\gamma$. Then $x\in\ker\Phi$ if and only if
	its numerator belongs to $\mathcal{P}$.
\end{enumerate}
\end{corollary}

\begin{proof}
Suppose that $(A_0\, A_\bullet)$ is the numerator of $x$.

(1) By Lemma~\ref{lem:Cramersrule}(1), $x$ is invertible if and only if $(A_0\, A_\bullet)^\varphi$
is invertible over $K$. That is, if and only if $(A_0\, A_\bullet)$ belongs to $\mathcal{P}$.

(1) By Lemma~\ref{lem:Cramersrule}(1), $x$ is invertible if and only if $(A_0\, A_\bullet)^\lambda$
is invertible over $R_\mathcal{P}$. Since $R_\mathcal{P}$
is a local ring with residue class graded division ring $R$-isomorphic to $K$, $x$ is invertible if and only if $(A_0\, A_\bullet)^{\Phi\lambda}$
is invertible over $K$.  That is, $x\in\ker \Phi$ if and only if $(A_0\, A_\bullet)$  belongs to $\mathcal{P}$.
\end{proof}

\begin{corollary}
 Let $R$ and $R'$ be $\Gamma$-graded rings with gr-prime matrix ideals
$\mathcal{P}$ and $\mathcal{P}'$, respectively, with corresponding graded epic $R$-division rings
$(K,\varphi)$ and $(K',\varphi')$ respectively. Let $f\colon R\rightarrow R'$
be a homomorphism of $\Gamma$-graded rings. The following assertions hold true.
\begin{enumerate}[\rm(1)]
	\item  $f$ extends to a gr-specialization if, and only if, $\mathcal{P}^f\subseteq\mathcal{P}'$.
	\item $f$ extends to a homomorphism $K\rightarrow K'$ if, and only if, 
	$\mathcal{P}^f\subseteq \mathcal{P}'$ and  $\Sigma^f\subseteq \Sigma'$, where $\Sigma=\mathfrak{M}(R)\setminus \mathcal{P}$ and 
	$\Sigma'=\mathfrak{M}(R')\setminus\mathcal{P}'$.
\end{enumerate} 
\end{corollary}

\begin{proof}
(1) First note that the set ${\mathcal{P}''}=\{A\in\mathfrak{M}(R)
\colon A^{f}\in\mathcal{P}'\}$ is a gr-prime matrix ideal whose corresponding graded epic $R$-division ring is
$\varphi'f\colon R\rightarrow \DC(\varphi'f)$. 

By Corollary~\ref{coro:specialization}, there exists a specialization from $(K,\varphi)$
to $(\DC(\varphi'f))$ if, and only if, $\mathcal{P}\subseteq \mathcal{P}''$.

(2) If $\mathcal{P}^f\subseteq \mathcal{P}'$ and $\Sigma^f\subseteq \Sigma'$, then $\mathcal{P}=\mathcal{P}''$, and therefore the gr-specialization
of (1) is in fact an isomorphism by Corollary~\ref{coro:specialization}. 
\end{proof}


\section{Gr-matrix ideals}\label{sec:grmatrixideals}

In this section, the concepts, arguments and proofs are an adaptation
of the ones in \cite[Section~7.3]{Cohnfreeeidealringslocalization} to the graded context.

\medskip

\emph{Throughout this section, let $\Gamma$ be a group}. 

\medskip

Let $R$ be a $\Gamma$-graded ring. 
A subset $\mathcal{I}$ of $\mathfrak{M}(R)$ is a \emph{gr-matrix pre-ideal} if the following
conditions are satisfied.

\begin{enumerate}[({I}1)]
	\item $\mathcal{I}$ contains all the homogeneous matrices that are not gr-full;
	\item If $A,B\in\mathcal{I}$ and their determinantal sum (with respect to
	a row or column) exists, then
	$A\nabla B\in\mathcal{I}$;
	\item If $A\in\mathcal{I}$, then $A\oplus B\in\mathcal{I}$ for all
	$B\in\mathfrak{M}(R)$;
	\item If $A\in\mathcal{I}$ and $E,F$ are permutation matrices of appropriate size, then 
	$EAF\in\mathcal{I}$.
\end{enumerate}

\medskip

If, moreover, we have
\begin{enumerate}[({I}5)]
	\item For $A\in\mathfrak{M}(R)$, if $A\oplus 1\in\mathcal{I}$, then $A\in\mathcal{I}$,
\end{enumerate}
we call $\mathcal{I}$ a \emph{gr-matrix ideal}.

\bigskip

Clearly, $\mathfrak{M}(R)$ is a gr-matrix ideal. A \emph{proper} gr-matrix ideal
is a gr-matrix ideal different from $\mathfrak{M}(R)$.

\begin{lemma}\label{lem:grmatrixideal}
Let  $R$ be a $\Gamma$-graded ring and $\mathcal{I}$ be a gr-matrix pre-ideal. 
 Let $A,B\in \mathfrak{M}(R)$. The following assertions hold true.
\begin{enumerate}[\rm(1)]
\item  If $A$ and $B$ are such that
$C=A\nabla B$ exists and $B$ is not gr-full. Then $A\in\mathcal{I}$ if and only if $C\in\mathcal{I}$.
\item Let $A\in\mathcal{I}$. The result of adding a suitable  right multiple of one column of $A$ to another column again
lies in $\mathcal{I}$. More precisely, if $A\in M_n(R)[\overline{\alpha}][\overline{\beta}]$ and $a\in 
R_{\beta_i\beta_j^{-1}}$,
then $(A_1\,\dotsc\, A_{j-1}\ A_j+A_ia\ A_{j+1}\,\dotsc\, A_n)\,$ belongs to $\mathcal{I}$.
\item If $A\oplus B\in\mathcal{I}$, then $B\oplus A\in\mathcal{I}$.
\item Suppose that $A\in M_m(R)[\overline{\alpha}][\overline{\beta}]$ and $B\in M_n(R)[\overline{\delta}][\overline{\varepsilon}]$.
For $C\in M_{n\times m}(R)[\overline{\delta}][\overline{\beta}]$, 
\[
\begin{pmatrix}
A&0\\C&B
\end{pmatrix}\in\mathcal{I}\ \textrm{ if and only if } \
\begin{pmatrix}
A&0\\0&B
\end{pmatrix}\in\mathcal{I}
\]
Similarly, for $C\in M_{m\times n}(R)[\overline{\beta}][\overline{\varepsilon}]$,
\[
\begin{pmatrix}
A&C\\0&B
\end{pmatrix}\in\mathcal{I}\ \textrm{ if and only if } \
\begin{pmatrix}
A&0\\0&B
\end{pmatrix}\in\mathcal{I}
\]
\end{enumerate}
If, moreover, $\mathcal{I}$ is a gr-matrix ideal, then the following assertions hold true.
\begin{enumerate}[\rm(1)]
\setcounter{enumi}{4}
\item Suppose that $A\in M_n(R)[\overline{\alpha}][\overline{\beta}]$, $B\in M_n(R)[\overline{\beta}][\overline{\delta}]$. Then $AB\in\mathcal{I}$ if and only if $A\oplus B\in\mathcal{I}$.
\item If $A$ and $B$ are such that
 $C=A\nabla B$ exists and $B\in\mathcal{I}$. Then $A\in\mathcal{I}$ if and only if $C\in\mathcal{P}$.
\item If an identity matrix $I_n$, $n\geq1$, belongs to $\mathcal{I}$, then $\mathcal{I}=\mathfrak{M}(R)$
\end{enumerate}
\end{lemma}

\begin{proof}
Note that (I1),(I2),(I3), (I4) are the same as (PM1),(PM2),(PM3),(PM6). Hence (1)--(6)
follow in exactly the same way as in Lemma~\ref{lem:grprimematrixideal}.

To prove (7), note that if $I_n\in\mathcal{I}$, for some $n\geq 1$, an application of (I5), shows that
$1\times 1$ matrix $1\in\mathcal{I}$. By (I3), any identity matrix $I_m$, $m\geq 1$, belongs to
$\mathcal{I}$. Again, using (I3), $I_m\oplus A\in\mathcal{I}$ for any positive integer $m$ and
matrix $A\in\mathfrak{M}(R)$. By (5), any $A\in\mathfrak{M}(R)$ belongs to $\mathcal{I}$,
as desired.
\end{proof}

One could think of defining  a gr-prime matrix ideal as a gr-matrix ideal $\mathcal{I}$ 
such that the following two
conditions are satisfied.
\begin{enumerate}[({I}1)]
\setcounter{enumi}{5}
\item $\mathcal{I}$ is a proper gr-matrix ideal. 
\item $\mathcal{I}$ satisfies (PM4).
\end{enumerate}
We proceed to show that both definitions are equivalent. Let $\mathcal{P}$
be a gr-prime matrix ideal, i.e. (PM1)--(PM6) in page~\pageref{def:grprimematrixideal} are satisfied.  Clearly, $\mathcal{P}$
satisfies (I1)--(I4) and (I7). By (PM5), $1\notin\mathcal{P}$. Therefore, by (PM4), 
if $A\oplus 1\in\mathcal{P}$, then $A\in\mathcal{P}$ for any $A\in\mathfrak{M}(R)$. Hence
(I5) is satisfied. Again by (PM5), $\mathcal{P}$ is a proper gr-matrix ideal.
Conversely, suppose that $\mathcal{I}$ satisfies (I1)--(I7). Clearly (PM1)--(PM4), (PM6)
are satisfied. By Lemma~\ref{lem:grmatrixideal}(7) and (I6), (PM5) is satisfied, as desired.

\bigskip

It is easy to prove that any intersection of gr-matrix (pre-)ideals is again a gr-matrix (pre-)ideal.
Thus, given a subset $\mathcal{S}\subseteq \mathfrak{M}(R)$, we define the \emph{gr-matrix (pre-)ideal
generated by $\mathcal{S}$} as the intersection of gr-matrix (pre-)ideals $\mathcal{I}$ that contain $\mathcal{S}$.
That is, $\bigcap_{\mathcal{S}\subseteq\mathcal{I}}\mathcal{I}$. Note that this gr-matrix (pre)-ideal
is contained in any gr-matrix (pre-)ideal that contains $\mathcal{S}$.

Now we fix some notation that will be used in what follows. 

Let $\mathcal{W}\subseteq\mathfrak{M}(R)$. We say that a matrix
$C\in\mathfrak{M}(R)$ is a determinantal sum of elements of $\mathcal{W}$ if there exist $A_1,\dotsc,A_m\in\mathcal{W}$,
$m\geq 1$, such that $A_1\nabla A_2\nabla\dotsc\nabla A_m$ exists for some choice of parenthesis
and equals $C$.

We will write $\mathcal{N}$ to denote the subset of $\mathfrak{M}(R)$ 
consisting of the matrices which are not gr-full.

We will denote the set of all identity matrices by $\mathfrak{I}$.

If $\mathcal{X}\subseteq\mathfrak{M}(R)$, we denote by $\mathcal{D}(\mathcal{X})$ the set of 
all  matrices in $\mathfrak{M}(R)$ which are of the form $E(X\oplus A)F$ where $X\in\mathcal{X}$,
$A\in\mathfrak{M}(R)$ and $E,F$ are permutation matrices of appropriate sizes. We remark that
we allow $A$ to be the empty matrix $\mathbb{O}$.

\begin{lemma}\label{lem:idealgenerated}
Let $R$ be a $\Gamma$-graded ring and $\mathcal{A}$ be a gr-matrix pre-ideal. Suppose that $\Sigma\subset 
\mathfrak{M}(R)$ satisfies the following two conditions
\begin{enumerate}[\rm(i)]
\item $1\in\Sigma$;
\item if $P,Q\in\Sigma$, then $P\oplus Q\in\Sigma$.
\end{enumerate}
Then the following assertions hold true
\begin{enumerate}[\rm(1)]
\item The set $\mathcal{A}/\Sigma\coloneqq\{A\in\mathfrak{M}(R)\colon A\oplus P\in\mathcal{A}
\textrm{ for some } P\in\Sigma\}$
is a gr-matrix ideal containing $\mathcal{A}$.
\item The gr-matrix ideal $\mathcal{A}/\Sigma$ is proper if and only if $\mathcal{A}\cap\Sigma=\emptyset$.
\item The gr-matrix ideal $\mathcal{A}/\mathfrak{I}$ is the gr-matrix ideal generated by $\mathcal{A}$.
\end{enumerate}
\end{lemma}

\begin{proof}
(1) Let $A\in\mathcal{A}$. By (I3), $A\oplus 1\in\mathcal{A}$. Since $1\in\Sigma$, $A\in\mathcal{A}/\Sigma$.
Hence $\mathcal{A}\subseteq \mathcal{A}/\Sigma$ and, by (I1), all non gr-full matrices belong to
$\mathcal{A}$. Therefore $\mathcal{A}/\Sigma$ satisfies (I1).

Let $A,B\in\mathcal{A}/\Sigma$ be such that $A\nabla B$ is well defined. There exist $P,Q\in\Sigma$
such that $A\oplus P,  B\oplus Q\in \mathcal{A}.$ By (I3), $A\oplus P\oplus Q$ and
$B\oplus Q\oplus P$ belong to $\mathcal{A}$. By (I4), $B\oplus P\oplus Q\in\mathcal{A}$. Now
$(A\nabla B)\oplus P\oplus Q=(A\oplus P\oplus Q)\nabla(B\oplus P\oplus Q)\in\mathcal{A}$ by
(I2). Hence $A\nabla B\in\mathcal{A}/\Sigma$ and $\mathcal{A}/\Sigma$ satisfies (I2).

Let $A\in\mathcal{A}/\Sigma$ and $B\in\mathfrak{M}(R)$. There exists $P\in\Sigma$ such that
$A\oplus P\in\mathcal{A}$. By (I3), $A\oplus P\oplus B\in\mathcal{A}$. Now (I4) implies that
$A\oplus B\oplus P\in\mathcal{A}$. Hence $A\oplus B\in\mathcal{A}/\Sigma$, and
$\mathcal{A}/\Sigma$ satisfies (I3).

Let $A\in\mathcal{A}/\Sigma$ and $E,F$ be permutation matrices of the same size as $A$.
There exists $P\in\Sigma$ such that $A\oplus P\in\mathcal{A}$. 
Since $E\oplus I$ and $F\oplus I$ are also permutation matrices, (I4) implies that
$(E\oplus I)(A\oplus P)(F\oplus I)=EAF\oplus P\in\mathcal{A}$. Hence $EAF\in\mathcal{A}/\Sigma$ and
(I4) is satisfied.

Let now $A\in\mathfrak{M}(R)$ such that $A\oplus 1\in\mathcal{A}/\Sigma$. Thus there exists
$P\in \Sigma$ such that $A\oplus 1\oplus P\in\mathcal{A}$. Since $1\oplus P\in\Sigma$,
then $A\in\mathcal{A}/\Sigma$ and $A/\Sigma$ satisfies (I5).

(2) Suppose that $\mathcal{A}\cap\Sigma\neq\emptyset$. Let $P\in\mathcal{A}\cap\Sigma$
and $M\in\mathfrak{M}(R)$. Then $P\oplus M\in\mathcal{A}$ by (I3). By (I4),
$M\oplus P\in\mathcal{A}$. Hence $M\in\mathcal{A}/\Sigma$. Therefore,
$\mathcal{A}/\Sigma=\mathfrak{M}(R)$.

Conversely, suppose that $\mathcal{A}/\Sigma=\mathfrak{M}(R)$. Thus, $1\in\mathcal{A}/\Sigma$
and there exists $P\in\Sigma$ such that $1\oplus P\in\mathcal{A}$.
Notice that $1\oplus P\in\Sigma$, by (i) and (ii). Therefore $\mathcal{A}\cap\Sigma\neq \emptyset$.

(3) Clearly $\mathfrak{I}$ satisfies conditions (i) and (ii). Thus $\mathcal{A}/\mathfrak{I}$
is a gr-matrix ideal that contains $\mathcal{A}$ by (1). Let now $\mathcal{B}$ be a
gr-matrix ideal such that $\mathcal{A}\subseteq \mathcal{B}$. If $A\in\mathcal{A}/\mathfrak{I}$,
then there exists $n\geq 1$ such that $A\oplus I_n\in\mathcal{A}\subseteq \mathcal{B}$.
By applying (I5) repeteadly, we obtain that $A\in\mathcal{B}$, as desired. 
\end{proof}

\begin{lemma}\label{lem:idealgenerated2}
Let $R$ be a $\Gamma$-graded ring and let $\mathcal{X}\subseteq\mathfrak{M}(R)$.
Let $\mathcal{A}(\mathcal{X})$ be the subset of $\mathfrak{M}(R)$ consisting of all the matrices
that can be expressed as determinantal sum of elements of $\mathcal{N}\cup\mathcal{D}(\mathcal{X})$.
The following assertions hold true.
\begin{enumerate}[\rm(1)]
	\item $\mathcal{A}(\mathcal{X})$ is the  gr-matrix pre-ideal generated by $\mathcal{X}$.
	\item $\mathcal{A}(\mathcal{X})/\mathfrak{I}$ is the  gr-matrix ideal generated by $\mathcal{X}$.
	\item The gr-matrix ideal generated by $\mathcal{X}$ is proper if and only if 
	$ \mathcal{A}(\mathcal{X})\cap\mathfrak{I}=\emptyset$.
\end{enumerate}
\end{lemma}

\begin{proof}
(1) $\mathcal{X}\subseteq \mathcal{A}(\mathcal{X})$ because $X=I(X\oplus\mathbb{O})I$
for all $X\in\mathcal{X}$. 
By definition of $\mathcal{A}(\mathcal{X})$, every homogeneous matrix that is not gr-full
belongs to $\mathcal{A}(\mathcal{X})$. By the same reason, if $A,B\in\mathcal{A}(\mathcal{X})$ and
$A\nabla B$ is defined, then $A\nabla B\in\mathcal{A}(\mathcal{X})$.

Let $A\in \mathcal{A}(\mathcal{X})$ and $B\in\mathfrak{M}(R)$. 
That $A\oplus B\in\mathcal{A}(\mathcal{X})$ follows from the following three facts. First,
for any $U,V\in\mathfrak{M}(R)$, when defined
$(U\nabla V)\oplus M=(U\oplus M)\nabla (V\oplus M)$. Second, for $X\in\mathcal{X}$ and
$U,M\in\mathfrak{M}(R)$ and permutation matrices $E,F$ of suitable size
$E(X\oplus U)F\oplus M=(E\oplus I)(X\oplus U\oplus M)(F\oplus I)$. Third,
if $U$ is not gr-full, then, for all $M\in\mathfrak{M}(R)$, $U\oplus M$ is not 
full for all $M\in\mathfrak{M}(R)$. Indeed, if $U=U_1U_2$, then
$U\oplus M=(U_1\oplus M)(U_2\oplus I)$.

If $A\in\mathcal{A}(\mathcal{X})$ and $E,F$ are permutation matrices of appropriate size, then
$EAF\in\mathcal{A}(\mathcal{X})$. This follows from the following facts. First, if
$A,B\in\mathfrak{M}(R)$ and $E,F$ are permutation matrices such that
$E(A\nabla B)F$ is defined, then $E(A\nabla B)F=EAF\nabla EBF$. Second, for $X\in\mathcal{X}$,
$U\in\mathfrak{M}(R)$ and permutation matrices $E,F,P,Q$ of appropriate sizes then
$P(E(X\oplus U)F)Q=(PE)(X\oplus U)(FQ)$. Third, if $U\in\mathfrak{M}(R)$ is not gr-full, and
$E,F$ are permutation matrices of appropriate size, then $EUF$ is not gr-full. Indeed,
if $U=U_1U_2$, then $EUF=(EU_1)(U_2F)$. 

Therefore, $\mathcal{A}(\mathcal{X})$ is a gr-matrix pre-ideal that contains $\mathcal{X}$.

Let now $\mathcal{B}$ be a gr-matrix pre-ideal such that $\mathcal{X}\subseteq\mathcal{B}$.
By (I1), $\mathcal{N}\subseteq\mathcal{B}$. By (I3) and (I4), $E(X\oplus A)F\in\mathcal{B}$
for all $X\in\mathcal{X}$, $A\in\mathfrak{M}(R)$ and permutation matrices $E,F$ of appropriate size.
By (I2), $\mathcal{A}(\mathcal{X})\subseteq\mathcal{B}$.

(2) Any gr-matrix ideal containing $\mathcal{X}$, must contain $\mathcal{A}(\mathcal{X})$.
By Lemma~\ref{lem:idealgenerated}(3), the result follows.

(3) By (2), the gr-matrix ideal generated by $\mathcal{X}$ equals
$\mathcal{A}(\mathcal{X})/\mathfrak{I}$. By Lemma~\ref{lem:idealgenerated}(2), 
$\mathcal{A}(\mathcal{X})/\mathfrak{I}$ is proper if and only if 
$\mathcal{A}(\mathcal{X})\cap\mathfrak{I}=\emptyset$.
\end{proof}

\begin{corollary}\label{coro:leastideal}
Let $R$ be a $\Gamma$ graded ring. The set
$\mathcal{A}(\mathcal{N})/\mathfrak{I}$ is the least gr-matrix ideal. Hence 
 $R$ has proper gr-matrix ideals if and only if no matrix of
$\mathfrak{I}$ can be expressed as a determinantal sum of matrices of $\mathcal{N}$.
\end{corollary}

\begin{proof}
The set $\mathcal{N}$ is contained in each gr-matrix ideal. By Lemma~\ref{lem:idealgenerated2}(2),
$\mathcal{A}(\mathcal{N})/\mathfrak{I}$ is the  gr-matrix ideal generated by $\mathcal{N}$.
Thus all gr-matrix ideals contain the gr-matrix ideal $\mathcal{A}(\mathcal{N})/\mathfrak{I}$.

Since any matrix in $\mathfrak{M}(R)$ of the form   
 $E(X\oplus A)F$ where $X\in\mathcal{\mathcal{N}}$,
$A\in\mathfrak{M}(R)$ and $E,F$ are permutation matrices of appropriate sizes,
again belongs to $\mathcal{N}$, then $\mathcal{D}(\mathcal{N})=\mathcal{N}$. Thus,
$\mathcal{A}(\mathcal{N})$ consists of the  matrices in $\mathfrak{M}(R)$ that can be expressed as
a determinantal sum of matrices from $\mathcal{N}$.

Now $R$ has proper gr-matrix ideals if and only if $\mathcal{A}(\mathcal{N})/\mathfrak{I}$
is proper. By Lemma~\ref{lem:idealgenerated}(3), this is equivalent to 
$\mathcal{A}(\mathcal{N})\cap\mathfrak{I}=\emptyset$. In other words, no matrix of
$\mathfrak{I}$ can be expressed as a determinantal sum of matrices of $\mathcal{N}$. 
\end{proof}

\begin{lemma}\label{lem:quotientideal}
Let $R$ be a $\Gamma$-graded ring, $\mathcal{I}$ be a gr-matrix ideal and 
$\mathcal{Z}\subseteq\mathfrak{M}(R)$.
Then the set $\mathcal{I}_\mathcal{Z}=\{A\in\mathfrak{M}(R)\colon A\oplus Z\in\mathcal{I}
\textrm{ for all }Z\in\mathcal{Z}\}$ is
a gr-matrix ideal.
\end{lemma}

\begin{proof}
Let $A\in\mathfrak{M}(R)$ and suppose it is not gr-full. If $A=BC$, then
$A\oplus Z=(B\oplus Z)(C\oplus I)$ for all
$Z\in\mathcal{Z}$. Thus $A\in\mathcal{I}_\mathcal{Z}$ and (I1) is satisfied.

Let $A,B\in\mathcal{I}_\mathcal{Z}$ and suppose that $A\nabla B$ exists. Then
$(A\nabla B)\oplus Z=(A\oplus Z)\nabla(B\oplus Z)$
for all $Z\in\mathcal{Z}$. Since $A\oplus Z,B\oplus Z\in\mathcal{I}$,
then $(A\nabla B)\oplus Z\in\mathcal{I}$ for all $Z\in\mathcal{Z}$. 
Hence $A\nabla B\in\mathcal{I}_\mathcal{Z}$, and
(I2) is satisfied.

Let $A\in\mathcal{I}_\mathcal{Z}$ and $B\in\mathfrak{M}(R)$. Since $A\oplus Z\in\mathcal{I}$
for all $Z\in\mathcal{Z}$ and $\mathcal{I}$ is a gr-matrix ideal, then
$A\oplus Z\oplus B\in\mathcal{I}$ for all $Z\in\mathcal{Z}$. By (I4),
$A\oplus B\oplus Z\in\mathcal{I}$ for all $Z\in\mathcal{Z}$. Therefore
$A\oplus B\in\mathcal{I}_\mathcal{Z}$ and (I3) is satisfied.

If $A\in\mathcal{I}_\mathcal{Z}$, $Z\in\mathcal{Z}$ and
$E,F$ are permutation matrices of appropriate size, then
$EAF\oplus Z=(E\oplus I)(A\oplus Z)(F\oplus I)$. It shows that
$EAF\in\mathcal{I}_\mathcal{Z}$ and (I4) is satisfied.

Suppose now that $A\in\mathfrak{M}(R)$ and that $A\oplus 1\in\mathcal{I}_\mathcal{Z}$.
Hence $A\oplus 1\oplus Z\in\mathcal{I}$ for all $Z\in\mathcal{Z}$. 
By (I4), $A\oplus Z\oplus 1\in\mathcal{I}$ for all $Z\in\mathcal{Z}$.
Now, by (I5), $A\oplus Z\in\mathcal{I}$ for all $Z\in\mathcal{Z}$,
which shows that $A\in\mathcal{I}_\mathcal{Z}$. Therefore (I5) is satisfied.
\end{proof}

Let $\mathcal{A}_1,\mathcal{A}_2$ be two gr-matrix ideals of a $\Gamma$-graded ring $R$. 
The \emph{product of $\mathcal{A}_1$ and $\mathcal{A}_2$}, denoted by $\mathcal{A}_1\mathcal{A}_2$,
is the gr-matrix ideal generated by the set 
$$\{A_1\oplus A_2\colon A_1\in\mathcal{A}_1, A_2\in\mathcal{A}_2\}.$$

A helpful description of $\mathcal{A}_1\mathcal{A}_2$ is given in the following lemma.

\begin{lemma}
Let $R$ be  $\Gamma$-graded ring and
$\mathcal{X}_1,\mathcal{X}_2\subseteq\mathfrak{M}(R)$. Set 
$$\mathcal{X}=\{X_1\oplus X_2\colon X_1\in\mathcal{X}_1, X_2\in\mathcal{X}_2\}.$$
Let $\mathcal{A}_1$ be the gr-matrix ideal generated by $\mathcal{X}_1$,
$\mathcal{A}_2$ be the gr-matrix ideal generated by $\mathcal{X}_2$ and
$\mathcal{A}$ be the gr-matrix ideal generated by $\mathcal{X}$. Then
$\mathcal{A}=\mathcal{A}_1\mathcal{A}_2$.

As a consequence, for any $A,B\in\mathfrak{M}(R)$, 
$\langle A\rangle\langle B\rangle=\langle A\oplus B\rangle$,
where $\langle A\rangle$ denotes the gr-matrix ideal generated by 
$\{A\}$.
\end{lemma}

\begin{proof}
First, $\mathcal{A}\subseteq \mathcal{A}_1\mathcal{A}_2$ because 
$X_1\oplus X_2\in\mathcal{A}_1\mathcal{A}_2$ for all $X_1\in\mathcal{X}_1,X_2\in\mathcal{X}_2$.

Now observe that $X_1\oplus X_2\in\mathcal{X}\subseteq \mathcal{A}$
for all $X_1\in\mathcal{X}_1$, $X_2\in\mathcal{X}_2$. By (I4),
$X_2\oplus X_1\in\mathcal{X}\subseteq \mathcal{A}$
for all $X_1\in\mathcal{X}_1$, $X_2\in\mathcal{X}_2$. Hence 
$\mathcal{X}_2$ is contained in the gr-matrix ideal $\mathcal{A}_{\mathcal{X}_1}$.
Thus, $\mathcal{A}_2\subseteq \mathcal{A}_{\mathcal{X}_1}$. It implies
that $A_2\oplus X_1\in\mathcal{A}$ for all $A_2\in\mathcal{A}_2$ and $X_1\in\mathcal{X}_1$.
Again by (I4), $X_1\oplus A_2\in\mathcal{A}$ 
for all $A_2\in\mathcal{A}_2$ and $X_1\in\mathcal{X}_1$. Therefore 
$\mathcal{X}_1$ is contained in the gr-matrix ideal $\mathcal{A}_{\mathcal{A}_2}$.
Thus $\mathcal{A}_1\subseteq \mathcal{A}_{\mathcal{A}_2}$. This means that
$A_1\oplus A_2\in\mathcal{A}$ for all $A_1\in\mathcal{A}_1$ and
$A_2\in\mathcal{A}_2$. Therefore $\mathcal{A}_1\mathcal{A}_2\subseteq \mathcal{A}$.
\end{proof}

Now we show that gr-prime matrix ideals behave like graded prime ideals of graded rings.

\begin{proposition}
Let $R$ be a $\Gamma$-graded ring. For a proper gr-matrix ideal $\mathcal{P}$, the following
are equivalent
\begin{enumerate}[\rm(1)]
	\item $\mathcal{P}$ is a gr-prime matrix ideal.
	\item For gr-matrix ideals $\mathcal{A}_1,\mathcal{A}_2$, if 
	$\mathcal{A}_1\mathcal{A}_2\subseteq \mathcal{P}$, then $\mathcal{A}_1\subseteq \mathcal{P}$
	or $\mathcal{A}_2\subseteq \mathcal{P}$.
	\item For gr-matrix ideals  $\mathcal{A}_1,\mathcal{A}_2$ that
	contain $\mathcal{P}$, if 
	$\mathcal{A}_1\mathcal{A}_2\subseteq \mathcal{P}$, then $\mathcal{A}_1=\mathcal{P}$
	or $\mathcal{A}_2=\mathcal{P}$.
\end{enumerate}
\end{proposition}

\begin{proof}
Suppose (1) holds true. Let $\mathcal{A}_1,\mathcal{A}_2$ be gr-matrix ideals such that
$\mathcal{A}_1\nsubseteq\mathcal{P}$ and $\mathcal{A}_2\nsubseteq\mathcal{P}$.
Hence there exist $A_1\in\mathcal{A}_1\setminus\mathcal{P}$ and 
$A_2\in\mathcal{A}_2\setminus\mathcal{P}$. Hence $A_1\oplus A_2\notin\mathcal{P}$.
It implies that $\mathcal{A}_1\mathcal{A}_2\nsubseteq\mathcal{P}$. Therefore (2) holds true.

Clearly (2) implies (3). 

Suppose (3) holds true and let $A_1,A_2\in\mathfrak{M}(R)$ be such that 
$A_1\oplus A_2\in\mathcal{P}$. Let $\mathcal{A}_1,\mathcal{A}_2$ be the gr-matrix ideals 
generated by $\mathcal{P}\cup\{A_1\}$ and $\mathcal{P}\cup\{A_2\}$, respectively.
Notice that  $X_1\oplus X_2\in\mathcal{P}$ for $X_1\in \mathcal{A}_1$,
$X_2\in\mathcal{A}_2$. Hence $\mathcal{A}_1\mathcal{A}_2\subseteq \mathcal{P}$.
By (3), either $\mathcal{A}_1=\mathcal{P}$ or $\mathcal{A}_2=\mathcal{P}$.
Hence $A_1\in\mathcal{P}$ or $A_2\in\mathcal{P}$, and (1) is satisfied.
\end{proof}

Let $\mathcal{A}$ be a gr-matrix ideal. The \emph{radical of $\mathcal{A}$}
is defined as the set
$$\sqrt{\mathcal{A}}=\{A\in\mathfrak{M}(R)\colon \oplus^r\!\! A\in\mathcal{A}
\textrm{ for some positive integer }r\}.$$

We say that a proper gr-matrix ideal $\mathcal{A}$ is \emph{gr-semiprime}
if $\sqrt{\mathcal{A}}=\mathcal{A}$.

\begin{lemma}
Let $R$ be a $\Gamma$-graded ring and let $\mathcal{A}$ be a gr-matrix ideal. The following
assertions hold true.
\begin{enumerate}[\rm(1)]
	\item $\sqrt{\mathcal{A}}$ is a gr-matrix ideal that contains $\mathcal{A}$.
	\item $\sqrt{\sqrt{\mathcal{A}}}=\sqrt{\mathcal{A}}$.
	\item If $\mathcal{A}$ is a gr-prime matrix ideal, then $\sqrt{\mathcal{A}}=\mathcal{A}$.
\end{enumerate}
\end{lemma}

\begin{proof}
(1) If $A\in\mathcal{A}$, then, for $r=1$, we obtain that 
$A=\oplus^1 A\in\mathcal{A}$. Hence $\mathcal{A}\subseteq \sqrt{\mathcal{A}}$.
In particular, all homogeneous matrices which are not gr-full belong to $\sqrt{\mathcal{A}}$.
Thus $\sqrt{\mathcal{A}}$ satisfies (I1).

Let $A,B\in\sqrt{\mathcal{A}}$ such that $A\nabla B$ exists. There exist $r,s\geq 1$
such that $\oplus^r A$, $\oplus^s B\in\mathcal{A}$. Set $n=r+s+1$. To prove that
$\sqrt{\mathcal{A}}$ satisfies (I2), it is enough to show that
$\oplus^n(A\nabla B)\in\mathcal{A}$.
For that aim, using $(A\nabla B)\oplus P=(A\oplus P)\nabla(B\oplus P)$,
one can prove by induction on $n$ that 
$\oplus^n(A\nabla B)$ is a determinantal sum of elements of the form
\begin{equation}\label{eq:radical}
C_1\oplus C_2\oplus\dotsb\oplus C_n
\end{equation}
where each $C_i$ equals $A$ or $B$. By the choice of $n$, there are at least $r$
$C_i$'s equal to $A$ or at least $s$ $C_i$'s equal to $B$. Either case,
there exist permutation matrices $E,F$ of appropriate size such that
$$C_1\oplus C_2\oplus\dotsb\oplus C_n=\left\{ \begin{array}{l} 
E((\oplus^rA)\oplus C_{r+1}'\oplus \dotsb \oplus C_n')F \\
E((\oplus^s B)\oplus C_{s+1}'\oplus \dotsb \oplus C_n')F  \end{array}\right. .$$
It implies that the elements in \eqref{eq:radical} belong to $\mathcal{A}$
by (I3). Now (I2) implies that $\oplus^n(A\nabla B)\in\mathcal{A}$, as desired.

Let now $A\in\sqrt{\mathcal{A}}$ and $B\in\mathfrak{M}(R)$. There exists $r\geq 1$
such that $\oplus^rA\in\mathcal{A}$. The equality
$\oplus^r(A\oplus B)=E((\oplus^rA)\oplus(\oplus^r B))F$ holds for some permutation matrices
$E,F$. Hence $\oplus^r(A\oplus B)\in\mathcal{A}$. Thus $A\oplus B\in\sqrt{\mathcal{A}}$
and $\sqrt{\mathcal{A}}$ satisfies (I3).

Let $A\in\sqrt{\mathcal{A}}$ such that $\oplus^rA\in\mathcal{A}$. For permutation matrices
$E,F$ of appropriate size
$$\oplus^r(EAF)=(\oplus^rE)(\oplus^rA)(\oplus^rF)\in\mathcal{A}.$$
Therefore $EAF\in\sqrt{\mathcal{A}}$ and $\sqrt{\mathcal{A}}$ satisfies (I4).

If $X\in\mathfrak{M}(R)$ is such that $X\oplus 1\in\sqrt{\mathcal{A}}$, then there
exists $t\geq 1$ such that $\oplus^t(X\oplus 1)\in\mathcal{A}$. But now
$(\oplus^t X)\oplus I_t=E(\oplus^t(X\oplus 1))F\in\mathcal{A}$.
Applying (I5), we get that $\oplus^tX\in\mathcal{A}$, and therefore $X\in\sqrt{\mathcal{A}}$.
Hence $\sqrt{\mathcal{A}}$ satisfies (I5).

(2) By (1), $\sqrt{\mathcal{A}}\subseteq \sqrt{\sqrt{\mathcal{A}}}$. Let
now $A\in\sqrt{\sqrt{\mathcal{A}}}$. It means that 
$\oplus^rA\in\sqrt{\mathcal{A}}$ for some positive integer $r$. Hence there exists a positive
integer $s$ such that $\oplus^s(\oplus^r A)\in\mathcal{A}$. Thus,
$\oplus^{rs}A=\oplus^s(\oplus^r A)\in\mathcal{A}$. Therefore $A\in\sqrt{\mathcal{A}}$, as desired.

(3) Suppose $\mathcal{A}$ is a gr-prime matrix ideal and let $A\in\sqrt{\mathcal{A}}$.
Hence $\oplus^rA\in\mathcal{A}$. By (PM4), $A\in\mathcal{A}$, as desired.
\end{proof}

\begin{proposition}\label{prop:maximalisprime}
Let $R$ be a $\Gamma$-graded ring. Suppose that the nonempty subset $\Sigma$ of
$\mathfrak{M}(R)$ and the gr-matrix ideal $\mathcal{A}$ satisfy the following two conditions.
\begin{enumerate}[\rm(i)]
	\item $A\oplus B\in\Sigma$ for all $A,B\in\Sigma$;
	\item $\mathcal{A}\cap\Sigma=\emptyset$.
\end{enumerate}
Then the set $W$ of gr-matrix ideals $\mathcal{B}$ such that
$\mathcal{A}\subseteq \mathcal{B}$ and $\mathcal{B}\cap\Sigma=\emptyset$
has maximal elements and each such maximal element is a gr-prime matrix ideal.
\end{proposition}

\begin{proof}
Let $(\mathcal{C}_i)_{i\in I}$ be a nonempty chain in $W$. Set 
$\mathcal{C}=\bigcup_{i\in I}\mathcal{C}_i$. It is not difficult to show that
$\mathcal{C}$ is a gr-matrix ideal. Then clearly
$\mathcal{A}\subseteq\mathcal{C}_i\subseteq\mathcal{C}$ and
$\mathcal{C}\cap\Sigma=(\bigcup_{i\in I}\mathcal{C}_i)\cap\Sigma=
\bigcup_{i\in I}(\mathcal{C}_i\cap\Sigma)=\emptyset$. By Zorn's lemma,
$W$ has maximal elements. Suppose that $\mathcal{P}$ is a maximal element of $W$.
Since $\mathcal{P}\cap\Sigma=\emptyset$, $\mathcal{P}$ is a proper gr-matrix ideal.
Let $\mathcal{A}_1,\mathcal{A}_2$ be gr-matrix ideals such that
$\mathcal{P}\subsetneq\mathcal{A}_1$, $\mathcal{P}\subsetneq\mathcal{A}_2$.
Since $\mathcal{P}$ is maximal in $W$, there exist $A_1\in\mathcal{A}_1\cap\Sigma$,
$A_2\in\mathcal{A}_2\cap\Sigma$. Then $A_1\oplus A_2\in\Sigma$ and
$A_1\oplus A_2 \notin\mathcal{P}$. Therefore $\mathcal{A}_1\mathcal{A}_2\neq \mathcal{P}$.
\end{proof}

\begin{corollary}\label{coro:primeidealsexistence}
Let $R$ be a $\Gamma$-graded ring. Let $\mathcal{A}$ be a proper gr-matrix ideal. Then
there exist  maximal gr-matrix ideals $\mathcal{P}$ with $\mathcal{A}\subseteq \mathcal{P}$, and
such maximal gr-matrix ideals are gr-prime matrix ideals. 
In particular, if there are proper gr-matrix ideals, then gr-prime matrix ideals exist.
\end{corollary}

\begin{proof}
By Lemma~\ref{lem:grmatrixideal}(7), no identity matrix belongs to $\mathcal{A}$. Apply now Proposition~\ref{prop:maximalisprime} to $\mathcal{A}$ and $\Sigma=\mathfrak{I}$.
\end{proof}

\begin{proposition}\label{prop:radicalofideal}
Let $R$ be a $\Gamma$-graded ring. For each proper gr-matrix ideal $\mathcal{A}$, the radical
$\sqrt{\mathcal{A}}$ is the intersection of all gr-prime matrix ideals that contain $\mathcal{A}$.
\end{proposition}

\begin{proof}
Let $\mathcal{P}$ be a prime matrix ideal such that $\mathcal{A}\subseteq\mathcal{P}$.
If $A\in\sqrt{\mathcal{A}}$, then $\oplus^rA\in\mathcal{A}\subseteq\mathcal{P}$ for some
positive integer $r$. By (PM4), $A\in\mathcal{P}$. Thus $\sqrt{\mathcal{A}}\subseteq\mathcal{P}$.

Let now $A\in\mathfrak{M}(R)\setminus\sqrt{\mathcal{A}}$. 
Notice that such $A$ exists because $\sqrt{\mathcal{A}}\subseteq\mathcal{P}$.
If we apply Proposition~\ref{prop:maximalisprime} to
$\mathcal{A}$ and $\Sigma=\{\oplus^rA\colon r \textrm{ positive integer }\}$,
we obtain a gr-prime matrix ideal $\mathcal{P}$ such that
$\mathcal{A}\subseteq\mathcal{P}$, $\mathcal{P}\cap\Sigma=\emptyset$. Therefore
$A$ does not belong to the intersection of the gr-prime matrix ideals that contain $\mathcal{A}$.
\end{proof}

\begin{corollary}
Let $R$ be a $\Gamma$-graded ring. A proper gr-matrix ideal is gr-semiprime if and only if it is
the intersection of gr-prime matrix ideals. \qed
\end{corollary}

Let $R$ be a $\Gamma$-graded ring. By Corollary~\ref{coro:leastideal}, 
$\mathcal{A}(\mathcal{N})/\mathfrak{I}$ is the least gr-matrix ideal.
We define the \emph{gr-matrix nilradical} of $R$ 
as the gr-matrix ideal $\mathfrak{N}=\sqrt{\mathcal{A}(\mathcal{N})/\mathfrak{I}}$.

\begin{theorem}
Let $R$ be a $\Gamma$-graded ring. The following assertions are equivalent.
\begin{enumerate}[\rm(1)]
	\item There exists a $\Gamma$-graded epic $R$-division ring $(K,\varphi)$.
	\item There exists a homomorphism of $\Gamma$-almost graded  rings from $R$
	to a $\Gamma$-almost graded division ring.
	\item The gr-matrix nilradical is a proper gr-matrix ideal.
	\item No identity matrix can be expressed as a determinantal sum
	of elements of $\mathcal{N}$.
\end{enumerate}

\end{theorem}

\begin{proof}
(1) is equivalent to (2) by Theorem~\ref{theo:gradedlocal}(2)(b). One could also argue as follows.
By Proposition~\ref{prop:viceversa}, (2) implies the existence of gr-prime matrix ideals,
and therefore of $\Gamma$-graded epic $R$-division rings by Theorem~\ref{theo:primematrixequalsdivisionring}.

If (1) holds, the gr-singular kernel of $\varphi$ is a gr-prime matrix ideal by 
Theorem~\ref{theo:primematrixequalsdivisionring}.
Thus (3) holds.

If (3) holds, then $\mathcal{A}(\mathcal{N})/\mathfrak{I}$ is a proper gr-matrix ideal. By
Corollary~\ref{coro:leastideal}, (4) holds.

Suppose that (4) holds true. Again by Corollary~\ref{coro:leastideal}, there exist proper
gr-matrix ideals. By Corollary~\ref{coro:primeidealsexistence}, gr-prime matrix ideal exists.
Now Theorem~\ref{theo:primematrixequalsdivisionring} implies (1).
\end{proof}

\begin{theorem}
Let $R$ be a $\Gamma$-graded ring. There exists a universal $\Gamma$-graded epic $R$-division ring
if and only if the gr-matrix nilradical is a gr-prime matrix ideal.
\end{theorem}

\begin{proof}
By Corollary~\ref{coro:specialization}, the existence of a universal $\Gamma$-graded
epic $R$-division ring is equivalent to the existence of a least gr-prime matrix ideal
$\mathcal{P}$. Hence the least gr-matrix ideal $\mathcal{A}(\mathcal{N})/\mathfrak{I}
\subseteq \mathcal{P}$ is proper. By Proposition~\ref{prop:radicalofideal}, 
$\mathfrak{N}$ is the intersection of all gr-prime
matrix ideals. Hence $\mathfrak{N}=\mathcal{P}$.

Conversely, if $\mathfrak{N}$ is a gr-prime matrix ideal, then $\mathcal{A}(\mathcal{N})/\mathfrak{I}$
is proper and, by Proposition~\ref{prop:radicalofideal}, $\mathfrak{N}$ is the intersection
of all gr-prime matrix ideals. Therefore $\mathfrak{N}$ is the least gr-prime matrix ideal.
\end{proof}

\begin{proposition}\label{prop:invertiblematrices}
Let $R$ be a $\Gamma$-graded ring and let $P,Q\in\mathfrak{M}(R)$. There exists a homomorphism of 
$\Gamma$-graded rings $\varphi\colon R\rightarrow K$ to a $\Gamma$-graded division ring $K$ such that
$P^\varphi$ is invertible over $K$ and $Q^\varphi$ is not invertible over $K$ if and only if
no matrix of the form $I\oplus (\oplus^r P)$ can be expressed as a determinantal sum of matrices
of $\mathcal{N}\cup \mathcal{D}(\{Q\})$.
\end{proposition}

\begin{proof}
The existence of such $(K,\varphi)$ is equivalent to the existence of  gr-prime
matrix ideals $\mathcal{P}$ such that $Q\in\mathcal{P}$ and $P\notin\mathcal{P}$.
The existence of such gr-prime matrix ideals is equivalent to the condition
$P\notin\sqrt{\langle Q\rangle}$, where $\langle Q\rangle$ denotes the gr-matrix
ideal generated by $Q$. Hence it is equivalent to the condition
that no matrix of the form $\oplus^rP\in\langle Q\rangle$. By Lemma~\ref{lem:idealgenerated2}(2),
$\langle Q\rangle$ is of the form $\mathcal{A}(\{Q\})/\mathfrak{I}$.
Therefore, by Lemmas~\ref{lem:idealgenerated} and \ref{lem:idealgenerated2}, 
everything is equivalent to the condition that no matrix
of the form $I\oplus(\oplus^rP)$ can be expressed as a determinantal sum of matrices
of $\mathcal{N}\cup\mathcal{D}(\{Q\})$, as desired.
\end{proof}

\begin{corollary}
Let $R$ be a $\Gamma$-graded ring and let $P,Q\in\mathfrak{M}(R)$. The following assertions hold true.
\begin{enumerate}[\rm(1)]
	\item There exists a $\Gamma$-graded epic $R$-division ring $(K,\varphi)$ such that
	$P^\varphi$ is invertible over $K$ if and only if no matrix of the form
	$I\oplus(\oplus^rP)$ can be expressed as a determinantal sum of matrices of $\mathcal{N}$.
	\item There exists a $\Gamma$-graded epic $R$-division ring $(K,\varphi)$ such that
	$Q^\varphi$ is not invertible over $K$ if and only if no identity matrix 
	can be expressed as a determinantal of matrices of $\mathcal{N}\cup\mathcal{D}(\{Q\})$.
\end{enumerate}
\end{corollary}

\begin{proof}
(1) In Proposition~\ref{prop:invertiblematrices}, let $Q=0$.

(2) In Proposition~\ref{prop:invertiblematrices}, let $P=1$.
\end{proof}

\begin{theorem}
Let $R$ be a $\Gamma$-graded ring. The following assertions are equivalent.
\begin{enumerate}[\rm(1)]
	\item There exists a $\Gamma$-graded epic $R$-division ring of fractions $(K,\varphi)$.
	\item There exists a homomorphism of $\Gamma$-almost graded rings
	$\varphi\colon R\rightarrow K$ with $K$ a $\Gamma$-almost graded division ring
	such that $\varphi(x)\neq 0$ for each $x\in\h(R)\setminus\{0\}$.
	\item $R$ is a $\Gamma$-graded domain and no matrix of the form $aI$
	with $a\in\h(R)\setminus\{0\}$ can be expressed as a determinantal sum of matrices of
	$\mathcal{N}$.
	\item No diagonal matrix with nonzero homogeneous elements on the main diagonal can 
	be expressed as a determinantal sum of matrices of $\mathcal{N}$.
\end{enumerate}
\end{theorem}

\begin{proof}
(1) and (2) are equivalent by Theorem~\ref{theo:gradedlocal}(b).

Suppose that (1) holds true. Then, for each diagonal matrix $A$ as in (4),
$A^\varphi$ is invertible. Thus, $A\notin\mathcal{P}$, the gr-prime matrix
ideal given as the gr-singular kernel of $\varphi$. In particular, $A$ cannot
be expressed as the determinantal sum of matrices in $\mathcal{N}$. Thus (4) holds.

Suppose (4) holds. Clearly no matrix of the form $aI$ with 
$a\in\h(R)\setminus\{0\}$ can be expressed as a determinantal sum of matrices of
	$\mathcal{N}$. Thus to prove (3), it remains to show that $R$ is a $\Gamma$-graded
	domain. Thus, let $a,b\in\h(R)$ of degrees $\gamma,\delta\in\Gamma$, respectively.
	If $ab=0$, then $\left(\begin{smallmatrix} a & 0 \\ 0 & b \end{smallmatrix} \right)
	\in M_2(R)[(\gamma,e)][(e,\delta^{-1})]$. Then 
	we can express $\left(\begin{smallmatrix} a & 0 \\ 0 & b \end{smallmatrix} \right)=
	\left(\begin{smallmatrix} a & 0 \\ 1 & b \end{smallmatrix} \right)\nabla 
	\left(\begin{smallmatrix} 0 & 0 \\ -1 & b \end{smallmatrix} \right)$ 
	as a determinantal sum of  matrices in $M_2(R)[(\gamma,e)][(e,\delta^{-1})]$.
	Note that $\left(\begin{smallmatrix} 0 & 0 \\ -1 & b \end{smallmatrix} \right)$ is
	hollow, and hence it is not gr-full.
	Furthermore, $\left(\begin{smallmatrix} a & 0 \\ 1 & b \end{smallmatrix} \right)=
	\left(\begin{smallmatrix} a  \\ 1  \end{smallmatrix} \right)
	\left(\begin{smallmatrix} 1 & b \end{smallmatrix} \right)$, where the factors
	belong to $M_{2\times 1}(R)[(\gamma,e)][e]$ and $M_{1\times 2}(R)[e][(e,\delta^{-1})]$, respectively.
	Hence $\left(\begin{smallmatrix} a & 0 \\ 0 & b \end{smallmatrix} \right)$
	can be expressed as a determinantal sum of matrices from $\mathcal{N}$.
	By (4), either $a=0$ or $b=0$. Hence $R$ is a $\Gamma$-graded domain and (3) holds.

Suppose now that (3) holds.  
If there does not exist a $\Gamma$-graded epic $R$-division ring
of fractions, then, by Corollary~\ref{coro:ultraproductoffractions},
there exists nonzero $a\in\h(R)$ such that $a^\varphi$ is not invertible
for every homomorphism of $\Gamma$-graded rings $\varphi\colon R\rightarrow K$
with $K$ a $\Gamma$-graded division ring. Hence the $1\times 1$ homogeneous
matrix $(a)$ belongs to the intersection of all gr-prime matrix ideals, i.e. 
$(a)\in\mathfrak{N}$. Hence $\oplus^r(a)\in\mathcal{A}(\mathcal{N})/\mathfrak{I}$.
Thus $I_s\oplus(\oplus^r(a))=I_s\oplus aI_r$ can be written as a determinantal sum of matrices
of $\mathcal{N}$. Then, since  $aI_s\oplus I_r\in \mathfrak{M}(R)$ and
it is diagonal, $aI_{r+s}=(aI_s\oplus I_r)(I_s\oplus aI_r)$
is a determinantal sum of matrices of $\mathcal{N}$, a contradiction. Therefore (1) holds.
\end{proof}



\section{gr-Sylvester rank functions}\label{sec:grSylvesterrankfunctions}

Throughout this section, let $\Gamma$ be a group. 

\medskip

The aim of this section is to show that the different definitions of gr-Sylvester rank functions (with values in $\mathbb{N}$) given below are equivalent between them and with the definition of a gr-prime matrix ideal, and thus they uniquely determine homomorphisms to graded division rings. We will adapt the definitions, results and proofs of \cite[Sections~1 and 3]{Malcolmsondetermining} and \cite[p.94--98]{Schofieldbook} to the graded situation. In defining gr-Sylvester rank functions, the main difference with the ungraded case stems from the  fact that, in the graded case, the same matrix
$A\in\mathfrak{M}_\bullet(R)$ can define more than one homomorphism between $\Gamma$-graded free modules.
This is reflected in properties  (MatRF4), (ModRF4) and (MapRF5) below.

We begin this section providing the different definitions
of gr-Sylvester rank functions for a $\Gamma$-graded ring (with values in $\mathbb{N}$), together with some of their basic properties.

\medskip

Let $R$ be a $\Gamma$-graded ring.
A \emph{gr-Sylvester matrix rank function} for $R$ is a map 
$\rank\colon \mathfrak{M}_{\bullet}(R)\rightarrow \mathbb{N}$ that satisfies the following conditions
\begin{enumerate}[({MatRF}1)]
\item $\rank((1))=1$, where $(1)$ is the identity matrix of size $1\times1$.
\item $\rank(AB)\leq \min\{\rank(A),\rank(B)\}$ for all compatible matrices $A,B\in\mathfrak{M}_\bullet(R)$.
\item $\rank\left(\begin{smallmatrix} A & 0 \\ 0 & B\end{smallmatrix}\right)=\rank(A)+\rank(B)$ for all
$A,B\in\mathfrak{M}_\bullet(R)$. 
\item $\rank\left(\begin{smallmatrix} A & C \\ 0 & B\end{smallmatrix}\right)\geq\rank(A)+\rank(B)$ for all
$A,B,C\in\mathfrak{M}_\bullet(R)$ such that $A$ has distribution $(\overline{\alpha},\overline{\beta})$, $B$ has distribution $(\overline{\delta},\overline{\varepsilon})$ and $C$ has distribution $(\overline{\alpha},
\overline{\varepsilon})$ for some finite sequences 
$\overline{\alpha}$, $\overline{\beta}$, $\overline{\delta}$, $\overline{\varepsilon}$ of elements of $\Gamma$ .  
\end{enumerate} 

Let $\rank_1,\rank_2$ be two gr-Sylvester matrix rank functions for $R$. We say that $\rank_1\leq \rank_2$ if 
$\rank_1(A)\leq \rank_2(A)$ for all $A\in\mathfrak{M}_\bullet(R)$. In this way, there is
defined a partial order in the set of gr-Sylvester matrix rank functions for $R$.

The following lemma describes some useful properties of gr-Sylvester matrix rank functions.

\begin{lemma}\label{lem:grrankproperties} 
Let $R$ be a $\Gamma$-graded ring and $\rank\colon\mathfrak{M}_\bullet(R)\rightarrow\mathbb{N}$
be a gr-Sylvester matrix rank function. Let $A\in M_{m\times n}(R)[\overline{\alpha}][\overline{\beta}]$. The following assertions hold true.
\begin{enumerate}[\rm(1)]
	\item $\rank(I_n)=n$ for all positive integers $n$.
	\item $\rank(Z)=0$ for all zero matrices of any size.
	\item The condition $\rank(A)\geq 0$ follows from {\rm(MatRF1)--(MatRF4)}.
	
	\item $\rank(A)=\rank(PA)=\rank(AQ)$ for all 
	invertible matrices $P\in M_m(R)[\overline{\alpha'}][\overline{\alpha}]$ and 
	$Q\in M_{n}(R)[\overline{\beta}][\overline{\beta'}]$.
	\item If $P\in\mathfrak{M}(R)$ is invertible of size $n\times n$, then $\rank(P)=n$.
	\item 	$\rank(A)\leq\min(m,n)$.
	\item $\rank(A)\leq \rank\left(\begin{smallmatrix} A\\ B \end{smallmatrix}\right)$ and
	$\rank(A)\leq \rank (A\ C)$ for all $B\in M_{m'\times n}(R)[\overline{\alpha'}][\overline{\beta}]$
		and $C\in M_{m\times n'}(R) [\overline{\alpha}][\overline{\beta'}]$.

\end{enumerate}

\end{lemma}

\begin{proof}
(1) follows from (MatRF1) and (MatRF3).

(2) We denote the zero matrix of size $m\times n$ by $0_{m\times n}$. From
$$(1\ 0_{1\times m})\left(\begin{array}{ll} 1 & 0 \\
0 & 0_{m\times n} \end{array}\right) \left(\begin{array}{l}
1\\ 0_{n\times 1} \end{array}\right)=(1),\quad 
\left(\begin{array}{l}1 \\ 0_{m\times 1} \end{array}\right)(1)(1\ 0_{1\times n})
= \left(\begin{array}{ll} 1 & 0 \\
0 & 0_{m\times n} \end{array}\right),$$
applying (MatRF1)--(MatRF3), we obtain that $1+\rank(0_{m\times n})\geq 1$ and
$1\geq 1+\rank(0_{m\times n})$, respectively.

(3) Let $Z$ be a zero matrix of appropriate size. By (2) and (MatRF2),
$0=\rank(Z)=\rank(AZ)\leq \min\{\rank(A),\rank(Z)\}\leq \rank(A).$

(4) By (MatRF2), $$\rank(PA)\leq\min\{\rank(P),\rank(A)\}\leq \rank(A),$$ 
$$\rank(A)=\rank(P^{-1}PA)\leq \min\{\rank(P^{-1}),\rank(PA)\}\leq\rank(PA).$$
Hence $\rank(A)=\rank(PA)$. The other case is shown in the same way.

(5) By (4), $\rank(I)=\rank(PI)=\rank(P)$.

(6) Since $A$ is $m\times n$, using (1), we obtain
$\rank(A)=\rank(I_mA)\leq \min\{\rank(I_m),\rank(A)\}\leq m$ and 
$\rank(A)=\rank(AI_n)=\min\{\rank(A),\rank(I_n)\}\leq n$.

(7) follows from
$$\rank(A)=\rank\left(\Big(I_m\ \ 0_{m\times m'}\Big) \left(\begin{array}{c} A\\ B\end{array}\right)\right)
\leq \rank\left(\left(\begin{array}{c} A\\ B\end{array}\right)\right),$$

$$\rank(A)=\rank\left(\Big(A\ \ C\Big) \left(\begin{array}{c} I_n\\ 0_{n'\times n}\end{array}\right)\right)
\leq \rank\left( \Big(A\ \ C \Big)\right).$$
\end{proof}

\medskip

Let $R$ be a $\Gamma$-graded ring. We will denote  the forgetful functor from the category
of finitely presented $\Gamma$-graded $R$-modules to the category of finitely presented $R$-modules by 
$\mathcal{F}$.

A \emph{gr-Sylvester module rank function} for $R$ is a function $\di$ on the class of finitely presented 
$\Gamma$-graded (right) $R$-modules with values on $\mathbb{N}$ such that
\begin{enumerate}[({ModRF}1)]
	\item $\di(R)=1$, where $R$ is the $\Gamma$-graded ring $R$ considered as a (right) $R$-module
	in the natural way.
	\item $\di(M_1\oplus M_2)=\di(M_1)+\di(M_2)$.
	\item For any exact sequence $M_1\rightarrow M_2\rightarrow M_3\rightarrow 0$ of 
	graded homomorphisms between finitely presented
	$\Gamma$-graded 	$R$-modules, $$\di(M_3)\leq \di(M_2)\leq \di(M_1)+\di(M_3).$$
	\item Let $f\colon R^n(\overline{\beta})\rightarrow R^m(\overline{\alpha})$ and
	$f'\colon  R^n(\overline{\beta'})\rightarrow R^m(\overline{\alpha'})$ be homomorphisms of $\Gamma$-graded
	$R$-modules. If $\mathcal{F}(f)=\mathcal{F}(f')$, then $\di(\coker f)=\di(\coker f')$.
\end{enumerate}
Let $\di_1,\di_2$ be two gr-Sylvester module rank functions for $R$ with values in $\mathbb{N}$. We say that $\di_1\leq \di_2$ if 
$\di_1(M)\leq \di_2(M)$ for all finitely presented $\Gamma$-graded $R$-modules $M$. In this way, there is
defined a partial order in the set of gr-Sylvester module rank functions for $R$.

The following easy but important remarks are in order.

\begin{lemma}
Let $R$ be a $\Gamma$-graded ring and $\di$ be a gr-Sylvester module rank function. 
\begin{enumerate}[{\rm (1)}]
\item $\di(0)=0$, where $0$ denotes the zero module.

\item $\di(M)\geq 0$ for all finitely presented $\Gamma$-graded $R$-modules
follows from {\rm(ModRF1)--(ModRF4)}.

\item If $M$ and $N$ are isomorphic as $\Gamma$-graded $R$-modules, then $\di(M)=\di(N)$. 

\item  $\di(R(\theta))=1$ for all $\theta\in\Gamma$. 

\item In {\rm(ModRF4)}, the condition $\mathcal{F}(f)=\mathcal{F}(f')$ implies that the lengths of the sequences
$\overline{\alpha}$ and $\overline{\alpha'}$ (respectively $\overline{\beta}$ and $\overline{\beta'}$)
coincide. 

\item $\di(R^n(\overline{\delta}))=\di(R^n(\overline{\delta'}))$ for any 
finite sequences $\overline{\delta}$, $\overline{\delta'}$ of the same length.
\end{enumerate}
\end{lemma}

\begin{proof}
(1) follows from (ModRF3) applied to
$0\rightarrow 0\oplus 0\oplus 0\rightarrow 0\rightarrow 0$.

(2)  follows
applying (ModRF3) to $M\rightarrow 0\rightarrow 0\rightarrow 0$.

(3) Make $M_1=0$, $M_2=M$ and $M_3=N$ in
(ModRF3).

(4)  By (ModRF4),  the natural inclusions of $R$ in the first component
$R\rightarrow R\oplus R(\theta)$, $R\rightarrow R\oplus R$, implies that $\di(R(\theta))=\di(R)$.

(5) Trivial.

(6) holds true by (ModRF2) and the foregoing.
\end{proof}
 
Let $R$ be a $\Gamma$-graded ring. 
A \emph{gr-Sylvester map rank function} for $R$ is a function $\rho$ on the class of
all homomorphisms of  
$\Gamma$-graded (right) $R$-modules between finitely generated $\Gamma$-graded projective
$R$-modules  with values on $\mathbb{N}$ such that
\begin{enumerate}[({MapRF}1)]
	\item $\rho(1_R)=1$, where $1_R$ denotes the identity map on $R$. 
	\item $\rho(gf)\leq \min\{\rho(f),\rho(g)\}$.
	\item $\rho\left(\begin{smallmatrix}f & 0 \\ 0 & g \end{smallmatrix}\right)=\rho(f)+\rho(g)$.
   \item $\rho\left(\begin{smallmatrix}f & h \\ 0 & g \end{smallmatrix}\right)\geq \rho(f)+\rho(g)$.
	\item Let $f\colon R^n(\overline{\beta})\rightarrow R^m(\overline{\alpha})$ and
	$f'\colon  R^n(\overline{\beta'})\rightarrow R^m(\overline{\alpha'})$ be homomorphisms of $\Gamma$-graded
	$R$-modules. If $\mathcal{F}(f)=\mathcal{F}(f')$, then $\rho(f)=\rho(f')$.
\end{enumerate}
Let $\rho_1,\rho_2$ be two gr-Sylvester map rank functions for $R$. We say that $\rho_1\leq \rho_2$ if 
$\rho_1(f)\leq \rho_2(f)$ for all homomorphism $f$ of $\Gamma$-graded modules  between
finitely generated $\Gamma$-graded projective modules. In this way, there is
defined a partial order in the set of gr-Sylvester map rank functions for $R$.

The proof of the following remarks can be proved very much as in Lemma~\ref{lem:grrankproperties}.

\begin{lemma}\label{lem:grmapproperties}
Let $R$ be a $\Gamma$-graded ring and $\rho$ be a gr-Sylvester map rank function. Let $f\colon P\rightarrow Q$ be a homomorphism of $\Gamma$-graded rings between the finitely generated 
$\Gamma$-graded projective modules $P$ and $Q$. The following assertions hold true.
\begin{enumerate}[\rm(1)]
	\item $\rho(1_{R^n(\overline\theta)})=n$  where $1_{R^n(\overline\theta)}$ denotes the identity homomorphism of the $\Gamma$-graded free $R$-module ${R^n(\overline\theta)}$.
	\item $\rho(0)=0$ where $0$ denotes the zero homomorphism between any finitely generated
	$\Gamma$-graded projective $R$-modules.
	\item The condition $\rho(f)\geq 0$ follows from {\rm(MapRF1)--(MapRF4)}.
	
	\item $\rho(f)=\rho(gf)=\rho(fh)$ for all isomorphisms of $\Gamma$-graded $R$-modules between
	finitely generated $\Gamma$-graded projective $R$-modules
	$g\colon Q\rightarrow Q'$, $h\colon P'\rightarrow P$. 
	\item If $f$ is invertible,  then $\rho(f)=\rho(1_P)=\rho(1_Q).$
	\item 	$\rho(f)\leq\min(\rho(1_P),\rho(1_Q))$.
	\item $\rho(f)\leq \rho\left(\begin{smallmatrix} f\\ g \end{smallmatrix}\right)$ and
	$\rho(f)\leq \rho (f\ h)$ for all homomorphisms of
	$\Gamma$-graded $R$-modules $g\colon P\rightarrow Q'$ and $h\colon P'\rightarrow Q$
	with $P',Q'$ being finitely generated $\Gamma$-graded projective $R$-modules. 	
\end{enumerate}
\end{lemma}

\begin{proof}
(1) Notice that the identity maps $1_{R(\theta)}\colon R(\theta)\rightarrow R(\theta)$ and
$1_R\colon R\rightarrow R$ are such that $\mathcal{F}(1_R)=\mathcal{F}(1_{R(\theta)})$. Hence
$\rho(1_{R(\theta)})=1$ for all $\theta\in\Gamma$. Now apply (MapRF3) to $1_{R^n(\overline{\theta})}$
for any $\overline{\theta}\in\Gamma^n$.

(2)  We denote by $0_{MN}$ the zero homomorphism $M\rightarrow N$ between the $\Gamma$-graded
$R$-modules $M$ and $N$. From the equalities $$(1_R\ 0_{QR})\left(\begin{array}{ll} 1_R & 0_{PR} \\
0_{RQ} & 0_{PQ} \end{array}\right) \left(\begin{array}{l}
1_R\\ 0_{RP} \end{array}\right)=(1_R),$$ 
$$\left(\begin{array}{l}1_R \\ 0_{RQ} \end{array}\right)(1_R)(1_R\ 0_{PR})
= \left(\begin{array}{ll} 1_R & 0_{PR} \\
0_{RQ} & 0_{PQ} \end{array}\right),$$
we obtain that $1+\rho(0_{PQ})=\rho\left(\begin{smallmatrix} 1_R & 0_{PR} \\
0_{RQ} & 0_{PQ} \end{smallmatrix}\right)\geq \rho(1_R)=1$ and
$1=\rho(1_R)\geq \rho\left(\begin{smallmatrix} 1_R & 0_{PR} \\
0_{RQ} & 0_{PQ} \end{smallmatrix}\right)=1+\rho(0_{PQ})$. Hence $\rho(0_{PQ})=0$.

(3) By (2) and (MapRF2), $0=\rho(0_{PQ})=\rho(0_{QQ}f)\leq \min\{\rho(0_{QQ}),\rho(f)\}\leq \rho(f)$.

(4) By (MapRF2), $$\rho(gf)\leq\min\{\rho(g),\rho(f)\}\leq \rho(f),$$ 
$$\rho(f)=\rho(g^{-1}gf)\leq \min\{\rho(g^{-1}),\rho(gf)\}\leq\rho(gf).$$
Hence $\rho(f)=\rho(gf)$. The other case is shown in the same way.

(5) By (4), $\rho(1_P)=\rho(f1_P)=\rho(f)$ and similarly $\rho(f)=\rho(1_Q)$.

(6) Using (MapRF2), we obtain
$\rho(f)=\rho(1_Qf)\leq \min\{\rho(1_Q),\rho(f)\}\leq \rho(1_Q)$ and 
$\rho(f)=\rho(f1_P)=\min\{\rho(f),\rho(1_P)\}\leq \rho(1_P)$.

(7) follows from
$$\rho(f)=\rho\left(\Big(1_Q\ \ 0_{Q'Q}\Big) \left(\begin{array}{c} f\\ g\end{array}\right)\right)
\leq \rho\left(\left(\begin{array}{c} f\\ g\end{array}\right)\right),$$

$$\rho(f)=\rho\left(\Big(f\ \ h\Big) \left(\begin{array}{c} 1_P\\ 0_{PP'}\end{array}\right)\right)
\leq \rho\left( \Big(f\ \ h \Big)\right).$$

\end{proof}

\subsection{Equivalence between gr-prime matrix ideals and gr-Sylvester matrix rank functions}
\label{subsec:gr-prime_gr-Sylvester_matrix}

The proofs and arguments contained in this subsection are the easy and natural adaptation of the ones in 
\cite[Section~3]{Malcolmsondetermining}.

\medskip

Let $D$ be a $\Gamma$-graded division ring. It is not difficult to show that $\mathfrak{M}_\bullet(D)\rightarrow \mathbb{N}$, $A\mapsto \drank(A)$,
satisfies (MatRF1)--(MatRF4).

Let $R$ be a $\Gamma$-graded ring and let $(K,\varphi)$ be a $\Gamma$-graded epic $R$-division ring.
One can induce a gr-rank function $\rank_\varphi$  for $R$ defining 
$\rank_\varphi(A)=\drank(A^\varphi)$ for all $A\in\mathfrak{M}_\bullet(R)$. Note that non-isomorphic $\Gamma$-graded
epic $R$-division rings induce different gr-rank functions because the gr-singular kernels
do not coincide, by Theorem~\ref{theo:primematrixequalsdivisionring}.
The aim of this section is to show that there are no other gr-rank functions for $R$.

\begin{lemma}\label{lem:grrankalternative}
Let $R$ be a $\Gamma$-graded ring. Suppose that $\rank\colon\mathfrak{M}_\bullet(R)\rightarrow \mathbb{N}$
is a gr-rank function and that $\rank(A)=n$ for some $A\in\mathfrak{M}_\bullet(R)$.
The following assertions hold true.
\begin{enumerate}[\rm(1)]
	\item Suppose that $A'$ is obtained as the result of eliminating a column (a row) of $A$, then
	$$\rank(A)-1\leq \rank(A')\leq \rank(A).$$
	\item If the result $A'$ of eliminating any of the columns (rows) of $A$ is such that
	$\rank(A')<\rank(A)$, then $A$ has exactly $n$ columns (rows).
	\item $n$ is the largest natural number such that there is a square submatrix $B$ of $A$ with
	$\rank(B)=\textrm{size of } B$.
\end{enumerate}
\end{lemma}

\begin{proof}
(1) Suppose that $A=(A'\ a)$, where $a$ is the last column of $A$. By Lemma~\ref{lem:grrankproperties}(7), $\rank(A')\leq\rank(A)$. Moreover
$$\rank(A)=\rank(A'\ a)\leq\rank\left(\begin{array}{cc} A'&a\\ 0 & 1 \end{array}\right)=
\rank\left(\begin{array}{cc} A'&0\\ 0 & 1 \end{array}\right)=\rank(A')+1,$$
where we have used Lemma~\ref{lem:grrankproperties}(4) on the second equality.
When the column we eliminate is not the last one, the result follows by Lemma~\ref{lem:grrankproperties}(4).

When $A'$ is obtained by eliminating some row, the result can be proved analogously.

(2) We claim that $A'$ obtained as the result of eliminating $m$ of the columns of $A$ satisfies
$\rank(A')=n-m$. We prove the claim by induction on $m$. 

If $m=1$, the result follows from (1). Suppose that the claim holds true for $m\leq k-1$.
Let $A'$ be the result of eliminating any $k$ columns and $a_1,a_2$ be two different columns of $A$ from
among those eliminated in obtaining $A'$. By induction hypothesis, 
$$\rank(A'\ a_1)=n-(k-1),\quad \rank(A'\ a_2)=n-(k-1),\quad \rank(A'\ a_1\ a_2)=n-(k-2).$$
Now
\begin{eqnarray*}\label{eq:}
\rank(A')+\rank(A'\ a_1\ a_2) & \leq & \rank\begin{pmatrix} A'&0&a_1&0\\ 0&A'&a_1&a_2 \end{pmatrix}
=\rank\begin{pmatrix} A'&a_1&0&0\\ 0&a_1&A'&a_2 \end{pmatrix}\\
&=& \rank\begin{pmatrix} A'&a_1&0&0\\ A'&a_1&A'&a_2 \end{pmatrix} =  
\rank\begin{pmatrix} A'&a_1&0&0\\ 0&0&A'&a_2 \end{pmatrix}\\ 
& = & \rank(A'\ a_1)+\rank(A'\ a_2),
\end{eqnarray*}
where we have used  Lemma~\ref{lem:grrankproperties}(4) in all but the last equality.
Thus $\rank(A')\leq n-k$. On the other hand, applying (1) several times, we get that $n-k\leq \rank(A')$, and
the claim is proved.

Since $\rank(A)\geq 0$ for all $A\in\mathfrak{M}_\bullet(R)$, the claim proves the result.

When we eliminate rows instead of columns, the result follows analogously.

(3) By (2), we can eliminate rows and/or columns of $A$ until we reach a square submatrix of $A$
of size exactly $n$.
\end{proof}

Now we are ready to prove the main result of this subsection.

\begin{theorem}\label{theo:grprimematrix_grSylvestermatrix}
Let $R$ be a $\Gamma$-graded ring. There is an anti-isomorphism of partially ordered sets
\begin{eqnarray*}
\{\textrm{gr-prime matrix ideals of }R \}&\longrightarrow& \{\textrm{gr-Sylvester matrix rank functions for }R\} \\
\mathcal{P}&\longmapsto& \rank_{\mathcal{P}} \\
\mathcal{P}_{\rank{}} &\longmapsfrom& \rank 
\end{eqnarray*}
defined as follows. 
\begin{enumerate}[\rm(a)]
\item If $\mathcal{P}$ is a gr-prime matrix ideal of $R$ and  
$A\in\mathfrak{M}_\bullet(R)$, then $\rank_\mathcal{P}(A)$ is the size of the largest square submatrix
of $A$ which is not in $\mathcal{P}$.  Equivalently, if $(K,\varphi_\mathcal{P})$
is a $\Gamma$-graded epic $R$-division ring with singular kernel $\mathcal{P}$,
then $\rank_\mathcal{P}(A)=\drank(A^{\varphi_\mathcal{P}})$.

\item Conversely, given a gr-Sylvester matrix rank function $\rank\colon \mathfrak{M}_\bullet(R)\rightarrow \mathbb{N}$,
then the set $$\mathcal{P}_{\rank}=\{A\in\mathfrak{M}(R)\colon\rank(A)<\textrm{size of }A\}$$
is a gr-prime matrix ideal.
\end{enumerate}
\end{theorem}

\begin{proof}
Theorem~\ref{theo:primematrixequalsdivisionring} proves that the two ways of describing the correspondence
$\mathcal{P}\mapsto \rank_{\mathcal{P}}{}$ are equivalent. By the comment at the beginning of 
Section~\ref{subsec:gr-prime_gr-Sylvester_matrix}, $\rank_{\mathcal{P}}{}$ is a gr-rank function for $R$.
Moreover, if $\mathcal{P}\subseteq\mathcal{Q}$ are gr-prime matrix ideals, then
$\rank_{\mathcal{P}{}}\geq \rank_{\mathcal{Q}{}}$ because there are more possible square
submatrices which are not in $\mathcal{P}$.
Conversely, suppose that $\rank{},\rank'{}$ are gr-Sylvester matrix rank functions with $\rank'{}\leq 
\rank{}$.  For $A \in \mathcal{P}_{\rank{}}$,  $\rank(A) < \textrm{size of }A$.  Since $\rank'(A) \leq \rank(A)$, 
$\rank'(A)\leq \textrm{size of } A$, and  $A \in \mathcal{P}_{\rank'}$. Therefore $\mathcal{P}_{\rank{}} \subseteq \mathcal{P}_{\rank'{}}$.

(a) Let now $\rank{}$ be a gr-Sylvester matrix rank function for $R$. Let $\mathcal{P}_{\rank{}}$ be defined
as (b) of the statement of the theorem. We have to show that $\mathcal{P}_{\rank{}}$ is
a gr-prime matrix ideal.

Let $A\in\mathfrak{M}(R)$ be an $n\times n$ non gr-full matrix. There exist
$\overline{\alpha},\overline{\beta}\in \Gamma^n$ such that $A\in M_n(R)[\overline{\alpha}][\overline{\beta}]$,  $\overline{\delta}\in \Gamma^{n-1}$ and matrices $B\in M_{n\times (n-1)}(R)[\overline{\alpha}][\overline{\delta}]$, $C\in M_{(n-1)\times n}(R)[\overline{\delta}][\overline{\beta}]$ such that 
$A=BC$. By Lemma~\ref{lem:grrankproperties}(6), $\rank(B)\leq n-1$ and $\rank(C)\leq n-1$. 
By (MatRF2), $\rank(A)=\rank(BC)\leq \min\{\rank(B),\rank(C)\}<n$. Thus $A\in\mathcal{P}_{\rank{}}$ and
$\mathcal{P}_{\rank{}}$ satisfies (PM1). 

Let $A\in\mathcal{P}_{\rank{}}$, $B\in\mathfrak{M}(R)$. Since $\rank(A)<\textrm{size of }A$, then
$\rank(A\oplus B)=\rank(A)+\rank(B)<\textrm{size of } A\oplus B$. Thus (PM3) follows. 

Let $A,B\in\mathfrak{M}(R)$ and supppose that $A\oplus B\in \mathcal{P}_{\rank{}}$. Then
$$\rank(A)+\rank(B)=\rank(A\oplus B)< \textrm{size of } A\oplus B=\textrm{size of } A +
\textrm{size of }B.$$  Hence $\rank(A)<\textrm{size of } A$ or $\rank(B)<\textrm{size of } B$.
 Thus (PM4) follows.

Let $A\in\mathfrak{M}(R)$ and suppose that $A\oplus 1\in\mathcal{P}_{\rank{}}$. It means that
$\rank(A)+1=\rank(A\oplus 1)<1+\textrm{size of } A$. Hence $\rank(A)<\textrm{size of }A$. Therefore
$A\in\mathcal{P}_{\rank{}}$ and (PM5) is satisfied.

By definition of gr-Sylvester matrix rank function, 
$(1)\notin\mathcal{P}_{\rank{}}$ and thus (PM5) holds.

By Lemma~\ref{lem:grrankproperties}(4), (PM6) follows. It remains to show (PM2).

Let $A,A'\in\mathcal{P}_r$ such that $A\nabla A'$ exists with respect to the last column.
Suppose that $A=(B\ c)$, $A'=(B\ c')$. Then $A\nabla A'=(B\ c+c')$. We claim that
$\rank(A\nabla A')\leq \max\{\rank(A),\rank(A')\}$. This claim implies that 
$A\nabla A'\in\mathcal{P}_{\rank{}}$. Then 
the case of the  determinantal sum with respect to any other column  follows from Lemma~\ref{lem:grrankproperties}(4) and the claim.

Now we prove the claim. 
\begin{eqnarray}
\rank(A)+\rank(A')=\rank(B\ c)+\rank(B\ c') & = & \rank\begin{pmatrix}
	B & c & 0 & 0 \\ 0 & 0 & B & c' 
\end{pmatrix} = \rank\begin{pmatrix}
	B & c & 0 & 0 \\ -B & 0 & B & c' 
\end{pmatrix} \nonumber \\
& = & \rank\begin{pmatrix}
	B & c & 0 & 0 \\ 0 & c & B & c' 
\end{pmatrix}= \rank\begin{pmatrix}
	B & c & 0 & c \\ 0 & c & B & c+c' 
\end{pmatrix} \nonumber\\
&\geq & \rank(B) +\rank(c\ B\ c+c').\label{eq:grrank}
\end{eqnarray}
If $\rank(A)\leq\rank(B)$, then $\rank(A)=\rank(B)$ by Lemma~\ref{lem:grrankproperties}(7).
By \eqref{eq:grrank}, $$\rank(A')\geq \rank(c\ B\ c+c')\geq \rank(B\ c+c')=\rank(A\nabla A'),$$ as desired.

If $\rank(A)>\rank(B)$, then $$\rank(A)\geq \rank(B)+1=\rank\begin{pmatrix}
	B & 0 \\ 0 & 1
\end{pmatrix}=\rank \begin{pmatrix}
	B & c+c' \\ 0 & 1
\end{pmatrix}\geq \rank(B\ c+c')=\rank(A\nabla A').$$
Interchanging the roles of $A$ and $A'$, we get
that if $\rank(A')\leq \rank(B)$, then $\rank(A)\geq \rank(A\nabla A')$ and
that if $\rank(A')>\rank(B)$, then $\rank(A')\geq \rank(A\nabla A')$. Thus, the claim is proved.

It remains to show that the maps $\mathcal{P}\mapsto\rank_{\mathcal{P}}{}$ and
$\rank{}\mapsto \mathcal{P}_{\rank{}}$ are inverse one of the other.

If $\mathcal{P}$ is a gr-prime matrix ideal, the gr-prime matrix ideal that corresponds to 
$\rank_\mathcal{P}{}$ is the set of matrices $A\in\mathfrak{M}(R)$ such that
$\rank_\mathcal{P}(A)<\textrm{size of } A$. That is, the set of matrices  
$A\in\mathfrak{M}(R)$ whose largest square submatrix that is not in $\mathcal{P}$
is less than the size of $A$. In other words, the matrices $A\in\mathcal{P}$. Therefore 
$\rank_\mathcal{P}{}\mapsto\mathcal{P}$. On the other hand, let now $\rank{}$ be
 a gr-rank function for $R$. Let $\rank_{\mathcal{P}_{\rank{}}}{}$ be the associated gr-rank function
associated to $\mathcal{P}_{\rank{}}$. If $A\in\mathfrak{M}_\bullet(R)$, then $\rank_{\mathcal{P}_{\rank{}}}(A)$
equals the size of a largest square submatrix of $A$ which is not in $\mathcal{P}_{\rank{}}$.
That is, the size of a largest square submatrix $B$ of $A$ such that $\rank(B)=\textrm{size of }B$.
Hence we have to show that $\rank(A)=n$ if and only if $n$ is the size of a largest
square submatrix of $A$ such that $\rank(B)=n$. But this now follows from 
Lemma~\ref{lem:grrankalternative}(3).
\end{proof}

\subsection{Equivalence between gr-Sylvester rank functions}

The following result can be proved in exactly the same way as in \cite[Lemma~2]{Malcolmsondetermining} where
the ungraded case is shown.

\begin{lemma}\label{lem:gradedSchanuel}
Let $R$ be a $\Gamma$-graded ring.
If $0\rightarrow K\rightarrow Q\rightarrow M\rightarrow 0$ and $0\rightarrow K'\rightarrow Q'\rightarrow M\rightarrow 0$ are exact sequences of $\Gamma$-graded $R$-modules with $Q$ and $Q'$
$\Gamma$-graded projective $R$-modules and with $K\subseteq Q$ and $K'\subseteq Q'$, then there is
an automorphism of $\Gamma$-graded $R$-modules of $Q\oplus Q'$ which maps
$K\oplus Q'$ onto $Q\oplus K'$. \qed
\end{lemma}

The next result was first stated in \cite[p.~97]{Schofieldbook} for the ungraded case. Our proof follows
the one of \cite[Theorem~4]{Malcolmsondetermining}. We would like to remark that the fact that $\mathbb{N}$ is the set of values of Sylvester rank functions is not used in the proof.

\begin{theorem}\label{theo:module_map_equivalent}
Let $R$ be a $\Gamma$-graded ring. There is an anti-isomorphism of partially ordered sets
\begin{eqnarray*}
\left\{\begin{array}{c}
\textrm{gr-Sylvester module rank } \\
\textrm{functions for }R 
\end{array}\right\}&\longrightarrow & \left\{\begin{array}{c}
\textrm{gr-Sylvester map rank } \\ 
\textrm{functions for }R
\end{array}\right\} \\
\di & \longmapsto & \rho_{\di} \\
\di_\rho & \longmapsfrom & \rho
\end{eqnarray*}
defined as follows. 
\begin{enumerate}[\rm(a)]
\item If $\di$ is a gr-Sylvester module rank function for $R$ and  
$f\colon P\rightarrow Q$ is a homomorphism of $\Gamma$-graded $R$-modules with
$P,Q$ finitely generated $\Gamma$-graded projective $R$-modules, then  
$\rho_{\di}(f)=\di(Q)-\di(\coker f)$.  

\item Conversely, let  $\rho$ be a gr-Sylvester map rank function and suppose that \linebreak
$f\colon P\rightarrow Q$ is a homomorphism of $\Gamma$-graded $R$-modules with
$P,Q$ finitely generated $\Gamma$-graded projective $R$-modules such that
$\coker f=M$. Then \linebreak $\di_\rho(M)=\rho(1_Q)-\rho(f)$.
\end{enumerate}
\end{theorem}

\begin{proof}
First we show that the correspondence is an anti-isomorphism of partially ordered sets.
Let $\rho_1\leq\rho_2$ be two gr-Sylvester map rank functions. Let $M$ be a finitely presented
$\Gamma$-graded $R$-module and suppose that $R^n(\overline{\beta})\stackrel{f}{\rightarrow} R^m(\overline{\alpha})\rightarrow M\rightarrow 0$ is
a presentation  of $M$ as $\Gamma$-graded $R$-module. Then $\rho_1(1_{R^m(\overline{\alpha})})=
\rho_2(1_{R^m(\overline{\alpha})})=m$ and $\rho_1(f)\leq \rho_2(f)$. Hence 
$\di_{\rho_1}(M)=\rho_1(1_{R^m(\overline{\alpha})})-\rho_1(f)\geq \rho_2(1_{R^m(\overline{\alpha})})-\rho_2(f)=\di_{\rho_2}(M)$.
Conversely, let $\di_1\leq\di_2$ be gr-Sylvester module rank functions for $R$. Let
$f\colon P\rightarrow Q$ be a homomorphism of $\Gamma$-graded $R$-modules with
$P,Q$ finitely generated $\Gamma$-graded projective $R$-modules. Let $P',Q'$ be
finitely generated $\Gamma$-graded projective $R$-modules such that
$P\oplus P'\cong R^n(\overline{\beta})$ and $Q\oplus Q'\cong R^m(\overline{\alpha})$ where
$m,n$ are positive integers and $\overline{\beta}\in \Gamma^n$, $\overline{\alpha}\in \Gamma^m$.
Consider $\left(\begin{smallmatrix} f & 0 \\ 0 & 0 \end{smallmatrix}\right)\colon P\oplus P'\rightarrow
Q\oplus Q'$. Then
\begin{eqnarray*}
\rho_{\di_1}(f)&=&\rho_{\di_1}\left(\begin{smallmatrix} f & 0 \\ 0 & 0 \end{smallmatrix}\right)=
\di_1(Q\oplus Q')-\di_1(\coker\left(\begin{smallmatrix} f & 0 \\ 0 & 0 \end{smallmatrix}\right))=
m-\di_1(\coker\left(\begin{smallmatrix} f & 0 \\ 0 & 0 \end{smallmatrix}\right)) \\
&\geq& m -\di_2(\coker\left(\begin{smallmatrix} f & 0 \\ 0 & 0 \end{smallmatrix}\right))=
\di_2(Q\oplus Q')- \di_2(\coker\left(\begin{smallmatrix} f & 0 \\ 0 & 0 \end{smallmatrix}\right))=
\rho_{\di_2}\left(\begin{smallmatrix} f & 0 \\ 0 & 0 \end{smallmatrix}\right)\\ &=&\rho_{\di_2}(f).
\end{eqnarray*}

Secondly, we show that the correspendences are one inverse of the other. Let $\di$ be a gr-Sylvester module rank
function. Let $M$ be a finitely presented $\Gamma$-graded $R$-module. Then, given a graded presentation
of $M$, $P\stackrel{f}{\rightarrow} Q\rightarrow M\rightarrow 0$, 
$\di_{\rho_{\di}}(M)=\rho_{\di}(1_Q)-\rho_{\di}(f)=\di(Q)-0-(\di(Q)-\di(M))=\di(M)$. Conversely,
let $\rho$ be a gr-Sylvester map rank function. Let $f\colon P\rightarrow Q$ be a homomorphism
of $\Gamma$-graded $R$-modules with $P,Q$ finitely generated $\Gamma$-graded projective $R$-modules.
Then $\rho_{\di_{\rho}}(f)=\di_{\rho}(Q)-\di_{\rho}(\coker f)=\rho(1_Q)-0-(\rho(1_Q)-\rho(f))=\rho(f)$,
as desired. 

(a) Suppose that $\di$ is a gr-Sylvester module rank function. For a homomorphism of $\Gamma$-graded
$R$-modules $f\colon P\rightarrow Q$ between finitely generated $\Gamma$-graded projective $R$-modules define $\rho_{\di}(f)=\di(Q)-\di(\coker f)$. Notice it is clearly well defined. We must show that $\rho_{\di}$ satisfies
(MapRF1)--(MapRF5).

Clearly $\rho_{\di}(1_R)=\di(R)-\di(0)=\di(R)=1$. Thus (MatRF1) is satisfied.

Let $g\colon Q\rightarrow T$ be a graded homomorphism between finitely generated $\Gamma$-graded $R$-modules.
By definition, $\rho_{\di}(gf)=\di(T)-\di(\coker(gf))$. From the natural homomorphism of graded modules
$T/\im(gf)\rightarrow T/\im g\rightarrow 0\rightarrow 0$, by (ModRF3), we get that
$\di(\coker g)\leq \di(\coker(gf))$. Thus $\rho_{\di}(gf)\leq\di(T)-\di(\coker g)=\rho_{\di}(g)$. Let
$\pi_1\colon Q\rightarrow Q/\im f$ and $\pi_2\colon T\rightarrow T/\im (gf)$ be the natural
homomorphisms (of $\Gamma$-graded $R$-modules), and consider $\overline{g}\colon Q/\im f\rightarrow T/\im(gf)$ induced from $g\colon Q\rightarrow T$. It is not difficult to show that the following sequence of
homomorphisms of $\Gamma$-graded $R$-modules is exact
$$Q\stackrel{\left(\begin{smallmatrix} g\\ \pi_1\end{smallmatrix}\right)}{\longrightarrow}T\oplus Q/\im f\stackrel{(-\pi_2\ \overline{g})}{\longrightarrow} T/\im(gf)\longrightarrow 0.$$
Hence, by (ModRF3), $\di(\coker(gf))\leq \di(T)+\di(\coker f)\leq \di(\coker (gf))+\di(Q)$. Now,
subtracting $\di(\coker(gf))+\di(\coker f)$ on both sides of the second inequality we obtain
$\rho_{\di}(gf)=\di(T)-\di(\coker gf)\leq \di(Q)-\di(\coker f)=\rho_{\di}(f)$. Thus (MapRF2) is satisfied.

Let $f\colon P\rightarrow Q$, $f'\colon P'\rightarrow Q'$ be homomorphisms of $\Gamma$-graded $R$-modules
where $P,Q,P',Q'$ are $\Gamma$-graded projective $R$-modules. Now
\begin{eqnarray*}
\rho_{\di}\left(\left(\begin{smallmatrix}f & 0 \\ 0 & f' \end{smallmatrix}\right)\right) & = &
\di(Q\oplus Q')-\di\left(\coker \left(\begin{smallmatrix}f & 0 \\ 0 & f' \end{smallmatrix}\right)\right) \\
& = & \di(Q)+\di(Q')-\di(\coker f\oplus \coker f')\\
& = & \di(Q)-\di(\coker f)+\di(Q')-\di(\coker f') \\
& = & \rho_{\di}(f)+\rho_{\di}(f').
\end{eqnarray*}
Thus (MapRF3) is satisfied.

Let $h\colon P'\rightarrow Q$ be a homomorphism of $\Gamma$-graded $R$-modules and consider the 
homomorphism of $\Gamma$-graded $R$-modules
$\left(\begin{smallmatrix} f & h \\ 0 & g \end{smallmatrix}\right)\colon
P\oplus P'\rightarrow Q\oplus Q'$. Let $\overline{\iota}\colon Q/\im f\rightarrow 
\frac{Q\oplus Q'}{\im \left(\begin{smallmatrix} f & h \\ 0 & g \end{smallmatrix}\right)}$
be induced from the natural inclusion $Q\rightarrow Q\oplus Q'$. Let
$\overline{\pi}\colon \frac{Q\oplus Q'}{\im \left(\begin{smallmatrix} f & h \\ 0 & g \end{smallmatrix}\right)}
\rightarrow Q'/\im f'$ be induced from the natural projection $Q\oplus Q'\rightarrow Q'$. It is not
difficult to show that 
$$\frac{Q}{\im f}\longrightarrow \frac{Q\oplus Q'}{\im\left(\begin{smallmatrix} f & h \\ 0 & g \end{smallmatrix}\right) }\longrightarrow \frac{Q'}{\im f'}\longrightarrow 0$$
is an exact sequence of finitely presented $\Gamma$-graded $R$-modules. From (ModRF3), one obtains that
$\di(\coker \left(\begin{smallmatrix} f & h \\ 0 & g \end{smallmatrix}\right))\leq \di(\coker f)+\di(\coker f')$. Therefore
\begin{eqnarray*}
\rho_{\di}\left( \left(\begin{smallmatrix} f & h \\ 0 & g \end{smallmatrix}\right)\right) & = & 
\di(Q\oplus Q')-\di\left(\coker \left(\begin{smallmatrix} f & h \\ 0 & g \end{smallmatrix}\right)\right) \\
& \geq & \di(Q)+\di(Q')-\di(\coker f)-\di(\coker g) \\
& = & \rho_{\di}(f)+\rho_{\di}(g),
\end{eqnarray*}
and (MapRF4) is satisfied. 

Let now $f\colon R^n(\overline{\beta})\rightarrow R^m(\overline{\alpha})$ and
$f'\colon R^n(\overline{\beta'})\rightarrow R^m(\overline{\alpha'})$ be homomorphisms of
$\Gamma$-graded $R$-modules such that $\mathcal{F}(f)=\mathcal{F}(f')$. By (ModRF4), 
$\di(\coker f)=\di(\coker f')$ and $\di(R^m(\overline{\alpha}))=\di(R^m(\overline{\alpha'}))$.
Hence $\rho_{\di}(f)=\rho_{\di}(f')$ and (MapRF5) is satisfied.
Therefore $\rho_{\di}$ is a gr-Sylvester map rank function.

\medskip

(b) Suppose that $\rho$ is a gr-Sylvester map rank function. If $M$ is a finitely presented
$\Gamma$-graded $R$-module and $f\colon P\rightarrow Q$ is a homomorphism of
$\Gamma$-graded $R$-modules with $P,Q$ finitely generated $\Gamma$-graded projective $R$-modules
such that $\coker f=M$, then we define $\di_\rho(M)=\rho(1_Q)-\rho(f)$. We must show
that $\di$ is well defined and satisfies (ModRF1)--(ModRF4).

We begin showing that $\di_\rho$ is well defined. Suppose that $P\stackrel{f}{\rightarrow}Q\rightarrow M
\rightarrow 0$ and $P'\stackrel{f'}{\rightarrow}Q'\rightarrow M
\rightarrow 0$ are two graded presentations with $P,P',Q,Q'$ finitely generated
$\Gamma$-graded projective $R$-modules. By Lemma~\ref{lem:gradedSchanuel}, there exists
an automorphism $h$ of $\Gamma$-graded $R$-modules of $Q\oplus Q'$ which maps
$\im f\oplus Q'$ onto $Q\oplus \im f'$. Since $Q\oplus Q'$ is a $\Gamma$-projective $R$-module,
we obtain the following diagram of homomorphisms of $\Gamma$-graded $R$-modules
$$\xymatrix{P\oplus Q'\ar[r]\ar@<-1ex>[d]_u\ar@/^1pc/[rr]^{\left(\begin{smallmatrix} f & 0 \\
0 & 1_{Q'} \end{smallmatrix}\right)}
 & \im f \oplus Q'\ar@{^{(}->}[r]\ar@<-1ex>[d]_h & Q\oplus Q'\ar@<-1ex>[d]_h   \\
 Q\oplus P'\ar[r]\ar@<-1ex>[u]_{u'}\ar@/_1pc/[rr]_{\left(\begin{smallmatrix} 1_Q & 0 \\
0 & f' \end{smallmatrix}\right)}  & Q \oplus \im f'\ar@{^{(}->}[r]\ar@<-1ex>[u]_{h^{-1}} & Q\oplus Q'\ar@<-1ex>[u]_{h^{-1}} },$$
Hence we obtain $h\left(\begin{smallmatrix}f & 0 \\ 0 & 1_{Q'} \end{smallmatrix}\right)
=\left(\begin{smallmatrix}1_Q & 0 \\ 0 & f' \end{smallmatrix}\right)u$ and
$h^{-1}\left(\begin{smallmatrix}1_Q & 0 \\ 0 & f' \end{smallmatrix}\right)=\left(\begin{smallmatrix}f & 0 \\ 0 & 1_{Q'} \end{smallmatrix}\right) u'$ and, by (MapRF2),
$$\rho(f)+\rho(1_{Q'})=\rho\left(h\left(\begin{smallmatrix}f & 0 \\ 0 & 1_{Q'} \end{smallmatrix}\right)\right)\leq \rho\left(\left(\begin{smallmatrix}1_Q & 0 \\ 0 & f' \end{smallmatrix}\right)\right)
= \rho(f')+\rho(1_Q),$$
$$\rho(f')+\rho(1_Q)=\rho\left(h^{-1}\left(\begin{smallmatrix}1_Q & 0 \\ 0 & f' \end{smallmatrix}\right) \right)\leq \rho\left(\left(\begin{smallmatrix}f & 0 \\ 0 & 1_{Q'} \end{smallmatrix}\right)\right)=
\rho(f)+\rho(1_{Q'}).$$
It implies $\rho(1_{Q'})-\rho(f')\leq \rho(1_Q)-\rho(f)$ and
$\rho(1_Q)-\rho(f)\leq \rho(1_{Q'})-\rho(f')$. Therefore
$\di_{\rho}(M)=\rho(1_Q)-\rho(f)$ is well defined.

Consider the exact sequence $0\stackrel{0}{\rightarrow} R\stackrel{1_R}{\rightarrow} R\rightarrow 0.$ By definition,
$\di_\rho(R)=\rho(1_R)-\rho(0)=1-0=1$. Thus (ModRF1) is satisfied.

Let $M_1,M_2$ be finitely presented $\Gamma$-graded $R$-modules. Let 
$P_1\stackrel{f_1}{\rightarrow} Q_1\rightarrow M_1\rightarrow 0$ and
$P_2\stackrel{f_2}{\rightarrow} Q_2\rightarrow M_2\rightarrow 0$ be exact sequences of 
homomorphisms of $\Gamma$-graded $R$-modules
with $P_1,P_2,Q_1,Q_2$ finitely generated $\Gamma$-graded projective $R$-modules. Then
$P_1\oplus P_2\stackrel{\left(\begin{smallmatrix}f_1 & 0 \\ 0 & f_2 \end{smallmatrix}\right)}{\longrightarrow} Q_1\oplus Q_2\longrightarrow M_1\oplus M_2\rightarrow 0$ is a $\Gamma$-graded
presentation of $M_1\oplus M_2$. By definition and (MapRF3),
$\di_\rho(M_1\oplus M_2)=\rho(1_{Q_1\oplus Q_2})-\rho\left(\left(\begin{smallmatrix}f_1 & 0 \\ 0 & f_2 \end{smallmatrix}\right)\right)=\rho(1_{Q_1})-\rho(f_1)+\rho(1_{Q_2})-\rho(f_2)=\di_{\rho}(M_1)+\di_\rho(M_2)$.
Hence (ModRF2) is satisfied. 

Let $M_1\stackrel{u}{\rightarrow}M_2\stackrel{v}{\rightarrow}M_3\rightarrow 0$ be an exact sequence of
homomorphisms of $\Gamma$-graded $R$-modules with $M_1,M_2,M_3$ finitely presented $\Gamma$-graded $R$-modules. Let $P_1\stackrel{f_1}{\rightarrow}Q_1\stackrel{g_1}{\rightarrow}M_1\rightarrow 0$ and
$P_3\stackrel{f_3}{\rightarrow}Q_3\stackrel{g_3}{\rightarrow}M_3\rightarrow 0$ be exact sequences of homomorphisms of $\Gamma$-graded $R$-modules with $P_1,Q_1,P_3,Q_3$ finitely generated $\Gamma$-graded
projective $R$-modules. By diagram chasing, it is easy to obtain the following commutative
diagram, with exact rows and columns, of homomorphisms of $\Gamma$-graded $R$-modules where
$\iota$ and $\pi$ are the natural inclusion and projection, respectively
$$
\xymatrix{ & & 0\ar[d] & 0\ar[d] & \\
 & &\ker g\ar[r]^{\pi_{|\ker g}}\ar[d] & \ker g_3\ar[r]\ar[d] & 0 \\
0\ar[r] & Q_1\ar[r]^\iota\ar[d]_{g_1} & Q_1\oplus Q_3\ar[r]^\pi\ar[d]^g & Q_3 \ar[r]\ar[d]^{g_3} & 0 \\
 & M_1\ar[d]\ar[r]^u & M_2\ar[d]\ar[r]^v & M_3\ar[d]\ar[r] & 0 \\
& 0 & 0 & 0 & }$$
Since $\ker g$ is a finitely generated $\Gamma$-graded $R$-module, there exist a finitely generated
$\Gamma$-graded projective $R$-module and a surjective homomorphism of $\Gamma$-graded $R$-modules
$P\stackrel{f}{\rightarrow}\ker g$. In this way, we obtain a commutative diagram, with exact rows and columns, of homomorphisms of $\Gamma$-graded $R$-modules
$$\xymatrix{ & & P\ar[d]^f\ar@{=}[r] & P\ar[d]^{\pi f=f_3'} & \\
0\ar[r] & Q_1\ar[d]^{g_1}\ar[r]^\iota & Q_1\oplus Q_3\ar[d]^g\ar[r]^\pi & Q_3\ar[r]\ar[d]^{g_3} & 0 \\
& M_1\ar[r]^u\ar[d] & M_2\ar[d]\ar[r]^v & M_3\ar[r]\ar[d] & 0 \\
& 0 & 0 & 0 &   } $$
Note that $\iota f_1(P_1)\subseteq \ker g$ and that $f\colon P\rightarrow Q_1\oplus Q_3$ is of the form
$\left(\begin{smallmatrix}\lambda \\ f_3' \end{smallmatrix}\right)$ for some
homomorphism of $\Gamma$-graded $R$-modules $\lambda\colon P\rightarrow Q_1$. Thus we can modify the foregoing diagram
to obtain $$\xymatrix{ & & P_1\oplus P\ar[d]^{\scriptsize\left(\begin{smallmatrix}
f_1 & \lambda \\ 0 & f_3' \end{smallmatrix}\right)}\ar[r]^\pi & P\ar[d]^{f_3'} & \\
0\ar[r] & Q_1\ar[d]^{g_1}\ar[r]^\iota & Q_1\oplus Q_3\ar[d]^g\ar[r]^\pi & Q_3\ar[r]\ar[d]^{g_3} & 0 \\
& M_1\ar[r]^u\ar[d] & M_2\ar[d]\ar[r]^v & M_3\ar[r]\ar[d] & 0 \\
& 0 & 0 & 0 &   } $$
Then, by (MapRF4), 
\begin{eqnarray*}
\di_{\rho}(M_2) & = & \rho(1_{Q_1\oplus Q_3})-\rho\left(\left(\begin{smallmatrix}
f_1 & \lambda \\ 0 & f_3' \end{smallmatrix}\right)\right) \\  & \leq & \rho(1_{Q_1})+\rho(1_{Q_3})
-\rho(f_1)-\rho(f_3') \\ & = & \di_{\rho}(M_1)+\di_{\rho}(M_3).
\end{eqnarray*}
Moreover,
\begin{eqnarray*}
\di_{\rho}(M_2) & = & \rho(1_{Q_1\oplus Q_3})-\rho\left(\left(\begin{smallmatrix}
f_1 & \lambda \\ 0 & f_3' \end{smallmatrix}\right)\right) \\
& \geq & \rho(1_{Q_1}) + \rho(1_{Q_3})-\rho(f'_3)-\rho((f_1\ \lambda)) \\
& = & \di_{\rho}(M_3)+\rho(1_{Q_1})-\rho((f_1\ \lambda)) \\
& \geq & \di_{\rho}(M_3),
\end{eqnarray*}
where we have used the fact that $\left(\begin{smallmatrix}
f_1 & \lambda \\ 0 & f_3' \end{smallmatrix}\right)=\left(\begin{smallmatrix}
f_1 & \lambda & 0 \\ 0& 0 & f_3' \end{smallmatrix}\right)\left(\begin{smallmatrix}
1_{P_1} & 0 \\ 0 & 1_P \\ 0 & 1_P \end{smallmatrix}\right)$, and properties
(MapRF2), (MapRF3)  on the first inequality, and Lemma~\ref{lem:grmapproperties}(6) on the second inequality. Therefore (ModRF3) is satisfied.

(ModRF4) follows easily. Indeed, let $f\colon R^n(\overline{\beta})\rightarrow R^m(\overline{\alpha})$,
$f'\colon R^n(\overline{\beta'})\rightarrow R^m(\overline{\alpha'})$ be homomorphisms
of $\Gamma$-graded $R$-modules such that $\mathcal{F}(f)=\mathcal{F}(f')$. By (MapRF5),
$\rho(f)=\rho(f')$ and therefore $\di_{\rho}(\coker f)=\rho(1_{R^m(\overline{\alpha})})-\rho(f)=
\rho(1_{R^m(\overline{\alpha'})}) -\rho(f')=\di_{\rho}(\coker f')$.
\end{proof}

The next result was given in \cite[Theorem~4]{Malcolmsondetermining}
in the ungraded context.

\begin{theorem}\label{theo:module_matrix_equivalent}
Let $R$ be a $\Gamma$-graded ring. There is an anti-isomorphism of partially ordered sets
\begin{eqnarray*}
\left\{\begin{array}{c}
\textrm{gr-Sylvester matrix rank } \\
\textrm{functions for }R 
\end{array}\right\}&\longrightarrow & \left\{\begin{array}{c}
\textrm{gr-Sylvester module rank } \\ 
\textrm{functions for }R
\end{array}\right\} \\
\rank & \longmapsto & \di_{\rank} \\
\rank_{\di} & \longmapsfrom & \di
\end{eqnarray*}
defined as follows. 
\begin{enumerate}[\rm(a)]
\item If $\rank$ is a gr-Sylvester matrix rank function for $R$
and $M$ is a finitely presented $\Gamma$-graded $R$-module with presentation
$R^n(\overline{\beta})\stackrel{A}{\rightarrow} R^m(\overline{\alpha})\rightarrow M\rightarrow 0$,
where $A\in M_{m\times n}(R)[\overline{\alpha}][\overline{\beta}]$, we define $\di_{\rank}(M)=m-\rank(A)$.

\item Conversely, let $\di$ is a gr-Sylvester module rank function for $R$. If  
$A\in M_{m\times n}(R)[\overline{\alpha}][\overline{\beta}]$,  we consider
$A$ as a homomorphism of $\Gamma$-graded $R$-modules $R^n(\overline{\beta})\rightarrow
R^m(\overline{\alpha})$ and define
$\rank(A)=m-\di(R^m(\overline{\alpha})/A(R^n(\overline{\beta})))$.
\end{enumerate}
\end{theorem}

\begin{proof}
First we show that the correspondence is an anti-isomorphism of partially ordered sets.
Let $\rank_1\leq\rank_2$ be two gr-Sylvester matrix rank functions. Let $M$ be a finitely presented
$\Gamma$-graded $R$-module and suppose that $R^n(\overline{\beta})\stackrel{A}{\rightarrow} R^m(\overline{\alpha})\rightarrow M\rightarrow 0$ is
a presentation  of $M$ as $\Gamma$-graded $R$-module. Then $\rank_1(1_{R^m(\overline{\alpha})})=
\rank_2(1_{R^m(\overline{\alpha})})=m$ and $\rank_1(A)\leq \rank_2(A)$. Hence 
$\di_{\rank_1}(M)=\rank_1(1_{R^m(\overline{\alpha})})-\rank_1(A)\geq \rank_2(1_{R^m(\overline{\alpha})})-\rank_2(f)=\di_{\rank_2}(M)$.
Conversely, let $\di_1\leq \di_2$ be two gr-Sylvester module rank functions for $R$. Let $A\in M_{m\times n}(R)[\overline{\alpha}][\overline{\beta}]$ and consider $A$ as a homomorphism of $\Gamma$-graded $R$-modules 
$R^n(\overline{\beta})\rightarrow R^m(\overline{\alpha})$. Now $\rank_{\di_1}(A)= m-\di_1(R^m(\overline{\alpha})/A(R^n(\overline{\beta}))) \leq m-\di_2(R^m(\overline{\alpha})/A(R^n(\overline{\beta}))) =\rank_{\di_2}(A)$. Therefore $\rank_{\di_1} \leq \rank_{\di_2}.$

Secondly, we show that the correspondences are one inverse of the other. Let $\di$ be a gr-Sylvester module rank
function. Let $M$ be a finitely presented $\Gamma$-graded $R$-module. Then, given a graded presentation
of $M$, $R^n(\overline{\alpha})\stackrel{A}{\rightarrow} R^m(\overline{\beta})\rightarrow M\rightarrow 0$, 
$\di_{\rank_{\di}}(M)=\rank_{\di}(1_{R^m(\overline{\alpha})})-\rank_{\di}(A)=\di(R^m(\overline{\alpha}))-0-(\di(R^m(\overline{\alpha}))-\di(M))=\di(M)$. Conversely,
let $\rank$ be a gr-Sylvester matrix rank function. Let $A\in M_{m\times n}(R)[\overline{\alpha}][\overline{\beta}]$. Then $\rank_{\di_{\rho}}(A)=\di_{\rank}(R^m({\overline{\alpha}}))-\di_{\rank}(\coker A)=\rank(1_{R^m(\overline{\alpha})})-0-(\rank(1_{R^m(\overline{\alpha})})-\rank(A))=\rank(A)$,
as desired. 

(a) Suppose that $\di$ is a Sylvester module rank function for $R$. Let $A\in M_{m\times n}(R)[\overline{\alpha}][\overline{\beta}]$ and consider $A$ as a homomorphism of $\Gamma$-graded $R$-modules 
$R^n(\overline{\beta})\rightarrow R^m(\overline{\alpha})$. Define $\rank_{\di}(A)=\di(R^m(\overline{\alpha}))-\di(\coker A)=m-\di(R^m(\overline{\alpha})/A(R^n(\overline{\beta}))).$
The fact that $\rank_{\di}$ is well defined follows from (ModRF4). That is, the matrix $A$ may define different homomorphisms of $\Gamma$-graded modules, but the value of $\rank_{\di}(A)$ is the same. 
The proof that $\rank_{\di}$ satisfies (MatRF1)--(MatRF4) follows in the same way as the proof that
$\rho_{\di}$ satisfies (MapRF1)--(MapRF5) in Theorem~\ref{theo:module_map_equivalent}.

(b) Suppose now that $\rank$ is a gr-Sylvester matrix rank function for $R$. Let 
$M$ be a finitely presented $\Gamma$-graded $R$-module with presentation
$R^n(\overline{\beta})\stackrel{A}{\rightarrow} R^m(\overline{\alpha})\rightarrow M\rightarrow 0$,
where $A\in M_{m\times n}(R)[\overline{\alpha}][\overline{\beta}]$, we define $\di_{\rank}(M)=m-\rank(A)$.
One can show that $\di_{\rank}$ is well defined and satisfies (ModRF1)--(ModRF4) in the same way
that one proves that $\di_{\rho}$ is well defined and satisfies (ModRF1)--(ModRF4).
There is another way to prove that $\di_{\rank}$ is well defined and satisfies (ModRF1)--(ModRF4). 
By Theorem~\ref{theo:grprimematrix_grSylvestermatrix}, let $(K,\varphi)$ be the corresponding
 $\Gamma$-graded epic $R$-division ring with $\rank$. Then,  for a finitely presented $\Gamma$-graded $R$-module $M$, the
$K$-module $M\otimes_R K$ is a $\Gamma$-graded free $K$-module. One can define $\di(M)=\dim_K(M\otimes_R K)$.
It is not difficult to show that $\di$ satisfies (ModRF1)--(ModRF4). 
Now let 
$M$ be a finitely presented $\Gamma$-graded $R$-module with presentation
$R^n(\overline{\beta})\stackrel{A}{\rightarrow} R^m(\overline{\alpha})\rightarrow M\rightarrow 0$,
where $A\in M_{m\times n}(R)[\overline{\alpha}][\overline{\beta}]$. Then, by Theorem~\ref{theo:grprimematrix_grSylvestermatrix}, $\dim (M\otimes_R K)$ equals
$m$ minus the number of columns of $A^\varphi$ which are right linearly independent over $K$ and that is
exactly $m-\rank(A)=\di_{\rho}(M)$.
\end{proof}

It is worth noting the  following corollary.  It is just a re-writing of parts of Theorems~\ref{theo:specialization}, 
\ref{theo:module_map_equivalent},  
\ref{theo:module_matrix_equivalent} and Corollary~\ref{coro:specialization}. 

\begin{corollary}
Let $R$ be a $\Gamma$-graded ring. Suppose that $(K_1,\varphi_1)$ and
$(K_2,\varphi_2)$ are $\Gamma$-graded epic $R$-division rings with corresponding
gr-prime matrix ideals $\mathcal{P}_1$, $\mathcal{P}_2$, respectively, and
gr-Sylvester rank functions $\rank_1$, $\rank_2$, $\di_1$, $\di_2$, $\rho_1$, $\rho_2$, respectively.
The following assertions are equivalent
\begin{enumerate}[\rm(1)]
	\item There exists a gr-specialization from $(K_1,\varphi_1)$ to $(K_2, \varphi_2)$.
	\item $\mathcal{P}_1\subseteq \mathcal{P}_2$.
	\item $\rank_1\geq \rank_2$.
	\item $\di_1\leq \di_2$.
	\item $\rho_1\geq \rho_2$.
\end{enumerate}
\end{corollary}



\section{gr-prime spectrum}\label{sec:grprimespectrum}

\emph{Throughout this section, 
let $\Gamma$ be a group and $\Omega\subseteq \Omega'$ be normal subgroups of $\Gamma$}.

\medskip

Let $R$ be a $\Gamma$-graded ring. It can be considered as a $\Gamma/\Omega$-graded ring too.
Now we introduce some notation in order to clarify which structure of graded object is being considered. 
We will denote by $\mathfrak{M}^\Gamma(R)$, $\mathfrak{M}^\Gamma_\bullet(R)$ and by $\mathfrak{M}^{\Gamma/\Omega}(R)$, $\mathfrak{M}^{\Gamma/\Omega}_\bullet(R)$ the corresponding sets of matrices. Notice
that $\mathfrak{M}^\Gamma(R)\subseteq \mathfrak{M}^{\Gamma/\Omega}(R)$ and 
 $\mathfrak{M}^\Gamma_\bullet(R)\subseteq \mathfrak{M}^{\Gamma/\Omega}_\bullet(R)$.
 We denote by $\Spec_\Gamma(R)$ the set of all $\Gamma$-gr-prime
matrix ideals and by $\Spec_{\Gamma/\Omega}(R)$ the set of all $\Gamma/\Omega$-gr-prime matrix ideals.
If $\Omega=\Gamma$, we will write $\Spec(R)$ instead of $\Spec_{\Gamma/\Gamma}(R)$. Note that
$\Spec(R)$ is the usual set of prime matrix ideals.

It follows directly from the definition that if $\mathcal{P}$ is a $\Gamma/\Omega$-gr-prime matrix ideal,
then $\mathcal{P}\cap \mathfrak{M}^\Gamma(R)$ is a $\Gamma$-gr-prime matrix ideal. Hence, there exists a map
$$\Spec_{\Gamma/\Omega}(R)\rightarrow \Spec_{\Gamma}(R),\quad \mathcal{P}\mapsto 
\mathcal{P}\cap\mathfrak{M}^\Gamma(R).$$
Suppose now that $\rank^{\Gamma/\Omega}\colon \mathfrak{M}_\bullet^{\Gamma/\Omega}(R)\rightarrow
\mathbb{N}$ is a Sylvester $\Gamma/\Omega$-gr-matrix rank function. It follows from the definition that the restriction of $\rank^{\Gamma/\Omega}$ to $\mathfrak{M}_\bullet^\Gamma(R)$
induces a $\Gamma$-gr-matrix rank function $\rank^\Gamma\colon \mathfrak{M}_\bullet^\Gamma(R)\rightarrow
\mathbb{N}$. In this way, there is a function
\begin{align*}
\left\{\begin{array}{c}
\Gamma/\Omega\textrm{-gr-Sylvester matrix} \\
\textrm{rank functions for }R 
\end{array}\right\} & & \longrightarrow & & \left\{\begin{array}{c}
\Gamma\textrm{-gr-Sylvester matrix} \\ 
\textrm{rank functions for }R
\end{array}\right\}, \quad \rank^{\Gamma/\Omega}\mapsto \rank^\Gamma.
\end{align*}
Similarly, since any $\Gamma$-graded (finitely generated, finitely presented, projective) $R$-module is also a $\Gamma/\Omega$-graded (finitely generated, finitely presented, projective)
$R$-module and any homomorphism between $\Gamma$-graded $R$-modules is also 
a homomorphism of $\Gamma/\Omega$-graded $R$-modules, by restriction, we obtain functions
 \begin{align*}
\left\{\begin{array}{c}
\Gamma/\Omega\textrm{-gr-Sylvester module} \\
\textrm{rank functions for }R 
\end{array}\right\} & & \longrightarrow & & \left\{\begin{array}{c}
\Gamma\textrm{-gr-Sylvester module} \\ 
\textrm{rank functions for }R
\end{array}\right\}, \quad \di^{\Gamma/\Omega}\mapsto \di^\Gamma.
\end{align*}
\begin{align*}
\left\{\begin{array}{c}
\Gamma/\Omega\textrm{-gr-Sylvester map} \\
\textrm{rank functions for }R 
\end{array}\right\} & & \longrightarrow & & \left\{\begin{array}{c}
\Gamma\textrm{-gr-Sylvester map} \\ 
\textrm{rank functions for }R
\end{array}\right\}, \quad \rho^{\Gamma/\Omega}\mapsto \rho^\Gamma.
\end{align*}

Considering $R$ as a $\Gamma/\Omega$-graded ring
and $\Omega'/\Omega$ as a normal subgroup of $\Gamma/\Omega$, we obtain maps
$\Spec_{\Gamma/\Omega'}(R)\rightarrow\Spec_{\Gamma/\Omega}(R)$ for each pair of such normal subgroups of $\Gamma$.
Hence if $\Spec_{\Gamma/\Omega}(R)$ is empty, then $\Spec_{\Gamma/\Omega'}(R)$ is also empty. In other words,
if there does not exist a $\Gamma/\Omega$-graded epic $R$-division ring, then there does not
exist a $\Gamma/\Omega'$-graded epic $R$-division ring.
In particular, for $\Omega'=\Gamma$ we obtain maps $\Spec(R)\rightarrow\Spec_{\Gamma/\Omega}(R)$,
 $\mathcal{Q}\mapsto\mathcal{Q}\cap\mathfrak{M}^{\Gamma/\Omega}(R)$,
for each normal subgroup $\Omega$ of $\Gamma$. Therefore, if there exist a normal subgroup
$\Omega$ of $\Gamma$ such that there does not exist a $\Gamma/\Omega$-graded epic $R$-division ring, then
there does not exist an epic $R$-division ring.

Let $\mathcal{Q}'\in\Spec_{\Gamma/\Omega'}(R)$ and let $\mathcal{Q}=\mathcal{Q}'\cap\mathfrak{M}^{\Gamma/\Omega}(R)$ be the corresponding element in $\Spec_{\Gamma/\Omega}(R)$. Let $(K_\mathcal{Q'},\varphi_\mathcal{Q'})$
be the $\Gamma/\Omega'$-graded epic $R$-division ring determined by $\mathcal{Q'}$, and let
$(K_\mathcal{Q},\varphi_\mathcal{Q})$ be the $\Gamma/\Omega$-graded epic $R$-division ring
determined by $\mathcal{Q}$. 
Let $x$ be a
homogeneous element of $R$ considered as a $\Gamma/\Omega$-graded ring. Notice it is also
a homogeneous element of $R$ considered as a $\Gamma/\Omega'$-graded ring.
If $x\notin\ker \varphi_{\mathcal{Q}'}$, then 
$x\in\mathfrak{M}^{\Gamma/\Omega'}(R)\setminus\mathcal{Q}'$. Thus
$x\in\mathfrak{M}^{\Gamma/\Omega}(R)\setminus\mathcal{Q}$, and therefore
$x\notin\ker\varphi_\mathcal{Q}$. Hence if $\varphi_\mathcal{Q'}$ is injective, then
$\varphi_\mathcal{Q}$ is also injective. In other words,
if $(K_{\mathcal{Q}'},\varphi_{Q'})$ is a $\Gamma/\Omega'$-graded epic $R$-division ring of 
fractions, then $(K_{\mathcal{Q}},\varphi_{Q})$ is also a $\Gamma/\Omega$-graded epic $R$-division ring of 
fractions. Therefore, if there exists a normal subgroup $\Omega$ of $\Gamma$ such that
there does not exist a $\Gamma/\Omega$-graded epic $R$-division ring of fractions, then there does not
exist an epic $R$-division ring of fractions.

Let $\mathcal{P'}\in\Spec_{\Gamma/\Omega'}(R)$ and set $\mathcal{P}=\mathcal{P}'\cap\mathfrak{M}^{\Gamma/\Omega}(R)\in
\Spec_{\Gamma/\Omega}(R)$. If $\mathcal{P}'\subseteq\mathcal{Q}'$, then
$\mathcal{P}\subseteq\mathcal{Q}$. Hence a specialization
from $(K_{\mathcal{P}'},\varphi_{{P}'})$ to $(K_{\mathcal{Q}'},\varphi_{{Q}'})$, implies
the existence of a specialization from $(K_{\mathcal{P}},\varphi_\mathcal{P})$
to $(K_\mathcal{Q},\varphi_\mathcal{Q})$ by Corollary~\ref{coro:specialization}. Notice
that it could happen that $\mathcal{Q}=\mathcal{P}$.

Also, if the map $\Spec_{\Gamma/\Omega'}(R)\rightarrow \Spec_{\Gamma/\Omega}(R)$ is surjective
and $R$ has a universal $\Gamma/\Omega'$-graded  epic $R$-division ring (of fractions),
then $R$ has a universal $\Gamma/\Omega$-graded epic $R$-division ring (of fractions).

Suppose that for each $\Gamma$-graded epic $R$-division ring $D$ 
there exist ring homomorphisms  to division rings. Then 
$\Spec_{\Gamma/\Omega}(R)\rightarrow\Spec_\Gamma(R)$ is surjective for each $\Omega\lhd \Gamma$.
Let $(D,\varphi)$ be a $\Gamma$-graded epic $R$-division ring with $\Gamma$-singular 
kernel $\mathcal{P}$. Let 
$\phi\colon D\rightarrow E$ be a ring homomorphism with $E$ a division ring.
Consider the composition $\phi\circ\varphi\colon R\rightarrow E$. It is
a homomorphism of $\Gamma/\Omega$-almost graded rings with $E$ a $\Gamma/\Omega$-almost graded division ring.
By Theorem~\ref{theo:gradedlocal}(2)(b), there exists $\psi\colon R\rightarrow D'$ a
$\Gamma/\Omega$-graded epic $R$-division ring, and a homomorphism $\rho\colon D'\rightarrow E$
such that $\phi\varphi=\rho\psi$. 
By Proposition~\ref{prop:almostdivisionring}, 
$$\{A\in\mathfrak{M}^{\Gamma}(R)\colon A^{(\phi\varphi)} \textrm{ invertible over } E\}=\{A\in\mathfrak{M}^{\Gamma}(R)\colon A^\varphi\textrm{ is invertible over } D\},$$
$$\{A\in\mathfrak{M}^{\Gamma/\Omega}(R)\colon A \textrm{ invertible over } D'\}=
\{A\in\mathfrak{M}^{\Gamma/\Omega}(R)\colon A^{(\rho\psi)} \textrm{ invertible over } E\}.$$
Now, since $\mathfrak{M}^{\Gamma}(R)\subseteq\mathfrak{M}^{\Gamma/\Omega}(R)$, we get that
$$\{A\in\mathfrak{M}^{\Gamma}(R)\colon A^\varphi\textrm{ inverts over } D\}=
\{A\in\mathfrak{M}^{\Gamma/\Omega}(R)\colon A \textrm{ inverts over } D'\}\cap \mathfrak{M}^\Gamma(R).$$
Hence, if $\mathcal{P'}$ is the $\Gamma/\Omega$-singular kernel of $(D',\psi)$, then
$\mathcal{P}=\mathcal{P'}\cap\mathfrak{M}^\Gamma(R)$.

We gather together what we have just proved in the following result.

\begin{theorem}\label{theo:relationsbetweenspectra}
Let $R$ be a $\Gamma$-graded ring.  The following assertions hold true.
\begin{enumerate}[\rm(1)]
	\item If there does not exist a $\Gamma/\Omega$-graded epic $R$-division ring (of fractions), then there does not
	exist a $\Gamma/\Omega'$-graded epic $R$-division ring (of fractions). Therefore, if there exists a normal subgroup
	$\Omega$ of $\Gamma$ such that there does not exist a $\Gamma/\Omega$-graded epic $R$-division ring (of fractions), then there does not 	exist an epic $R$-division ring (of fractions).
	\item Let $(K_\mathcal{P'},\varphi_\mathcal{P'})$, $(K_\mathcal{Q'},\varphi_\mathcal{Q'})$ be $\Gamma/\Omega'$-epic $R$-division rings, such that there exists a specialization from $(K_\mathcal{P'},\varphi_\mathcal{P'})$ to $(K_\mathcal{Q'},\varphi_\mathcal{Q'})$, then there exists a specialization
	between the corresponding $\Gamma/\Omega$-graded epic $R$-division rings. 
	\item If the map $\Spec_{\Gamma/\Omega'}(R)\rightarrow \Spec_{\Gamma/\Omega}(R)$,
	$\mathcal{Q'}\mapsto \mathcal{Q'}\cap\mathfrak{M}^{\Gamma/\Omega}(R)$, is surjective, then
	the existence of a universal $\Gamma/\Omega'$-graded epic $R$-division ring implies
	the existence of a universal $\Gamma/\Omega$-graded epic $R$-division ring. Therefore,
	if $\Spec(R)\rightarrow \Spec_{\Gamma/\Omega}(R)$,
	$\mathcal{Q'}\mapsto \mathcal{Q'}\cap\mathfrak{M}^{\Gamma/\Omega}(R)$, is surjective,
	the existence of a universal $R$-division ring, implies the existence of a universal
	$\Gamma/\Omega$-graded epic $R$-division ring.
	\item If for each $\Gamma$-graded epic $R$-division ring 
there exist ring homomorphisms  to division rings, then 
$\Spec_{\Gamma/\Omega}(R)\rightarrow\Spec_\Gamma(R)$, $\mathcal{Q}\mapsto \mathcal{Q}\cap
\mathfrak{M}^\Gamma(R)$, is surjective.
	\qed
\end{enumerate}
\end{theorem}

Note that there is a corresponding statement of Theorem~\ref{theo:relationsbetweenspectra}  in terms of the different gr-Sylvester rank functions instead of gr-prime matrix ideals.

\medskip

Let  $R=\bigoplus_{\gamma\in\Gamma}R_\gamma$ be a $\Gamma$-graded ring. 
In the foregoing, we gave a correspondence from the set of $\Gamma/\Omega$-graded epic
$R$-division rings to the set of $\Gamma$-graded epic $R$-division rings. We proceed to
give a more down to earth description of such correspondence.
Recall that
$R$ can be regarded as a $\Gamma/\Omega$-graded
ring making $R=\bigoplus_{\Delta\in\Gamma/\Omega}R_\alpha$ where 
$R_\alpha=\bigoplus_{\gamma\in\alpha}R_\gamma$ for each $\alpha\in\Gamma/\Omega$.

Let $E=\bigoplus_{\alpha\in\Gamma/\Omega}E_\alpha$ be a $\Gamma/\Omega$-graded division ring. 
Consider the group ring
$E[\Gamma]=\bigoplus_{\gamma\in\Gamma}E\gamma$. We construct a $\Gamma$-graded division ring $D=\bigoplus_{\gamma\in\Gamma}D_\gamma$
which is a $\Gamma$-graded subring of $E[\Gamma]$ in the same way as in 
\cite[Proposition~1.2.2]{NastasescuvanOystaeyenMethodsgraded}.
For each $\gamma\in\Gamma$, there exists a unique $\alpha\in\Gamma/\Omega$ such that $\gamma\in\alpha$.
Set $D_\gamma=E_\alpha \gamma\subseteq E\gamma$. Note that
$$D_\gamma D_{\gamma'}=E_\alpha\gamma E_{\alpha'}\gamma=E_\alpha E_{\alpha'}\gamma\gamma'\subseteq
E_{\alpha\alpha'}\gamma\gamma'=D_{\gamma\gamma'}.$$ 
Hence $D$ is a $\Gamma$-graded ring. Since $E$ is a $\Gamma/\Omega$-graded
division ring, any nonzero homogeneous element of $D$ is invertible.
Thus $D$ is a $\Gamma$-graded division ring.

Suppose that $(E,\varphi)$ is a $\Gamma/\Omega$-graded epic $R$-division ring. Let
$\gamma\in\Gamma$ and $\alpha\in\Gamma/\Omega$ such that $\gamma\in\alpha$.
For each $a_\gamma\in R_\gamma$, $\varphi(a_\gamma)\in E_\alpha$. Then define
$\psi(a_\gamma)=\varphi(a_\gamma)\gamma\in D_\gamma$. In this way, we
obtain a homomorphism of $\Gamma$-graded rings $\psi\colon R\rightarrow D$.

Let $A=(a_{ij})\in M_n(R)[\overline{\delta}][\overline{\varepsilon}].$ We claim that $A^\varphi$ is invertible in $E$
if and only if $A^\psi$ is invertible in $D$. Indeed, let $\alpha_i,\beta_j\in\Gamma/\Omega$
be such that $\delta_i\in\alpha_i$, $\varepsilon_j\in\beta_j$. Then
$A^\varphi=(b_{ij})$ with $b_{ij}\in E_{\alpha_i\beta_j^{-1}}$ and
$(A^\varphi)^{-1}=(c_{ij})$ with $c_{ij}\in E_{\beta_i\alpha_j^{-1}}$. 
Then $A^{\psi}=(b_{ij}\alpha_i\beta_j^{-1})$ is invertible in $D$ with inverse
$(A^{\psi})^{-1}=(c_{ij}\beta_i\alpha_j^{-1})$. Conversely, if $A^\psi$ is invertible
with inverse $(A^\psi)^{-1}=(d_{ij}\beta_i\alpha_j^{-1})$ where $d_{ij}\in E_{\beta_i\alpha_j^{-1}}$,
then $(A^\varphi)^{-1}=(d_{ij})$. 
Hence, let $\mathcal{P}\in\Spec_{\Gamma/\Omega}(R)$. If $(E,\varphi)$ is the 
$\Gamma/\Omega$-graded epic $R$-division ring associated to $\mathcal{P}$, then the 
$\Gamma$-graded epic $R$-division ring associated to $\mathcal{P}\cap\mathfrak{M}^\Gamma(R)$
is determined by the $\Gamma$-graded division ring $\psi\colon R\rightarrow D$.
That is, the $\Gamma$-graded epic $R$-division ring $\psi\colon R\rightarrow D'$ where
$D'$ the graded division ring generated by $\im\psi$.

\bigskip

Now we proceed to give an important family of examples of Theorem~\ref{theo:relationsbetweenspectra}(4).

Let $(\Gamma,<)$ be an ordered group. Let $D=\bigoplus_{\gamma\in\Gamma}D_\gamma$ be a
$\Gamma$-graded division ring. Given a map $f\colon\Gamma\rightarrow D$, let
$\supp f=\{\gamma\in\Gamma\colon f(\gamma)\neq 0\}$. We will write $f$ as a series. 
Thus, $f=\sum_{\gamma\in\Gamma}a_\gamma$ means that $f(\gamma)=a_\gamma\in D$ for each
$\gamma\in\Gamma$. Consider the set 
$$D((\Gamma;<))=\left\{f=\sum_{\gamma\in\Gamma}a_\gamma\colon a_\gamma\in D_\gamma
\textrm{ for all  }\gamma\in\Gamma,\ \supp f\textrm{ is well
ordered}\right\}.$$
$D((\Gamma;<))$ is an abelian group under the natural sum. That is,
for $f=\sum_{\gamma\in\Gamma}a_\gamma$, $f'=\sum_{\gamma\in\Gamma}a'_\gamma$, then
$$f+f'=\sum_{\gamma\in\Gamma}(a_\gamma+a'_\gamma).$$
One can then define the product in $D((\Gamma;<))$ as
$$ff'=\sum_{\gamma\in\Gamma}\left(\sum_{\delta\varepsilon=\gamma}a_\delta a'_{\varepsilon}\right).$$
These operations endow $D((\Gamma;<))$ with a ring structure. We regard
$D$ as a subring of $D((\Gamma;<))$ identifying $D$ with the series of $D((\Gamma;<))$ of finite support.
Malcev and Neumann independently showed that $D((\Gamma;<))$ is in fact a division ring 
\cite{Malcev}, \cite{Neumann}. Hence, we have just shown that for every $\Gamma$-graded division ring
there exists a homomorphism of rings to a division ring.

Now we proceed to show that every $D((\Gamma;<))$ contains a $\Gamma/\Omega$-graded division ring and that
it corresponds to $D$ via $\Spec_{\Gamma/\Omega}(R)\rightarrow \Spec_{\Gamma}(R)$.
Let $\Omega$ be a normal subgroup of $\Gamma$. Consider $D$ as a $\Gamma/\Omega$-graded ring.
For each $\alpha\in\Gamma/\Omega$, define the subset of $D((\Gamma;<))$
$$E_\alpha=\left\{f=\sum_{\gamma\in\Gamma}a_\gamma\in D((\Gamma;<))\colon
\supp f\subseteq \alpha  \right\}.$$
Note that $E_\alpha$ is an additive subgroup of $D((\Gamma;<))$.
Let $\alpha,\beta\in\Gamma/\Omega$. Suppose that $f=\sum_{\gamma\in\Gamma}a_\gamma\in E_\alpha$
and $f'=\sum_{\gamma\in\Gamma}a_\gamma'\in E_\beta$. Then
$$ff'=\sum_{\gamma\in\Gamma}\left(\sum_{\delta\varepsilon=\gamma}a_\delta a'_{\varepsilon}\right)\in
E_{\alpha\beta}.$$ Hence $E_\alpha E_\beta\subseteq E_{\alpha\beta}$.
Moreover, if $\alpha\in\Gamma/\Omega$, $$E_\alpha\cap \left(\sum_{\beta\in\Gamma/\Omega,\, 
\beta\neq \alpha} E_\beta \right)=\{0\},$$
because $\Gamma$ is the disjoint union $\Gamma=\bigcup_{\beta\in\Gamma/\Omega}\beta.$
Hence $E(\Omega)=\bigoplus_{\alpha\in\Gamma/\Omega}E_\alpha$ is a $\Gamma/\Omega$-graded ring. 
Furthermore, let $f=\sum_{\gamma\in\Gamma}a_\gamma\in E_\alpha$, $f\neq 0$.
Then $f$ is invertible in $D((\Gamma;<))$ with inverse
$$f^{-1}=\left(\sum_{n\geq0}(-1)^ng^n\right)a_{\gamma_0}^{-1}$$
where $\gamma_0=\min\supp f$ and $g=\sum_{\gamma\in\Gamma} a_{\gamma_0}^{-1}a_\gamma$.
Since $\supp\gamma\subseteq \alpha$, $\gamma_0\in\alpha$ and $\gamma_0^{-1}\in\alpha^{-1}$,
then $\supp g\subseteq E_e$ where $e$ denotes the identity element in $\Gamma/\Omega$.
Hence $\supp g^n\subseteq E_e$ for each integer $n\geq 0$
and $\supp(\sum_{n\geq 0}(-1)^ng^n)\subseteq E_e$. Thus, $\supp f^{-1}\subseteq \alpha^{-1}$.
Therefore $E(\Omega)$ is a $\Gamma/\Omega$-graded division ring and the embedding
$\phi_\Omega\colon D\hookrightarrow E(\Omega)$ is a homomorphism of $\Gamma/\Omega$-graded rings.
Let $D(\Omega)$ be the $\Gamma/\Omega$ graded division subring of $E(\Omega)$ generated by $D$.
Then $(D(\Omega),\phi_\Omega\colon D\hookrightarrow D(\Omega))$ is a $\Gamma/\Omega$-graded
epic $D$-division ring. 

Let $R=\bigoplus_{\gamma\in\Gamma} R_\gamma$ be a $\Gamma$-graded ring,
where $(\Gamma,<)$ is an ordered group. Let $\mathcal{P}\in\Spec_\Gamma(R)$
with corresponding epic $R$-division ring $(K,\varphi)$. Consider
$K((\Gamma,<))$. Then for each $\Omega\lhd \Gamma$, we get that
$\Spec_{\Gamma/\Omega}(R)\rightarrow \Spec_\Gamma(R)$ is surjective. 
Indeed, if $\mathcal{Q}\in\Spec_{\Gamma/\Omega}(R)$ is the corresponding $\Gamma/\Omega$-graded
prime matrix ideal to the $\Gamma/\Omega$-graded epic $R$-division ring
$(K(\Omega),\phi_\Omega\varphi)$, then $\mathcal{Q}\mapsto\mathcal{P}$ by Proposition~\ref{prop:almostdivisionring}.

We would like to remark that $D(\Gamma)$, the division subring of
$D((\Gamma;<))$ generated by $D$, does not depend on the order
$<$ of $\Gamma$ by \cite{Hughes} or \cite{DicksHerberaSanchez}. Hence,
since $D(\Omega)$ is just $\DC(\phi_\Omega)$, then $D(\Omega)$ does not depend on the order $<$ of
$\Gamma$. 

\medskip

We end this section with a concrete application of the results  in this section.  
Let $K$ be a field, $X$ be a nonempty set and $K\langle X\rangle$ be the free $K$-algebra on $X$. 
It is well known that $K\langle X\rangle$ has a universal division ring of fractions \cite[Section~7.5]{Cohnfreeeidealringslocalization}. 
Let now $\Gamma$ be a group and $X\rightarrow \Gamma$, $x\mapsto \hat{x}$, be a map. 
Then $K\langle X\rangle =\bigoplus_{\gamma\in\Gamma} K\langle X\rangle_\gamma$ is
a $\Gamma$-graded ring where $K\langle X\rangle_\gamma$ is the $K$-vector space
spanned by the monomials $x_1x_2\dotsc x_r$ such that $\hat{x}_1\hat{x}_2\cdots \hat{x}_r=\gamma$.
If $(\Gamma,<)$ is an ordered group, then $K\langle X\rangle$ has a $\Gamma$-graded
universal division ring of fractions by the foregoing example and Theorem~\ref{theo:relationsbetweenspectra}(3),(4).


\section{Inverse limits and ultraproducts in the category of graded epic $R$-division rings}
\label{sec:inverselimitsandultraproducts}

For details on filters and ultrafilters we refer the reader to \cite{BourbakiGeneralTobology}.

Let $I$ be a nonempty set. A \emph{filter} on $I$ is a set $\mathfrak{F}$ of subsets of $I$
which has the following properties
\begin{enumerate}[(F1)]
	\item Every subset of $I$ that contains a set of $\mathfrak{F}$ belongs to $\mathfrak{F}$.
	\item Every finite intersection of sets of $\mathfrak{F}$ belongs to $\mathfrak{F}$.
	\item The empty set is not in $\mathfrak{F}$.
\end{enumerate}
The set of filters on $I$ is partially ordered by inclusion. 
An \emph{ultrafilter} on $I$ is a maximal filter.
By \cite[Theorem~1, p.60]{BourbakiGeneralTobology}, each filter is contained in an ultrafilter. An ultrafilter $\mathfrak{U}$ on $I$ has the following property: if $J,K$ are subsets of $I$ such that
$J\cup K=I$, then either $J\in\mathfrak{U}$ or $K\in \mathfrak{U}$.

The concrete ultrafilters we will be dealing with are constructed as follows. Let $(I,\leq)$ be a directed
preordered set. For each $i\in I$, the set $\mathfrak{S}(i)=\{j\in I\colon i\leq j\}$
is called a \emph{section} of $I$ relative to $i$. The set $\mathfrak{S}$ consisting
of all sections relative to elements of $I$ is a filter base and there exists a filter containing 
$\mathcal{S}$ \cite[Proposition~2, p.59]{BourbakiGeneralTobology}. Therefore there exists 
an ultrafilter of $I$ containing $\mathfrak{S}$. 

\bigskip

Let $\Gamma$ be a group.

Let $I$ be a set and $\mathfrak{U}$ be an ultrafilter on $I$. For each $i\in I$,
let $R_i=\bigoplus_{i\in I}R_{i\gamma}$ be a $\Gamma$-graded ring.
We proceed to define the graded ultraproduct of the family $\{R_i\}_{i\in I}$ following \cite{IonNita}.
Consider the ring $P=\prod_{i\in I} R_i$ and consider the following subset $S$ of $P$
$$S=\bigoplus_{\gamma\in\Gamma}\left(\prod_{i\in I} R_{i\gamma} \right).$$
Note that $S$ is a subring of $P$ which is $\Gamma$-graded with $S_\gamma=\prod_{i\in I}R_{i\gamma}$.
For each $\gamma\in\Gamma$, if $x=(x_{i\gamma})_{i\in I}\in S_\gamma$, let
$z(x)=\{i\in I\colon x_{i\gamma}=0\}$. The set $Z_\gamma=\{x\in S_\gamma\colon z(x)\in\mathfrak{U}\}$
is an additive subgroup of $S_\gamma$. Moreover, if $y\in S_\delta$ and $x\in Z_\gamma$,
then $yx\in Z_{\delta\gamma}$ and $xy\in Z_{\gamma\delta}$. Therefore $Z=\bigoplus_{\gamma\in\Gamma}Z_\gamma$
is a graded ideal of $S$. Then the $\Gamma$-graded ring $U=S/Z$ is called the \emph{graded ultraproduct}
of the family of $\Gamma$-graded rings $\{R_i\}_{i\in I}$.

A homogeneous element element $x\in U_\gamma$ is the class of an element $(x_i)_{i\in I}\in S_\gamma$,
where each $x_i\in R_{i\gamma}$. We will write $x=[(x_i)_{i\in I}]_\mathfrak{U}$.
Observe that if $x=[(x_i)_{i\in I}]_\mathfrak{U}$ and $y=[(y_i)_{i\in I}]_\mathfrak{U}$,
then $x=y$ if and only if the set $\{i\in I\colon x_i=y_i\}\in\mathfrak{U}$.

Suppose that $(R_i,\varphi_i)$ is a $\Gamma$-graded $R$-ring for each $i\in I$. 
Hence $\varphi_i\colon R\rightarrow R_i$ is a homomorphism of $\Gamma$-graded rings. Then
there exists a unique homomorphism of rings $\varphi'\colon R\rightarrow \prod_{i\in I}R_i$
such that $\pi_i\varphi'=\varphi_i$ for each $i\in I$. Observe that $\im \varphi'\subseteq S$. Composing
with the natural homomorphism $S\rightarrow S/Z=U$, we obtain a homomorphism
of $\Gamma$-graded rings $\varphi\colon R\rightarrow U$. Hence $U$ is a $\Gamma$-graded $R$-ring in
a natural way. This fact and the following lemma will be very useful in this section.

\begin{lemma}\label{lem:ultraproduct}
Let $\Gamma$ be a group. Let $I$ be a nonempty set and $\mathfrak{U}$ be an ultrafilter on $I$.
\begin{enumerate}[\rm(1)]
	\item If $R_i$ is a $\Gamma$-graded division ring for each $i\in I$, then the ultraproduct $U$
	of the family $\{R_i\}_{i\in I}$ is a $\Gamma$-graded division ring.
	\item If $R_i$ is a $\Gamma$-graded local ring with graded maximal ideal $\mathfrak{m}_i$ for each $i\in I$,
	then the ultraproduct $U$ of the family $\{R_i\}_{i\in I}$ is a $\Gamma$-graded local ring
	with residue $\Gamma$-graded division ring $V$, the ultraproduct of the family of $\Gamma$-graded division rings $\{R_i/\mathfrak{m}_i\}_{i\in I}$.
\end{enumerate}
\end{lemma}

\begin{proof}
(1) Let $x\in U_\gamma$. Then $x=[(x_i)_{i\in I}]_\mathfrak{U}$ for some $x_{i}\in R_{i\gamma}$.
If $x$ is nonzero, then $J=\{i\in I\colon x_i\neq 0\}\in \mathfrak{U}$. For each $i\in I$, define
$${x_i}'=\left\{\begin{array}{ll} x_i^{-1} & \textrm{if } i\in J 
\\ 0 & \textrm{if } i\notin J\end{array}\right. .$$
Notice that ${x_i}'\in R_{i\gamma^{-1}}$ for each $i \in I$.
Then $x'=[({x'}_i)_{i\in I}]_\mathfrak{U}\in U_{\gamma^{-1}}$ and $xx'=x'x=1$, as desired.

(2) A homogeneous element $x=[(x_i)_{i\in I}]_\mathfrak{U}$ is invertible in $U$ if and only if
the set $\{i\in I\colon x_i \textrm{ is invertible in }U\}\in\mathfrak{U}$ if and only if
$\{i\in I\colon x_i\notin \mathfrak{m}_i\}\in \mathfrak{U}$. Therefore, the ideal $\mathfrak{m}$ 
generated by the homogeneous noninvertible elements, that is, the set 
$\{[(x_i)_{i\in I}]_{\mathfrak{U}}\in \h(U)\colon \{i\in I\colon x_i\in\mathfrak{m}_i\}\in\mathfrak{U}\}$,
is a proper ideal of $U$. Hence $U$ is a graded local ring.

It is not difficult to prove that the projections $R_i\rightarrow R_i/\mathfrak{m}_i$, 
$a\mapsto\overline{a}$, induce a surjective homomorphism of $\Gamma$-graded rings
$U\rightarrow V$, $[(x_i)_{i\in I}]_\mathfrak{U}\mapsto [(\overline{x_i})_{i\in I}]_\mathfrak{U}$.
Note that a homogeneous element $[(x_i)_{i\in I}]_\mathfrak{U}$ is in the kernel if and only if
it belongs to $\mathfrak{m}$.
\end{proof}

The following corollary was used in Section~\ref{sec:grmatrixideals}

\begin{corollary}\label{coro:ultraproductoffractions}
Let $R$ be a $\Gamma$-graded domain. Suppose that, for each $a\in\h(R)\setminus\{0\}$, there
exists a homomorphism of $\Gamma$-graded rings $\varphi_a\colon R\rightarrow K_a$,
where $K_a$ is a $\Gamma$-graded division ring such that $\varphi_a(a)\neq 0$. Then
there exists a $\Gamma$-graded epic $R$-division ring of fractions.
\end{corollary}

\begin{proof}
Let $I=\h(R)\setminus\{0\}$. For each $a\in I$, let $I_a=\{\lambda\in I\colon \varphi_\lambda(a)\neq 0\}$.
Let $E=\{a_1,\dotsc,a_n\}$ be a finite subset of $I$.
Then $\bigcap_{i=1}^n I_{a_i}\neq \emptyset$, because
$\varphi_{a_1\dotsb a_n}(a_i)\neq 0$ for each $i=1,\dotsc,n$. Hence the set
$\mathfrak{B}=\mathcal\{I_a\colon a\in I\}$ is a set of subsets of $I$ such that
no finite subset of $\mathfrak{B}$ has empty intersection. By 
\cite[Proposition~1,p.58]{BourbakiGeneralTobology}, there exists a filter on $I$ containing
$\mathfrak{B}$. By \cite[Theorem~1, p.60]{BourbakiGeneralTobology}, there exists
an ultrafilter $\mathfrak{U}$ on $I$ containing $\mathfrak{B}$. By Lemma~\ref{lem:ultraproduct}(1),
the ultraproduct $U$ of the family $\{K_a\}_{a\in I}$ is a $\Gamma$-graded
division ring and there exists a homomorphism of $\Gamma$-graded rings
$\varphi\colon R\rightarrow U$, defined by 
$\varphi(x)=[(\varphi_a(x))_{a\in I}]_{\mathfrak{U}}$. Since the set
$I_x\in\mathfrak{U}$, then $\varphi(x)\neq 0$ for each $x\in\h(R)\setminus\{0\}$.
Therefore $\varphi$ is injective. 
\end{proof}

Let $R$ be a $\Gamma$-graded ring. Consider the category $\mathcal{E}_R$ of $\Gamma$-graded epic $R$-divison
rings with specializations as morphisms defined in Section~\ref{sec:categoryspecializations}.

First we look at how  inverse systems are in this category. An inverse system in $\mathcal{E}_R$
is a pair $((K_i,\varphi_i)_{i\in I}, (\psi_{i,j})_{i\geq j})$ where
$(I,\leq)$ is a directed preordered set, $(K_i,\varphi_i)$ is a $\Gamma$-graded
epic $R$-division ring for each $i\in I$, and $\psi_{i,j}$ is a specialization
from $(K_i,\varphi_i)$ to $(K_j,\varphi_j)$, $i\geq j$, such that 
\begin{equation}\label{eq:inversesystem}
\psi_{j,k}\circ\psi_{i,j}=\psi_{i,k} \textrm{ for all } i,j,k\in I, i\geq j\geq k.
\end{equation}
Observe that \eqref{eq:inversesystem} is superfluous because since the specializations already exist, and
there is at most one specialization between graded epic $R$-division rings, the equality in 
\eqref{eq:inversesystem} holds trivially.

Now we look at inverse limits. An inverse limit of the inverse system \linebreak
$((K_i,\varphi_i)_{i\in I}, (\psi_{i,j})_{i\geq j})$ in $\mathcal{E}_R$
is a pair $((K,\varphi),(\psi_i)_{i\in I})$ where $(K,\varphi)$ is a
$\Gamma$-graded epic $R$-division ring and $\psi_i$ is a specialization
from $(K,\varphi)$ to $(K_i,\varphi_i)$ for each $i\in I$ such that the following
properties are satisfied: 
\begin{enumerate}[(i)]
	\item $\psi_{i,j}\circ \psi_i=\psi_j$ for all $i,j\in I$, $i\geq j$.
	\item For each pair $((K',\psi'),(\psi_i')_{i\in I})$ that satisfies (i), i.e.
	$\psi_{i,j}\circ \psi_i'=\psi'_j$ for $i\geq j$, then there is a specialization
	$\phi\colon L\rightarrow K$ such that $\psi_i\circ \phi=\psi'$ for all $i\in I$.
\end{enumerate}
Again note that (i) and the equality of specializations in (ii) are superfluous.

\begin{theorem}
Let $R$ be a $\Gamma$-graded ring. Let $((K_i,\varphi_i)_{i\in I}, (\psi_{i,j})_{i\geq j})$
be an inverse system in $\mathcal{E}_R$ indexed on the directed nonempty  preordered set $(I,\leq)$.
Consider an ultrafilter $\mathfrak{U}$ on $I$ that contains all the sections $\mathfrak{S}(i)$, 
	$i\in I$, of $I$. Set $$\Sigma_i=\{A\in\mathfrak{M}(R)\colon A^{\varphi_i} \textrm{ is
	invertible over } K_i\},\quad i\in I.$$
The following assertions hold true.
\begin{enumerate}[\rm(1)]
	\item There exists a $\Gamma$-graded epic $R$-division ring $(K,\varphi)$ which is the inverse limit
	of $((K_i,\varphi_i)_{i\in I}, (\psi_{i,j})_{i\geq j})$.
	\item If $\Sigma=\{A\in\mathfrak{M}(R)\colon A^\varphi \textrm{ is invertible over }K\}$,
	then $\Sigma=\bigcup_{i\in I}\Sigma_i$.
	\item  The $\Gamma$-graded epic $R$-division ring determined by the
	ultraproduct $U$ of the family $\{(K_i,\varphi_i)\}_{i\in I}$ equals $(K,\varphi)$.
  
	\item The ultraproduct $V$ of the family $\{R_{\Sigma_i}\}_{i\in I}$ is a $\Gamma$-graded local
	ring with residue $\Gamma$-graded division ring equal to $U$.
	\item $R_\Sigma$ embeds in $V$.
\end{enumerate}
\end{theorem}

\begin{proof}
For each $i\in I$, let $\mathcal{P}_i$ be the singular kernel of $\varphi_i\colon R\rightarrow K_i$.

(1) Since there exists a specialization $\psi_{i,j}\colon K_i\rightarrow K_j$ for $i,j\in I$, $i\geq j$,
then $\mathcal{P}_i\subseteq\mathcal{P}_j$ by Corollary~\ref{coro:specialization}. Thus,
the family of prime matrix idelas $\{\mathcal{P}_i\}_{i\in I}$ is directed from below. 
Set $\mathcal{P}=\bigcap_{i\in I}\mathcal{P}_i$. It is not difficult to prove that $\mathcal{P}$
satisfies (PM1)--(PM3), (PM5) and (PM6) in the definition of gr-prime matrix ideal. To show
that $\mathcal{P}$ satisfies (PM4), let $A,B\in\mathfrak{M}(R)\setminus\mathcal{P}$. There exist
$i,j\in I$ such that $A\notin \mathcal{P}_i$ and $B\notin \mathcal{P}_j$. Since $I$ is directed, there
exists $k\in I$ such that $i\leq k$, $j\leq k$. Thus $\mathcal{P}_k\subseteq\mathcal{P}_i$ and
$\mathcal{P}_k\subseteq \mathcal{P}_j$, and both $A$ and $B$ do not belong to the 
gr-prime matrix ideal $\mathcal{P}_k$. Hence $A\oplus B\notin\mathcal{P}_k$, and therefore
$A\oplus B\notin\mathcal{P}$.

Let $(K,\varphi)$ be the $\Gamma$-graded epic $R$-division ring corresponding to $\mathcal{P}$.
Since $\mathcal{P}\subseteq \mathcal{P}_i$ for all $i\in I$, there exists a unique specialization
$\psi_i\colon K\rightarrow K_i.$ Consider now a pair $((L,\varphi'),(\psi_i')_{i\in I})$ where
$L$ is a $\Gamma$-graded epic $R$-division ring and $\psi_i'$ is a gr-specialization from
$(L,\varphi')$ to $(K_i,\varphi_i)$ for each $i\in I$. Let $\mathcal{Q}$ be the singular kernel
of $\varphi'\colon R\rightarrow L$. Then $\mathcal{Q}\subseteq\mathcal{P}_i$ for each $i\in I$. 
Hence $\mathcal{Q}\subseteq \mathcal{P}$. Therefore, there exists a specialization 
$\phi\colon L\rightarrow K$ by Corollary~\ref{coro:specialization}.

(2) By (1), the singular kernel of $(K,\varphi)$ equals $\mathcal{P}=\bigcap_{i\in I}\mathcal{P}_i$.
Hence, $\Sigma=\mathfrak{M}(R)\setminus\mathcal{P}=\bigcup_{i\in I}\Sigma_i$.

(3) Let $U$ be the graded ultraproduct of the family $\{(K_i,\varphi_i)\}_{i\in I}$ of $\Gamma$-graded
epic $R$-division rings, and let
$\varphi\colon R\rightarrow U$ be the canonical homomorphism of $\Gamma$-graded rings.
Let $L$ be the $\Gamma$-graded division subring of $U$ generated by the image of $\varphi$.
Consider  the $\Gamma$-graded epic $R$-division ring $(L,\varphi)$.

Let $A\in\bigcup_{i\in I}\Sigma_i$ be such that $A\in M_n(R)[\overline{\alpha}][\overline{\beta}]$for some
$\overline{\alpha},\overline{\beta}\in \Gamma^n$.
By Theorem~\ref{theo:specialization}, note that $\Sigma_j\subseteq\Sigma_i$ for all $i,j$, $i\geq j$.
Let $t\in I$ be such that $A\in\Sigma_t$. Then $A\in\Sigma_i$ for all $i\in\mathfrak{S}(t)$.
Define $M_i=\left\{\begin{array}{ll} (A^{\varphi_i})^{-1} & \textrm{if } i\in\mathfrak{S}(t) \\
0 & \textrm{if } i\notin\mathfrak{S}(t). 
 \end{array}\right.$
If $M_i$ has $(u,v)$-entry $m_{uv}^i\in (K_i)_{\beta_u\alpha_v^{-1}}$, let
$M=(m_{uv})\in M_n(L)$ where $m_{uv}=[(m^i_{uv})_{i\in I}]_\mathfrak{U}$. Note now that
$A^\varphi M=MA^{\varphi}=I$. Indeed, if $AM=([c^i_{uv}]_\mathfrak{U})$, then
$[c^i_{uv}]_\mathfrak{U}\in U_e$, $\mathfrak{S}(t)\subseteq\{i\in I\colon c_{uu}^i=1\}$
and $\mathfrak{S}(t)\subseteq\{i\in I\colon c^i_{uv}=0,\, u\neq v\}$. Similarly one can show that
$MA=I$.

Conversely, suppose that $A\in\Sigma_n[\overline{\alpha}][\overline{\beta}]$. Let $(A^{\varphi})^{-1}=(m_{uv})$, $u,v=1,\dotsc,n$,
where $m_{u,v}=[(m_{uv}^i)_{i\in I}]_\mathfrak{U}$.
For each $i\in I$, let $M_i\in M_n(K_i)[\overline{\alpha}][\overline{\beta}]$ be the matrix
with $(u,v)$-entry $m^i_{uv}$. Since $A^\varphi(A^\varphi)^{-1}$ and $(A^\varphi)^{-1}A^\varphi$
equal the identity matrix, the set 
$$\mathcal{U}=\{i\in I\colon M_i \textrm{ is the inverse of } A^{\varphi_i} \textrm{ over } K_i\}\in \mathfrak{U}.$$
In particular $\mathcal{U}\neq \emptyset$. If $i\in\mathcal{U}$, then $A\in\Sigma_i$, and
therefore $A\in\bigcup_{i\in I}\Sigma_i$. Hence, the singular kernel of $(L,\varphi)$ is
$\mathcal{P}=\mathfrak{M}(R)\setminus(\bigcup_{i\in I}\Sigma_i)$.

(4) is Lemma~\ref{lem:ultraproduct}(2).

(5) First observe that $R_\Sigma=\varinjlim R_{\Sigma_i}$. Let 
$\tau_i\colon R_{\Sigma_i}\rightarrow R_\Sigma$ and
$\lambda_i\colon R\rightarrow R_{\Sigma_i}$ be the natural homomorphism of $\Gamma$-graded rings.
By (3), there exists a natural homomorphism of $\Gamma$-graded rings $\rho\colon R_\Sigma\rightarrow U$.
Since $V$ is $\Gamma$-graded local with residue $\Gamma$-graded division ring equal to $U$,
the universal property of $R_\Sigma$, induces a  homomorphism of $\Gamma$-graded rings $\rho\colon R_\Sigma\rightarrow V$.
We must prove that $\rho$ is injective. Suppose that $x\in\ker\rho$ is homogeneous.
Then $x$ is the $(u,v)$-entry of the 
inverse of a matrix $A^\varphi$ with $A\in\Sigma=\bigcup_{i\in I}\Sigma_i.$
Let $i\in I$ be such that $A\in\Sigma_i$. Note that $A\in\Sigma_l$ for all $l\in I$, $i\leq l$.
Let $x^l_{uv}$ be the $(u,v)$-entry of $A^{\lambda_l}$, $i\leq l$.
Then $\tau_i(x^l_{uv})=x$. Then $\rho(x)=[(z_{uv}^l)_{l\in I}]_\mathfrak{U}$
where $z_{uv}^l=\left\{\begin{array}{ll}x_{u,v}^l & \textrm{if } i\leq l\\
0 & \textrm{otherwise} \end{array}\right.$. 
Since $\rho(x)=0$, then $[(z_{uv}^l)_{l\in I}]_\mathfrak{U}=0$ and
the set $\{j\in I\colon z_{uv}^j=0\}\in\mathfrak{U}$. It implies
that $\tau_i(x_{uv}^j)=0$. Therefore $x=0$.
\end{proof}

Note that the  proof of  (3)  shows also (1) in a more elementary way.



\bibliographystyle{amsplain}
\bibliography{grupitosbuenos}

\end{document}